\documentclass{article}
\usepackage[all]{xy}\usepackage[latin1]{inputenc}       
\usepackage[dvips]{graphics,graphicx}
\usepackage{amsfonts, amssymb, amsmath,xcolor,mathrsfs, amstext, comment, mathtools}
\usepackage{amsbsy, amsopn, amscd, amsxtra, amsthm, authblk, enumitem, dsfont}
\usepackage{upref}
\usepackage{geometry}
\usepackage[colorlinks,
            linkcolor=blue,
            anchorcolor=green,
            citecolor=blue
            ]{hyperref}
\usepackage[nameinlink,capitalize]{cleveref}
\usepackage{microtype}
\usepackage[displaymath]{lineno}
\usepackage{float}
\usepackage{bbm}
\usepackage{tikz, tikz-3dplot}
\usetikzlibrary{calc}
\usepackage{subcaption}
\usepackage{bm}
\usetikzlibrary{arrows.meta,calc,patterns}
\definecolor{xcol}{RGB}{28,91,156}
\definecolor{ycol}{RGB}{194,84,38}

\allowdisplaybreaks

\numberwithin{equation}{section}
\newtheorem{theorem}{Theorem}[section]
\newtheorem{lemma}[theorem]{Lemma}
\newtheorem{proposition}[theorem]{Proposition}
\newtheorem{corollary}[theorem]{Corollary}
\newtheorem{definition}[theorem]{Definition}
\newtheorem{remark}[theorem]{Remark}

\renewcommand{\P}{\mathbb{P}}
\newcommand{\E}{\mathbb{E}}
\newcommand{\Z}{\mathbb{Z}}
\newcommand{\R}{\mathbb{R}}

\renewcommand{\1}{\mathbf 1}

\usepackage{graphicx} 
\newcommand{\diam}{\mathrm{diam}}
\newcommand{\fC}{\mathsf{C}}
\newcommand{\fD}{\mathsf{D}}

\newcommand{\dist}{\mathrm{dist}}
\newcommand{\fB}{\mathsf{B}}

\newcommand{\xx}{\mathsf{x}}
\renewcommand{\sf}[1]{\mathsf{#1}}
\newcommand{\yy}{\mathsf{y}}
\newcommand{\zz}{\mathsf{z}}
\newcommand{\ii}{\mathsf{i}}
\newcommand{\NN}{\mathsf{N}}
\newcommand{\aux}{\mathrm{AUX}}

\newcommand{\bC}{\bar{C}}
\newcommand{\bc}{\bar{c}}
\newcommand{\con}[1]{\stackrel{#1}{\longleftrightarrow}}
\newcommand{\ncon}[2][1]{
\mathrel{%
    \ooalign{%
      $\overset{#2}{\xleftrightarrow{\hspace{#1cm}\vspace{1cm}}}$\cr
      \hidewidth\raisebox{0.1ex}{{$/$}}\hidewidth\cr 
    }%
  }%
}

\title{Exponential decay for the two-dimensional Heisenberg model}
\author{Jian Ding\thanks{New Cornerstone Science Laboratory, School of Mathematical Sciences, Peking University} \qquad\qquad Fenglin Huang$^*$ \qquad\qquad Peng Yang$^*$}
\date{\today}

\begin{document}

\maketitle

\begin{abstract}
    In this work, we prove exponential decay for the spin correlation of the spin $O(\NN)$ model for $\NN\geq3$ at every inverse temperature $ \beta>0 $. In particular, our work confirms Polyakov's prediction in 1975 on the Heisenberg model, i.e., the case of $\NN=3$.
    To this end, we investigate a generalized $O(\NN)$ model where $ \NN $ is a positive \emph{real number}, as follows: we sample a loop soup on the metric graph (say of a box) with intensity $ \tfrac{\NN}{2} $ subject to the conditioning  that the local time at each lattice vertex is $\beta$, and then for each loop we independently assign it to copies $1\leq \ii\leq \lceil\NN\rceil$ such that it is assigned to copy $\ii$ with probability $ \tfrac{1}{\NN} $ for $ 1\leq \ii\leq \lfloor\NN\rfloor $. Then, by Lupu's isomorphism theorem, we can construct a spin field by assigning a random sign to each loop soup cluster in each copy independently, such that the field coincides with the spin $O(\NN)$ model when $\NN$ is an integer. We then establish exponential decay for every real-valued $\NN>2$. Our result reveals an underlying mechanism for exponential decay, which emerges as soon as $\NN>2$ (i.e., as soon as beyond the range of the XY model) and does not necessarily rely on the rotation symmetry in the sphere $ \mathbb{S}^{2} $.
\end{abstract}


\section{Introduction}
\label{sec:intro}

For $\NN\geq1$, we study the spin $O(\NN)$, in particular, the Heisenberg model (where $\NN=3$) on the square lattice. Consider a finite graph $ G=(V, E)\subset\Z^2 $, and let $ \partial V:=\{u\in\Z^2\setminus V: v\sim u\text{ for some }v\in V \} $ be the external boundary of $V$ (here $ u\sim v $ means $ |u-v|_2=1 $). Let $ \mathbb{S}^{\mathsf{N}-1}\subset\R^{\NN} $ be the $ (\mathsf{N}-1) $-sphere. For a boundary condition $\zeta\in(\mathbb{S}^{\mathsf{N}-1})^{\partial V}$, define the spin $O(\NN)$ Hamiltonian
\begin{equation}
    \label{eq:def-Hamiltonian}
    H^{\zeta}_G(\sigma):=-\sum_{u\sim v}\sigma_u\cdot\sigma_v - \sum_{(u,v): u\in V, v\in\partial V,u\sim v}\sigma_u\cdot\zeta_v
\end{equation}
for a spin configuration $ \sigma=(\sigma_v)_{v\in V}\in (\mathbb{S}^{\NN-1})^{V} $.
For an inverse temperature $\beta>0$, define the Gibbs measure $ \langle \cdot \rangle^{\zeta}_{G,\beta} $ by
\begin{equation}
    \label{eq:def-Heisenberg}
    \langle F(\sigma) \rangle^{\zeta}_{G,\beta}:=\frac{1}{Z^{\zeta}_{G,\beta}}\int F(\sigma)e^{-\beta H^{\zeta}_{G}(\sigma)}d\sigma,\qquad Z^{\zeta}_{G,\beta}:=\int e^{-\beta H^{\zeta}_G(\sigma)}d\sigma
\end{equation}
for all measurable $ F $, where $ d\sigma=\prod_{v\in V}d\sigma_v $ is the product of Haar measures.

The main result of this paper is to prove Polyakov's prediction~\cite{Polyakov75} on the exponential decay in two dimensions.
Let $ B_L:=B(o, L)\subset\Z^2 $ be a two-dimensional box centered at $o$ with side length $L\geq1$.
\begin{theorem}
    \label{thm-1.1}
    Fix $\mathsf N \geq 3$.
    For all $\beta>0$, there exists a constant $ c_{\beta}>0 $ such that for all $L\geq1$,
    \begin{equation}
        \label{eq:thm1.1}
        \langle \sigma_o \rangle^{\zeta}_{B_L, \beta} \leq c^{-1}_{\beta}e^{-c_{\beta}L}
    \end{equation}
    where $\zeta$ is a boundary condition such that the following holds for all $ 1\leq \ii\leq \NN $, $ \zeta^{\ii}_{\cdot} $ are the same for all $ v\in\partial B_L $. In particular, in the free boundary condition $\mathsf{f}$, (formally) corresponding to $ \zeta\equiv0 $  and $u,v\in B_{L/2}$, we have that
    \begin{equation}
        \label{eq:thm1.1-1}
        \langle\sigma_u\sigma_v\rangle_{B_L, \beta}\leq c_{\beta}^{-1}e^{-c_{\beta}|u-v|_2}.
    \end{equation}
\end{theorem}
\begin{remark}
    \label{rmk: for thm1.1, with coupling constants}
    The same proof applies verbatim to prove the same result when, in the definition of~\eqref{eq:def-Hamiltonian}, the dot product of spins $ \sigma_u\cdot\sigma_v $ is multiplied by $ \sqrt{J_u J_v} $ as long as $ 0\leq J_u\leq1 $ for all $u$.
\end{remark}
\begin{remark}
    \label{rmk: 1.3}
    It is natural to expect that~\eqref{eq:thm1.1} holds for all boundary conditions $\xi$, but it does not seem to be that immediate. On Oct 6, 2026, OpenAI released a proof of~\eqref{eq:thm1.1-1}, and our feeling is that it may not yield~\eqref{eq:thm1.1} with all boundary conditions in a straightforward way either.
\end{remark}
Our proof approach for Theorem~\ref{thm-1.1} is via a \emph{generalized} spin $O(\NN)$ model, which is defined for all positive real numbers $ \NN $ from the perspective of loop soups.
To this end, we let $ \widetilde{G} $ be the metric graph of $G$ (we may introduce killing so that the random walk thereon is transient) and sample a loop soup $ \mathcal{L}^{\NN/2}_G $ on $\widetilde{G}$ with intensity $ \NN $. For notation convenience, we will denote by $ \overline{\NN}=\lceil\NN\rceil $ and $ \underline{\NN}=\lfloor\NN\rfloor $ (as the minimal integer great equal to $\mathsf N$ and the maximal integer less equal to $\mathsf N$ respectively), and we will construct an $ \mathbb{R}^{\overline{\NN}} $-valued spin field from this sampled loop soup (where each coordinate corresponds to a copy). For each sampled loop, we independently assign it to all the $ \overline{\NN} $ copies such that it  is assigned to the $ \ii $-th copy with probability $\tfrac{1}{\NN}$ for $1\leq\ii\leq\underline{\NN}$ (thus, except for the $ \overline{\NN} $-th copy, each other copy corresponds to a loop soup with intensity $\tfrac{1}{2}$, referred to as a critical loop soup). Further, we denote by $ (L^{\ii}_v)_{1\leq\ii\leq\overline{\NN}, v\in V} $ the local time field where $ L^{\ii}_v $ is the combined local time at $v$ for all loops assigned to the $\ii$-th copy. Then for any non-negative $ \beta=(\beta_v)_{v\in V} $, we may impose the following local time conditioning:
\begin{equation}
    \label{eq: def loop soup local time condition}
    L^{ *}_v = \beta_v,\qquad\forall v\in V.
\end{equation}
where $ L^{ *}_v=\sum^{\NN}_{\ii=1}L^{\ii}_v $ is the local time of the soup $ \mathcal{L}^{\NN/2}_G $ at $v$.
Then, by Lupu's isomorphism theorem (see~\cite[Theorem~1]{Lupu14}), for the special case when $ \NN $ is an integer, there exists a coupling with the spin $O(\NN)$ model such that
\begin{equation}
    \label{eq: local-time --- spin}
    L^{\ii}_v=\beta_v\cdot(\sigma^{\ii}_v)^2,\qquad\forall v\in V,
\end{equation}
and that the signs of $ \sigma $ may be sampled by assigning  each loop soup cluster in each copy a random sign. We will further elaborate on this definition and this coupling in Section~\ref{subsec: conditional loop soup representation}, and for now take it for granted that we have justified our extension of the spin $O(\NN)$ model to a model with continuous parameterization of $ \NN $.
For exposition convenience in the introduction, we will focus only on the case when the loop soup is taken with free boundary condition, and we denote by $ \nu^{\mathsf{f}, \mathsf{N}}_{\widetilde{G},\beta} $ as the law for the loop soups in copies $ 1, \dots, \overline{\NN} $ subject to the local time conditioning~\eqref{eq: def loop soup local time condition}. In the rest of the paper, we will assume that $ 2<\NN\leq3 $ and the more general case of $ \NN>3 $ follows verbatim. With this assumption, we then identify the three copies as the $ \xx,\yy,\zz $ copies, and we note that the intensity for the loop soup in the $ \zz $-copy is not necessarily $1/2$.
That said, the $\xx$ and $\yy$ copies are still loop soups with intensity $\frac{1}{2}$ (of course subject to the conditioning~\eqref{eq: def loop soup local time condition}), and thus by conditioning on the local time field of the $\zz$-copy, we can still find an XY model, see~\eqref{eq:intro-xy-projection} below for further details.
\begin{theorem}
    \label{thm-1.2}
    Fix $ 2<\NN\leq3 $.
    For all nonnegative $ \beta $ with $ \max_v\beta_v\leq\beta_0 $, there exists a constant $ c>0 $ depending only on $\beta_0, \NN$ such that for all $L\geq1$  and $u,v\in B_{L/2}$,
    \begin{equation}
        \label{eq:thm1.2}
        \nu^{\mathsf{f},\NN}_{\widetilde{B}_L,\beta}\Big(\big\{
        u\text{ is connected to } v \text{ in the }\xx\text{-copy}
        \big\}\Big) \leq c^{-1}e^{-c|u-v|_{\infty}},
    \end{equation}
where the connection above is by a path in the metric graph with all $\xx$-copy local times being positive (see \eqref{eq: local time definition}).
\end{theorem}

\begin{remark}
    \label{rmk: 1.4}
    As hinted earlier, Theorem~\ref{thm-1.2} then implies  exponential decay in Theorem~\ref{thm-1.1} as desired (modulo consideration of boundary conditions). From what we understand, the proof of OpenAI for Theorem~\ref{thm-1.1} relies crucially on the non-commutativity of rotations on $ \mathbb{S}^2 $ and on a magically (in our opinion) constructed family of spin waves on $ \mathbb{S}^2 $. In a way, one may argue that this method  is more from the ``quantum'' side of the $O(3)$ model. Our approach, on the other hand, is from the perspective of loop soups and (as it will be seen) employs a number of carefully chosen spin waves on $\mathbb{S}^1$, and is perhaps more from the ``classic'' side of the $O(3)$ model. The core ideas of the two approaches, encapsulated in the proof of Theorem~\ref{thm-annulus-local-time-crossing}, are drastically different. While we admire the magical nature of the OpenAI construction and admit that our proof is longer, we do feel that our proof has the advantage of explaining the underlying mechanism by revealing an interesting phenomenon that a slight perturbation on the XY model from the loop soup perspective produces exponential decay. In a recent work of Garban and Montgolfier~\cite{GarbandeMontgolfier} a polynomial lower bound on spin correlation for large $\beta$ was shown for a spin model with Hamiltonian $ H(\sigma)=-\sum_{u\sim v}(\sigma^{\xx}_u\sigma^{\xx}_v+\sigma^{\yy}_u\sigma^{\yy}_v+\lambda\sigma^{\zz}_u\sigma^{\zz}_v) $ for $ \lambda\in(-1,1) $. Thus, this demonstrates that a slight perturbation of the spin Hamiltonian of the Heisenberg model destroys the exponential decay. The combination of~\cite{GarbandeMontgolfier} and this work seems to reveal some interesting phenomena on the exponential decay of continuous spin systems.
\end{remark}

In what follows, we outline our proof ideas at a very high level, and in the meantime, we emphasize key ingredients and inspirations we borrow from previous works.
We refer readers to Peled and Spinka~\cite{Peled2019_On_Notes} for an excellent account on the topic of $O(\NN)$ models.

We omit discussions and literature review on the case for $\mathsf{N}=1$ (i.e., the Ising model), since it is rather different from the case of continuous spins, whereas $\mathsf{N}\geq2$. In particular, the case of $\NN=2$ is the XY model, with spins $(\cos\theta_v,\sin\theta_v)\in\mathbb S^1$.
In two dimensions, the Berezinskii--Kosterlitz--Thouless (BKT) transition of the XY model was predicted in~\cite{Berezinskii71,KT73,Kosterlitz74}, where vortex--antivortex pairs bind at low temperature, producing a phase with polynomially decaying correlations rather than spontaneous magnetization as for the Ising model. The absence of spontaneous magnetization at positive temperature is a consequence of the Mermin--Wagner theorem in its classical formulation~\cite{MW66,Mermin67}.
A key quantity of interest is the infinite-volume two-point function, defined by the free boundary limit as follows:
\[
g_\beta(x) :=\lim_{\Lambda\uparrow\mathbb Z^2} \big\langle\theta_0-\theta_x)\big\rangle^{\mathrm{XY}}_{\Lambda,\beta}.
\]
Then, the complex-rotation argument of McBryan and Spencer~\cite{MS77}, implementing the spin-wave mechanism, gives that for a constant $C_{\beta}>0$
\[
0\leq g_\beta(x)\leq C^{-1}_\beta(1+|x|)^{-C_\beta}, \qquad 0<\beta<\infty.
\]
The complementary power-law lower bound at sufficiently low temperature was first proved by Fr\"ohlich and Spencer~\cite{FS81} (see also~\cite{Peled_BKT_Notes} for a streamlined presentation), through a multiscale analysis of the dual Coulomb gas, with a lower-bound exponent that tends to zero as $\beta\to\infty$. More recently, Lammers' delocalization theorem for integer-valued height functions on cubic planar graphs~\cite{Lammers22} led to alternative proofs: van Engelenburg and Lis~\cite{vEL23} used a loop representation of spin correlations, while Aizenman, Harel, Peled, and Shapiro~\cite{AHPS21} related dual height level sets to spin correlations in the Villain and XY models. In particular,~\cite[Theorem~1]{vEL23} gives the sharp dichotomy
\begin{equation}
    \label{eq:intro-xy-bkt}
    \begin{aligned}
        \beta<\beta_{\mathrm{BKT}}
        &\quad\Longrightarrow\quad g_\beta(x)\leq C^{-1}_\beta e^{-C_\beta|x|},\\
        \beta\geq\beta_{\mathrm{BKT}}
        &\quad\Longrightarrow\quad g_\beta(x)\geq c(1+|x|)^{-1},
    \end{aligned}
    \qquad \beta_{\mathrm{BKT}}\in(0,\infty).
\end{equation}
The height--spin duality has been developed further in~\cite{Lammers23,vEL25}.
In contrast, on $\mathbb Z^d$ with $d\geq3$, the XY model (and spin $O(\NN)$ model in general, for all $ \NN\geq1 $) has long-range order at sufficiently low temperature: its two-point function remains bounded away from zero at large distances. This was established in~\cite{FSS76} by Fr\"ohlich, Simon and Spencer using infrared bounds and reflection positivity; Garban and Spencer later obtained an alternative proof~\cite{GS22} without reflection positivity, through a disordered model on the Nishimori line and correlation inequalities.

For the two-dimensional Heisenberg model, the nonlinear sigma-model analysis of Polyakov~\cite{Polyakov75} predicted a different picture, where exponential decay of correlations at every finite inverse temperature, with a correlation length that becomes very large at low temperature, is expected.
Related foundational work includes the finite-energy continuum configurations of Belavin and Polyakov~\cite{BP75} and the renormalization-group analysis of Br\'ezin and Zinn-Justin~\cite{BZ76}.
Patrascioiu and Seiler~\cite{PS93,PS02} also developed a percolation approach to possible massless behavior, whose numerical evidence concerns a constrained version of the interaction.
A natural way to compare the Heisenberg and XY models is to condition on the norms of the projections onto the X and Y coordinates. Write
\begin{equation}
    \label{eq:intro-xy-projection}
    \begin{gathered}
        r_v=\sqrt{(\sigma_v^\xx)^2+(\sigma_v^\yy)^2} =\sqrt{1-(\sigma_v^\zz)^2},\qquad (\sigma^{\xx}_v, \sigma^{\yy}_v) = r_v(\cos \theta_v, \sin \theta_v),\\
        \mathbb P_{G,\beta}(d\theta\mid r) \propto \exp\!\left\{\sum_{uv\in E}\beta r_ur_v \cos(\theta_u-\theta_v)\right\} \prod_{v\in V}\frac{d\theta_v}{2\pi}.
    \end{gathered}
\end{equation}
Indeed, at fixed $r$ the interaction of the third coordinates is independent of the angles $\{\theta_v\}$, and integrating their signs yields precisely the XY model with nonnegative random couplings $K_{uv}=\beta r_ur_v$.
Aru, Garban, and Sep\'ulveda~\cite{AGS25} analyzed these couplings. In torus limits, after rotating the spin at the origin to $(0,0,1)$,~\cite[Theorem~1.3]{AGS25} showed that $\sqrt\beta(\sigma^x,\sigma^y)$ converges locally to a two-component Gaussian free field pinned at the origin as $\beta\to\infty$. Their estimates for the small-norm sets of this field imply exponential connectivity bounds for fixed-radius (denoted as $k$) neighbourhoods of edges with $K_{uv}\leq R$, for every prescribed $k$ and $R$ at sufficiently large $\beta$, see~\cite[Theorem~1.2 and Corollary~1.4]{AGS25}.
Thus, large projected couplings have strong connectivity.
In addition, in~\cite[Theorem~1.5]{AGS25} they constructed strongly mixing random subgraphs approaching $\mathbb Z^2$, with exponentially decaying connectivity in their complements, where the XY model still has exponential decay at an arbitrarily prescribed inverse temperature.
In a different direction, Dario and Garban~\cite[Theorem~1.3]{DG25} proved a BKT phase for the XY model on any supercritical Bernoulli percolation cluster. More precisely, if $\omega_v$ are independent Bernoulli variables of parameter $p>p_c^{\mathrm{site}}(\mathbb Z^2)$ and the XY couplings are $\lambda\omega_u\omega_v$, then, for $\lambda$ sufficiently large,
\begin{equation}
    \label{eq:intro-bernoulli-xy}
    \mathbb E_p\!\left[\big\langle\cos(\theta_0-\theta_x)\big\rangle^{\mathrm{XY}}_{\lambda,\omega}\right] \geq c_{p,\lambda}(1+|x|)^{-C_{p,\lambda}}, \qquad c_{p,\lambda},C_{p,\lambda}>0.
\end{equation}
The average over $\omega$ is part of this statement.
These results manifest the delicacy of the problem strongly,  where the geometry and dependence of the projected couplings matter; this provides very valuable philosophical guidance in our work.
In addition, on Oct 6, 2026, OpenAI released a proof for Theorem~\ref{thm-1.1}; in a very recent work~\cite{GarbandeMontgolfier}, a BKT transition was proved for the Heisenberg model with small perturbation, as discussed already in Remark~\ref{rmk: 1.4}. These works appeared during the completion stage of this manuscript, and to our best knowledge did not influence the mathematical arguments in this paper. 

The key challenge in proving Theorem~\ref{thm-1.2} is to detect the difference between $\NN = 2$ and $ \NN>2 $.
Despite the delicacy implied by the works of~\cite{GS22, AGS25}, we nevertheless pursue the approach of projecting to the XY model.
One natural thought would be that if, in the $\zz$-copy, we see ``many loops'', then it ``should'' reduce the crossing probability in the remaining XY copies.
Then, even if ``many loops'' in $\zz$-copy happen with a small (but fixed) probability at each scale (which will be referred to as defected region for the XY model), overall, the non-defective region forms a fractal percolation process and, in particular, has a small volume.
If the defective region indeed has XY local times 0, then indeed a simple spin-wave argument would imply that the crossing probability for the XY model in the non-defective region is very small, plausibly leading to an exponential decay with a standard coarse-graining step.
The philosophy here is close to the proof for the decay of the 2D random field Ising model~\cite{DingXia21, AHP20}, with a key difference that tortuosity~\cite{AizenmanBurchard} plays a key role in~\cite{DingXia21, AHP20}, whereas for us, a reduced volume suffices to improve the spin-wave argument.

There are two important challenges in the aforementioned roadmap: on the one hand, we need to show indeed ``many loops'' in the $\zz$-copy indeed leads to a small crossing probability in the $\xx$ and $\yy$ copy; on the other hand, it seems impossible to actually treat defective regions as with 0 local time for the projected XY model.  
In order to address the first challenge, we first prove that the total number of macroscopic loops has a Gaussian decay, and thus, upon observing many macroscopic loops in the $\zz$-copy, the conditional probability of seeing a macroscopic loop in the $\xx$ or $\yy$ copy is small. Then, we show that the XY model arising from the conditioning of the local time field of $ \zz $-copy (which has a good conditional probability of containing many macroscopic loops) typically has no macroscopic loop. Based on this, one can then derive that the probability for the existence of a loop soup cluster crossing through a macroscopic annulus is small. The attempt to prove the Gaussian decay is reminiscent of~\cite{AizenmanPercolationSuperExpDecay}, which proved the Gaussian decay on the number of disjoint crossing clusters for two-dimensional critical Bernoulli percolation. However, it seems that the strong and somewhat mysterious correlation for traces in different spatial locations by the same loop forms a serious obstacle for adapting the argument in~\cite{AizenmanPercolationSuperExpDecay}. In the end, we dealt  with this (see Section~\ref{sec: loop Gaussian decay}) using a spin wave argument using complex rotation, taking advantage of various geometric representations of the loop soup model, including the random path model in~\cite{TaggiLRO}. Provided with this, we then show (see Proposition~\ref{prop: small loop crossing implies small cluster crossing}) a small crossing probability for loop soup clusters arising from a small probability for the existence of macroscopic loops.

The second challenge in improving the spin wave argument is also serious.
In order to gain from the defective region, we have to employ a rather novel extension of the original spin wave argument as in \cite{MW66, Mermin67}, see also~\cite[Chapter~9]{Velenik}.
In particular, based on the formulation of the defective region, it seems inevitable to work with different representations of the XY model (in terms of angles, loop soups, etc).
However, the even more important challenge lies in the necessity to keep tracking the multiplicity of the spin wave map: this was straightforward in the original argument, but seemingly rather challenging in our setting since ``no loop soup cluster crossing'', which is the main feature we wish to use for the defective region, is an event with small but strictly positive probability.
Thus, one cannot simply treat them as if they never happened. See Section~\ref{sec:fractality} for a more elaborative discussion on how to address this difficulty.

\paragraph{Statement:}
We have been working intensively on this project since June 2025, and were completing the draft when, on October 6, 2026, OpenAI released a proof of the result stated here as Theorem~\ref{thm-1.1}. The two approaches are drastically different; see Remark~\ref{rmk: 1.4}. To the best of our knowledge, the mathematical content of our work has not been influenced by OpenAI's proof.

The authors developed the mathematical arguments and wrote the manuscript, with limited, task-specific AI assistance as detailed below.

1. For Lemma~\ref{lem: update--domination}, we checked the degree-three case ourselves and then used AI (GPT-5.2 thinking) to check the degree-four case. The proof presented here was written by the authors and may have benefited from insights suggested by those computations. At that stage, we were not aware of~\cite{Chayes98} or~\cite{MC1998}.

2. We used AI (GPT-5.6 Sol)  to identify a reference in optimization that was important to an earlier approach to proving a Gaussian tail bound for the number of loops. Neither that approach nor the reference is used in the present proof.

3. We used AI (GPT-5.6 Sol) to document some of our blackboard discussions for record-keeping. To the best of our knowledge, we did not incorporate AI-generated mathematical content from these records. After many revisions of the manuscript, however, we cannot reliably identify any remaining textual influence of this assistance. In addition, we used GPT-6 Astra to help producing  both Figures \ref{fig: type of good box} and \ref{fig: local time crossing exploration} after a hand-drawing version.

4. At various stages, we asked AI systems (both GPT-5.6 Sol and GPT-6 Astra)  whether they could solve the problem and how likely a solution was to appear in OpenAI's release. Their responses tended to reinforce our impression that solving the problem would require substantial computational resources. The October 6 release shaped our assessment of AI's capabilities.

5. On the day of submission to arXiv, we used GPT-6 Astra for a preliminary check of the manuscript and made some edits in response. The mathematical edit clarified a measure decomposition in Proposition~\ref{prop: symmetric spin wave}.

\newcommand{\f}{\mathsf{f}}
\newcommand{\w}{\mathsf{w}}
\newcommand{\rd}{\mathrm{d}}
\section{Preliminaries}
\label{sec: preliminaries}

In this section, we give precise definitions for employed terminologies, and collect a few preliminary facts on the spin $O(\NN)$ model.

To begin with, we first explain the convention of constants throughout the paper. Let $ \beta_0>0 $ be a large constant serving as the upper bound for local times (or in the language of Theorem \ref{thm-1.1}, serving as the upper bound for the inverse-temperature). In general, constants chosen to be large are written in the uppercase ($C$, for example), while small constants are in the lowercase. The constants $c_1, c_2, \dots,C_1,C_2,\dots$ have fixed values within each  section, though they may be redefined from one section to the next. In contrast, constants with overlines $\bar{c}_1, \bar{c}_2, \dots, \bar C_1, \bar C_2,\dots$ may take different values in different places, even within the same subsection or proof; that is, their meaning is local and may change from one occurrence to another.

\subsection{Formal definition of spin models}
\label{subsec: notations}

First we give the definition of the spin $O(\NN)$ model, for an integer-valued $\mathsf N$.
Let us first define the  partition functions for the spin models.
We will use plain capitals (see e.g. \eqref{eq: def--spin O(N) partition function} below) for spin partition functions throughout this paper.
(An alternative partition function used in the loop soup representation will be introduced in Section~\ref{subsec: conditional loop soup representation}.)
Fix a finite graph $G=(V,E)$, nonnegative coupling constants $(J_e)_{e\in E}$, and a local time condition $\beta=(\beta_v)_{v\in V}$ with $\beta_v\geq0$ (a special case is when all $ \beta_v=\beta $, as mentioned in the introduction).
Throughout this subsection, $\sigma_v$ is a unit spin and $\phi_v$ is the corresponding field, with conventions
\begin{equation}
    \label{eq: def--spin O(N) spin}
    \sigma_v\in\mathbb S^{\NN-1},\qquad \phi_v=\sqrt{2\beta_v}\,\sigma_v.
\end{equation}
For $A\subset V$, let $\zeta\in (\mathbb S^{\mathsf N-1})^A$ be a spin boundary condition on $A.$
Writing $d\sigma$ for the normalized Haar measure on $\mathbb S^{\NN-1}$, we define the partition function (note that $\delta$ below is a Dirac mass)
\begin{equation}
    \label{eq: def--spin O(N) partition function}
    \begin{aligned}
        Z^\zeta_{G,\beta}
        &= \int_{(\mathbb S^{\NN-1})^V}\exp\big(2\sum_{uv\in E}J_{uv}\sqrt{\beta_u\beta_v}\,\sigma_u\cdot\sigma_v\big)\prod_{v\in A} \delta(\sigma_v-\zeta_v)\prod_{v\in V\setminus A} \rd\sigma_v.
    \end{aligned}
\end{equation}
The measure is then defined as 
\begin{equation}
    \label{eq: def--spin O(N) measure}
    \begin{aligned}
        \mu^\zeta_{G,\beta}(\rd\sigma)
        &= \frac{1}{Z^\zeta_{G,\beta}}\cdot \exp\big(2\sum_{uv\in E}J_{uv}\sqrt{\beta_u\beta_v}\,\sigma_u\cdot\sigma_v\big)\prod_{v\in A} \delta(\sigma_v-\zeta_v)\prod_{v\in V\setminus A} \rd\sigma_v.
    \end{aligned}
\end{equation}
Throughout this paper, we assume that $ \max_{e\in E}J_e\leq1 $ without loss of generality.

\subsection{Loop soup representation}
\label{subsec: conditional loop soup representation}

In this subsection, we introduce the metric graph loop soup (following \cite{Lupu14}) and introduce the spin $O(\sf N)$ model as a loop soup with local time condition. Let $\bar\NN=\lceil\NN\rceil$. Let $\mathsf{c}_i=1$ for $1\le i\le \bar\NN-1$ and let $\mathsf{c}_{\bar \NN}=\NN+1-\bar \NN$.
Throughout the loop/excursion representation, loops and excursions are considered modulo time reversal, i.e., we consider the undirected loops/excursions (we switched to work with the undirected loops/excursions, i.e., modulo time reversal, since this is more convenient for our analysis later). All discrete skeletons are understood under this convention.

\subsubsection{The metric graph loop soup}
Let $D=(V,E)$ be a finite subgraph of $\mathbb Z^2$.
Define coupling constants $J_{uv}>0$ (zero-coupling edges are omitted) on its edges, and choose killing rates $\kappa_v\geq0$. Thus, the matrix
\[
A_{vv}=a_v=\kappa_v+\sum_{w:vw\in E}J_{vw},\qquad A_{uv}=-J_{uv}\ (u\ne v)
\]
is positive definite, where $J_{uv}=0$ off $E$.
We also consider its metric graph $\widetilde{D}$, where we replace each edge $uv$ by a cable (i.e., an interval) of length $1$.
Define $X$ to be the Brownian motion on $\widetilde D$ (see~\cite[Section~2]{Lupu14}), with generator $\frac12\frac{d^2}{ds^2}$ on each cable and killing at rate $\kappa_v$ per unit vertex local time and with skeleton generator $-A$, as follows:
\[
H_t=\sum_{v\in V}L_t^v(X),\qquad \tau_s=\inf\{t:H_t>s\},\qquad Y_s=X_{\tau_s},\qquad \operatorname{gen}(Y)=-A.
\]
Here, vertex local times use the continuous cable-length normalization of~\cite[Section~2]{Lupu14}.
Write $\mathsf m_D$ for the Brownian loop  measure of this killed process, modulo circular time-shift and in the usual $dt/t$ convention (as in~\cite[Section~1]{Lupu14}), and write $\mathsf{M}^\f_D$ for the law of independent Poisson loop soups $\mathfrak{L}=(\mathfrak L^1,\ldots,\mathfrak L^{\bar\NN})$, whereas $\mathfrak L^{\ii}$ is of intensity $ \mathsf{c}_\ii\mathsf m_D $ (note that in the undirected case $ \mathsf{c}_{\ii}=1 $ corresponds to the directed case where $ \alpha_{\ii}=1/2 $). The superscript $\f$ denotes the free boundary condition, and thus $\mathfrak L$ does not contain excursions. 
All loops and connectivity below are on $\widetilde D$, including loops contained entirely inside a cable or loops that visit only one lattice vertex.
The vector of the local time fields is $ \mathbf{L}(\mathfrak{L})=(L^1(\mathfrak{L}),\ldots,L^{\bar\NN}(\mathfrak{L})) $, where $L^{\ii}$ is the local time field of the $\ii$-th copy, given by
\begin{equation}\label{eq: local time definition}
    L_v^{\mathsf{i}}(\mathfrak{L})=\sum_{\gamma\in\mathfrak L^{\mathsf{i}}}L_v(\gamma),\qquad \forall v\in\widetilde D.
\end{equation}
Here $L_v(\gamma)$ is the local time accumulated at $v$ by the loop $\gamma$. In addition, let $L^*_v(\mathfrak{L})=\sum_{i=1}^{\bar\NN}L^i_v(\mathfrak{L})$.
For notation convenience, for two sets $A, B$, we use $A\stackrel{\mathsf i}{\longleftrightarrow} B$ to denote the local time connecting event, i.e., $A$ is connected to $B$ by a metric graph path where all local times are strictly positive in the $\mathsf i$-copy. When we put $\stackrel{\mathsf{xy}}{\longleftrightarrow}$ it means connection in either the $\mathsf{x}$- or the $\mathsf{y}$-copy.

Let $F$ be a nonnegative measurable weight on the soups.
Now we define the weighted partition function with respect to the local time condition $\beta$ as follows:
\begin{equation}
    \label{eq:weighted-density}
    \mathcal{Z}^{\f}_{D,\beta}(F) :=\int F(\mathfrak L) \prod_{v\in D}\delta(L_v^{*}(\mathfrak{L})-\beta_v)\cdot \,\mathsf M^\f_D(\rd\mathfrak L).
\end{equation}  
For notation convenience, we write $ \mathcal{Z}^{\f}_{D,\beta}= \mathcal{Z}^{\f}_{D,\beta}(1)$ as the partition function on $D$ with free boundary condition and local time condition $\beta$.
We then define the probability density of a loop configuration $\mathfrak{L}$ by
\begin{equation}
    \label{eq: loop free law}
    \nu^{\f}_{D,\beta}(\rd\mathfrak L) :=\frac{1}{\mathcal{Z}^{\f}_D(\boldsymbol{\beta})} \prod_{v\in D}\delta(L_v^{*}(\mathfrak{L})-\beta_v) \,\mathsf M^{\f}_D(\rd\mathfrak L).
\end{equation} We remark that although the choice of killing rate $\kappa_v$ changes  $\mathsf M_D$ by a factor $\exp(-\sum_{v\in V}\kappa_v L^*_v)$, it does not influence the measure $\nu^{\f}_{D,\beta}$ since the factor cancels in the right-hand side of \eqref{eq: loop free law}.

\subsubsection{Two copy measure with local time boundary conditions}\label{sec: local time boundary condition}
In this part, we focus on the case $\NN=2$ and refer to the two copies as $\xx$ and $\yy$. 
We now define partition functions under local time boundary conditions. Consider a subset $A\subset D\subset V $.
We define the (copywise) local time boundary condition {$ \xi=(\xi^{\xx},\xi^{\yy},\mathcal{P}^{\xx}_A,\mathcal{P}^{\yy}_A) $} on $A$ (with the default choice of $A = \partial D$, in which case the partition $ \mathcal{P}^{\ii}_A $ will serve as the connectivity of vertices in $A$ induced by the $\ii$-copy local time cluster restricted to $ D^c $, as shown in Lemma~\ref{lem: block-parity}) where $\xi^\mathsf{x}, \xi^{\yy} \in [0, \infty)^A$ satisfying that $\xi^{\xx}_v + \xi^{\yy}_v = \beta_v$ for all $v\in A$, and $ \mathcal{P}^{\ii}=\{B^{\ii}_1,\dots,B^{\ii}_{m_{\ii}}\} $ (each $B^{\ii}_j$ may be referred to as a block) is a partition of $A$ for $ \ii=\xx,\yy $.
In contrast to the free boundary measure, we now allow excursions from $B^{\ii}_k$ to $B^{\ii}_j$. For $k\neq j$, let $N_{k,j}^{\ii}$ denote the the total number of excursions in $\ii$-copy between $B^{\ii}_k$ and $B^{\ii}_j$ for $\ii\in \{\xx, \yy\}$ (recall that excursions are undirected).
Let $\1_{\xi}(\mathfrak{L})$ indicate the event such that \begin{enumerate}
    \item For any $\ii\in\{\xx,\yy\}$ and any $v\in A$, we have $L_v^\ii=\xi^\ii_v$.
    \item For any $\ii\in\{\xx,\yy\}$ and $1\le k\le m_{\ii}$, $\sum_{j\neq k}N^\ii_{k,j}$ is even.
\end{enumerate}  Let $F$ be a nonnegative measurable weight on the soups. For $u,v\in A$, let $\mathsf{m}^{A}_D(v,u)$ denote the mass under the (undirected) excursion measure in $A$ with endpoints $u,v\in A$.
In addition, we define $\mathsf{m}^{\xi,\ii}_D(v,u)=\mathsf{m}^{A}_D(v,u)\cdot 2\sqrt{\xi^\ii_u\cdot \xi^\ii_v}$ for $\ii\in\{\xx,\yy\}.$
Then we can define the measure $\mathsf{M}_D^\xi$ to be the law of an independent Poisson loop soup with intensity $\frac{\mathsf{m}_D}{2}$ composed with an independent Poisson excursion with intensity $\frac{\mathsf{m}_D^{\xi,\ii}}{2}$ for each copy $\ii\in \{\xx,\yy\}$.
Then the weighted partition function with local time boundary condition $\xi$ is defined as
\begin{equation}
    \label{eq:weighted-density with local time bc}
    \mathcal{Z}^{\xi}_{D,\beta}(F) :=\int F(\mathfrak L) \prod_{v\in D\setminus A}\delta(L_v^{*}-\beta_v)\cdot\1_{\xi}(\mathfrak L)\cdot \,\mathsf M^\xi_D(\rd\mathfrak L).
\end{equation} 
For notation convenience, we write $ \mathcal{Z}^{\xi}_{D,\beta}= \mathcal{Z}^{\xi}_{D,\beta}(1)$ as the partition function on $D$ with the local time boundary condition $\xi$ and local time condition $\beta$.
We then define the corresponding probability density of a loop configuration $\mathfrak{L}$ as
\begin{equation}
    \label{eq: loop law with local time bc}
    \mu^{\xi}_{D,\beta}(\rd\mathfrak L) :=\frac{1}{\mathcal{Z}^{\xi}_{D,\beta} } \prod_{v\in D\setminus A}\delta(L_v^{*}-\beta_v){\cdot\1_{\xi}(\mathfrak L)\cdot} \,\mathsf M^{\xi}_D(\rd\mathfrak L).
\end{equation} 

For notation convenience, when we  consider the local time boundary condition, we may omit the word ``local time" if it is clear in the context. When we  consider the spin boundary condition, we decide to keep the word ``spin" to emphasize. Since we will frequently consider both the XY loop measure (i.e., with intensity $1$) and the $O(\NN)$ loop measure (i.e., with intensity $\NN/2$), we will use the notation $\mu$ to denote the XY loop measure and use the notation $\nu$ to denote the $O(\NN)$ loop measure. We sometimes need to consider the $O(\NN)$ loop measure with some local time boundary condition $\xi=(\xi^\xx,\xi^\yy,\xi^\zz,\mathcal{P}_A^\xx,\mathcal{P}_A^\yy,\mathcal{P}_A^\zz)$  where $\xi^\mathsf{x}, \xi^{\yy} \in [0, \infty)^A$ satisfying that $\xi^{\xx}_v + \xi^{\yy}_v = \beta_v$ for all $v\in A$, and $\xi_v^\zz=0$ for $v\in A$.  In this case, we call the $\xi$ an $\xx\yy$ local time boundary condition and denote the measure on $D$ with the boundary condition $\xi$ as $\nu^\xi_{D,\beta}$. 

\subsubsection{Spatial restriction}
We record the following restriction property.
\begin{lemma}
    \label{lem:cpf-spatial-restriction}
    Fix $D, H\subset V$. Fix a boundary condition $\xi=(\xi^\xx,\xi^\yy,\mathcal{P}^\xx,\mathcal{P}^\yy)$ on $\partial D$. For any configuration $\mathfrak{L}=(\mathfrak{L}^\xx,\mathfrak{L}^\yy)$, it can be decomposed into    $\mathfrak{L}_H=(\mathfrak{L}_H^\xx,\mathfrak{L}_H^\yy)$ and $\mathfrak{L}_{D\setminus H}=(\mathfrak{L}_{D\setminus H}^\xx,\mathfrak{L}_{D\setminus H}^\yy)$ such that $\mathfrak L^{\mathsf  i}_H$ consists of all loops (no boundary excursions here) in $\mathfrak L^\ii$ intersecting $H$ for $\ii\in \{\xx, \yy\}$ (and naturally $\mathfrak L^\ii = \mathfrak L^{\ii}_H \cup \mathfrak L^{\ii}_{D\setminus H}$). (We note that $\mathfrak L_{D\setminus H}$ may contain an excursion {intersecting} $H$). We can also decompose $\mathsf{M}^\xi_D$ into the product of two measure $\bar{\mathsf{M}}^\xi_{H}$ and $\mathsf{M}^\xi_{D\setminus H}$ such that 
    $$\mathsf{M}^\xi_D(\mathfrak{L})= \bar{\mathsf{M}}^\xi_{H}(\mathfrak{L}_H)\cdot \mathsf{M}^\xi_{D\setminus H}(\mathfrak{L}_{D\setminus H}).
    $$ For any $\mathfrak{L}_H$, we define the remaining local time condition as $\beta'_v=\beta_v-L^*_v(\mathfrak{L}_H)$. Since no loop in $\mathcal L_H$ intersects $\partial D$, we have
    \begin{equation}
        \label{eq:cpf-domain-partition}
        \mathcal Z^{\xi}_{D,\beta} =\int \bar{\mathsf M}^\xi_H(\rd\mathfrak{L})\, \mathcal Z^{\xi}_{D,\beta'}.
    \end{equation}
    Consequently,
    given $\mathfrak{L}_H$, the remaining loops have the loop soup law in $D$ with the local time condition $\beta'$ and boundary condition $\xi$.
\end{lemma}
\begin{proof}
By \eqref{eq:weighted-density with local time bc}, we calculate that 
\begin{align*}
    \mathcal{Z}^{\xi}_{D,\beta} &=\int_{\mathfrak{L}}\prod_{v\in V}\delta(L_v^{*}(\mathfrak{L})-\beta_v)\cdot\1_{\xi}(\mathfrak L)\cdot \,\mathsf M^\xi_D(\rd\mathfrak L)\\&=\int_{\mathfrak{L}_H,\mathfrak{L}_{D\setminus H}}\prod_{v\in V}\delta(L_v^{*}(\mathfrak{L}_H)+L^*_v(\mathfrak{L}_{D\setminus H})-\beta_v)\cdot\1_{\xi}(\mathfrak L)\cdot \bar{\mathsf{M}}^\xi_{H}(\rd\mathfrak{L}_H)\cdot \mathsf{M}^\xi_{D\setminus H}(\rd\mathfrak{L}_{D\setminus H})\\&=\int_{\mathfrak{L}_H}\bar{\mathsf{M}}^\xi_{H}(\rd\mathfrak{L}_H)\cdot\int_{\mathfrak{L}_{D\setminus H}}\prod_{v\in V}\delta(L^*_v(\mathfrak{L}_{D\setminus H})-\beta'_v)\cdot\1_{\xi}(\mathfrak L)\cdot  \mathsf{M}^\xi_{D\setminus H}(\rd\mathfrak{L}_{D\setminus H})\\&=\int \bar{\mathsf M}^\xi_H(\rd\mathfrak{L})\, \mathcal Z^{\xi}_{D,\beta'}.
\end{align*} 
   Thus, we finish the proof of Lemma \ref{lem:cpf-spatial-restriction}.
\end{proof}
\begin{remark}\label{rmk: spatial restriction for excursion}
    In the special case $\partial D\subset H$, we can decompose the loop configuration $\mathfrak{L}=(\mathfrak{L}^\xx,\mathfrak{L}^\yy)$ into    $\mathfrak{L}_H=(\mathfrak{L}_H^\xx,\mathfrak{L}_H^\yy)$ and $\mathfrak{L}_{D\setminus H}=(\mathfrak{L}_{D\setminus H}^\xx,\mathfrak{L}_{D\setminus H}^\yy)$ such that $\mathfrak L^{\mathsf  i}_H$ consists of all loops and boundary excursions  in $\mathfrak L^\ii$ intersecting $H$ for $\ii\in \{\xx, \yy\}$ (and naturally $\mathfrak L^\ii = \mathfrak L^{\ii}_H \cup \mathfrak L^{\ii}_{D\setminus H}$). Then we have $\mathfrak L_{D\setminus H}$ does not contain any excursion intersecting $H$. In this case, \eqref{eq:cpf-domain-partition} also holds with the same proof.
\end{remark}
\subsubsection{Domain Markov Property}
We record the domain Markov property for the critical loop soup. Fix a region $D$ and let $\Phi$ denote the measure $\mathsf{M}_D$ with $\NN=1$. Then $\Phi_D$ is the critical loop soup measure on $D$.  
\begin{lemma}\label{lem: DMP for loop on critical loop soup}
Fix $H\subset D$ and fix a local time configuration $\omega$ on $H$. Let $\xi$ be the boundary condition on $\partial H$ induced by $\omega$. That is, $\xi_v = \omega_v$ for $v\in \partial H$, and we partition $\partial H$ into blocks where two vertices are in the same block if they are connected by a path in $H$ with all $\omega$-values positive. (Note that this echoes the local time condition.) Then conditioned on $\omega$, the measure for the loops and excursions in $D\setminus H$ is $\mu^\xi_{D\setminus H,\beta}$.
\end{lemma}
\begin{proof}
    The proof of this lemma follows from combining \cite[Proposition 7]{Wer16} and \cite[Theorem 1.3]{Aid20}. To this end, one needs to repeat some computations as in \cite[Theorem 4.5]{Aid20}. Indeed, \cite[Equation (4.6)]{Aid20} gives the conditional probability for the excursions, and following the calculation at the end of the proof of \cite[Theorem 4.5]{Aid20} yields the desired result.
\end{proof}
Applying Lemma \ref{lem: DMP for loop on critical loop soup} to each copy, we arrive at the following corollary.
\begin{corollary}
    \label{cor: DMP for loop}
    Let $ \NN\geq2 $ be an integer. 
    Consider two regions $ B_1 \subset B_2$ and a boundary condition $\eta$ on $\partial B_2$. Let $\Xi$ be the collection of boundary conditions on $\partial B_1$.
    We write $ \mu^\eta_{B_2,\beta}(\xi) $ for the distribution of the boundary condition $\xi$ induced by the local time configuration on $B_2\setminus B_1$ {and the boundary condition $\eta$}.
    Then, for all local time conditions $\beta$ and every event $\mathcal A$ measurable with loops and excursions on $B_1$ (by excursion on $B_1$ we mean segments of loop or excursion from larger boxes cut at $\partial B_1$), we have
    \begin{equation}
        \int_{\Xi}\mu^\xi_{\fB_1,\beta}(\mathcal{A}) \mu^\eta_{\fB_2, \beta}(\rd\xi) = \mu^\eta_{\fB_2,\beta}(\mathcal{A}).\nonumber
    \end{equation}
\end{corollary}
\begin{proof}
Let $\eta=(\eta^1,\cdots,\eta^\NN,\mathcal{P}^1,\cdots,\mathcal{P}^\NN)$.
    Let $\omega=(\omega^1,\cdots,\omega^\NN)$ be a local time configuration on $B_2\setminus B_1$. Note that $\mu^\eta_{\fB_2,\beta}$ can be interpreted as the product measure of $\mu^{\eta^\ii,\mathcal{P}^\ii}_{B_2}$ $(\ii\in\{1,\cdots,\NN\})$, subject to the local time condition $\sum_{\ii\in\{1,\cdots,\NN\}}L^\ii_v=\beta_v$ for $v\in B_2$. Let $\xi=(\xi^1,\cdots,\xi^\NN,\mathcal{P}^{1'},\cdots,\mathcal{P}^{\NN'})$ be the local time boundary condition induced by $\omega$ and $\eta$. Then by Lemma \ref{lem: DMP for loop on critical loop soup}, we get that $\mu^{\eta^\ii,\mathcal{P}^\ii}_{B_2}(\cdot\mid \omega^\ii)=\mu^{\xi^\ii,\mathcal{P}^{\ii'}}_{B_2}$. Therefore, $\mu^\eta_{\fB_2, \beta}(\cdot\mid\omega)$ is the product measure of $\mu^{\xi^\ii,\mathcal{P}^{\ii'}}_{B_2}$ $(\ii\in\{1,\cdots,\NN\})$, subject to the same local time condition $\sum_{\ii\in\{1,\cdots,\NN\}}L^\ii_v=\beta_v$ for $v\in B_1$. In conclusion, we get that $\mu^\eta_{\fB_2, \beta}(\cdot\mid\omega)=\mu^\xi_{\fB_1,\beta}$.
\end{proof}
\subsubsection{Copy reassignment}
We record the following reassignment property for the loops with intensity $ \NN\in\mathbb{R}^+ $.

Let $\mathfrak L^*$ be a loop soup configuration in $D$ without label.
We independently assign to each loop $\ell\in\mathfrak L^*$ a uniform label $U_\ell\in(0,\NN]$, and define its copy by $\lceil U_\ell\rceil$. Recall that $\mathsf{c}_i=1$ for $1\le i\le \bar\NN-1$ and let $\mathsf{c}_{\bar \NN}=1+\NN-\bar \NN$. We record the following lemma, which follows immediately from the Poisson splitting property.

\begin{lemma}
    \label{lem:yyy}
    Let $\ell_1,\ldots,\ell_m$ be distinct loops selected in a way that is measurable to using only the unlabelled loop collection.
For each assignment $\mathbf a={(\mathsf{a}_1,\ldots,\mathsf{a}_m)}\in\{\mathsf{x},\mathsf{y},\mathsf{z}\}^m$, we have
    \begin{equation*}
        \nu_{D,\beta}^{(\NN)}(\ell_i\in \mathfrak{L}^{{\mathsf{a}_i}},~ 1\le i\le m\mid \{\ell_1, \ldots, \ell_m\})
        =\prod_{i=1}^{m}\frac{\mathsf{c}_{{\mathsf{a}_i}}}{\NN}.
    \end{equation*}
\end{lemma}

\subsection{Relation between spin model and loop model}
Now we discuss the relation between the spin model defined in \eqref{eq: def--spin O(N) partition function} and the loop model defined in \eqref{eq: loop free law} and \eqref{eq: loop law with local time bc}. Since the spin model can only be defined when $\NN\in\mathbb Z$, we focus on the case that $\NN=2$.

\paragraph{Free boundary condition case}
By the isomorphism theorem (see, e.g.~\cite[Theorem~1]{Lupu14}), we have the following coupling between loop soup and GFF.
For every loop configuration $\mathfrak{L}=(\mathfrak{L}^{\xx}, \mathfrak{L}^{\yy})$, we assign independent signs to each of the local time clusters in $\mathfrak{L}^i$ uniformly at random. Then we get a spin configuration $\sigma$ (the absolute value is determined via \eqref{eq: 2.9} below). For any spin configuration $\sigma$, we define the local time configuration $L=(L^\xx,L^\yy)$ on vertices by $L^\ii_v=\beta_v(\sigma_v^\ii)^2$ (we can then also interpolate each metric edge with independent Brownian bridges whose absolute values will then give local times thereon; see Lemma \ref{lem: FK-cluster-loop} and discussions preceding it). We say a loop configuration $\mathfrak{L}$ is compatible with the spin configuration $\sigma$ (denoted as $\mathfrak{L}\sim \sigma$) if for any $1\le i\le \NN$  and $v\in D$ \begin{equation}
    L_v(\mathfrak{L}^\ii)=\beta_v(\sigma_v^\ii)^2.\label{eq: 2.9}
\end{equation}
By~\cite[Theorem~1 and Proposition~2.1]{Lupu14}, there exists a constant $\mathfrak{Z}^\f_{D,\beta}$ depending only on $D$ and $\beta$ such that  
\begin{equation}
    \int_{\mathfrak{L}\sim \sigma}\mathsf{M}^\f(\rd\mathfrak{L})=\frac{1}{\mathfrak{Z}^\f_{D,\beta}}\cdot \exp\big({\,2\sum_{uv\in E}J_{uv}\sqrt{\beta_u\beta_v}\, \sigma_u\cdot\sigma_v}\Big) \prod_v\rd\sigma_v.
\end{equation} 
\paragraph{General boundary condition case}
Let $\zeta$ be a spin boundary condition on $ B\subset D $ as
\begin{equation}
    \label{eq: def spin b.c.}
    \zeta_v=(\zeta^{\xx}_v,\zeta^{\yy}_v)\in\mathbb{S}^{1},\qquad\forall v\in B.
\end{equation}
Let $\xi_0$ be the local time boundary condition with boundary partition $\mathcal{P}=(\mathcal{P}^{\xx}, \mathcal{P}^{\yy})$ such that for any $v\in B$ and $ \ii\in \{\xx,\yy\}$, $\xi^\ii_0(v)={\beta_v}(\zeta_v^\ii)^2$ and that vertices of opposite signs are in different partitioning sets copywise (we refer to such $\mathcal{P}$ as a valid boundary partition).
Identities in below follow from the isomorphism theorem and the independent cluster sign construction in~\cite[Proposition~2.1 and Lemma~3.2]{Lupu14}, and a Bayes calculation. First
we can couple the spin configuration $\sigma$ under the boundary condition $\zeta$ with the local time configuration under the local time boundary condition $\xi_0$. Let $\omega$ be a local time configuration. For $\ii\in \{\xx, \yy\}$, {we write $\kappa^\ii(\omega,\mathcal{P})$ to be the number of boundary local time clusters in the $\ii$-copy under the boundary partition $\mathcal{P}$}, and let $V_{\ii}^+, V_\ii^-$ denote the vertices on the boundary $B$ with plus and minus $\ii$-copy values.  On the one hand, for any local time configuration $\omega$ under $\mu^{\zeta}_{\fD,\gamma}$, it must satisfy $V_{\ii}^+\ncon[0.5]{\ii} V_\ii^-$ for $\ii\in\{\xx,\yy\}$; we refer such $\omega$ as a valid local time configuration under $\mu^{\zeta}_{\fD,\gamma}$. On the other hand, for any valid local time configuration $\omega$ under $\mu^{\zeta}_{\fD,\gamma}$, to get the desired spin on $B$ as in $\zeta$, the probability is $ 2^{-\kappa^\xx(\omega)-\kappa^\yy(\omega)}$. Thus, we get that for any local time configuration $\omega$ \begin{equation}\label{eq: relation between spin boundary condition and local time boundary condition}
    \mu^{\zeta}_{\fD,\gamma}(\omega)\propto \mu^{\xi_0}_{\fD,\gamma}(\omega)\cdot 2^{-\kappa^\xx(\omega{,\mathcal{P}})-\kappa^\yy(\omega{,\mathcal{P}})} \cdot \1_{V_{\xx}^+\ncon[0.5]{\xx} V_\xx^-}(\omega)\cdot \1_{V_{\yy}^+\ncon[0.5]{\yy} V_\yy^-}(\omega).
\end{equation}

For a fixed spin boundary condition $ \zeta $, the law of the local time field does not depend on the (valid) boundary partition $ \mathcal{P} $. This follows from the Edwards-Sokal coupling (see~\cite[Fact~3.4 and Corollary~3.5]{LubetzkySly12}), as we will elaborate in what follows. Fix $ \ii\in\{\xx,\yy\} $. {Note that given $\{L_v^\ii\}_{v\in D}$, the local time percolation through edges in $\ii$-copy follows from the law of FK-Ising model with coupling constant $\sqrt{2L_u^\ii L_v^\ii}$.}
Let $ \mathsf{e}^{\ii} $ denote the FK edge configuration in copy $\ii$. Write $ k^{\ii}(\mathsf{e}, \mathcal{P}) $ for the number of FK clusters after wiring blocks in $ \mathcal{P}^{\ii} $, and let $ \kappa^{\ii}(\mathsf{e}, \mathcal{P}) $ count those intersecting the boundary. Then in the joint density of vertex local times and FK edges, fixing the boundary signs changes the cluster factor to $ 2^{k^{\ii}_0} := 2^{ k^{\ii} }\cdot 2^{-\kappa^{\ii}} $, where $ k^{\ii}_0 $ is the number of FK clusters disjoint from the boundary (before applying the prescribed wiring). Since boundary partitioning only joins clusters intersecting the boundary, the preceding factor $ 2^{ k^{\ii}-\kappa^{\ii} } $ is independent of the valid partitions $ \mathcal{P}^{\ii} $, thus so does~\eqref{eq: relation between spin boundary condition and local time boundary condition}.
However, only knowing the local time $\xi^\ii_0(v)$ on $v\in B$ is not enough to determine the law of the loop and excursion configuration. Thus, we introduce the spin boundary condition with partition $(\zeta,\mathcal{P})$.  The corresponding local time boundary condition is thus $\xi=(\xi^\xx_0,\xi^\yy_0,\mathcal{P}^\xx,\mathcal{P}^\yy)$. As before, we  write $\kappa^\ii(\omega,\mathcal{P})$ to be the number of boundary local time clusters in $\ii$-copy under the boundary partition $\mathcal{P}$. Similar to \eqref{eq: relation between spin boundary condition and local time boundary condition}, we have \begin{equation}\label{eq: loop relation between spin boundary condition and local time boundary condition}
    \mu^{\zeta,\mathcal{P}}_{\fD,\gamma}(\rd\mathfrak{L})\propto \mu^{\xi}_{\fD,\gamma}(\rd\mathfrak{L})\cdot 2^{-\kappa^\xx(\mathfrak{L},\mathcal{P})-\kappa^\yy(\mathfrak{L}),\mathcal{P}} \cdot \1_{V_{\xx}^+\ncon[0.5]{\xx} V_\xx^-}(\mathfrak{L})\cdot \1_{V_{\yy}^+\ncon[0.5]{\yy} V_\yy^-}(\mathfrak{L}).
\end{equation} 
In the special case that $\zeta^\ii$ has the same sign over $B$ and $\mathcal{P}_\w^\ii=\{ B\}$ for all $ \ii\in\{\xx,\yy\} $, the last three terms in \eqref{eq: loop relation between spin boundary condition and local time boundary condition} become constant, and thus we get that \begin{equation}\label{eq: special case for spin bc to local time bc}
    \mu^{\zeta,\mathcal{P}_\w}_{D,\beta}(\rd\mathfrak{L})\propto \mu^{\xi}_{D,\beta}(\rd\mathfrak{L}).
\end{equation}
Therefore, without further announcement, we will use the notation $\mu^{\zeta}_{D,\beta}$ to denote both the measure for the spin configuration under the spin boundary condition $\zeta$ and also the measure for the loop/excursion configuration  under the local time boundary condition $\xi=({\beta_v}(\zeta_v^\xx)^2,{\beta_v}(\zeta_v^\yy)^2,\mathcal{P}_\w^\xx,\mathcal{P}_\w^\yy)$. We use the notation  $\mathcal{Z}^{\zeta}_{D,\beta}$ to denote the local time partition function as follows: \begin{equation}\label{eq: spin boundary local time boundary condition}
 \mathcal{Z}^{\zeta}_{D,\beta} :=\int  \prod_{v\in D\setminus B}\delta(L_v^{*}-\beta_v)\cdot\1_{\xi}(\mathfrak L)\cdot \,\mathsf M^\xi_D(\rd\mathfrak L)\cdot \1_{V_{\xx}^+\ncon[0.5]{\xx} V_\xx^-}(\mathfrak{L})\cdot \1_{V_{\yy}^+\ncon[0.5]{\yy} V_\yy^-}(\mathfrak{L}).    
\end{equation}
\begin{lemma}
\label{lem: spin local time factorization}
Fix $B\subset  D$. Recall from \eqref{eq: def--spin O(N) partition function} that $Z^{\zeta}_{D,\beta}$ is the spin partition function. 
There exist strictly positive functions $F_D$ of the spin boundary condition $\zeta$ and $G_D$ of the interior local time condition such that for any spin boundary condition $\zeta$ on $B$, we have
\begin{equation}
\label{eq: spin local time factorization}
    \frac{\mathcal{Z}^{\zeta}_{D,\beta}}{Z^{\zeta}_{D,\beta}}
    =F_D(\zeta)\,G_D(\beta|_{D\setminus B}).
\end{equation}
\end{lemma}
\begin{proof}
For any spin boundary condition $\zeta$, let $\mathbb P^\zeta$ denote the law of two independent metric graph GFFs with boundary values $\sqrt{2\beta}\zeta$: (note that here we do not require $|\sigma_v|=1$)
\begin{equation}
\label{eq: boundary gaussian density}
    \mathbb P^\zeta(d\sigma)
    =c_D(\zeta)\exp\!\left[
       -\sum_{v\in D\setminus B}a_v\beta_v
       +2\sum_{u\sim v}J_{uv}\sqrt{\beta_u\beta_v}\sigma_u\cdot\sigma_v
    \right]
    \prod_{v\in D\setminus B}\rd \sigma_v\prod_{v\in B}\delta(\sigma_v-\zeta_v).
\end{equation}
For any loop configuration $\mathfrak{L}$, we say it is valid for a spin configuration $\sigma$ (denoted as $\mathfrak{L}\sim \sigma$) if there exists a sign assignment to each local time cluster of $\mathfrak{L}$ such that the resulting spin configuration is $\sigma$.
Let $K(\mathfrak{L},\mathcal{P}_\w)$ denote the number of clusters (we only consider the clusters that contain at least one vertex) in $\mathfrak{L}$ with respect to the partition $\mathcal{P}_\w=(\mathcal{P}_\w^\xx,\mathcal{P}_\w^\yy)$ (recall that $\mathcal{P}_\w^\ii=\{ B\}$). In the case that both $\zeta^\xx$ and $\zeta^\yy$ have the same sign on $B$, by \cite[Propositions~2.3, 2.4]{ALS20} we get that  \begin{equation}
\int_{\mathfrak{L}\sim \sigma}\1_{\xi}(\mathfrak L)\cdot\mathsf{M}_D^\xi(\rd\mathfrak{L})\cdot 2^{-K(\mathfrak{L},\mathcal{P}_\w)}\propto\mathbb P^\zeta(\rd\sigma).\label{eq: isomorphism with bc}
\end{equation}
For a general spin boundary condition $\zeta$, 
by the same reason as in the paragraph before \eqref{eq: relation between spin boundary condition and local time boundary condition}, we can extend \eqref{eq: isomorphism with bc} to $\zeta$ and get that \begin{equation}
\int_{\mathfrak{L}\sim \sigma} \1_{\xi}(\mathfrak L)\cdot\mathsf{M}_D^\xi(\rd\mathfrak{L})\cdot 2^{-K(\mathfrak{L},\mathcal{P}_\w)}    =\bC_1(\zeta)\mathbb P^\zeta(\rd\sigma).\label{eq: isomorphism with general bc}
\end{equation} Let $\mathcal{E}=\{V_{\xx}^+\ncon[0.5]{\xx} V_\xx^-\}\cap \{V_{\yy}^+\ncon[0.5]{\yy} V_\yy^-\}$.
Since each $\mathfrak{L}$ is valid for $2^{K(\mathfrak{L})}$ spin configurations, we get from \eqref{eq:weighted-density with local time bc} that
\begin{align}
    \label{eq: gaussian radial density}
    \mathcal{Z}^\zeta_{D,\beta}
    &=\int_{\sigma\in(\mathbb{S}^1)^{D\setminus B}}\int_{\mathfrak{L}\sim \sigma} \prod_{v\in I}\delta(L_v^{*}-\beta_v)\cdot\1_{\xi}(\mathfrak L)\cdot \,\mathsf M^\xi_D(\rd\mathfrak L)\cdot 2^{-K(\mathfrak{L},\mathcal{P}_\w)}\rd\sigma
\end{align}
Since for any $\mathfrak{L}\sim \sigma$, it satisfies $L_v^\ii(\mathfrak{L})=\beta_v(\sigma^\ii_v)^2$ for $\ii\in\{\xx,\yy\}$ and thus $L_v^*=\beta_v$. Combining~\eqref{eq: isomorphism with general bc} with~\eqref{eq:weighted-density}, we get that
\begin{equation}
    \label{eq: gaussian radial density 2}
    \begin{aligned}
        \mathcal{Z}^\zeta_{D,\beta}
        &\stackrel{\eqref{eq: gaussian radial density}}{=} \int_{\sigma\in(\mathbb{S}^1)^{D\setminus B}}\rd\sigma\int_{\mathfrak{L}\sim \sigma} \1_{\xi}(\mathfrak L)\cdot \,\mathsf M^\xi_D(\rd\mathfrak L)\cdot 2^{-K(\mathfrak{L},\mathcal{P}_\w)}\nonumber
        \\&\ \ =\int_{\sigma\in(\mathbb{S}^1)^{D\setminus B}}c_1(\zeta)\cdot \mathbb P^\zeta(\rd\sigma)
        = c_1(\zeta)\cdot c_D(\zeta)\exp\Big(-\sum_{v\in D\setminus B}a_v\beta_v\Big) Z_{D,\beta}^\zeta,
    \end{aligned}
\end{equation}
proving Lemma~\ref{lem: spin local time factorization} with $F_D=c_1(\zeta)\cdot c_D(\zeta)$ and $G_D=\exp\Big(-\sum_{v\in D\setminus B}a_v\beta_v\Big)$.
\end{proof}

For different  partitions $\mathcal{P}_1$ and $\mathcal{P}_2$, the measures $\mu^{\zeta,\mathcal{P}_1}$ and $\mu^{\zeta,\mathcal{P}_2}$ are different. However, when considering the loops disjoint from the boundary, different partitions will give the same law given some disconnection event.
\begin{lemma}\label{lem: spin boundary condition plus different partitions}
    Let $D_1\subset D_2$ be two regions and let $\zeta$ be a spin boundary condition on $\partial D_2$. Let $\mathcal{P}_1$ and $\mathcal{P}_2$ be two partitions and let $\mathcal{A} $ denote the event that the local time cluster does not connect $\partial D_2$ to $D_1$ in either copy. Let $\mathfrak{L}_1$ be the collection of loops intersecting $D_1$. Then conditioned on $\mathcal{A}$, the law of $\mathfrak{L}_1$ under $\mu^{\zeta,\mathcal{P}_1}_{D_2,\beta}$ is the same to the law under $\mu^{\zeta,\mathcal{P}_2}_{D_2,\beta}$. 
\end{lemma}
\begin{proof}
    Conditioned on $\mathcal{A}$, we can sample the local time configuration $\omega$ under the law $\mu^{\zeta}_{D_2,\beta}$. Note that the law of local time configurations does not depend on the boundary partitions. Fix a local time configuration $\omega$. Let $C_1,\dots, C_m$ be the local time clusters intersecting $D_1$. Since $C_i\cap \partial D_2=\emptyset$ {on $A$}, we can sample the loop configuration $\mathfrak{L}_1$ on the clusters $C_1,\dots, C_m$, and this does not depend on the boundary partitions.
\end{proof}

\paragraph{From spin configuration to local time configuration}
By Lupu's isomorphism theorem~\cite{Lupu14}, given the spin configuration $\sigma$, to recover the $\ii$-copy local times on an edge $e=\{u,v\}$, we can sample an independent Brownian bridge with endpoint values $\sigma^\ii_u$ and $\sigma^\ii_v$. Consequently, we define the edge mark to be an independent random variable that uniquely determines the behavior of the Brownian bridge on the  edge $e=\{u,v\}$ once given $\sigma^\ii_u$ and $\sigma^\ii_v$. To summarize, we have the following lemma.
\begin{lemma}
    \label{lem: FK-cluster-loop}
    Let $ \mathfrak{e}=(\mathfrak{e}^{\ii}_e)_{e\in E, \ii=\xx,\yy} $ be the independent edge marks. There exists a deterministic bijection function $\Phi = (\Phi(\mathfrak{e}^i_e, \sigma^{\ii}_{e^+}, \sigma^\ii_{e^-}))_{e\in E, \ii\in\{\xx,\yy\}}$  (where $e^\pm$ are the two endpoints of the edge $e$) such that for any $\mathfrak{e}$ and $\sigma$, $\Phi(\mathfrak{e},\sigma)$ is a local time configuration on the metric graph (the explicit form of this function is not relevant to us).
\end{lemma}

\subsection{The random current representation and XY spin correlations}
\label{subsec:random path model}
In this subsection, we study the XY model via its random current representation, and a loop model (also known as the random path model, see e.g.~\cite{TaggiLRO}) induced by the random current representation.
To connect this to the loop soup model, we consider the skeleton of a metric graph loop $\ell$: it is a sequence of neighboring lattice vertices $ v_0,\dots, v_m=v_0 $ such that $ v_{i+1} $ is the next lattice point which is different from $ v_i $ and is visited by $ \ell $ after $ v_i $ for $ 0\leq i\leq m-1 $. Thus,  a skeleton is a random walk loop. We will show in Lemma~\ref{lem: rpm equivalence} that the metric graph loops have the same skeletons as in the random path model.
Lemma~\ref{lem: cluster-loop comparison, all-x local time b.c.} then gives an upper bound of the spin correlation in terms of the loop connecting probability, under certain specific boundary conditions.
Fix $ \NN=2 $, and recall from~\eqref{eq: def--spin O(N) partition function} the measure of the XY model. For a subset $B\subset V$, recall from Section~\ref{sec: local time boundary condition} the local time boundary condition $\xi$.
For each component $ \mathsf{C}\in\mathcal{P}^{\ii}_B $, let $ \eta^{\ii}_{\mathsf{C}}=\eta^{\ii}_{u}:={\rm sgn}(\sigma^{\ii}_u) $ be the sign assigned to $ \mathsf{C} $ in copy $\ii$ for all $u\in\mathsf{C}$ and $\ii=\xx,\yy$ (by definition each component in the boundary partition is assigned with a common sign, for each copy $\ii=\xx,\yy$).
Throughout this subsection, we abbreviate the expectation with respect to the XY measure on $G$ by $ \langle\cdot\rangle_{G,\beta} $, and omit the subscript occasionally.

We begin with the random current representation. Recalling from Section \ref{sec: local time boundary condition} the XY partition function with local time boundary condition $\xi$ on $ B\subset V $ (note that the restriction of the normalized Haar measure $ \rd\sigma_u $ on $ \beta^{\ii}_u(\sigma^{\ii}_u)^2=\xi^{\ii}_u $ is the normalized counting measure on $ \{-1,1\} $ denoted by $ d\eta^{\ii}_u $), we have (in below let $ \varepsilon\odot\xi $ denote the spin boundary condition $ \{(\varepsilon^{\xx}_u\sqrt{\frac{\xi^{\xx}_u}{\beta_u}}, \varepsilon^{\yy}_u\sqrt{\frac{\xi^{\yy}_u}{\beta_u}})\}_{u\in B} $)
\begin{equation}
    \label{eq: Taylor expansion of Z}
    \begin{aligned}
        &\quad\ Z^{\xi}_{G, \beta}
        := 2^{-|\mathcal{P}^{\xx}_B|-|\mathcal{P}^{\yy}_B|} \sum_{\varepsilon^{\xx}\in \{-1,1\}^{\mathcal{P}^{\xx}_B}}\sum_{\varepsilon^{\yy}\in \{-1,1\}^{\mathcal{P}^{\yy}_B}} Z^{\varepsilon\odot\xi}_{G, \beta}
        \\&= 2^{-|\mathcal{P}^{\xx}_B|-|\mathcal{P}^{\yy}_B|}\sum_{\substack{\varepsilon^{\xx}\in \{-1,1\}^{\mathcal{P}^{\xx}_B},\\ \varepsilon^{\yy}\in \{-1,1\}^{\mathcal{P}^{\yy}_B} }}
        \int_{(\mathbb{S}^{1})^{V\setminus B}}\exp(\sum_{u\sim v}c^{\beta}_{uv}\sigma_u\cdot\sigma_v)
        \prod_{u\in B}\mathbf{1}_{\eta^{\xx}_u=\varepsilon^{\xx}_u,\eta^{\yy}_u=\varepsilon^{\yy}_u}
        \prod_{v\in V\setminus B}\rd\sigma_v\\
        &= \int_{\sigma}\prod_{u\sim v}\exp(c^{\beta}_{uv}\sigma^{\xx}_u\sigma^{\xx}_v+c^{\beta}_{uv}\sigma^{\yy}_u\sigma^{\yy}_v)\rd\theta d\eta,
    \end{aligned}
\end{equation}
where $ c^\beta_{uv}:=2J_{uv}\sqrt{\beta_u\beta_v} $, and $\rd\theta$ is the product of $ \rd\theta_v $ over $ v\in V\setminus B $ and $ d\eta $ is the product of $ d\eta^{\ii}_{\mathsf{C}^{\ii}} $ over $ \mathsf{C}^{\ii}\in \mathcal{P}^{\ii}_B $ and $ \ii\in\{\xx,\yy\} $ (recall that $ \sigma_v=(\cos\theta_v, \sin\theta_v) $).
Now we analyze~\eqref{eq: Taylor expansion of Z} via Taylor expansion. We call a nonnegative integer function $ n: E\to\mathbb{Z}_{\geq0} $ a current.
For a current $ (n_e)_{e\in E}\in\mathbb{Z}^E_{\geq0} $, let
\begin{equation}
    \label{eq: def random current weight}
    w_{\beta}(n):=\prod_{u\sim v}\frac{(c^{\beta}_{uv})^{n_{uv}}}{n_{uv}!},\qquad  d_v(n):=\sum_{e\ni v}n_e,\qquad\partial n:=\{v\in V: d_v(n)\text{ is odd}\}.
\end{equation}
Expanding the two exponentials in~\eqref{eq: Taylor expansion of Z} on every edge introduces two currents $ n^{\xx} $ and $ n^{\yy} $. Write $ d^{\sf{i}}_v=d_v(n^{\sf{i}}) $ with $ \sf{i}\in\{\xx,\yy\} $ for brevity, and let $ m_v:=\frac{d^{\xx}_v+d^{\yy}_v}{2} $.
Then we have
\begin{equation}
    \label{eq: Taylor expansion of Z-1}
    \begin{aligned}
        Z^{\xi}_{G, \beta}
        &=\int\sum_{n^{\xx}, n^{\yy}}\prod_{u\sim v}\frac{(c^{\beta}_{uv})^{n^{\xx}_{uv}}}{n^{\xx}_{uv}!}\cdot \frac{(c^{\beta}_{uv})^{n^{\yy}_{uv}}}{n^{\yy}_{uv}!} (\sigma^{\xx}_u\sigma^{\xx}_{v})^{n^{\xx}_{uv}}(\sigma^{\yy}_u\sigma^{\yy}_{v})^{n^{\yy}_{uv}}\rd\theta d\eta \\
        &= \sum_{n^{\xx}, n^{\yy}}w_{\beta}(n^{\xx})w_{\beta}(n^{\yy})\prod_{v\in V\setminus B}\int(\sigma^{\xx}_v)^{d^{\xx}_v}(\sigma^{\yy}_v)^{d^{\yy}_v}\rd\theta_v\prod_{\substack{\ii\in\{\xx,\yy\},\\ \mathsf{C}^{\ii}\in\mathcal{P}^{\ii}_B } }\int\prod_{u\in\mathsf{C}^{\ii}}(\sigma^{\ii}_u)^{d^{\ii}_u}d\eta_{\mathsf{C}^{\ii}}.
    \end{aligned}
\end{equation}
At each $ v\in V\setminus B $, we have (below $ 2|a $ means that $a$ is a multiple of $2$)
\[
\frac{1}{2\pi}\int^{2\pi}_0 (\cos\theta_v)^{d^{\xx}_v}(\sin\theta_v)^{d^{\yy}_v}\rd\theta_v = \frac{(d^{\xx}_v-1)!!(d^{\yy}_v-1)!!}{2^{m_v} m_v!}\cdot \mathbf{1}_{2\mid d^{\xx}_v}\mathbf{1}_{2\mid d^{\yy}_v}
\]
with the convention $ (-1)!!=1 $.
For each $ \ii\in\{\xx,\yy\} $ and each $ \mathsf{C}^{\ii}\in \mathcal{P}^{\ii}_B $, we have
\[
    \int\prod_{u\in\mathsf{C}^{\ii}}(\sigma^{\ii}_u)^{d^{\ii}_u}d\eta_{\mathsf{C}^{\ii}}
    = \frac{1}{2}\sum_{\eta^{\ii}_{\mathsf{C}^{\ii}}\in\{-1,1\}}\prod_{u\in\mathsf{C}^{\ii}}(\sigma^{\ii}_u)^{d^{\ii}_u}
    = \mathbf{1}_{2|\sum_{u\in\mathsf{C}^{\ii}}d^{\ii}_u} \prod_{u\in\mathsf{C}^{\ii}} \big(\frac{\xi^{\ii}_u}{\beta_u}\big)^{d^{\ii}_u/2}.
\]
Therefore, only currents satisfying
\begin{equation}
    \label{eq: sourceless currents}
    \mathcal{N}_{\mathcal{P}}:=
    \left\{
    (n^{\xx}, n^{\yy}) \in\mathbb{Z}^{E}_{\geq0}\times\mathbb{Z}^{E}_{\geq0}\;\middle|\;
    \begin{aligned}
        &2 \mid d^{\ii}_v  &&\text{for every }v\in V\setminus B, \\
        &2 \mid \sum_{u\in\mathsf{C}^{\ii}}d^{\ii}_u  &&\text{for every }\mathsf{C}^{\ii}\in\mathcal{P}^{\ii}_B,
    \end{aligned}
    \qquad\text{for }\ii\in\{\xx,\yy\} \right\}.
\end{equation}
makes non-zero contribution to~\eqref{eq: Taylor expansion of Z-1}, and thus we have
\begin{equation}
    \label{eq: Taylor expansion of Z-2}
    Z^{\xi}_{G,\beta}
    = \sum_{(n^{\xx}, n^{\yy})\in\mathcal{N}_{\mathcal{P}}}
    w_{\beta}(n^{\xx})w_{\beta}(n^{\yy}) \prod_{v\in V\setminus B}\frac{(d^{\xx}_v-1)!!(d^{\yy}_v-1)!!}{2^{m_v} m_v!}\prod_{u\in B}\big(\frac{\xi^{\xx}_u}{\beta_u}\big)^{d^{\xx}_u/2}\big(\frac{\xi^{\yy}_u}{\beta_u}\big)^{d^{\yy}_u/2}.
\end{equation}

The random current representation~\eqref{eq: Taylor expansion of Z-2} can be interpreted via uniform pairings at each vertex.
For each edge $ e\in E $ and each $ n^{\ii} $, $ \ii\in\{\xx,\yy\} $, consider $ n^{\ii}_e $ replicas of $e$ labeled by $ 1,\dots,n^{\ii}_e $ for the copy $ \ii $.
At each vertex $v\in V\setminus B$, the factor $ (d^{\xx}_v-1)!!(d^{\yy}_v-1)!! $ in~\eqref{eq: Taylor expansion of Z-2} counts the number of pairings of its adjacent replicate edges within each $\xx$- and $\yy$-copy.
At each boundary vertex $ u\in B $ there is only a factor $ \big(\frac{\xi^{\xx}_u}{\beta_u}\big)^{d^{\xx}_u/2}\big(\frac{\xi^{\yy}_u}{\beta_u}\big)^{d^{\yy}_u/2} $ in~\eqref{eq: Taylor expansion of Z-2}, and thus we do not introduce any pairing at these vertices.
Letting $ \pi\in \Pi(n^{\xx},n^{\yy}) $ be a pairing configuration given by such pairings at vertices $ v\in V\setminus B $, we can assign a weight
\[
W^{\xi}_{\beta}(n^{\xx}, n^{\yy}):=w_{\beta}(n^{\xx})w_{\beta}(n^{\yy})\prod_{v\in V\setminus B}\frac{1}{2^{m_v}m_v!}\prod_{u\in B}\big(\frac{\xi^{\xx}_u}{\beta_u}\big)^{d^{\xx}_u/2}\big(\frac{\xi^{\yy}_u}{\beta_u}\big)^{d^{\yy}_u/2}
\]
to the configuration $ (n^{\xx},n^{\yy}; \pi) $, so that
\begin{equation*}
    \begin{aligned}
        Z^{\xi}_{G,\beta}
        &= \sum_{(n^{\xx}, n^{\yy})\in\mathcal{N}_{\mathcal{P}}}
        w_{\beta}(n^{\xx})w_{\beta}(n^{\yy}) \prod_{v\in V\setminus B}\frac{(d^{\xx}_v-1)!!(d^{\yy}_v-1)!!}{2^{m_v} m_v!}\prod_{u\in B}\big(\frac{\xi^{\xx}_u}{\beta_u}\big)^{d^{\xx}_u/2}\big(\frac{\xi^{\yy}_u}{\beta_u}\big)^{d^{\yy}_u/2}\\
        &= \sum_{(n^{\xx}, n^{\yy})\in\mathcal{N}_{\mathcal{P}}} \sum_{\pi\in\Pi(n^{\xx}, n^{\yy})} W^{\xi}_{\beta}(n^{\xx}, n^{\yy}).
    \end{aligned}
\end{equation*}
Observe that each $ (n^{\xx},n^{\yy};\pi) $ induces a (copywise) loop/excursion configuration $ \mathfrak{L}(n^{\xx},n^{\yy};\pi) $ where each loop or each excursion in copy $\ii$ consists of replicated edges in copy $\ii$ paired at vertices in $ V\setminus B$ (thus the endpoints of excursions are in $B$) for $\ii=\xx,\yy$. For brevity we call $ \mathfrak{L}(n^{\xx},n^{\yy};\pi) $ a (copywise) random current loop soup (with local time boundary condition $\xi$). Note that by setting $ B=\varnothing $, the random current loop soup coincides with the (free boundary) random path model (see e.g.~\cite{TaggiLRO}).
\begin{lemma}
    \label{lem: rpm equivalence}
    The previously defined random current loop soup $ \mathfrak{L}(n^{\xx},n^{\yy};\pi) $ with local time boundary condition $\xi$ has the law of the (undirected) discrete skeletons (without considering the skeletons of point loops/excursions or loops contained in one edge) of two independent metric graph loop/excursion soups of intensity $1/2$ with local time condition $\xi$ and local time condition that $ L^{\xx}_v+L^{\yy}_v=\beta_v $ for $ v\in V\setminus B $.
\end{lemma}
\begin{proof}
    In the preceding discussion, we have already defined the law of $ \mathfrak{L}(n^{\xx},n^{\yy};\pi) $, with the weight $ W^{\xi}_{\beta}(n^{\xx},n^{\yy}) $.
    We plug $ c^{\beta}_{uv}=2J_{uv}\sqrt{\beta_u\beta_v} $ back into $ W^{\xi}_{\beta}(n) $, and get that (recall that we leave all the replicated edges unpaired at vertices in $B$)
    \begin{equation}
        \label{eq: Weight for random current loop}
        W^{\xi}_{\beta}(n) = \prod_{u\sim v}\prod_{\ii\in\{\xx,\yy\}}\frac{J^{n^{\ii}_{uv}}_{uv}}{n^{\ii}_{uv}!}
        \prod_{v\in V\setminus B}\frac{\beta^{m_v}_v}{\Gamma(m_v+1)}
        \prod_{u\in B}\prod_{\ii\in\{\xx,\yy\}}(2\xi^{\ii}_u)^{d^{\ii}_u/2}.
    \end{equation}
    For brevity, let $ m^{\ii}_v=d^{\ii}_v/2 $ for $ v\in V $ and $ \ii\in\{\xx,\yy\} $. We now compute the weight of the discrete skeleton of a continuous time random walk loop/excursion soup with boundary condition $\xi$ and local time conditioning.
    Recall that the simple random walk on $G$ has jump rates $J_{uv}$ (from $u$ to $v$) and total rates $ a_v=\kappa_v+\sum_{v\sim u}J_{uv} $ (where $ \kappa_v $ is the killing rate). Then its transition probabilities are $ p_{uv}=J_{uv}/a_u $.
    In each copy, take the (unconditioned) intensity-$\frac{1}{2}$ Poisson loop soup killed on $B$, together with the independent Poisson excursions from $B$ to $B$ with intensity $ 2\sqrt{\xi^{\ii}_a\xi^{\ii}_b} $ for each excursion from $a$ to $b$.
    In addition, writing $ \mathcal{E}^{\ii} $ as the collection of all such excursions for each copy $\ii=\xx,\yy$, we can define the parity event
    \[
    {\rm Even}(\mathcal{P}^{\ii}_B):=\{2\big|\sum_{s:s\neq r}N^{\ii}_{rs}\text{ for every }r=1,\dots,m^{\ii}\},
    \]
    where $ N^{\ii}_{rs} $ counts the excursions in $ \mathcal{E}^{\ii} $ between $ B^{\ii}_r\in\mathcal{P}^{\ii}_B $ and $ B^{\ii}_s\in\mathcal{P}^{\ii}_B $.
    Without the local time conditioning in $ V\setminus B $, the weight of configuration $ \mathfrak{L}(n^{\xx},n^{\yy};\pi) $ is proportional to (up to a constant irrelevant to the skeleton and local time)
    \begin{equation}
        \label{eq: free skeleton weight}
        \prod_{e\in E}\prod_{\ii=\xx,\yy}\frac{J^{n^{\ii}_e}_e}{n^{\ii}_e!}\prod_{v\in V\setminus B}\prod_{\ii=\xx,\yy}a^{-m^{\ii}_v}_v\prod_{u\in B}\prod_{\ii=\xx,\yy}(2\xi^{\ii}_u)^{m^{\ii}_u}.
    \end{equation}
    This is because each passing at a vertex $v\in V\setminus B$ of a trajectory (loop or excursion) contributes a factor $ a^{-1}_v $ (and each traversing to an edge $e\in E$ contributes a factor $J_e$), and each endpoint (of an excursion in copy $\ii$) at $ u\in B $ contributes a factor $ \sqrt{2\xi^{\ii}_u} $, and the labeling of replicated edges in copy $\ii$ contributes a factor $\frac{1}{n^{\ii}_e!}$.
    
    Conditional on these skeletons, we first analyze the local time on each (interior) vertex $ v\in V\setminus B $. Each visit to such $ v $ contributes an independent $ {\rm Exp}(a_v) $ waiting time. The zero-jump loops (i.e., point loops) at $v$ in each copy contribute a local time with density $ \frac{1}{2t}e^{-a_vt}dt $, and therefore their total contribution to the local time at $v$ is an independent $ {\rm Gamma}(1/2, a_v) $ variable where $ a_v $ is its rate.
    Consequently the total local time at $ v\in V\setminus B $ of the intensity-$\tfrac{1}{2}$ continuous time random walk loop soup, has distribution $ {\rm Gamma}(m^{\ii}_v+\frac{1}{2},a_v) $.
    Multiplying these Gamma densities by the skeleton weight~\eqref{eq: free skeleton weight} cancels the factors $ a^{-m^{\ii}_v} $, and thus the joint density of the skeleton and the (interior) local time at $ (\mathfrak{L}(n^{\xx},n^{\yy};\pi), b=(b^{\ii}_v)_{v\in V\setminus B}) $ is proportional to (up to a constant irrelevant to the skeleton and local time)
    \begin{equation}
        \label{eq: joint density}
        \prod_{e\in E}\prod_{\ii=\xx,\yy}\frac{J^{n^{\ii}_e}_e}{n^{\ii}_e!}\prod_{v\in V\setminus B}\prod_{\ii=\xx,\yy}\frac{e^{-a_v b^{\ii}_v}(b^{\ii}_v)^{m^{\ii}_v-\tfrac{1}{2}}}{\Gamma(m^{\ii}_v + \tfrac{1}{2})}\prod_{u\in B}\prod_{\ii=\xx,\yy}(2\xi^{\ii}_u)^{m^{\ii}_u}.
    \end{equation}
    Then at each $ v\in V\setminus B $, integrating the density~\eqref{eq: joint density} over the local time conditioning $ b^{\xx}_v+b^{\yy}_v=\beta_v $ implies that the density at $ (\mathfrak{L}(n^{\xx},n^{\yy};\pi), \beta=(\beta_v)_{v\in V\setminus B}) $ is proportional to
    \begin{equation}
        \label{eq: joint density-2}
        \begin{aligned}
            &\quad\ \int\prod_{e\in E}\prod_{\ii=\xx,\yy}\frac{J^{n^{\ii}_e}_e}{n^{\ii}_e!}
            \prod_{u\in B}\prod_{\ii=\xx,\yy}(2\xi^{\ii}_u)^{m^{\ii}_u}
            \prod_{v\in V\setminus B}\prod_{\ii=\xx,\yy}\frac{e^{-a_v b^{\ii}_v}(b^{\ii}_v)^{m^{\ii}_v-\tfrac{1}{2}}}{\Gamma(m^{\ii}_v + \tfrac{1}{2})}\mathbf{1}_{b^{\xx}_v + b^{\yy}_v=\beta_v,\,\forall v\in V\setminus B}\rd b \\
            &=e^{-\sum_{v\in V\setminus B}a_v\beta_v}\prod_{e\in E}\prod_{\ii=\xx,\yy}\frac{J^{n^{\ii}_e}_e}{n^{\ii}_e!}
            \prod_{v\in V\setminus B}\frac{\beta^{m_v}_v}{\Gamma(m_v+1)}
            \prod_{u\in B}\prod_{\ii=\xx,\yy}(2\xi^{\ii}_u)^{m^{\ii}_u}\,.
        \end{aligned}
    \end{equation}
    The right hand side of~\eqref{eq: joint density-2} is exactly the random current loop soup weight~\eqref{eq: Weight for random current loop} times a factor $ e^{-\sum_{v\in V\setminus B}a_v\beta_v} $, which depends only on the underlying graph (via jump rates) and the prescribed copywise boundary local time $ \beta $.
    Therefore, the continuous time random walk loop soup with boundary condition $\xi$ and the local time conditioning has a discrete skeleton exactly the same as the previously defined random current loop soup.
    Finally, the lemma follows by the coupling (applied before our local time conditioning) between the continuous time random walk loop soup and the metric graph loop soup (see e.g.~\cite[Section~2]{Lupu14}), where the same time change (stopped on $B$) maps the metric graph excursions to the aforementioned continuous time random walk excursions and preserves their interior vertex local time and their endpoints. Therefore, the parity events $ {\rm Even}(\mathcal{P}^{\ii}_B) $ are also preserved, and thus by Lemma~\ref{lem: block-parity} we get the desired metric graph loop/excursion soups under boundary condition $\xi$ and local time conditioning.
\end{proof}

Now we give the following lemma that bounds the two-point correlation function in an XY model with prescribed local time boundary conditions in terms of the loop connecting probability.
For nonintersecting subsets $ A, B\subset V $ choose arbitrary partitions $ \mathcal{P}^{\xx}_A $ of $A$, $ \mathcal{P}^{\xx}_B $ of $B$, and $ \mathcal{P}^{\yy}_{A\cup B} $. Let $ \xi^{A, B} $ be a local time boundary condition on $A\cup B$ such that
\begin{equation}
    \label{eq: cluster-loop for subsets with b.c.}
    \xi^{A, B}:=\big((\xi^{\xx}|_{A\cup B}=\beta|_{A\cup B}, \xi^{\yy}|_{A\cup B}=0), (\mathcal{P}^{\xx}_A\cup\mathcal{P}^{\xx}_B, \mathcal{P}^{\yy}_{A\cup B}) \big).
\end{equation}
Further, we define
\[
P^{\xi^{A, B}}_{\beta}(A, B):=\mu^{\xi^{A, B}}_{G,\beta}[\exists\,\ell\in \mathfrak{L}^{\xx}: \ell\cap A\neq\varnothing,\,\ell\cap B\neq\varnothing].
\]
\begin{lemma}
    \label{lem: cluster-loop comparison, all-x local time b.c.}
    Let $ \NN=2 $ and fix an arbitrarily chosen local time condition $ \beta $. For all disjoint $ A, B\subset V $ and all $ (u,v)\in A\times B $, we have (where we abbreviate $\xi^{A, B}$ as $\xi$)
    \begin{equation}
        \label{eq: cluster-loop, all-x b.c.}
        0\leq\langle\sigma_u\cdot\sigma_v\rangle^{\xi}_{\beta} \leq 2^{3/4}P^{\xi}_{\beta}(A, B)^{1/4}.
    \end{equation}
\end{lemma}
\begin{proof}
    We will prove this lemma by expressing the quantities under consideration in terms of auxiliary XY models (with ${\rm U}(1)$ symmetry), and then applying the Ginibre inequality. Let $S=A\cup B$, and $ \mathcal{P}^{\xx}_S:=\mathcal{P}^{\xx}_A\cup\mathcal{P}^{\xx}_B $. Note that since on $ \xi $ we have $ \xi^{\yy}|_{S}\equiv0 $, the choice of boundary connectivity in $\yy$-copy does not affect the measure $ \mu^{\xi}_{G,\beta} $. Therefore without loss of generality we may assume $ \mathcal{P}^{\yy}_S:=\mathcal{P}^{\yy}_{A\cup B}=\mathcal{P}^{\xx}_A\cup\mathcal{P}^{\xx}_B $, and abbreviate $ \mathcal{P}^{\ii}_S $ as $ \mathcal{P}_S $.

    We begin with defining the auxiliary XY models. Let $ e_1=(1,0) $ and $ e_2=(0,1) $.
    For each sign assignment $ \varepsilon\in\{-1,1\}^{\mathcal{P}_S} $ and each $ (i,j)\in\{1,2\}^2 $, let $ Z^{\varepsilon}_{ij} $ be the partition function of the XY model with spin boundary condition (below $ \varepsilon_u:=\varepsilon_{\mathsf{C}} $ for all $ u\in\mathsf{C}\in\mathcal{P}_S $)
    \[
    \sigma_u = \varepsilon(u)e_i\quad\text{for $u\in A$}
    \qquad\text{and}\qquad
    \sigma_v = \varepsilon(v)e_j\quad\text{for $v\in B$}, 
    \]
    and define
    \begin{equation}
        \label{eq: def-Z_{ij}}
        Z_{ij}:=2^{-|\mathcal{P}_S|}\sum_{\varepsilon\in\{-1,1\}^{\mathcal{P}_S}} Z^{\varepsilon}_{ij}.
    \end{equation}
    Namely, $ Z_{ij} $ is the partition function of the XY model with spin boundary condition $ \sigma_u\in\{\pm e_i\}$ for all $ u\in A $ and $ \sigma_v\in\{\pm e_j\} $ for all $ v\in B $. By copy symmetry, we have
    \begin{equation}
        \label{eq: Z_11=Z_22}
        Z_{11}=Z_{22}\qquad\text{and}\qquad Z_{12}=Z_{21}.
    \end{equation}
    Given $ Z_{ij} $, we can now define another auxiliary XY measure with the partition function
    \begin{equation}
        \label{eq: def-Z_0}
        \widehat{Z}:=2^{-4}\sum_{i,j\in\{1,2\}}Z_{ij}
    \end{equation}
    and write $ \widehat{\E} $ for the expectation under this auxiliary XY measure. Letting $ \hat{q}_m:=\widehat{\E}[\cos(m(\theta_u-\theta_v))] $ (where $\theta_u$ is the angle of $ \sigma_u $ under the angle representation; note that $ \theta_u\in\{0,\pi/2,\pi,3\pi/2\} $ for all $ u\in A\cup B $) for $m\geq1$, we claim that~\eqref{eq: cluster-loop, all-x b.c.} follows from
    \begin{equation}
        \label{eq: cluster-loop x-b.c. --1}
        \hat{q}_1 = \frac{\langle\sigma_u\cdot\sigma_v\rangle^{\xi}_{\beta}}{2-P^{\xi}_{\beta}(A, B)},
        \quad \hat{q}_2=\frac{P^{\xi}_{\beta}(A, B)}{2-P^{\xi}_{\beta}(A, B)},
    \end{equation}
    and
    \begin{equation}
        \label{eq: cluster-loop x-b.c. --2}
        \hat{q}_2\geq\hat{q}^4_1.
    \end{equation}
    In fact, plugging~\eqref{eq: cluster-loop x-b.c. --1} in~\eqref{eq: cluster-loop x-b.c. --2} yields
    \[
    \big(\langle\sigma_u\cdot\sigma_v\rangle^{\xi}_{\beta}\big)^4\leq P^{\xi}_{\beta}(A, B)\big( 2-P^{\xi}_{\beta}(A, B) \big)^3\leq 8P^{\xi}_{\beta}(A, B),
    \]
    which combined with the nonnegative correlation proves~\eqref{eq: cluster-loop, all-x b.c.}.
    Therefore it suffices to prove~\eqref{eq: cluster-loop x-b.c. --1} and~\eqref{eq: cluster-loop x-b.c. --2}.

    First we prove~\eqref{eq: cluster-loop x-b.c. --1} by analyzing the auxiliary XY measures. For the model with partition function~\eqref{eq: def-Z_0}, since $ \theta_u\in\{0,\pi/2,\pi,3\pi/2\} $ for all $ u\in A\cup B $, we know that $ \cos(2(\theta_u-\theta_v))=1 $ when spins in $ A $ and $ B $ take non-zero values in the same copy (which corresponds to the $Z_{11}$ and $Z_{22}$ fractions), and $ \cos(2(\theta_u-\theta_v))=-1 $ when spins in $ A $ and $ B $ take non-zero values in the opposite copy (which corresponds to the $Z_{12}$ and $Z_{21}$ fractions). Therefore we have
    \begin{equation}
        \label{eq: hat{q}_2}
        \hat{q}_2=\frac{Z_{11}+Z_{22}-Z_{12}-Z_{21}}{Z_{11}+Z_{22}+Z_{12}+Z_{21}}\stackrel{\eqref{eq: Z_11=Z_22}}{=}\frac{1-\frac{Z_{12}}{Z_{11}}}{1+\frac{Z_{12}}{Z_{11}}}.
    \end{equation}
    For $ \hat{q}_1 $, first observe that by definition~\eqref{eq: def-Z_{ij}}, the partition function for the model under $\xi$ defined in the statement is actually $ Z_{11} $. Since in the auxiliary model corresponding to $\widehat{Z}$), we have $ \cos(\theta_u-\theta_v)=0 $ when spins in $ A $ and $ B $ take non-zero values in the opposite copy (which corresponds to the $Z_{12}$ and $Z_{21}$ fractions), by expanding the total probability we know that
    \begin{equation}
        \label{eq: hat{q}_1}
        \hat{q}_1=\frac{Z_{11}+Z_{22}}{Z_{11}+Z_{22}+Z_{12}+Z_{21}}\cdot \langle\sigma_u\cdot\sigma_v\rangle^{\xi}_{\beta}=\frac{\langle\sigma_u\cdot\sigma_v\rangle^{\xi}_{\beta}}{1+\frac{Z_{12}}{Z_{11}}}.
    \end{equation}
    Now it suffices to express $ P^{\xi}_{\beta}(A, B) $ in terms of $ Z_{11} $ and $ Z_{12} $. In fact, we will prove that (note that plugging~\eqref{eq: P^xi in Z_ij} below in~\eqref{eq: hat{q}_2} and~\eqref{eq: hat{q}_1} proves~\eqref{eq: cluster-loop x-b.c. --1} immediately)
    \begin{equation}
        \label{eq: P^xi in Z_ij}
        \frac{Z_{12}}{Z_{11}} = 1 - P^{\xi}_{\beta}(A, B).
    \end{equation}
    To this end, we use the metric graph loop/excursion soup representation.
    In $Z_{11}$, we restrict to the configurations where every loop intersects at most one of $A$ and $B$ (and thus such configurations have total weight $ Z_{11}(1-P^{\xi}_{\beta}(A, B)) $). Then we can reassign the copy of every loop intersecting $B$ from $\xx$-copy to $\yy$-copy, which gives a bijection to the configurations contributing to $Z_{12}$. Indeed, in each image configuration, loops intersecting $A$ are in $\xx$-copy, while loops intersecting $B$ are in $\yy$-copy, and the inverse map changes the copy of each loop intersecting $B$ from $\yy$ to $\xx$. By Lemma~\ref{lem:yyy} this bijection is also weight-preserving.
    Consequently, this bijection implies
    \[
    Z_{11}(1-P^{\xi}_{\beta}(A, B))=Z_{12},
    \]
    proving~\eqref{eq: P^xi in Z_ij}.

    Now it suffices to prove~\eqref{eq: cluster-loop x-b.c. --2}.
    To this end, choose arbitrarily $ a_0\in A $ and $b_0\in B$, and for each component $ \mathsf{C}\in\mathcal{P}_S $, choose an (arbitrary) $ u_{\mathsf{C}}\in\mathsf{C} $. Then, consider the XY model with Hamiltonian
    \begin{equation}
        \label{eq: Cor 2.8}
        \begin{aligned}
            -\sum_{u,v\in V:\,u\sim v}c^{\beta}_{u,v}\cos(\theta_u-\theta_v) - \lambda
            & \Big[\sum_{v\in A\cup B}\cos(4\theta_v) + \sum_{v\in A\setminus\{a_0\}}\cos(2(\theta_v-\theta_{a_0})) \\
            &\quad+ \sum_{v\in B\setminus\{b_0\}}\cos(2(\theta_v-\theta_{b_0})) + \sum_{\mathsf{C}\in\mathcal{P}_S}\sum_{v\in\mathsf{C}\setminus\{u_{\mathsf{C}}\}}\cos(\theta_v-\theta_{u_{\mathsf{C}}}) \Big].
        \end{aligned}
    \end{equation}
    Since Hamiltonian~\eqref{eq: Cor 2.8} has the form (below $ \alpha\cdot\theta=\sum_{v\in V}\alpha_v\theta_v $)
    \[
    -\sum_{\alpha\in I}t_{\alpha}\cos(\alpha\cdot\theta)\qquad\text{for some finite }I\subset\mathbb{Z}^V,\qquad t_{\alpha}\geq0,
    \]
    the Ginibre inequality holds for all $\lambda>0$. Letting $\lambda\to\infty$, we know that the XY model under Hamiltonian~\eqref{eq: Cor 2.8} converges to the XY model with partition function $ \widehat{Z} $, and thus the Ginibre inequality holds (for the latter model) as well.
    Now we consider two independent samples $ \theta $ and $ \theta^{\prime} $ from the product measure $\widehat{\mathbb{E}}^{\otimes2}$, and let $ U=\theta_u-\theta_v $, $ V=\theta^{\prime}_u-\theta^{\prime}_v $. Letting $\theta^{\prime}$ be sampled from a replica of $\widehat{\E}$, defined on a disjoint replicated graph $ G^{\prime}=(V^{\prime}, E^{\prime}) $ of $G$, we can view the joint law of $ (\theta, \theta^{\prime}) $ as an XY measure defined on the graph $ (V\cup V^{\prime}, E\cup E^{\prime}) $, with Hamiltonian defined as the summation of the Hamiltonians of $ \theta $ and $\theta^{\prime}$ under $\widehat{\E}$.
    Applying the Ginibre inequality~\cite{Ginibre70} with respect to this XY model to $ \cos(U+V) $ and $ \cos(U-V) $, we have (note that $ \widehat{\E}[\sin U]=\widehat{\E}[\sin V]=0 $ by symmetry)
    \begin{equation*}
        \begin{aligned}
            \hat{q}_2
            &=\widehat{\E}^{\otimes 2}[\cos(U+V)\cos(U-V)]\geq\widehat{\E}^{\otimes 2}[\cos(U+V)]\widehat{\E}^{\otimes 2}[\cos(U-V)]\\
            &= \big(\widehat{\E}^{\otimes2}[\cos U\cos V]\big)^2 - \big( \widehat{\E}^{\otimes2}[\sin U\sin V] \big)^2 
            =\hat{q}^4_1,
        \end{aligned}
    \end{equation*}
    which proves~\eqref{eq: cluster-loop x-b.c. --2}, and thus concludes the proof.
\end{proof}
\begin{remark}\label{rmk: boundary two point function}
Let $A_u$ denote the component containing $u$ in $\mathcal{P}_A$ and let $B_v$ denote the component containing $v$ in $\mathcal{P}_B$. 
    Then by \cite{Lupu14}, we get that \begin{equation}
        \langle\sigma_u\cdot\sigma_v\rangle^{\xi}_{\beta} =\mu^{\xi}_{\beta}(u\con{\xx}v)=\mu^{\xi}_{\beta}(A_u\con{\xx}B_v).\nonumber
    \end{equation} 
\end{remark}

\section{An FKG inequality for the XY model}
\label{sec:FKG-XY}

In this section, we incorporate a copywise FKG inequality (see Theorem~\ref{thm:fkg-xy}) for the XY model, and present an RSW type estimate for the (copywise) cluster crossing probability as a consequence (see Corollary~\ref{cor:fkg-boundary-domination}).
The FKG inequality for copywise local times for the XY model was proved in~\cite{Chayes98} (see also~\cite{MC1998} for an extension to general spin $O(\NN)$ model). In addition to the FKG inequality, we also need the monotonicity with respect to the boundary condition, which can be derived from the FKG inequality (with some work). Instead of deriving this from the FKG inequality, we choose to present a self-contained proof from the dynamics point of view.
Our proof uses the loop soup connection to the XY model, with the following underlying intuition: if local times are larger (for the $\xx$-copy), then the loop soup clusters are (stochastically) larger (this is by the FKG for the FK-Ising model, which by Lemma~\ref{lem: loop cluster is FK cluster} is exactly the edge configurations in the loop soup model conditioned on local times, where an edge is open if all local times are nonzero in the corresponding cable edge); in addition, conditioned on local times and edge configurations outside a certain vertex $u$, the conditional law (with respect to the loop soup) for the local time at $u$ is increasing.
Crucially, for the XY models, if local times for the $\xx$-copy is larger, then it necessarily implies that the local times for the $\yy$-copy is smaller, and thus the corresponding conditional law for the local time at $u$ for the $\yy$-copy is decreasing. We then show that the aforementioned stochastic monotonicity on local times at $u$ can be preserved after we condition on their sums  being $\beta_u$. Once this is proved, then the desired FKG inequality can be shown to hold by running a dynamics, whose equilibrium measure is the XY measure.
Our proof follows this road map closely: in Section~\ref{subsec: block-parity decomposition} we decompose the local time of a vertex (in a critical loop soup) in terms of local time contributed by loops and excursions (see Lemma~\ref{lem: block-parity}); in Section~\ref{subsec: FK representation} we recall the Edwards--Sokal coupling and prove that the loop clusters are the same as the FK clusters; then in Section~\ref{subsec: update domination} we prove the domination of Glauber updates (Lemma~\ref{lem: update--domination}), and thus prove the FKG inequality (Theorem~\ref{thm:fkg-xy}).

\subsection{Block-parity decomposition}
\label{subsec: block-parity decomposition}
In this subsection we consider the (unconditional) critical loop soup (i.e., has intensity $\frac{1}{2}$), and introduce the block-parity decomposition, which will be used in the one-site update domination.
Let $\mathfrak{L}$ be a critical loop soup on a finite transient metric graph.
Consider a subset $ D_1\subset V $ and let $ A:=\partial D_1 $ be its outer boundary,
and $ D_2 = V\setminus D_1 $.
Let $ \widetilde{D}_1 $ and $ \widetilde{D}_2 $ be metric graphs such that cable edges with at least one endpoint in $ D_1 $ belong to $ \widetilde{D}_1 $, and cable edges connecting two vertices in $D_2$ belong to $ \widetilde{D}_2 $.
In particular, local time configurations on $A$ and the partition of $A$ induced by loop clusters restricted to $ \widetilde{D}_2 $ serve as boundary conditions on $ \widetilde{D}_1 $.
Let $\mathcal{P}_A$ be the partition of $A$ induced by the local time cluster connectivity of $ \mathfrak{L} $ restricted to $ \widetilde{D}_2 $ (recall the prescribed boundary partition introduced in Section~\ref{sec: local time boundary condition}).
Fix a local time condition $ \lambda=(\lambda_a)_{a\in A}>0 $ and a (realization of the) partition $\mathcal P=\{B_1,\ldots,B_m\}$ of $A$, and put $ q_a=\sqrt{2\lambda_a} $ (here we use a different choice of notation other than $\xi$, to emphasize that we are dealing with one copy of critical loop soup without local time condition in this subsection).
For every pair of $ a, b\in A $, let $\mathsf{n}_{ab}$ denote the (undirected) Brownian excursion measure between $a$ and $b$ in $ \widetilde{D}_1 $ (supported on excursions that start at $ a $ and then only hit $A$ at $b$),
and set
\begin{equation}
    \label{eq:fkg-excursion-intensity}
    \Theta_{\lambda}=\sum_{\{a, b\}\subset A: a\neq b}q_aq_b\,\mathsf{n}_{ab} + \frac{1}{2}\sum_{a\in A}q^2_a\mathsf{n}_{aa}.
\end{equation}

We will study the excursion configurations (i.e., the collection of all excursions) in $ \widetilde{D}_1 $ denoted by $ \mathcal{E} $. Let $N_{rs}=N_{rs}(\mathcal{E})$ count the excursions between $B_r$ and $B_s$ in $ \mathcal{E} $ for $ 1\leq r, s\leq m $, and define the event (recall the similar definitions for random current excursions in the proof of Lemma~\ref{lem: rpm equivalence})
\[
\operatorname{Even}(\mathcal P)=\Big\{\sum_{s:s\ne r}N_{rs}\in2\mathbb Z \text{ for every }r\Big\}.
\]
Furthermore, let $\mathcal{E}_{\mathcal{P}, \lambda}$ be a Poisson point process of intensity $\Theta_{\lambda}$ conditioned on event $ \operatorname{Even}(\mathcal P) $, and let $L^{\prime}_0$ be the local time field of an independent critical loop soup in $\widetilde{D}_1$ killed on hitting $A$ (by isomorphism theorem equivalently, $ L^{\prime}_0 $ is the half-square of an independent GFF on $ \widetilde{D}_1 $ with zero boundary condition on $A$).

\begin{lemma}[Block-parity decomposition]
    \label{lem: block-parity}
    Using notations above, we have the following: for any non-negative $ \lambda $,
    \begin{equation}
        \label{eq: block-parity}
        (L_z)_{z\in \widetilde{D}_1}\,\big|\{L {|}_{A}=\lambda,\ \mathcal{P}_A=\mathcal P\}\ \overset{d}{=}\ \big(L^{\prime}_0(z)+\sum_{\gamma\in\mathcal{E}_{\mathcal{P}, \lambda}} L_z(\gamma)\big)_{z\in \widetilde{D}_1}
    \end{equation}
    where $ L_z(\gamma) $ is the local time accumulated at $z$ by $\gamma$ (we emphasize that in this lemma the original laws are critical loop soups, without conditioning from the XY model).
\end{lemma}

\begin{proof}
By the coupling (see~\cite[Section~2]{Lupu14}) between the metric graph loop soup and the continuous time random walk loop soup (or by fine-meshing the edges of $D_1$ and letting the mesh size converge to zero), it suffices to prove~\eqref{eq: block-parity} restricted to the discrete graph $ D_1 $, which is done by comparing the Laplace transforms against nonnegative, compactly supported test functions $f$ defined on $ D_1$.
Recall that $ \mathcal{E} $ is the excursion configuration.
Write
\[
T_f(\gamma) = \sum_{v\in D_1}L_v(\gamma)f(v)
\quad\text{for $ \gamma\in\mathcal{E} $,}
\qquad\text{and let}\quad
\mathsf{L}_0(f)=\mathbb E\exp\{-\langle f,L^{\prime}_0\rangle\}.
\]
For a boundary-sign profile (which takes a constant sign on each $B_i$ for $i=1,\dots, m$) $\varepsilon\in\{-1,1\}^m$, and an excursion $\gamma$ with endpoints in $B_r,B_s$, define $\chi_\varepsilon(\gamma)=\varepsilon_r\varepsilon_s$.
We will use the following parity observation, which allows us to represent the parity event $ {\rm Even}(\mathcal{P}) $ in terms of the $\chi$-function.

As in the Ising model random current representation (where we view each $ B_r $ as a big `node' with $ \varepsilon_r $ being its spin, and each $ \gamma $ being one count of current among these nodes), we claim that
\begin{equation}
    \label{eq: parity-identity}
    \mathbf 1_{\operatorname{Even}(\mathcal P)}(\mathcal E)=2^{-m}\sum_{\varepsilon\in\{-1,1\}^m}\prod_{\gamma\in\mathcal E}\chi_\varepsilon(\gamma).
\end{equation}
Indeed, only excursions between different blocks enter the right-hand side non-trivially (and the left-hand side also depends only on these excursions); summation over $ \epsilon_r\in\{-1,1\} $ at each block $ B_r $ with an odd number of excursions induces a factor of zero, while summation at each block with an even number of excursions induces a factor of two.
Thus, this verifies~\eqref{eq: parity-identity}.
The Poisson exponential formula, therefore, gives the Laplace transform of the right-hand side of~\eqref{eq: block-parity} as
\begin{equation}
    \label{eq: parity--laplace of RHS}
    \mathsf{L}_0(f)\cdot \frac{\displaystyle\sum_\varepsilon\exp\!\Big\{\int\big(\chi_\varepsilon(\gamma) e^{-T_f(\gamma)}-1\big)\,\Theta_{\lambda}(d\gamma)\Big\}}{\displaystyle\sum_\varepsilon\exp\!\Big\{\int\big(\chi_\varepsilon(\gamma)-1\big)\,\Theta_{\lambda}(d\gamma)\Big\}}.
\end{equation}
Note that all the displayed integrals are finite for these test measures.

We now compute the Laplace transform for the left-hand side of~\eqref{eq: block-parity} from the GFF perspective.
We begin with the following observation.
On the event $\{L_A=\lambda,\mathcal{P}_A=\mathcal P\}$ and for each boundary-sign profile $\varepsilon$, the boundary values have the form $q_v^\varepsilon=\varepsilon_rq_v$ for every $v\in B_r$.
Provided with a boundary-sign profile $\varepsilon$, the GFF domain Markov property (see~\cite[Section~3]{Lupu14} for example) gives
\begin{equation}
    \label{eq: parity-DMC}
    \varphi|_{\widetilde{D}_1}\ \overset d=\ \varphi_0+h^\varepsilon,
\end{equation}
where $\varphi_0$ is the zero-boundary GFF in $\widetilde{D}_1$ and $h^\varepsilon$ is the harmonic extension of $q^\varepsilon$ in $\widetilde{D}_1$.
For $ j\in\{1,2\} $, by the domain Markov property, we can write the partition function of $ \varphi|_{\widetilde{D}_j} $ given the boundary condition $ q^{\varepsilon} $ with partition $ \mathcal{P} $ on $A$ as $ \mathcal{Z}_j(q^{\varepsilon}) $.
Then conditioned on $q$ and $ \mathcal{P} $, we claim that  the partition function $ \mathcal{Z}_2(q^{\varepsilon}) $ does not depend on $\varepsilon$:
since the partition (of $A$) is induced from the cluster connectivity of the GFF configuration restricted to $D_2$, flipping the sign of each cluster does not change the partition; the local time and the FK-connectivity in $D_2$ also remain unchanged, and thus their probability density remains the same.
The conditional (on $q$ and $ \mathcal{P} $) weight of $\varepsilon$ is therefore
\begin{equation}
    \label{eq:fkg-boundary-sign-weights}
    \P[\varphi|_A=q^{\varepsilon} \big| L_A=\lambda,\mathcal{P}_{A}=\mathcal{P} ]
    =\frac{\mathcal{Z}_1(q^\varepsilon)}{\sum_{\delta\in\{-1,1\}^m}\mathcal{Z}_1(q^\delta)},
\end{equation}
and thus we have
\begin{equation}
    \label{eq: parity--Laplace transform of LHS}
    \begin{aligned}
        &\qquad\E\big[ \exp\{ -\langle f, L|_{D_1} \rangle \} \big| \{L|_A=\lambda, \mathcal{P}_{A}=\mathcal{P}\} \big]\\
        &= \sum_{\varepsilon\in\{-1,1\}^{m}}\P[\varphi|_A=q^{\varepsilon} \big| L_A=\lambda,\mathcal{P}_{A}=\mathcal{P} ] \cdot \E[ \exp\{-\langle f, L|_{D_1} \rangle\} \big| \varphi|_A=q^{\varepsilon}, \mathcal{P}_{A}=\mathcal{P} ]\\
        &= \frac{1}{\sum_\delta\mathcal{Z}_1(q^\delta)}\sum_{\varepsilon\in\{-1,1\}^{m}} \mathcal{Z}_1(q^\varepsilon) \cdot \E[ \exp\{-\langle f, L|_{D_1} \rangle\} \big| \varphi|_A=q^{\varepsilon}, \mathcal{P}_{A}=\mathcal{P} ].
    \end{aligned}
\end{equation}

Write
\[
I_f(\varepsilon) :=\E[ \exp\{-\langle f, L|_{D_1} \rangle\} \big| \varphi|_A=q^{\varepsilon}, \mathcal{P}_{A}=\mathcal{P} ],
\]
and let
\[
U(\varepsilon):=\int\big(\chi_\varepsilon(\gamma)-1\big)\,\Theta_{\lambda}(d\gamma),\qquad V_f(\varepsilon):=\int\chi_\varepsilon(\gamma)(e^{-T_f(\gamma)}-1)\,\Theta_{\lambda}(d\gamma).
\]
We will prove~\eqref{eq: block-parity} by showing that
\begin{equation}
    \label{eq: block-parity-STP}
    \begin{aligned}
        &\quad\ \text{Laplace transform of the LHS of~\eqref{eq: block-parity}}
        = \frac{\sum_{\varepsilon}\mathcal{Z}_1(q^{\varepsilon}) I_f(\varepsilon)}{\sum_{\delta}\mathcal{Z}_1(q^{\delta})}\\
        &= \mathsf{L}_0(f)\frac{\sum_{\varepsilon}\exp(U(\varepsilon)+V_f(\varepsilon))}{\sum_{\delta}\exp(U(\varepsilon))}
        = \text{Laplace transform of the RHS of~\eqref{eq: block-parity}}.
    \end{aligned}
\end{equation}
In fact, the first equality of~\eqref{eq: block-parity-STP} follows from~\eqref{eq: parity--Laplace transform of LHS} and the last equality of~\eqref{eq: block-parity-STP} was derived from~\eqref{eq: parity--laplace of RHS} (see also discussions preceding it).
Therefore, it suffices to prove the second equality of~\eqref{eq: block-parity-STP}.

To this end, we will show that
\begin{equation}
    \label{eq: parity--gaussian identities}
    \mathcal{Z}_1(q^\varepsilon) = \mathcal{Z}_1(q^1) \,e^{U(\varepsilon)},
    \qquad
    I_f(\varepsilon)=\mathbb E\exp\!\left\{-\tfrac12 \langle f,(\varphi_0+h^\varepsilon)^2\rangle\right\}=\mathsf{L}_0(f)\,e^{V_f(\varepsilon)}
\end{equation}
where $ q^1 $ is given by $ \varepsilon\equiv 1 $, namely the boundary condition with plus signs on every $B_r\in\mathcal{P}$.
Write the test function as $ f=\sum_{x\in D_1} f_x\mathbf{1}_x $, with non-negative entries $ f_x\geq0 $.
Let $M$ be the Laplacian matrix in $D_1$ (the diagonal of $M$ includes the killings on hitting $A$), and let $R_{av}$ be the conductance from $a\in A$ to $ v\in D_1 $.
Write
\[
G=M^{-1},
\qquad
F=\operatorname{diag}(f_x:x\in D_1),
\qquad
G_f=(M+F)^{-1}.
\]
Then we have (see~\cite[Section~5]{LeJan2007})
\begin{equation}
    \label{eq: parity--DMC}
    \varphi_0|_{D_1}\sim\mathcal N(0,G),
    \qquad
    h^\varepsilon|_{D_1}=GRq^\varepsilon.
\end{equation}

We first compute the excursions integrals $U(\varepsilon)$ and $V_f(\varepsilon)$.
Observe that every excursion from $a$ to $b\ne a$ visits $\widetilde{D}_1$ (by our metric graph convention of $ \widetilde{D}_1 $ and $\widetilde{D}_2$, edges with both endpoints in $A$ belong to $\widetilde{D}_2$), and excursions avoiding $\widetilde{D}_1$ have $T_f=0$.
An excursion from $a\in A$ entering $D_1$ at $x$ contributes an entrance factor $R_{xa}$; its subsequent path is a process killed on $A$, with an exit factor $R_{yb}$ on hitting $A$ at $b\in A$ from $D_1\ni y\sim b$.
Then by~\cite[Proposition~1]{LeJan2007} (and the discussion preceding it), for all $a,b\in A$,
\begin{equation*}
    \int_{\{\gamma:\gamma\cap D_1\ne\varnothing\}} e^{-T_f(\gamma)}\,\mathsf{n}_{ab}(d\gamma)
    =\sum_{x,y\in {D_1}}R_{xa}G_f(x,y)R_{yb}=(R^{\mathsf T}G_fR)_{ab}.
\end{equation*}
Consequently (letting $f=0$)
\begin{align*}
    \mathsf{n}_{ab}(1)
    &=(R^{\mathsf T}GR)_{ab},
    &&a\ne b,\\
    \int(e^{-T_f(\gamma)}-1)\,\mathsf{n}_{ab}(d\gamma)
    &=(R^{\mathsf T}(G_f-G)R)_{ab},
    &&a,b\in A.
\end{align*}
The left-hand side in the second equality is a finite integral even when $\mathsf n_{aa}(1)=\infty$, since only excursions visiting $\widetilde{D}_1$ contribute to it.
Since $\chi_\varepsilon(\gamma)q_aq_b =q_a^\varepsilon q_b^\varepsilon$ for an excursion $\gamma$ from $a$ to $b$, we obtain
\begin{equation}
    \label{eq: parity--U(vareps)}
    \begin{aligned}
        U(\varepsilon)
        &= \int\big(\chi_\varepsilon(\gamma)-1\big)\,\Theta_{\lambda}(d\gamma)
        = \frac{1}{2}\sum_{a,b\in A}q_aq_b\int(\chi_{\epsilon}(\gamma)-1)\mathsf{n}_{ab}(d\gamma) \\
        &= \frac12\sum_{a,b\in A: a\neq b} (q_a^\varepsilon q_b^\varepsilon-q_aq_b) (R^{\mathsf T}GR)_{ab}
        = \frac12\left[(q^\varepsilon)^{\mathsf T}R^{\mathsf T}GRq^\varepsilon -q^{\mathsf T}R^{\mathsf T}GRq\right],
    \end{aligned}
\end{equation}
and
\begin{equation}
    \label{eq: parity--V_f(vareps)}
    \begin{aligned}
        V_f(\varepsilon)
        &= \int\chi_\varepsilon(\gamma)(e^{-T_f(\gamma)}-1)\,\Theta_{\lambda}(d\gamma)
        = \frac{1}{2}\sum_{a,b\in A}q_aq_b \int\chi_{\epsilon}(\gamma) (e^{-T_f(\gamma)}-1) \mathsf{n}_{ab}(d\gamma)
        \\ &= \frac12\sum_{a,b\in A}q_a^\varepsilon q_b^\varepsilon (R^{\mathsf T}(G_f-G)R)_{ab}
        = \frac12(q^\varepsilon)^{\mathsf T} R^{\mathsf T}(G_f-G)Rq^\varepsilon.
    \end{aligned}
\end{equation}

We now compute the Gaussian integrals $\mathcal{Z}_1(q^{\varepsilon})$ and $ I_f(\varepsilon) $. For $ \mathcal{Z}_1(q^{\varepsilon}) $, the calculation
\[
\mathcal{Z}_1(z)=\int_{\mathbb R^{D_1}} e^{-\frac12 y^{\mathsf T}My+y^{\mathsf T}Rz}\,dy =(2\pi)^{|{D_1}|/2}(\det M)^{-1/2} e^{\frac12z^{\mathsf T}R^{\mathsf T}GRz},
\]
yields
\begin{equation}
    \label{eq: parity--Z_1(q^vareps)}
    \frac{\mathcal{Z}_1(q^\varepsilon)}{\mathcal{Z}_1(q^1)}
    =\frac{\displaystyle\int_{\mathbb R^{D_1}} e^{-\frac12y^{\mathsf T}My+y^{\mathsf T}Rq^\varepsilon}\,dy} {\displaystyle\int_{\mathbb R^{D_1}} e^{-\frac12y^{\mathsf T}My+y^{\mathsf T}Rq^1}\,dy}
    =\exp\!\Big\{\frac12\left[ (q^\varepsilon)^{\mathsf T}R^{\mathsf T}GRq^\varepsilon -(q^1)^{\mathsf T}R^{\mathsf T}GRq^1\right]\Big\} =e^{U(\varepsilon)},
\end{equation}
and thus proves the first equality of~\eqref{eq: parity--gaussian identities}.
For $I_f(\varepsilon)$, by~\eqref{eq: parity--DMC} the conditional (on $ q^{\varepsilon} $ and $ \mathcal{P} $) interior field $ \varphi|_{D_1} $ has mean $GRq^\varepsilon$ and covariance $G$. Thus
\begin{equation*}
    \begin{aligned}
        I_f(\varepsilon)
        &=\frac{\displaystyle\int_{\mathbb R^{D_1}} e^{-\frac12y^{\mathsf T}(M+F)y+y^{\mathsf T}Rq^\varepsilon}\,dy}{\displaystyle\int_{\mathbb R^{D_1}} e^{-\frac12y^{\mathsf T}My+y^{\mathsf T}Rq^\varepsilon}\,dy}
        =\left(\frac{\det M}{\det(M+F)}\right)^{1/2} \exp\!\left\{\frac12(q^\varepsilon)^{\mathsf T} R^{\mathsf T}(G_f-G)Rq^\varepsilon\right\}\\
        &= \mathsf{L}_0(f) \exp\!\left\{\frac12(q^\varepsilon)^{\mathsf T} R^{\mathsf T}(G_f-G)Rq^\varepsilon\right\}
        = \mathsf{L}_0(f) \cdot e^{V_f(\varepsilon)},
    \end{aligned}
\end{equation*}
where the first equality in the second row is by the (killed) loop-soup/GFF isomorphism theorem (see for example~\cite [Theorem~1]{Lupu14} or~\cite[Proposition~3]{LeJan2007}), and the second equality in the second row is derived from~\eqref{eq: parity--V_f(vareps)}.
Combining the preceding display and~\eqref{eq: parity--Z_1(q^vareps)} proves~\eqref{eq: parity--gaussian identities}, and thus concludes the proof of~\eqref{eq: block-parity}.
\end{proof}

\subsection{Edwards--Sokal coupling and copywise FK representation}
\label{subsec: FK representation}

We now give the (copywise) FK representation of the spin $O(\NN)$ model, as prerequisite for the XY model FKG inequality and comparison of boundary conditions.

Recall the spin convention from~\eqref{eq: def--spin O(N) spin}, and recall that $ L^{\xx} $, $ L^{\yy} $ are local time fields for the $\xx$, $\yy$ copies respectively.
Conditionally on $ L^{\xx}=\beta^{\xx} $ and $ L^{\yy}=\beta^{\yy} $, the two sign fields are independent Ising models, with coupling constants $2J_{uv}\sqrt{\beta_u^{\mathsf{i}}\beta^{\mathsf{i}}_v}$ in copy $\mathsf{i}=\mathsf{x}, \mathsf{y}$.
Consider their Edwards--Sokal edge configurations $\omega^{\mathsf{x}},\omega^{\mathsf{y}}$ as follows.
Given the signs, edges are \emph{open} independently indicated
\begin{equation}
    \label{eq: FK edge parameters}
    p^{\mathsf{i}}_{uv}=1-e^{-2J_{uv}|\phi_u^{\mathsf{i}}\phi^{\mathsf{i}}_v|}=1-e^{-4J_{uv}\sqrt{\beta_u^{\mathsf{i}}\beta^{\mathsf{i}}_v}},\qquad \mbox{for the edge $(u,v)$ in the $\mathsf{i}$ copy}
\end{equation}
when $ \eta^{\mathsf{i}}_u\eta^{\mathsf{i}}_v=1 $ (we use the convention that $ \omega^{\mathsf{i}}_{u,v}=1 $ when open, and $ \omega^{\mathsf{i}}_{u,v}=0 $ when closed), and are \emph{closed} always otherwise.
We denote by $\mu_\beta^\xi$ the marginal for $(\beta^{\mathsf{x}},\omega^{\mathsf{x}},\omega^{\mathsf{y}})$ (note that $\beta^{\mathsf{y}}=\beta-\beta^{\mathsf{x}}$ is determined by $\beta^{\mathsf{x}}$, and denote by $ \Omega^{\mathsf{x}} $ the collection of all possible realizations for $ (\beta^{\mathsf{x}}, \omega^{\mathsf{x}}, \omega^{\mathsf{y}}) $).
We define an $ \mathsf{x} $-preferred partial order on $\Omega^{\mathsf{x}}$ by (with coordinate-wise comparisons in below)
\begin{equation}
    \label{eq:fkg-order}
    (\overline{\beta}^{\mathsf{x}}, \overline{\omega}^{\mathsf{x}}, \overline{\omega}^{\mathsf{y}}) \preceq (\widetilde{\beta}^{\mathsf{x}},\widetilde\omega^{\mathsf{x}},\widetilde\omega^{\mathsf{y}})\ \Longleftrightarrow\ \overline{\beta}^{\mathsf{x}}\leq\widetilde{\beta}^{\mathsf{x}},\ \overline{\omega}^{\mathsf{x}}\leq\widetilde\omega^{\mathsf{x}},\ \overline{\omega}^{\mathsf{y}}\geq\widetilde\omega^{\mathsf{y}}
\end{equation}
for all $ (\overline{\beta}^{\mathsf{x}},\overline{\omega}^{\mathsf{x}},\overline{\omega}^{\mathsf{y}}), (\widetilde{\beta}^{\mathsf{x}},\widetilde\omega^{\mathsf{x}},\widetilde\omega^{\mathsf{y}})\in\Omega^{\mathsf{x}} $.
Thus this order increases with $\beta^{\mathsf{x}},\omega^{\mathsf{x}}$ and decreases with $\beta^{\mathsf{y}},\omega^{\mathsf{y}}$.
An \emph{$\mathsf{x}$-increasing} function is defined as an increasing function of $(\beta^{\mathsf{x}},\omega^{\mathsf{x}})$, and a \emph{$\mathsf{y}$-increasing} function is an increasing function of $(\beta^{\mathsf{y}},\omega^{\mathsf{y}})$.
The former is increasing, and the latter is decreasing with respect to the $\mathsf{x}$-preferred order as in~\eqref{eq:fkg-order}.

\begin{lemma}
    \label{lem: loop cluster is FK cluster}
    Consider the XY loop soup local time fields under the boundary condition $\xi$ defined in Section \ref{sec: local time boundary condition}, and their coupled (via the isomorphism theorem) sign fields $ L^{\ii}_z=\frac{1}{2}(\varphi^{\ii}_z)^2 $ with $ \ii=\xx,\yy $. Also define $ \widehat{\omega}^{\ii}_e:=\mathbf{1}_{L^{\ii}_z>0,\,\forall z\in\widetilde{e}} $ for each edge $e$ and its cable edge $\widetilde{e}$. Then:
    \begin{itemize}
        \item Conditional on the (copywise and signed) spin fields, these indicators are independent over edges and copies, with
        \[
        \mu^{\xi}_{G, \beta}(\widehat{\omega}^{\ii}_{uv}=1\ \big|\ \varphi|_V) = \mathbf{1}_{\varphi^{\ii}_u\varphi^{\ii}_v>0}(1-e^{-2J_{uv}|\varphi^{\ii}_u\varphi^{\ii}_v|}).
        \]
        Thus they have the exact same Edwards--Sokal law~\eqref{eq: FK edge parameters}.
        \item For every cable subgraph $ \widetilde{D} $ formed by vertices and complete cable edges, the partition of its positive local time vertices induced by $ \{z\in\widetilde{D}: L^{\ii}_z>0 \} $ equals the partition induced by the open edges $ \widehat{\omega}^{\ii}\big|_{E(D)} $.
    \end{itemize}
\end{lemma}
\begin{proof}
    Under the coupling~\cite[Theorem~1]{Lupu14} between metric loop soup and GFF, the local time field is the half-square of the spin field. Conditional on the spin fields on vertices, the values on cable edges are independent Brownian bridges, and thus the first assertion follows from the Brownian motion reflection principle.
    For the second assertion, observe that a path with positive signs between vertices traverses precisely those cable edges where the signed field $\varphi^{\ii}$ stays positive, which in turn are exactly the FK open edges defined as above, and thus proving the second assertion.
\end{proof}
\begin{remark}
    Lemma~\ref{lem: loop cluster is FK cluster} states that the (copywise) loop clusters are exactly the previously defined FK clusters, thus in the rest of this section we denote by $\omega^{\ii}_{e}$ for the edge openness in both of them.
    Consequently, the boundary connectivity used in Lemma~\ref{lem: block-parity} is equivalent to the boundary connectivity induced by the corresponding (copywise) FK clusters, thus allowing us to write the explicit update law at a vertex in Lemma~\ref{lem: update site density}.
\end{remark}

\subsection{FKG inequality via update domination}
\label{subsec: update domination}
Recall the vertex set $V$ and the (possibly empty) boundary set $B$.
For $u\in V\setminus B$, let $E_u$ be the edges incident to $u$ and $\mathcal{N}_u$ its neighbors.

We define the configuration outside $u$ by
\begin{equation}
    \label{eq: def--update environment}
    \tau_u=\big((\beta_v^{\mathsf{x}})_{v\ne u},\omega^{\mathsf{x}}|_{E\setminus E_u},\omega^{\mathsf{y}}|_{E\setminus E_u}\big),
\end{equation}
and we denote by $\Omega_u$ the collection of all possible configurations $\tau_u$.
Given $ \tau_u $, the local times in the $\mathsf{y}$-copy satisfies $\beta_v^{\mathsf{y}}=\beta_v-\beta_v^{\mathsf{x}}$ for $v\neq u$.
Let $\mathcal P_\tau^{\mathsf{x}},\mathcal P_\tau^{\mathsf{y}}$ be the partitions of $\mathcal{N}_u$ induced by connectivities generated by FK clusters on $ V\setminus\{u\} $ in the $\mathsf{x}$- and $\mathsf{y}$-copy respectively.
For copy label $\mathsf{i}=\mathsf{x,y}$ and $ \mathsf{C}\in\mathcal P_\tau^{\mathsf{i}}$, define
\[
A_\mathsf{C}^{\mathsf{i}}=\sum_{v\in \mathsf{C}}J_{uv}\sqrt{2\beta^{\mathsf{i}}_v}.
\]
The update at $u$ resamples $\beta_u^{\mathsf{x}}$ (and hence $\beta_u^{\mathsf{y}}$) together with (openness/closeness of) edges in $E_u$ in both copies, conditional on $\tau_u$. In what follows we denote the updating vertex by $u$ and abbreviate $ \tau_u $ as $\tau$.

\begin{lemma}[One-site update density]
    \label{lem: update site density}
    The conditional density of $\beta_u^{\mathsf{x}}$ is
    \begin{equation}
        \label{eq: update site density}
        \rho_\tau(t)=\frac{1}{\bar{C}_\tau\sqrt{t(\beta_u-t)}}\prod_{\mathsf{C}\in\mathcal P_\tau^{\mathsf{x}}}\cosh\bigl(\sqrt{2t}\,A_\mathsf{C}^{\mathsf{x}}\bigr)\prod_{\mathsf{C}\in\mathcal{P}_\tau^{\mathsf{y}}}\cosh\bigl(\sqrt{2(\beta_u-t)}\,A_\mathsf{C}^{\mathsf{y}}\bigr),\qquad 0<t<\beta_u,
    \end{equation}
    where $\bar{C}_\tau$ is a normalizing constant. Conditional on $\beta_u^{\mathsf{x}}=t$ and $\tau$, the two incident edge configurations are independent FK-Ising cluster configurations on the star (rooted at $u$ with leaves $\mathcal{N}_u$), with boundary connectivities $\mathcal P_\tau^{\mathsf{x}},\mathcal P_\tau^{\mathsf{y}}$ (induced by clusters off $u$) and parameters~\eqref{eq: FK edge parameters}. (Here, the conditioning fixes the copied local times and the exterior FK edge configurations. The signs assigned in the Edwards--Sokal coupling are summed out, with a common sign on each exterior FK cluster, for both $ \mathsf{x} $ and $\mathsf{y}$-copies.)
\end{lemma}
\begin{proof}
    Fix $i\in\{\mathsf{x}, \mathsf{y}\}$, and abbreviate $w^{\mathsf{i}}_v=J_{uv}\sqrt{2\beta^{\mathsf{i}}_v}$ and $\mathcal P=\mathcal P_\tau^{\mathsf{i}}$. Let $a_u>0$ be the diagonal entry of the Laplacian matrix at $u$.
    The local time accumulated at $u$ by loops disjoint from $ \mathcal{N}_u $ is $\operatorname{Gamma}(1/2,a_u)$, with the second parameter being the rate.
    In addition, for excursions between vertices in $\mathcal{N}_u$ that visit $u$, each of them will contribute a local time at $u$, as an independent exponential variable with rate $ a_u $.
    Without taking into account the block-parity conditioning $ {\rm Even}(\mathcal{P}) $ (let $\mathbb P_0$ denote the law without this conditioning), the number of excursions between $ v,w\in \mathcal{N}_u $ are independent Poisson variable with mean
    \begin{equation}
        \label{eq:fkg-star-intensities}
        \Lambda_{vw}=\frac{w_vw_w}{a_u}\quad(v\neq w),\qquad\Lambda_{vv}=\frac{w_v^2}{2a_u}.
    \end{equation}
    For brevity we introduce an apriori order on $\mathcal{N}_u$ and define the summation convention
    \[
    \sum_{v\leq w} := \sum_{v\in\mathcal{N}_u} + \sum_{\{v,w\}\subset\mathcal{N}_u:v\neq w}
    \]
    for functions defined on $ \mathcal{N}_u $.
    Let $ M_{vw}=M_{wv}\sim\operatorname{Poisson}(\Lambda_{vw}) $ be independent Poisson variables for all $ v, w\in\mathcal{N}_u $, and write
    \[
    N^0_u=\sum_{v\leq w}M_{vw},
    \qquad
    \Lambda_*:=\sum_{v\leq w}\Lambda_{vw}.
    \]
    Furthermore, excursions not hitting $u$ are irrelevant: on the one hand, they do not contribute to the local time at $u$; on the other hand, within the star, they do not contribute to connectivity between any $ v,w\in \mathcal{N}_u $.

    Let $N_u$ be the number of excursions hitting $u$.
    Then $N_u$ is distributed as $N^0_u$ given $ {\rm Even}(\mathcal{P}) $.
    For the sign-profile $\varepsilon\in\{-1,1\}^{\mathcal P}$, set (we slightly abuse the notation of $ \varepsilon_{\mathsf{C}}=\varepsilon_a $ for all $ a\in \mathsf{C}\in\mathcal{P} $ since the sign-profile $ \varepsilon $ assigns a common sign to all vertices in each cluster in the partition $\mathcal{P}$)
    \[
    H_\varepsilon=\sum_{\mathsf{C}\in\mathcal P}\varepsilon_\mathsf{C} \sum_{v\in \mathsf{C}}w_v.
    \]
    Since $ \sum_{v\leq w}\Lambda_{vw}\varepsilon_{\mathsf{C}(v)}\varepsilon_{\mathsf{C}(w)} = H_\varepsilon^2/(2a_u) $, we have
    \begin{equation}
        \label{eq: update-star-counts}
        \mathbb P(N_u=n)=\mathbb P_0\bigl(N_u^{0}=n \big| \operatorname{Even}(\mathcal P)\bigr)=\big(\sum_\varepsilon e^{H_\varepsilon^2/(2a_u)}\big)^{-1}\cdot \sum_\varepsilon \frac{H_\varepsilon^{2n}}{(2a_u)^n n!}
    \end{equation}
    In fact, we can prove~\eqref{eq: update-star-counts} via a generating function calculation. Let $\mathsf{C}(v)\in\mathcal P$ be the block containing $v$.
    For $z\in[0,1]$, the parity identity~\eqref{eq: parity-identity} and the Poisson generating function give
    \begin{equation*}
        \begin{aligned}
            &\quad\ \mathbb{E}_0\!\left[z^{N^0_u} \mathbf1_{\operatorname{Even}(\mathcal P)}\right]
            = 2^{-|\mathcal P|}\sum_\varepsilon \prod_{v\leq w} \mathbb E_0\!\left[ \bigl(z\varepsilon_{\mathsf{C}(v)}\varepsilon_{\mathsf{C}(w)}\bigr)^{M_{vw}} \right]
            \\&= 2^{-|\mathcal P|}e^{-\Lambda_*} \sum_\varepsilon \exp\!\big(z\sum_{v\leq w} \Lambda_{vw}\varepsilon_{\mathsf{C}(v)}\varepsilon_{\mathsf{C}(w)}\big)
            = 2^{-|\mathcal P|}e^{-\Lambda_*} \sum_\varepsilon \exp\!\big(\frac{zH_\varepsilon^2}{2a_u}\big).
        \end{aligned}
    \end{equation*}
    Taking the coefficient of $z^n$ gives
    \[
    \mathbb P_0\bigl(N^0_u=n,\,\operatorname{Even}(\mathcal P)\bigr) =2^{-|\mathcal P|}e^{-\Lambda_*} \sum_\varepsilon\frac{H_\varepsilon^{2n}}{(2a_u)^n n!}.
    \]
    Dividing by the same generating function evaluated at $z=1$ (which then gives the probability for $ \mathrm{Even}(\mathcal{P}) $) proves~\eqref{eq: update-star-counts}.
    
    Conditional on $N_u=n$, the local time $L_u$ is $\operatorname{Gamma}(n+\frac{1}{2},a_u)$. Using~\eqref{eq: update-star-counts} and using
    \[
    n!\,\Gamma(n+1/2)=\frac{\sqrt\pi\,(2n)!}{4^n},\qquad\sum_\varepsilon\cosh(rH_\varepsilon)=2^{|\mathcal P|}\prod_{C\in\mathcal P}\cosh\!\big(r\sum_{v\in C}w_v\big),
    \]
    we obtain the density of $L_u$ as (up to normalization)
    \begin{equation}
        \label{eq:fkg-one-copy-density}
        t^{-1/2}e^{-a_ut}\prod_{\mathsf{C}\in\mathcal P} \cosh\!\big(\sqrt{2t}\sum_{v\in \mathsf{C}}w_v\big).
    \end{equation}
    Before imposing the total local time conditioning at $u$,
    let $T^{\mathsf{x}}_u,T^{\mathsf{y}}_u$ represent its copywise local times. Conditional on $\tau$, by~\eqref{eq:fkg-one-copy-density} (applying it to the $\mathsf{x}$- and the $\mathsf{y}$-copy) they are independent and have densities
    \[
    f_i(s)=\bar{c}_i s^{-1/2}e^{-a_us} \prod_{\mathsf{C}\in\mathcal P_\tau^{\mathsf{i}}} \cosh\!\left(\sqrt{2s}\,A_\mathsf{C}^{\mathsf{i}}\right),\qquad s>0,\quad i=x,y,
    \]
    where $\bar{c}_i>0$ are normalizing constants. The change of variables $(T^{\mathsf{x}}_u,T^{\mathsf{y}}_u)=(t,b-t)$ has Jacobian one, so
    \[
    \mathbb P(T^{\mathsf{x}}_u\in dt,T^{\mathsf{x}}_u+T^{\mathsf{y}}_u\in db\mid\tau)
    =f_x(t)f_y(b-t)\mathbf1_{\{0<t<b\}}\,dt\,db.
    \]
    Therefore for $0<t<\beta_u$,
    \begin{equation*}
        \begin{aligned}
            \rho_\tau(t)
            &=\frac{f_x(t)f_y(\beta_u-t)}{\displaystyle\int_0^{\beta_u}f_x(s)f_y(\beta_u-s)\,ds}
            \\&= \frac{1}{\bar{C}_\tau\sqrt{t(\beta_u-t)}} \prod_{\mathsf{C}\in\mathcal P_\tau^{\mathsf{x}}} \cosh\!\left(\sqrt{2t}\,A_\mathsf{C}^{\mathsf{x}}\right) \prod_{\mathsf{C}\in\mathcal P_\tau^{\mathsf{y}}} \cosh\!\left(\sqrt{2(\beta_u-t)}\,A_\mathsf{C}^{\mathsf{y}}\right),
        \end{aligned}
    \end{equation*}
    since the factor $c_xc_y e^{-a_u\beta_u}$ does not depend on $t$ and cancels between the numerator and denominator, and thus proves~\eqref{eq: update site density}.

    One can also check the same formula directly in the Edwards--Sokal representation. Put $r_x=\sqrt{2t}$ and $r_y=\sqrt{2(\beta_u-t)}$, and set $K_{uv}^{\mathsf{i}}=J_{uv}r_i\sqrt{2\beta_v^{\mathsf{i}}}$ for $ i\in\{\mathsf{x}, \mathsf{y}\} $.
    The Edward--Sokal coupling writes (denote the sign of $i$-copy at $v$ by $\eta^{\mathsf{i}}_v$)
    \[
    e^{K_{uv}^{\mathsf{i}}\eta_u^{\mathsf{i}}\eta_v^{\mathsf{i}}}
    =e^{-K_{uv}^{\mathsf{i}}}\big[1+(e^{2K_{uv}^{\mathsf{i}}}-1) \mathbf1_{\{\eta_u^{\mathsf{i}}=\eta_v^{\mathsf{i}}\}}\big],
    \]
    which, after expanding over edges and summing over feasible sign choices (i.e., the sign on each cluster $ C\in\mathcal{P} $ is constant), gives
    \[
    \mathbb P(\omega^{\mathsf{i}}|_{E_u}=\alpha\mid t,\tau)
    \propto
    2^{k_{\mathcal P_\tau^{\mathsf{i}}}(\alpha)} \prod_{e\in E_u}(p_e^{\mathsf{i}})^{\alpha_e}(1-p_e^{\mathsf{i}})^{1-\alpha_e},
    \qquad p_e^{\mathsf{i}}=1-e^{-2K_e^{\mathsf{i}}}.
    \]
    Here $k_{\mathcal P_\tau^{\mathsf{i}}}(\alpha)$ counts the components of the star where connectivity is jointly given by the partition $\mathcal P_\tau^{\mathsf{i}}$ and the edge configuration $\alpha$ in the star, since the summation over each of such clusters gives a factor of two.
    The weights are (conditionally) independent between $\mathsf x$- and $\mathsf y$-copies, thus proving the asserted two FK-Ising laws and their independence.
\end{proof}

\begin{lemma}[XY update domination]
    \label{lem: update--domination}
    Suppose $ \tau^{\pm}=((\beta^{\mathsf{x},\pm}_v)_{v\ne u},\omega^{\mathsf{x},\pm}|_{E\setminus E_u},\omega^{\mathsf{y},\pm}|_{E\setminus E_u})\in\Omega_u $ satisfy that
    \[
    \beta_v^{\mathsf{x},+}\geq \beta_v^{\mathsf{x},-}\quad(v\ne u),\qquad\omega^{\mathsf{x},+}|_{E\setminus E_u}\geq\omega^{\mathsf{x},-}|_{E\setminus E_u},\qquad\omega^{\mathsf{y},+}|_{E\setminus E_u}\leq \omega^{\mathsf{y},-}|_{E\setminus E_u}.
    \]
    Then their updates at $u$ admit a coupling such that
    \begin{equation}
        \label{eq: update-domination}
        \beta_u^{\mathsf{x},+}\geq \beta_u^{\mathsf{x},-},\qquad \omega^{\mathsf{x},+}|_{E_u}\geq\omega^{\mathsf{x},-}|_{E_u},\qquad \omega^{\mathsf{y},+}|_{E_u}\leq\omega^{\mathsf{y},-}|_{E_u}.
    \end{equation}
    In fact, $\rho_{\tau^+}(t)/\rho_{\tau^-}(t)$ (recalling~\eqref{eq: update site density}) is non-decreasing in $t$.
\end{lemma}
\begin{proof}
For nonnegative weights $w=(w_v)$ on a finite set and a partition $\mathcal P$, define
\[
F_{\mathcal P,w}(r)=\prod_{\mathsf{C}\in\mathcal P}\cosh\!\big(r\sum_{v\in \mathsf{C}}w_v\big).
\]
Then by~\eqref{eq: update site density} we have
\begin{equation}
    \label{eq: update--F and rho}
    \rho_{\tau}(t) = \frac{1}{C_{\tau}\sqrt{t(\beta_u-t)}}\cdot F_{\mathcal{P}^{\mathsf{x}}_{\tau}, w^{\mathsf{x}}}(\sqrt{2t})\cdot F_{\mathcal{P}^{\mathsf{y}}_{\tau}, w^{\mathsf{y}}}(\sqrt{2(\beta_u-t)})
\end{equation}
where $ w_v^{\mathsf{i}}=J_{uv}\sqrt{2\beta_v^{\mathsf{i}}} $ for $ i\in\{\mathsf{x},\mathsf{y}\} $.
For two partitions $ \mathcal{P} $ and $ \mathcal{P}^{\prime} $, we say $ \mathcal{P} \leq \mathcal{P}^{\prime} $ if for each $ \mathsf{C}\in \mathcal{P} $ there exists $ \mathsf{C}^{\prime}\in \mathcal{P}^{\prime} $ such that $ \mathsf{C}\subset \mathsf{C}^{\prime} $. Define
\[
f(r):=\frac{F_{\mathcal{P}^{\prime}, w^{\prime}}(r)}{F_{\mathcal{P}, w}(r)}.
\]

Then we claim that when $ (w,\mathcal{P})\leq (w^{\prime}, \mathcal{P}^{\prime}) $, the function $f(r)$ is increasing in $r$. The monotonicity with respect to $r$ is obvious, when $ \mathcal{P} = \mathcal{P}^{\prime} $ and $ w \leq w^{\prime} $ (since $ \cosh (cr) $ is increasing when $ c\geq 0 $), and it remains to prove the monotonicity when $ w = w^{\prime} $ and $ \mathcal{P} < \mathcal{P}^{\prime} $. To this end, by induction it suffices to consider $ \mathcal{P} $ and $ \mathcal{P}^{\prime} $ such that $ \mathcal{P} \setminus \{ \mathsf{C}_1, \mathsf{C}_2\} = \mathcal{P}^{\prime}\setminus \{\mathsf{C}^{\prime}\} $ with $ \mathsf{C}_1\cup \mathsf{C}_2=\mathsf{C}^{\prime} $ (i.e., $ \mathcal{P}^{\prime} $ is obtained from $\mathcal{P}$ by merging two sets $ \mathsf{C}_1, \mathsf{C}_2 $).
Then (since the common clusters shared by the two partitions cancel each other and the weights are identical)
\[
f(r)=\frac{\cosh (r\sum_{v\in \mathsf{C}}w_v)}{\cosh (r\sum_{v\in \mathsf{C}_1}w_v)\cosh (r\sum_{u\in \mathsf{C}_2}w_u)}=1+\tanh (r\sum_{v\in \mathsf{C}_1}w_v)\tanh (r\sum_{u\in \mathsf{C}_2}w_u)\geq1,
\]
verifying the desired monotonicity.

For the $\mathsf{x}$-copy, $w^{\mathsf{x},+}\geq w^{\mathsf{x},-}$ and $\mathcal P^{\mathsf{x},+}\geq \mathcal P^{\mathsf{x},-}$. Thus (by the preceding claim on the monotonicity of $f$) the ratio $$ \frac{F_{\mathcal{P}^{\mathsf{x},+}_{\tau}, w^{\mathsf{x},+}}(\sqrt{2t})}{F_{\mathcal{P}^{\mathsf{x},-}_{\tau}, w^{\mathsf{x},-}}(\sqrt{2t})} $$ is non-decreasing in $\sqrt{2t}$.
For the $\mathsf{y}$-copy, $w^{\mathsf{y},+}\leq w^{\mathsf{y},-}$ and $ \mathcal P^{\mathsf{y},+}\leq \mathcal P^{\mathsf{y},-} $, so
\[
\frac{F_{\mathcal{P}^{\mathsf{y},+}_{\tau}, w^{\mathsf{y},+}}(\sqrt{2t})}{F_{\mathcal{P}^{\mathsf{y},-}_{\tau}, w^{\mathsf{y},-}}(\sqrt{2t})}
\]
is non-increasing similarly.
Combined with~\eqref{eq: update--F and rho}, it implies that $ \rho_{\tau^+}(t)/\rho_{\tau^-}(t) $ is increasing in $t$. This immediately implies the existence of the coupling such that $ \beta^{\mathsf{x},+}_u\geq \beta^{\mathsf{x},-}_u $.

It remains to couple the incident edges. This follows immediately from the fact that the FK-Ising model is monotone with respect to the boundary conditions (i.e., corresponding to  the partition $ \mathcal{P}^{\pm} $ here) and coupling constants (see~\cite{GrimmettRC2003} for instance).
\end{proof}

Now we are ready to give the FKG inequality for the XY model.
\begin{theorem}[FKG inequalities for the XY model]
    \label{thm:fkg-xy}
    Consider the $\mathrm{XY}$ model defined by~\eqref{eq: def--spin O(N) partition function}.
    For every copywise local time boundary condition $\xi$ as defined in Section \ref{sec: local time boundary condition} (in addition, we allow it to include the boundary partitioning conditions), the marginal $\mu_\beta^\xi$ on $ (\beta^{\mathsf{x}},\omega^{\mathsf{x}},\omega^{\mathsf{y}}) $ is positively associated with respect to the ordering given in~\eqref{eq:fkg-order}: for bounded increasing functions $F,H$ of $ (\beta^{\mathsf{x}},\omega^{\mathsf{x}},\omega^{\mathsf{y}}) $,
    \begin{equation}
        \label{eq:fkg-associated}
        \mu_\beta^\xi(FH)\geq \mu_\beta^\xi(F)\,\mu_\beta^\xi(H).
    \end{equation}
    In particular, two $\mathsf{x}$-increasing functions are positively correlated, as are two $\mathsf{y}$-increasing (namely, decreasing in ordering~\eqref{eq:fkg-order}) functions. If $F$ is $\mathsf{x}$-increasing and $G$ is $\mathsf{y}$-increasing, then
    \begin{equation}
        \label{eq:fkg-cross-suppression}
        \mu_\beta^\xi(FG)\leq \mu_\beta^\xi(F)\,\mu_\beta^\xi(G).
    \end{equation}
    Consequently, for an $\mathsf{x}$-increasing event $\mathcal A$ of positive probability,
    \[
    \mu_\beta^\xi(F\mid\mathcal A)\geq\mu_\beta^\xi(F),\qquad\mu_\beta^\xi(G\mid\mathcal A)\leq\mu_\beta^\xi(G).
    \]
\end{theorem}
\begin{proof}
    In order to prove the FKG inequality for the desired measure, we design a version of Glauber dynamics with update kernel $T$ such that (i) every update preserves positive association and (ii) the dynamics converges to $ \mu^{\xi}_{\beta} $ in equilibrium.

    The dynamics is defined as follows. We start from the same (local time and FK-edge) configuration in both the $\mathsf{x}$- and $ \mathsf{y} $-copies.
    Listing the vertices whose local times are not fixed by the boundary condition as $u_1,\ldots,u_n$, we update the vertices in the fixed cyclic order $ u_1,u_2,\ldots,u_n $ iteratively (we call each round of updates from $u_1$ to $u_n$ a sweep).
    At each selected vertex (in each sweep), perform the following two steps.
    \begin{enumerate}
        \item\label{item: Glauber step 1} update its local times $ \beta^{\mathsf{x}}_u $ (and thus $ \beta^{\mathsf{y}}_u $) according to the conditional distribution given $ \tau_u $ (recall~\eqref{eq: def--update environment});
        \item\label{item: Glauber step 2} conditioned on the local times at all vertices, we sample the edge connections according to the FK-Ising model for both $\mathsf{x}$- and $\mathsf{y}$-copies (with parameters given by~\eqref{eq: FK edge parameters}) independently.
    \end{enumerate}
    We first verify (i), the preservation of positive association.
    To this end, suppose that $ \vartheta $ is a law on the configuration space satisfying positive association. We next show that after Step $i$ (for $ i=1,2 $), the law $ \vartheta^{\prime} $ we obtain still satisfies the positive association.

    \textsf{Proof of Step~\ref{item: Glauber step 1}.}
For $s\in(0,1)$, define (recall~\eqref{eq: update site density})
\[
Q_\tau(s) :=\inf\left\{t\in[0,\beta_u]: \int_0^t\rho_\tau(r)\,dr\geq s\right\}.
\]
Lemma~\ref{lem: update--domination} gives
\[
 \tau\preceq\widetilde\tau
 \quad\Longrightarrow\quad
 Q_\tau(s)\leq Q_{\widetilde\tau}(s),
\]
and $Q_\tau$ is non-decreasing in $s$.
Given $(\beta^{\mathsf{x}},\omega^{\mathsf{x}},\omega^{\mathsf{y}})$, recall that
\[
 \tau_u
 =\left(
     (\beta_v^{\mathsf{x}})_{v\ne u},
     \omega^{\mathsf{x}}|_{E\setminus E_u},
     \omega^{\mathsf{y}}|_{E\setminus E_u}
   \right).
\]
For $s\in(0,1)$, define
\[
 \widehat\beta_v^{\mathsf{x}}(s)
 =\begin{cases}
    Q_{\tau_u}(s),&v=u,\\
    \beta_v^{\mathsf{x}},&v\ne u,
   \end{cases}
 \qquad
 \widehat\beta_v^{\mathsf{y}}(s)=\beta_v-\widehat\beta_v^{\mathsf{x}}(s).
\]
Then by definition, Step~\ref{item: Glauber step 1} is performed via replacing $(\beta^{\mathsf{x}},\omega^{\mathsf{x}},\omega^{\mathsf{y}})$ by
$(\widehat\beta^{\mathsf{x}}(s),\omega^{\mathsf{x}},\omega^{\mathsf{y}})$, with $s$ sampled by an independent uniform measure on $ (0,1) $.
Then the corresponding update kernel satisfies (since the local times at vertices other than $u$ and FK-edge configurations in both $\mathsf{x}$- and $\mathsf{y}$-copies remain unchanged)
\[
(K_uF)(\beta^{\mathsf{x}},\omega^{\mathsf{x}},\omega^{\mathsf{y}})
=\int_0^1 F(\widehat\beta^{\mathsf{x}}(s),\omega^{\mathsf{x}},\omega^{\mathsf{y}})\,ds.
\]
Consequently, $K_uF$ is increasing whenever $F$ is increasing.
For bounded increasing $F,H$ (of $\mathsf{x}$-copy), the scalar covariance identity gives
\begin{align*}
 &K_u(FH)-(K_uF)(K_uH)\\
 &=\frac12\int_0^1\int_0^1
   [F(\widehat\beta^{\mathsf{x}}(s),\omega^{\mathsf{x}},\omega^{\mathsf{y}})-F(\widehat\beta^{\mathsf{x}}(t),\omega^{\mathsf{x}},\omega^{\mathsf{y}})]
   [H(\widehat\beta^{\mathsf{x}}(s),\omega^{\mathsf{x}},\omega^{\mathsf{y}})-H(\widehat\beta^{\mathsf{x}}(t),\omega^{\mathsf{x}},\omega^{\mathsf{y}})]\,ds\,dt
 \geq0.
\end{align*}
If $\vartheta$ is positively associated, it follows that
\begin{align*}
 \operatorname{Cov}_{\vartheta K_u}(F,H)
 &=\vartheta\!\left[K_u(FH)-(K_uF)(K_uH)\right]
     +\operatorname{Cov}_\vartheta(K_uF,K_uH)
 \geq0,
\end{align*}
where $ \vartheta K_u $ is the law of $ \vartheta $ after an update $ K $ at $u$. Corresponding results also hold similarly for $ F $ being increasing in $ \mathsf{x} $-copy and $ H $ being increasing in $ \mathsf{y} $-copy.
Thus Step~~\ref{item: Glauber step 1} preserves positive association.

    \textsf{Proof of Step~\ref{item: Glauber step 2}.}
    Since by Step~\ref{item: Glauber step 1} we already have the domination of $ \mathsf{x} $- and $ \mathsf{y} $-copy local times and FK-edge configurations in the sense of~\eqref{eq: update-domination}, the desired positive association after Step~\ref{item: Glauber step 2} follows from the FKG inequality of the FK-Ising model (see~\cite{GrimmettRC2003} for instance).

    {Now we can verify (i). Let $T_u$ denote the complete two-step update kernel at $u$, and let $ T=T_{u_1}T_{u_2}\cdots T_{u_n} $ be the kernel of one full sweep. The preceding arguments show that each $T_u$ preserves positive association, and thus, starting from a deterministic admissible configuration $\omega_0$, the law $\delta_{\omega_0}T^m$ is positively associated for every $m\geq0$.

    Next we verify (ii). Note that each $T_u$ preserves $\mu_\beta^\xi$: the first step preserves its local-time marginal, and the second step preserves the conditional FK-edge distribution. Consequently, $\mu_\beta^\xi T=\mu_\beta^\xi$.
    By Lemma~\ref{lem: update site density} we have
    \[
    \rho_\tau(t)
    \geq
    \frac{\overline{c}_u}{\pi\sqrt{t(\beta_u-t)}}\geq\frac{2\overline{c}_u}{\pi\beta_u}>0
    \]
    for every vertex $u$, where $\overline{c}_u>0$ is a constant not depending on the update environment $\tau$. Therefore the sweep update kernel $T$ satisfies the minorisation condition of~\cite[Theorem~8]{RobertsRosenthal2004}, and thus verifying (ii).}
\end{proof}

\begin{remark}[Conditional extension to $O(\NN)$]
    \label{rem:fkg-on-conditional}
    For $\NN\geq3$, once we condition on the total local time for the first two copies at each vertex $v$ being $\beta_v$, the first two copies then form an $\text{XY}$ model with $\beta=(\beta_v)_{v\in V}$. Under the local time conditioning $\beta$, the FKG inequality for the first copy follows from Theorem~\ref{thm:fkg-xy}.
\end{remark}

\paragraph{Box-crossing probabilities.}
Now with the FKG inequalities, we can give the following estimate of the box-crossing probabilities in the XY model.
The key to bounding the crossing probability is to relate it to spin correlations, which can then be bounded by the spin wave argument. However, the challenge lies in the fact that the connecting probability between two sets cannot be easily written as spin correlations. In order to address this, we apply the FKG inequality, which essentially allows us to glue a set to a point, and thus, the aforementioned strategy can be carried out.
Recall the order $ \mathcal{P}\leq\mathcal{P}^{\prime} $ right below~\eqref{eq: update--F and rho}.
For further convenience, we give the following definition.
\begin{definition}
    \label{def: best x bc for crossing}
    For an annulus $A$, we define the all-$\mathsf x$-copy local time condition $\xi^*$ on $\partial A$ such that on $\partial A$ all the local times are on the $\mathsf x$-copy and that the outer and inner boundaries of $\partial A$ are wired into single points, respectively. 
\end{definition}
\begin{remark}
    \label{rmk: xi^* largest boundary condition}
    Note that $ \xi^* $ is the largest possible $\xx$-preferred boundary condition that can be induced by local time connectivities outside $A$. This is obvious by noting that the inner and outer boundaries cannot be connected by a path outside $A$.
\end{remark}
\begin{corollary}
    \label{cor:fkg-boundary-domination}
    Let $ D\subset \mathbb{Z}^2 $, and consider copywise local time boundary condition on $D$.
    Write $\xi\preceq\widetilde\xi$ if $\xi^{\mathsf{x}}\leq\widetilde\xi^{\mathsf{x}}$ coordinate-wise (or equivalently $ \xi^{\mathsf{y}} \geq \widetilde\xi^{\mathsf{y}} $) and their boundary partitions satisfy $ \mathcal{P}_{\xi^{\xx}}\leq\mathcal{P}_{\widetilde{\xi}^{\xx}} $ and $ \mathcal{P}_{\xi^{\yy}}\geq\mathcal{P}_{\widetilde{\xi}^{\yy}} $.
    Then
    \begin{equation}
        \label{eq:fkg-boundary-monotonicity}
        \xi\preceq\widetilde\xi\quad\Longrightarrow\quad \mu_\beta^\xi\preceq\mu_\beta^{\widetilde\xi} \quad\text{under the order~\eqref{eq:fkg-order}}.
    \end{equation}
    In addition, for the annulus $ A=\{v\in\Z^2: \ell \leq|v|_{\infty}\leq L \} $ and event
    \[
    \mathcal{C}^{\mathsf{x}}:=\{\mbox{there exists an $\xx$-cluster crossing the annulus }A\},
    \]
    there exists a constant $ C_{\beta_0}>0 $ depending only on $ \beta_0 $ such that
    \begin{equation}
        \label{eq:fkg-crossing-boundary-comparison-2}
        \mu_\beta^{\xi^*}(\mathcal C^{\mathsf{x}})\leq \frac{C_{\beta_0}}{\log(1+\frac{L}{\ell})}
    \end{equation}
    for all $ 1\leq\ell\leq L $ and all local time condition $ \beta\leq\beta_0 $.
\end{corollary}
\begin{proof}
    Use the same updates as in Theorem~\ref{thm:fkg-xy} with the two ordered boundary conditions $ \xi,\,\widetilde{\xi} $. Lemma~\ref{lem: update--domination} couples the chains in order at every step. Passing to equilibrium proves~\eqref{eq:fkg-boundary-monotonicity}.

    In addition, for the annulus case, suppose that under $\xi^*$ the inner boundary is wired into a single vertex $ v_{in} $ and the outer boundary is wired into a single vertex $v_{out}$ separately.
    Under the boundary condition $\xi^*$, we have that $ \sigma^{\xx}_{v_{in}} $ and $ \sigma^{\xx}_{v_{out}} $ are just the signs of the $ \xx $-spin values on $v_{in}$ and $v_{out}$, and therefore
    \[
    \mu^{\xi^*}_{\beta}(\mathcal C^{\mathsf{x}}) = \E_{\mu^{\xi^*}_{\beta}}[\sigma_{v_{in}}\sigma_{v_{out}}] \leq \frac{C_{\beta_0}}{\log(1+\frac{L}{\ell})},
    \]
    where the last inequality follows from a standard McBryan-Spencer argument (see the argument of~\cite[Theorem~9.12]{Velenik}).
\end{proof}

\begin{remark}
    \label{rmk: RSW, CBC, FKG}
    In the rest of the paper we will refer to~\eqref{eq:fkg-crossing-boundary-comparison-2} as $ \mathsf{RSW} $, to~\eqref{eq:fkg-boundary-monotonicity} as $ \mathsf{CBC} $,  and to Theorem~\ref{thm:fkg-xy} as $ \mathsf{FKG} $.
\end{remark}

\section{Gaussian decay of crossing loops}
\label{sec: loop Gaussian decay}

In this section, we consider the XY model and prove a Gaussian tail estimate for the number of loops crossing two subsets; our proof uses the random current representation (see Section~\ref{subsec:random path model}) and the McBryan-Spencer argument.

Recall the graph $G=(V,E)$ and the XY measure $\mu_{G,\beta}$. For brevity we will omit the subscript in the XY expectation $ \langle\cdot\rangle_{\mu_{G,\beta}}$ occasionally and put (possible) boundary conditions in the superscript. For $R\geq1$, consider the annulus $ D:=\Lambda_{2R}\setminus\Lambda_R $.
Let $ K_R $ be the number of $\xx$-loops crossing $D$, namely the number of loops intersecting both $\partial\Lambda_{R}$ and $\partial\Lambda_{2R}$. We have the following Gaussian tail estimate for $ K_R $ (note that sending $R\to\infty$ yields a Gaussian decay in $k$ in~\eqref{eq: apriori--Gaussian tail for number of crossing loops}).
\begin{theorem}
    \label{thm: apriori--Gaussian tail for number of crossing loops}
    For all $\beta_0>0$, there exists a constant $ C'>0 $ depending only on $\beta_0$, such that
    \begin{equation}
        \label{eq: M.G.F. of loop crossing number}
        \big\langle (\cosh t)^{K_R}\big\rangle_{\beta}\leq \exp\big( C'R^2\big[\cosh (\frac{t}{R}) -1\big] \big)
    \end{equation}
    for all $ t\geq0 $ and local time condition $ \beta $ with $\max_{v}\beta_v\leq \beta_0$.
    Consequently, there exists a constant $c=c(\beta_0)>0$ such that for all $ k, R\geq1 $
    \begin{equation}
        \label{eq: apriori--Gaussian tail for number of crossing loops}
        \mu_{\beta}[K_R\geq k]\leq c^{-1}\exp\{-c kR\log(1+\frac{k}{R})\}.
    \end{equation}
\end{theorem}
We will prove Theorem~\ref{thm: apriori--Gaussian tail for number of crossing loops} by the following more general result for the number of loops intersecting both of two (arbitrarily chosen) disjoint sets. Suppose $ A,B\subset V $ are two disjoint subsets, and let $ K_{A, B} $ be the number of loops in both copies that visit both $ A $ and $B$.
For every nonnegative $b=(b_v)$, let $ c^b_{uv}:=2J_{uv}\sqrt{b_ub_v} $.
Consider $ h:V\to\mathbb{R} $ as a function satisfying that
\begin{equation}
    \label{eq: def property of test function h}
    h_v\equiv1\quad\forall v\in A,\text{ and}\qquad h_u\equiv0\quad\forall u\in B,
\end{equation}
and in addition for all $ t\geq0 $ define
\begin{equation}
    \label{eq: def-Dirichlet energy in MS argument}
    \mathcal{E}_b(t,h):=\sum_{u\sim v}c^b_{uv}[\cosh\!\big(t(h_u-h_v)\big) -1].
\end{equation}
\begin{theorem}
    \label{thm: Gaussian tail for general loops passing two subsets}
    For all nonnegative $b$ and all $h$ satisfying~\eqref{eq: def property of test function h}, we have
    \begin{equation}
        \label{eq: M.G.F. for general loops passing two subsets}
        \big\langle (\cosh t)^{K_{A, B}} \big\rangle_{b}\leq e^{\mathcal{E}_b(t, h)},\qquad \forall t\geq0.
    \end{equation}
\end{theorem}
\begin{proof}[Proof of Theorem~\ref{thm: apriori--Gaussian tail for number of crossing loops} assuming Theorem~\ref{thm: Gaussian tail for general loops passing two subsets}]
Let $ A=\partial\Lambda_{R} $ and $ B=\partial\Lambda_{2R} $, and consider
\begin{equation*}
    h_*(v):=
    \begin{cases}
        1, & |v|_{\infty}\leq R,\\
        2-\frac{1}{R}|v|_{\infty}, & R<|v|_{\infty}\leq 2R,\\
        0, & |v|_{\infty}\geq 2R.
    \end{cases}
\end{equation*}
Note that $ |h_u-h_v|\leq 1/R $ on every edge $(u,v)$, and there are at most $ 100R^2 $ edges with non-zero gradient (for $h$). Therefore by $ c^{\beta}_{uv}=2J_{uv}\sqrt{\beta_u\beta_v}\leq 2\beta_0 $, the following holds for some $\overline{C}$ depending only on $ \beta_0 $
\[
\mathcal{E}_{\beta}(t,h_*)\leq \overline{C} R^2 (\cosh (t/R)-1).
\]
Combined with~\eqref{eq: M.G.F. for general loops passing two subsets}, this proves~\eqref{eq: M.G.F. of loop crossing number}.
The Gaussian tail bound~\eqref{eq: apriori--Gaussian tail for number of crossing loops} follows from Markov's inequality as follows (by taking $ t=R\,\mathrm{arcsinh}[ k/(\overline{C}R) ]$)
\[
\mu_{\beta}[K_R\geq k]\leq \inf_{t\geq0}\exp\{- k\log(\cosh t) + \overline{C} R^2 (\cosh (t/R)-1) \}.\qedhere
\]
\end{proof}
For further convenience, we state the following corollary, which considers macroscopic loops instead of crossing loops.
\begin{corollary}\label{cor: Gaussian tail for number of macroscopic loops}
    For all $\delta>0$, there exists a constant $c$ depending on $(\beta_0,\delta)$ such that for all $\beta$ with $ \max_{v}\beta_v\leq\beta_0 $, any integer $k\ge0$ and any box $D$ with diameter at least $k$, we have
    \begin{equation}
       \mu_{\beta}\Big(
       \exists\ \text{at least $k$ $\xx\yy$-loops in $D$ with diameter}\,\geq\,\delta\cdot \diam(D)
       \Big)\leq c^{-1}\exp(-c k^2).\label{eq: gauss decay for large diameter loop}
    \end{equation}
\end{corollary}
\begin{proof}
    Write $ L={\rm diam}(D) $, and without loss of generality assume $ \delta<1 $.
    Let $ E_k $ denote the event that there exist at least $k$ $\xx\yy$-loops in $D$ with diameter at least ${\delta\cdot \diam(D)}$.
    
    Choose boxes $ B_1,\dots,B_{m} $  covering $D$ with the same diameter $r$ satisfying $ \tfrac{\delta L}{8}\leq r\leq\tfrac{\delta L}{4} $, where $m=m(\delta)\geq1$ is some integer depending only on $\delta$. Consider the annulus $ A_j = 2B_j\setminus B_j $ for $1\leq j\leq m$, where $ 2B_j $ is the concentric enlargement of $B_j$ with a factor of two. Consequently, every loop in $D$ with diameter at least $ \delta L $ intersects some $B_j$ and exits $ 2B_j $, and thus crossing $A_j$.
    Therefore, if $E_k$ occurs, we can find some $A_j$ being crossed by at least $ \lceil k/(2m) \rceil $ distinct loops in either the $ \xx $- or the $\yy$-copy.
    Theorem~\ref{thm: apriori--Gaussian tail for number of crossing loops} and a union bound then give
    \begin{equation*}
        \mu_{\beta}[E_k]
        \leq 2\overline{C}m\exp\big( -\overline{C}^{-1} \frac{k}{2m} r\log\big(1+\frac{k}{2mr}\big) \big)
        \leq \overline{c}^{-1}\exp(-\overline{c}k^2)
    \end{equation*}
    for some constants $ \overline{c}, \overline{C}>0 $ depending only on $ (\beta_0,\delta) $, concluding the proof.
\end{proof}

Now we give the proof of Theorem~\ref{thm: Gaussian tail for general loops passing two subsets}, via a carefully chosen family of boundary conditions on $ S:=A\cup B $; this allows us to convert the moment generating function for the number of loop segments crossing $A$ and $B$ into ratios of partition functions, which can then be estimated by the McBryan-Spencer argument.

\subsection{Counting the crossing loop segments}
\label{subsec: counting the loop segments}
Let $ e_{\xx}=(1,0) $.
Let $ \Xi^* $ be the collection of boundary conditions $ \xi^* $ such that
\begin{equation}
    \label{eq: def all-x local time boundary condition}
    \xi^*_v\in\{e_{\xx}, -e_{\xx}\} \quad \forall v\in S
\end{equation}
and consider the Gibbs density (with local time condition $b=(b_v)$)
\begin{equation*}
    \frac{1}{Z^{\Xi^*}_{b}}\cdot 2^{-|S|}\sum_{\xi^*\in\Xi^*}\mathbf{1}_{\sigma|_S=\xi^*}\cdot\exp\big( \sum_{u\sim v:\,u,v\in V}c^b_{uv}\sigma_u\cdot\sigma_v \big),
\end{equation*}
where $Z^{\Xi^*}_{b}$ is the normalizing constant, i.e., the partition function (the notation $ \cdot|_S $ means the restriction on $S$).
We will denote the expectation under the preceding density by $ \langle\cdot\rangle^{\Xi^*}_b $ for brevity.
In addition, for $ \varphi\in\mathbb{C} $ we let $ \Xi^{\varphi} $ denote the collection of boundary conditions $ \xi^{\varphi} $ such that
\begin{equation}
    \label{eq: def rotated all-x b.c.}
    \begin{aligned}
        &\xi^{\varphi}_v\in\{(\cos\varphi, \sin\varphi), (-\cos\varphi, -\sin\varphi)\} &&\forall v\in A,\\
        &\xi^{\varphi}_u\in\{e_{\xx}, -e_{\xx}\} &&\forall u\in B,
    \end{aligned}
\end{equation}
and define
\begin{equation}
    \label{eq: def-Z_b(varphi)}
    Z_b(\varphi):=2^{-|S|}\sum_{\xi^{\varphi}\in\Xi^{\varphi}}\int_{\sigma}\mathbf{1}_{\sigma|_S=\xi^{\varphi}}\cdot \exp\big( \sum_{u\sim v}c^b_{uv}\sigma_u\cdot\sigma_v \big) {\mathrm d}\sigma.
\end{equation}
Here ${\mathrm d}\sigma$ means the counting measure for spins on $S$ and the Haar measure on $\mathbb{S}^1$ for spins on $S^c$.
Note that
\begin{equation}
    \label{eq: Z^Xi^*_b=Z_b(0)}
    Z^{\Xi^*}_{b}=Z_b(0).
\end{equation}

Now for each loop, list its visits to $S$ in the cyclic order. Count consecutive visits lying in different boundaries $ A $, $ B $ (so that each pair of such consecutive visits correspond to a segment between $A$ and $B$) and let $M_{A, B}$ be the sum of these segment counts over all loops.
Observe that $ M_{A,B}\geq2K_{A,B} $.
\begin{lemma}
    \label{lem: M.G.F. of the segment number}
    For every $\varphi\in\mathbb{C} $ and for all $ b=(b_v) $, we have
    \begin{equation}
        \label{eq: M.G.F. of the segment number}
        \big\langle (\cos\varphi)^{M_{A, B}} \big\rangle^{\Xi^*}_b = \frac{Z_b(\varphi)}{Z_b(0)}.
    \end{equation}
\end{lemma}
\begin{proof}
    We will use the random current representation (see~\cite{TaggiLRO, vEL23} for instance) of the XY model.
    We begin with the angle representation of the XY model. Let $ \varphi_u=\varphi \cdot \mathbf{1}_{u\in A} $. At each vertex $v\in V$, we write
    \begin{equation*}
        \begin{aligned}
            & \sigma_v=(\sigma^{\xx}_v, \sigma^{\yy}_v) = ( \cos\theta_v, \sin\theta_v ), && v\in V\setminus S;\\
            & \sigma_u=(\sigma^{\xx}_u, \sigma^{\yy}_u) = (\epsilon_u \cos\varphi_u, \epsilon_u\sin\varphi_u ), && u\in S
        \end{aligned}
    \end{equation*}
    where $ \theta_v \in [0,2\pi) $ and $ \epsilon_u \in \{-1,1\} $.
    Note that different from Lemma~\ref{lem: rpm equivalence}, the boundary condition here only has $ 2^{|S|} $ choices of signs (since by~\eqref{eq: def rotated all-x b.c.} the signs assigned to the $\xx$- and $\yy$-coordinates are the same at each $ u\in S $).
    That said, we can still apply the same Taylor expansion as in Section~\ref{subsec:random path model} and integrate the spins (on $ S^c $) and signs (on $S$) out. Namely, letting $ {\mathrm d}\theta $ stand for product of Haar measures on $ \mathbb{S}^1 $ and ${\mathrm d}\epsilon$ stand for uniform measures on $ \{-1,1\}^S $, we have (recall~\eqref{eq: def random current weight} and~\eqref{eq: Taylor expansion of Z} as well as notations introduced in between)
    \begin{equation}
        \label{eq: random current of xy partition function-1}
        \begin{aligned}
            Z_b(\varphi)
            &= \int\prod_{u\sim v}\exp(c^b_{uv}\sigma^{\xx}_u\sigma^{\xx}_v)\exp(c^b_{uv}\sigma^{\yy}_u\sigma^{\yy}_v) {\mathrm d}\theta {\mathrm d}\epsilon\\
            &= \sum_{n^{\xx}, n^{\yy}} w_b(n^{\xx})w_b(n^{\yy})\cdot\Big(\prod_{v\in S^c}\int(\sigma^{\xx}_v)^{d^{\xx}_v}(\sigma^{\yy}_v)^{d^{\yy}_v}{\mathrm d}\theta_v\Big)\cdot\Big( \prod_{u\in S}\int(\sigma^{\xx}_u)^{d^{\xx}_u}(\sigma^{\yy}_u)^{d^{\yy}_u}{\mathrm d}\epsilon_u \Big).
        \end{aligned}
    \end{equation}
    The integration at each $ v\in S^c $ remains the same as in Section~\ref{subsec:random path model}, while at each $u\in S$, we have (note that $ \varphi_u=0 $ for $u\in B$ by definition)
    \[
    \begin{aligned}
        &\quad\ \frac{1}{2}\sum_{\epsilon_u\in\{-1,1\}}(\epsilon_u\cos\varphi_u)^{d^{\xx}_u}(\epsilon_u\sin\varphi_u)^{d^{\yy}_u}
        = (\cos\varphi_u)^{d^{\xx}_u}(\sin\varphi_u)^{d^{\yy}_u}\mathbf{1}_{2\mid (d^{\xx}_u+d^{\yy}_u)}\\
        &=(\cos\varphi)^{d^{\xx}_u}(\sin\varphi)^{d^{\yy}_u}\mathbf{1}_{2\mid (d^{\xx}_u+d^{\yy}_u)}\mathbf{1}_{u\in A}+\mathbf{1}_{2|d^{\xx}_u}\mathbf{1}_{d^{\yy}_u=0}\mathbf{1}_{u\in B}.
    \end{aligned}
    \]
    Therefore only currents in
    \begin{equation}
        \label{eq: surviving currents}
        \mathcal{N}:=\big\{(n^{\xx}, n^{\yy})\in\mathbb{Z}^{E}_{\geq0}\times\mathbb{Z}^{E}_{\geq0}:\partial n^{\xx}=\partial n^{\yy}\subset A,\text{ and } d^{\yy}_u=0 \text{ for all } u\in B\big\}
    \end{equation}
    makes non-zero contribution to~\eqref{eq: random current of xy partition function-1}, and thus we have
    \begin{equation}
        \label{eq: random current of xy partition function-2}
        Z_b(\varphi)
        = \sum_{(n^{\xx}, n^{\yy})\in\mathcal{N}}
        w_b(n^{\xx})w_b(n^{\yy})\times \prod_{v\in S^c}\frac{(d^{\xx}_v-1)!!(d^{\yy}_v-1)!!}{2^{m_v} m_v!}\times\prod_{u\in A}(\cos\varphi_u)^{d^{\xx}_u}(\sin\varphi_u)^{d^{\yy}_u}.
    \end{equation}

    Similar to Section~\ref{subsec:random path model}, we can still pair the replicated edges induced by the random current representation~\eqref{eq: random current of xy partition function-2} at each vertex. However we now introduce the pairing on the combined current (and its induced replicated edges) for further convenience.
    Precisely, letting $ n=n^{\xx}+n^{\yy} $, we have $ \partial n=\varnothing $ and $ m_v=d_v(n)/2 $. For each edge $ e\in E $, consider $n_e$ replicas of it (we call them replicated edges) labeled by $ 1,\dots,n_e $. Assign each labeled replicated edge to $\xx$- or $\yy$-copy, subject to~\eqref{eq: surviving currents}.
    Let $ \mathcal{N}(n) $ denote the collection of such assignments satisfying that $ (n^{\xx}, n^{\yy})\in\mathcal{N} $. Observe that $ \mathcal{N}(n)\neq\varnothing $ if and only if $ \partial n=\varnothing $.
    The number of choices producing the prescribed assignment $ n_e=n^{\xx}_e+n^{\yy}_e $ is $ {n_e\choose n^{\xx}_e} $, and note that
    \[
    \frac{(c^b_e)^{n^{\xx}_e}}{n^{\xx}_e!}\frac{(c^b_e)^{n^{\yy}_e}}{n^{\yy}_e!} = \frac{(c^b_e)^{n_e}}{n_e!}{n_e\choose n^{\xx}_e}.
    \]
    Therefore we have
    \begin{equation}
        \label{eq: resummation of copywise current weights}
        \sum_{(n^{\xx},n^{\yy})\in \mathcal{N}}w_b(n^{\xx})w_b(n^{\yy})F(n^{\xx},n^{\yy})=\sum_{\partial n=\varnothing}w_b(n)\sum_{(n^{\xx}, n^{\yy})\in\mathcal{N}(n)}F(n^{\xx},n^{\yy})
    \end{equation}
    for all measurable function $F$ of $ (n^{\xx},n^{\yy}) $.
    At each vertex $v\in S^c$, we still use the factor $ (d^{\xx}_v-1)!!(d^{\yy}_v-1)!! $ to count the number of pairings of its adjacent edges within each ($\xx$ or $\yy$) copy.
    At each boundary vertex $ u\in S $, however, we introduce a uniform pairing of its adjacent replicated edges regardless of their copy assignments (namely, summing over all $ (2m_u-1)!! $ possible pairings and then dividing by $ (2m_u-1)!! $). Let $ \pi^{\xx\yy}\in \Pi(n^{\xx},n^{\yy}) $ be the product of the aforementioned pairings at all vertices, and define
    \begin{equation}
        \label{eq: def W_b(n)-Sec4}
        W_b(n):=w_b(n)\prod_{v\in S^c}\frac{1}{2^{m_v}m_v!}\prod_{u\in S}\frac{1}{(2m_u-1)!!}
    \end{equation}
    for current $n$.
    Then plugging
    \[
    F(n^{\xx}, n^{\yy}):=\sum_{\pi^{\xx\yy}\in \Pi(n^{\xx}, n^{\yy})}\prod_{v\in S^c}\frac{1}{2^{m_v}m_v!}\prod_{u\in S}\frac{1}{(2m_u-1)!!}\times\prod_{u\in A}(\cos\varphi_u)^{d^{\xx}_u}(\sin\varphi_u)^{d^{\yy}_u},
    \]
    into~\eqref{eq: resummation of copywise current weights}, we have (noting that summation over $ \pi^{\xx\yy} $ in $F$ contributes the double factorials in~\eqref{eq: random current of xy partition function-2})
    \begin{equation}
        \label{eq: Z_b(varphi)}
        \begin{aligned}
            Z_b(\varphi)
            &\stackrel{\eqref{eq: random current of xy partition function-2}}{=}
            \sum_{(n^{\xx}, n^{\yy})\in\mathcal{N}}w_b(n^{\xx})w_b(n^{\yy})F(n^{\xx}, n^{\yy}) = \sum_{\partial n=\varnothing}w_b(n)\sum_{(n^{\xx}, n^{\yy})\in\mathcal{N}(n)}F(n^{\xx}, n^{\yy})
            \\&\ \ = \sum_{\partial n=\varnothing}w_b(n) \prod_{v\in S^c}\frac{1}{2^{m_v}m_v!}\prod_{u\in S}\frac{1}{(2m_u-1)!!}  \sum_{\substack{(n^{\xx}, n^{\yy})\in\mathcal{N}(n),\\ \pi^{\xx\yy}\in \Pi(n^{\xx}, n^{\yy})}} \prod_{u\in A}(\cos\varphi_u)^{d^{\xx}_u}(\sin\varphi_u)^{d^{\yy}_u}
            \\&\stackrel{\eqref{eq: def W_b(n)-Sec4}}{=} \sum_{\partial n=\varnothing} W_b(n) \sum_{\substack{(n^{\xx}, n^{\yy})\in\mathcal{N}(n),\\ \pi^{\xx\yy}\in \Pi(n^{\xx}, n^{\yy})}} \prod_{u\in A}(\cos\varphi_u)^{d^{\xx}_u}(\sin\varphi_u)^{d^{\yy}_u}.
        \end{aligned}
    \end{equation}

    Now we analyze the preceding summation and convert it into a summation over $ \pi\in\Pi(n) $, where $ \Pi(n) $ is the collection of all possible pairings of replicated edge induced by $ n $ regardless of the copies, and then followed by a summation over all possible ($\xx$ or $\yy$) copy assignments.
    First note that each pairing $\pi$ induces a loop configuration $ \mathfrak{L}(n, \pi) $ where each loop consists of replicated edges paired at each vertex. For each configuration $ \mathfrak{L}(n, \pi) $, let $ K_0(n,\pi) $ count the number of loops avoiding $S$ under the pairing $\pi$. For each loop $\ell\in\mathfrak{L}(n,\pi)$ intersecting $S$, we list its successive visits to $S$ as $ u_1, u_2, \dots, u_{q} $ in the cyclic order. Then we know that the loop segment between each $ (u_j, u_{j+1}) $ is an excursion from $S$ to $S$, and we let $ M(\ell) $ be the number of these excursions ending in different sets $A$ and $B$ (we use the convention $ M(\ell)=0 $ when $\ell\cap S=\varnothing$). Letting $ M(n,\pi):=\sum_{\ell\in\mathfrak{L}(n,\pi)}M(\ell) $, we have $ M(n,\pi)=M_{A, B}(\mathfrak{L}(n,\pi)) $.
    We now prove that for every $ n $ with $ \partial n=\varnothing $,
    \begin{equation}
        \label{eq: sum over copy assignments}
        \sum_{\substack{(n^{\xx}, n^{\yy})\in\mathcal{N}(n),\\ \pi^{\xx\yy}\in\Pi(n^{\xx}, n^{\yy})}}\prod_{u\in A}(\cos\varphi_u)^{d^{\xx}_u}(\sin\varphi_u)^{d^{\yy}_u} = \sum_{\pi\in\Pi(n)}2^{K_0(n,\pi)}(\cos\varphi)^{M(n,\pi)}.
    \end{equation}
    In order to prove the preceding identity, we consider loops generated by a given pairing $ \pi\in\Pi(n) $. We then keep track of the legitimate copy assignments for each such  loop, which depends on the number of visits to $S$.
    For a loop avoiding $S$, its pairings are at vertices in $S^c$, forcing the same copy along the loop, which gives a factor of two. For a loop visiting $S$, recall its successive visits $ \{u_j\}_{1\leq j\leq q} $ to $S$.
    Between each two successive visits $ u_j $ and $ u_{j+1} $, the pairings at vertices in $ S^c $ again forces the same copy; summing over the two possible copies gives a factor
    \[
    \cos \varphi_{u_j}\cos \varphi_{u_{j+1}} + \sin \varphi_{u_j}\sin \varphi_{u_{j+1}} = \cos(\varphi_{u_j} - \varphi_{u_{j+1}}).
    \]
    Since the boundary factor $ \prod_{u\in A}(\cos\varphi_u)^{d^{\xx}_u}(\sin\varphi_u)^{d^{\yy}_u} $ is just the product of corresponding (with respect to copy $\xx$ or $\yy$) factor ($\cos\varphi_u$ or $\sin\varphi_u$ respectively) over all the adjacent replicated edges and all boundary vertices $u$ (there are $ d^{\xx}_u+d^{\yy}_u $ adjacent replicated edges at each $ u $, and each of them belongs to an excursion from $S$ to $S$), we can decompose the boundary factor by each loop visiting $S$ and each of its aforementioned excursion $ (u_j, u_{j+1}) $ such that each excursion carries its own factors $ \cos\varphi_{u_j}\cos\varphi_{u_{j+1}} $ (if it is assigned to $\xx$-copy) or $ \sin\varphi_{u_j}\sin\varphi_{u_{j+1}} $ (if assigned to $\yy$-copy).
    Therefore, the boundary factor is in fact equivalent to summing the factors of each excursion of each loop over their copy assignments.
    Consequently, each loop $\ell$ visiting $S$ contributes a factor
    \[
    \prod^q_{j=1} \cos(\varphi_{u_j} - \varphi_{u_{j+1}}) = (\cos \varphi)^{M(\ell)}
    \]
    to the whole boundary factor $ \prod_{u\in A}(\cos\varphi_u)^{d^{\xx}_u}(\sin\varphi_u)^{d^{\yy}_u} $.
    Note that $ M(\ell) $ does not necessarily equal $q$, since there might be excursions with both endpoints in $ A $ (or $B$); in addition, each of these excursions only contributes a factor of 1. Therefore~\eqref{eq: sum over copy assignments} follows from the preceding argument, by summing over $\pi$ (whereas the copy assignment then corresponds to a summation over $ (n^{\xx}, n^{\yy}) $ and $ \pi^{\xx\yy} $ on the left). Combined with~\eqref{eq: Z_b(varphi)}, it yields that
    \begin{equation}
        \label{eq: random current of xy partition function-3}
        Z_b(\varphi)=\sum_{\partial n=\varnothing}\sum_{\pi\in\Pi(n)}W_b(n)2^{K_0(n,\pi)}(\cos\varphi)^{M(n,\pi)}.
    \end{equation}
    The ratio then states (noting that $ M(n,\pi)=M_{A,B} $, and plugging in $ \varphi=0 $)
    \[
    \frac{Z_b(\varphi)}{Z_b(0)} = \frac{\sum_{\partial n=\varnothing}\sum_{\pi\in\Pi(n)}W_b(n)2^{K_0(n,\pi)}(\cos\varphi)^{M(n,\pi)}}{\sum_{\partial n=\varnothing}\sum_{\pi\in\Pi(n)}W_b(n)2^{K_0(n,\pi)}}
    =\big\langle(\cos\varphi)^{M_{A,B}}  \big\rangle^{\Xi^*}_b
    \]
    where the last equality is derived by noticing that $ (\cos\varphi)^{M(n,\pi)} $ is a measurable function of the loop configuration $ \mathfrak{L}(n,\pi) $, and by~\eqref{eq: Z_b(varphi)} and~\eqref{eq: sum over copy assignments} its weight under $ \mu^{\Xi^*}_b $ is exactly $W_b(n)2^{K_0(n,\pi)}$ (recall by~\eqref{eq: Z^Xi^*_b=Z_b(0)} the partition function defining $\mu^{\Xi^*}_b$ is exactly $ Z^{\Xi^*}_b=Z_b(0) $).
    This then concludes the proof by recalling Lemma~\ref{lem: rpm equivalence}.
\end{proof}

\subsection{The McBryan-Spencer argument}
\label{subsec: MS argument}
Provided with Lemma~\ref{lem: M.G.F. of the segment number}, we can control the moment generating function of $M_{A, B}$ by bounding the partition function ratio $ \frac{Z_b(\varphi)}{Z_b(0)} $, which is via the McBryan-Spencer argument.
\begin{lemma}
    \label{lem: McBryan-Spencer bound}
    Recall the definition of $ M_{A, B} $ right above Lemma~\ref{lem: M.G.F. of the segment number}.
    For all $h$ satisfying~\eqref{eq: def property of test function h}, all $t\geq0$ and all nonnegative $b$, we have
    \begin{equation}
        \label{eq: McBryan-Spencer bound}
        \big\langle (\cosh t)^{M_{A, B}}\big \rangle^{\Xi^*}_b
        \leq e^{\mathcal{E}_b(t,h)}.
    \end{equation}
\end{lemma}
\begin{proof}
    Letting $ \varphi=it $, we derive from Lemma~\ref{lem: M.G.F. of the segment number} that
    \begin{equation}
        \label{eq: MB bound-1}
        \big\langle (\cosh t)^{M_{A, B}}\big \rangle^{\Xi^*}_b
        =\frac{Z_b(it)}{Z_b(0)}.
    \end{equation}
    Now we estimate the ratio in~\eqref{eq: MB bound-1}.
    By definition~\eqref{eq: def-Z_b(varphi)}, we can write (note that the rotation $\varphi$ is defined only on $A$)
    \begin{equation}
        \label{eq: integration of partition funciton-1}
        Z_b(\varphi)=2^{-|S|}\sum_{\varepsilon\in\{-1,1\}^S}\int_{[0,2\pi]^{S^c} }\exp\big\{ \sum_{u\sim v} c^b_{uv}\cos\big(\theta_u-\theta_v +\varphi \cdot[\mathbf{1}_{A}(u)-\mathbf{1}_{A}(v)] \big) \big\}{\mathrm d}\theta,
    \end{equation}
    where we used the angle representation $ \sigma_u=(\cos\theta_u,\sin\theta_u) $ of the XY model, and in particular on $S$, we use  the following correspondence
    \[
    \theta_v=\begin{cases}
        0,&\varepsilon_v=1,\\
        \pi,&\varepsilon_v=-1.
    \end{cases}
    \]
    Now let $h$ be a function satisfying~\eqref{eq: def property of test function h}.
    By Cauchy's theorem, we can translate the integration path of $ Z_b(\varphi) $ with respect to $ \theta_v $ for each $ v\in S^c $ by $ ith_v $, using the rectangle with vertices $ 0,2\pi,2\pi +ith_v, ith_v $.
    Plugging in $ \varphi=it $, this yields
    \begin{equation}
        \label{eq: integration of partition funciton-2}
        Z_b(it)=2^{-|S|}\sum_{\varepsilon\in\{-1,1\}^S}\int_{[0,2\pi]^{S^c} }\exp\big\{ \sum_{u\sim v} c^b_{uv}\cos\big(\theta_u-\theta_v +it(h_u-h_v) \big) \big\}{\mathrm d}\theta.
    \end{equation}
    Since $ c^b_{uv}\geq0 $, by
    \[
    {\rm Re}\cos (\theta +i a)\leq \cos\theta +\cosh a -1,
    \]
    we have (recalling~\eqref{eq: def-Dirichlet energy in MS argument})
    \begin{align*}
        \big| \exp\big\{ \sum_{u\sim v} c^b_{uv}\cos\big(\theta_u-\theta_v +it(h_u-h_v) \big) \big\} \big|
        &\leq \exp\big\{{\rm Re}\big[ \sum_{u\sim v} c^b_{uv}\cos\big(\theta_u-\theta_v +it(h_u-h_v) \big) \big]\big\}\\
        &\leq \exp\big\{ \sum_{u\sim v} c^b_{uv}\cos(\theta_u-\theta_v)  +  \mathcal{E}_b(t,h) \big\}.
    \end{align*}
    Combined with~\eqref{eq: integration of partition funciton-2}, it gives that $ |Z_b(it)|\leq e^{\mathcal{E}_b(t,h)} Z_b(0) $.
    Combined with~\eqref{eq: MB bound-1} it then yields that
    \[
    \big\langle (\cosh t)^{M_{A, B}}\big \rangle^{\Xi^*}_b
    \leq e^{\mathcal{E}_b(t,h)},
    \]
    which concludes the proof.
\end{proof}
Now we can handle the free boundary condition via exploration of loops and prove Theorem~\ref{thm: Gaussian tail for general loops passing two subsets}.
\begin{proof}[Proof of Theorem~\ref{thm: Gaussian tail for general loops passing two subsets}]
    Let $ K^{\xx}_{A, B} $ and $K^{\yy}_{A, B}$ count the number of crossing loops in $\xx$ and $\yy$-copy respectively.
    Reveal the whole collection of $\yy$-loops that intersect $S$ (including point loops)
    \[
    \eta^{\yy}:=\{ \ell\in\mathfrak{L}^{\yy}: \ell\cap S\neq\varnothing \}
    \]
    and let
    \[
    b_v:=\beta_v - L_v(\eta^{\yy}):=\beta_v -\sum_{\ell\in\eta^{\yy}}L_v(\ell)
    \]
    be the local time condition for remaining loops.
    Note that without the local time conditioning, the Poisson restriction yields independence of $\eta^{\yy}$, the (collection of) $\yy$-loops avoiding $S$, and all the $\xx$-loops.  Therefore conditional on $\eta^{\yy}$ and on the local time condition $\beta$, the two remaining soups have exactly the XY model law with local time condition $b$ and an all-$\xx$ local time boundary condition $ \Xi^* $ (recall~\eqref{eq: def all-x local time boundary condition}). Consequently, every crossing loop in this new XY model law is in the $\xx$-copy, meanwhile every original $\xx$-loop is retained as well.
    Since $0\leq b\leq\beta$, by monotonicity of $ \mathcal{E}_b $ with respect to $b$ we have (recall~\eqref{eq: def-Dirichlet energy in MS argument})
    \begin{equation}
        \label{eq: M.G.F. of crossing loops under free b.c.}
        \big\langle(\cosh t)^{K_{A, B}}\big\rangle_{\beta}\leq \big( \big\langle(\cosh t)^{2K^{\xx}_{A, B}}\big\rangle_{\beta}\big\langle(\cosh t)^{2K^{\yy}_{A, B}}\big\rangle_{\beta} \big)^{1/2}\leq e^{\mathcal{E}_{\beta}(t,h)}
    \end{equation}
    for all $ h $ satisfying~\eqref{eq: def property of test function h}, where the first inequality is by $ K_{A, B}=K^{\xx}_{A, B}+K^{\yy}_{A, B} $ and the Cauchy-Schwarz inequality, and the second inequality is by Lemma~\ref{lem: McBryan-Spencer bound} and the fact that $ M_{A,B}\geq 2K_{A, B} $ and $ \cosh t\geq1 $ for all $t\geq0$.
\end{proof}

\newcommand{\mcc}{\mathcal}
\section{Spin wave empowered by perturbation from \texorpdfstring{$\zz$}{}-copy}
\label{sec:fractality}
In this section, we will fix $\NN=2+\epsilon$ for some $\epsilon\in (0,1]$ and omit it from the superscript. Recall that we have fixed $\beta_0>0$ and fixed $\{\beta_u\}$ to be a local time condition such that $\beta_u\le \beta_0$ for all vertex $u$.

We first prove the following slightly stronger version of Theorem  \ref{thm-1.2} using Theorem \ref{thm-annulus-local-time-crossing} below with a standard coarse-graining technique. Recall the notations after \eqref{eq: local time definition}.
\begin{theorem}
    \label{thm-1.2'}
    There exists a constant $ c_1>0 $ depending only on $\beta_0, \NN$ such that for all $L\geq1$, $u,v\in B_{L/2}$ and any $\xx\yy$ boundary condition $\xi$ on $\partial\widetilde{B}_L$, 
    \begin{equation}
        \label{eq:thm1.2'}
        \nu^{\xi}_{{B}_L,\beta}\Big(\big\{ u\con{\mathsf{xyz}} v  
        \big\}\Big) \leq c_1^{-1}e^{-c_1|u-v|_{\infty}}.
    \end{equation}
\end{theorem}

\begin{theorem}\label{thm-annulus-local-time-crossing}
     For any $c_2>0$, there exists a constant $M_0>0$ depending only on $\beta_0, \NN, c_2$ such that for any $M\ge M_0$ and any $\xx\yy$ boundary condition $\xi$ on $\partial\Lambda_{M,2M}$ (recall the end of Section~\ref{sec: local time boundary condition}) we have 
    \begin{equation}
        \nu_{\Lambda_{M,2M},\beta}^{\xi} (\partial\Lambda_{M}\con{\mathsf{xy}}\partial\Lambda_{2M})\le c_2.\label{eq: 5.2}
    \end{equation}
\end{theorem}

\begin{proof}[Proof of Theorem~\ref{thm-1.2'} assuming Theorem~\ref{thm-annulus-local-time-crossing}]
    All connections below are of edges in the domain under consideration, without using the additional connectivity in the boundary condition.
    Consider the event
    \[
    \mathcal C(u,v) :=\bigcup_{\mathsf{i}\in\{\mathsf{x,y}\}} \{u\stackrel{\mathsf{i}}{\longleftrightarrow} v\}.
    \] Since under $\xx\yy$ boundary condition the $\zz$-copy is dominated by the $\xx,\yy$-copy, we get that \begin{equation}
        \nu^{\xi}_{B_L,\beta}\Big(\big\{u\con{\mathsf{xyz}}v\big\}\Big)\le 2\nu^{\xi}_{B_L,\beta}\big(\mathcal C(u,v)\big).
    \end{equation}

    Let $\kappa=10^{-100}$ and choose an integer $M\geq\max\{1,M_0(\kappa)\}$ as in Theorem~\ref{thm-annulus-local-time-crossing}. Put
    \[
    \Lambda_t(w)=\{u\in\mathbb Z^2:|u-w|_\infty\leq t\},
    \qquad
    Q_v=Mv+\{0,\ldots,M-1\}^2
    \]
    for $ w\in\mathbb Z^2 $.
    Let $\mathcal C_v$ be the event that for either $\xx$- or $\mathsf  y$-copy there exists a local time cluster crossing through the annulus $ A_w:=\Lambda_{2M}(Mw)\setminus\Lambda_{M-1}(Mw) $.
    
    Let $\mathcal F_w^{\mathrm{ext}}$ be the $\sigma$-field generated by the $\xx\yy$ boundary condition $\xi_w$ for $A_w$ and the local time $\beta^\zz$ for $\zz$-copy  local times from loops intersecting $A_w^c$, where $\xi_w$ encodes the $\xx\yy$ local times on vertices
    outside the interior of $A_w$ (including those on its boundary) and (local time) connectivities outside $A_w$.
    Thus Theorem~\ref{thm-annulus-local-time-crossing} gives
    \begin{equation}\label{eq: dmp for three copies}
        \nu^{\xi}_{B_L,\beta}(\mathcal{C}_w\mid\mathcal F_w^{\mathrm{ext}}) =\nu_{A_w,\beta-\beta^\zz}^{\xi_w}(\mathcal{C}_w) \leq\kappa
    \end{equation}
    almost surely.
    For any deterministic $w_1,\ldots,w_s$ whose enlarged boxes $\Lambda_{2M}(Mw_j)$ are pairwise disjoint,  since for $j\neq i$, $\mathcal C_{w_j}$ is measurable with respect to the edge configuration inside $A_{w_j}$, we get  $\mathcal C_{w_j} $ is measurable in $ \mathcal F_i^{\mathrm{ext}}$. Therefore, we get from \eqref{eq: dmp for three copies} that \begin{equation}
        \nu^{\xi}_{\Lambda_{M,2M},\beta}(\mathcal{C}_{w_i}\mid \bigcap_{j=1}^{i-1}\mathcal{C}_{w_j}) \leq\kappa.\nonumber
    \end{equation}
    Iterating the preceding inequality gives that
    \begin{equation}
        \label{eq: lattice animal energy}
        \nu^{\xi}_{B_L,\beta}\Bigl(\bigcap_{j=1}^s\mathcal C_{w_j}\Bigr) \leq\kappa^s.
    \end{equation}

    Let $R=|u-v|_\infty$. Assume first that $R\geq16M$, and let
    \[
    r=\lfloor R/2\rfloor,\qquad n=\lfloor R/(4M)\rfloor.
    \]
    On $\mathcal C(u,v)$, choose a self-avoiding (nearest-neighbor) path $ \gamma $ from $u$ to $v$, and stop it when hitting a vertex $t$ with $|t-u|_\infty=r$.
    Then this path lies in $\Lambda_r(u)$.
    Let $\Gamma$ be the collection of $w$ such that $Q_w$ intersects $\gamma$. Then $\Gamma$ forms a lattice animal, i.e., a connected subset of $\mathbb Z^2$ where two vertices are neighboring each other if their $|\cdot|_\infty$-distance is at most 1. Since $\max_{w\in \Gamma}|w-u|_\infty \geq r/M-1$, we see that $m = |\Gamma| \geq n$. We next bound the number of lattice animals containing $u$ (clearly $u\in \Gamma$) of size $m$. To this end, it suffices to bound the number of spanning trees rooted at $u$ of size $m$, which in turn can be encoded by a contour rooted at $u$ of length $2m$ (the contour corresponds to the depth-first search of the tree). Thus, the number of such lattice animals is at most $8^{2m}$.

    Next, we analyze the probability for the crossing event. To this end, we divide $\Gamma$ into 25 subsets according to
    \[
    v\bmod5\in\{0,1,2,3,4\}^2.
    \]
    One of the $25$ subsets contains at least $m/25$ vertices. Distinct vertices in this class have $|\cdot|_\infty$-distance at least $5$, and thus their enlarged boxes of radius $2M$ are disjoint. Therefore, by~\eqref{eq: lattice animal energy}
    \[
    \nu^{\xi}_{B_L,\beta}\Bigl(\bigcap_{w\in\Gamma} \mathcal{C}_w\Bigr)
    \leq\kappa^{m/25}.
    \]
    Since the number of the aforementioned lattice animals is at most $8^{2m}$, taking a union bound over them yields
    \begin{equation*}
        \nu^{\xi}_{B_L,\beta}\big(\mathcal C(u       ,v)\big) \leq \sum^{\infty}_{m=n}8^{2m}\kappa^{m/25}\leq \bar{c}_1^{-1}\exp(-\bar{c}_1|u-v|_{\infty})
    \end{equation*}
    for some constant $ \bar{c}_1>0 $ depending on $M_0$ (which in turn depends on $\beta_0, \mathsf N$).
\end{proof}
\begin{proof}[Proof of Theorem \ref{thm-1.1}]
    The proof of Theorem \ref{thm-1.1} is similar to that of Theorem \ref{thm-1.2}. The only difference is that in Theorem \ref{thm-1.1}, the measure has a spin boundary condition. Recall \eqref{eq: relation between spin boundary condition and local time boundary condition}. For the spin boundary condition, at least one of $V_\xx^+$ and $V_\xx^-$ is the empty set, and thus $V^+_\mathsf x \ncon[0.5]{} V^-_{\mathsf x}$ holds automatically. The same is true for the $\yy$- and $\mathsf{z}$-copy. 
    Thus, when considering local time configurations, the spin boundary condition is equivalent to the local time boundary condition $\bar{\mathsf{b}}$, where $\bar{\mathsf{b}}$ is obtained from $\mathsf{b}$ by $\bar{ \mathsf{b}}_v^{\mathsf{i}} = \beta (\mathsf{b}_v^{\mathsf{i}})^2$ for $v\in \partial B_L$ and $\mathsf i\in \{\xx, \yy, \mathsf z\}$ and wiring all vertices together for all copies  (recall \eqref{eq: special case for spin bc to local time bc}). Then Theorem \ref{thm-1.1} follows from \eqref{eq: 5.2} by a simple union bound.
\end{proof}

The following theorem encapsulates all the difficulties for the proof of Theorem~\ref{thm-annulus-local-time-crossing}.
\begin{theorem}\label{thm: thin annulus local time crossing}
There exists a constant $c_3>0$ depending only on $\beta_0, \NN$ such that the following holds. For any $\iota>0$, there exists $M_0$ depending on $(\beta_0,\NN, c_3, \iota)$
 such that for any $M\ge M_0$ and any $\xx\yy$ boundary condition $\xi$ on $\partial\Lambda_{M,2M}$ we have 
    \begin{equation}\label{eq: 5.5}
        \nu_{\Lambda_{M,(1+\iota)M},\beta}^{\xi} (\partial\Lambda_{M}\con{\mathsf{xy}}\partial\Lambda_{(1+\iota)M})\le 1-c_3.
    \end{equation}
\end{theorem}
\begin{proof}[Proof of Theorem \ref{thm-annulus-local-time-crossing} assuming Theorem \ref{thm: thin annulus local time crossing}]
    Choose $\bar{C}_1>0$ large enough and let $\iota=2^{1/\bC_1}-1.$ Let $\Lambda_i=\Lambda_{(1+\iota)^{i-1}M,(1+\iota)^{i}M}$, and let $\mathcal{A}_i$ denote the event that $\partial\Lambda_{(1+\iota)^{i-1}M}\con{\mathsf{x}}\partial\Lambda_{(1+\iota)^{i}M}$.  
 Let $\mathcal F_i^{\mathrm{ext}}$ be the $\sigma$-field generated by the $\xx\yy$ boundary condition $\xi_i$ for $\Lambda_i$ and the local time $\beta_i^\zz$ for $\zz$-copy local times from loops intersecting $\Lambda_i^c$, where $\xi_i$ encodes  $\xx\yy$ local times on vertices
    outside the interior of $\Lambda_i$ (including those on its boundary) and (local time) connectivities outside $\Lambda_i$. 
    Thus, Theorem~\ref{thm: thin annulus local time crossing} gives
    \begin{equation}\label{eq: dmp for three copies 2}
        \nu^{\xi}_{\Lambda_{M,2M},\beta}(\mathcal{A}_i\mid\mathcal F_i^{\mathrm{ext}}) =\nu_{\Lambda_i,\beta-\beta^\zz_i}^{\xi_i}(\mathcal{A}_i) \leq1- c_3.
    \end{equation}
 Note that for $j\neq i$, $\mathcal{A}_j$ is measurable with respect to the edge configuration inside $\Lambda_{j}$ and thus $\mathcal{A}_j$ is measurable in $ \mathcal F_i^{\mathrm{ext}}$. Therefore, we get that \begin{equation}
        \nu^{\xi}_{\Lambda_{M,2M},\beta}(\mathcal{A}_i\mid \bigcap_{j=1}^{i-1}\mathcal{A}_j) \leq1- c_3.\nonumber
    \end{equation}
Iterating the preceding inequality gives that \begin{align}
            \nu_{\Lambda_{M,2M},\beta}^\xi (\partial\Lambda_{M}\con{\mathsf{x}}\partial\Lambda_{2M})
            \le (1-c_3)^{\bC}.\nonumber
        \end{align} 
    Choosing $\bC>0$ large enough such that $(1-c_3)^{\bC}\le c_2$ completes the proof of Theorem~\ref{thm-annulus-local-time-crossing}.
\end{proof}

We now discuss the main ideas underlying the proof of Theorem~\ref{thm: thin annulus local time crossing}. Our starting point is the classic spin wave argument, which gives a polynomial decay of the spin correlation for the XY model. With such an argument, one can prove a version of Theorem~\ref{thm: thin annulus local time crossing} with $c_3\to0$ as $\iota\to0$ or $\beta\to\infty$ (see Corollary \ref{cor:fkg-boundary-domination}). The key innovation in Theorem~\ref{thm: thin annulus local time crossing} is that $c_3$ is uniform for all $\iota$.

To achieve this, we exploit the fact that our XY model has random local times which arise from subtracting the local times of the $\mathsf{z}$-copy. Ideally, if in most of the region the $\mathsf{z}$-copy takes all the $\beta$ local times (and thus the XY model has 0 local time), then the original spin wave argument would have energy cost only in the region with non-zero local time and thus be much more efficient.

However, the above condition is too ideal to hold. Indeed, it was proved in~\cite{AGS25} that in the XY model the vertices with large local times percolate in a strong fashion. This forces us to work with a much weaker condition on the environment. As already mentioned in the introduction, one reasonable hope is that in a ``good'' box the probability of having a local time crossing in either the $\mathsf{x}$- or the $\mathsf{y}$-copy is very small. It occurs to us that the best hope is to show that a box is good with some small but fixed probability (depending on $\beta$). This will be incorporated in Lemma~\ref{lem: good box probability}.
Next, we focus the discussions on how to proceed, provided with Lemma~\ref{lem: good box probability}. We first aim to show that the union of the good boxes occupies almost all the volume, despite the fact that each has a small probability. This is possible since we may consider good boxes in many scales, and by controlling the correlations, we arrive at our desired volume estimate as in Lemma~\ref{lem: covering property for good boxes}.

The employment of the spin wave argument is nevertheless highly nontrivial. In an ideal scenario, we assume that in each good box there is no local time crossing in either the $\xx$- nor the $\yy$-copy. Then, we can sample all the loops intersecting the good boxes, and conditioning on this, we arrive at a new XY model, where the local times become 0 on the good boxes. 

The main issue in the above picture is that we only know that most good boxes have no crossings, but not all. In order to fix this, we may a priori acquire the information on which good boxes have crossings. Naturally, we wish to sample the loops intersecting these boxes and apply the spin wave argument on the new XY model. This then leads to another issue on whether such ``no-crossing'' property can be preserved by our spin wave argument---it cannot in general. In order to address this, we have to classify the boxes into two types, as in Definition~\ref{def: base spin wave}, whereas Type I preserves the no-crossing property after the spin wave and Type II does not. We then show in Lemma~\ref{lem: type probability} that Type II occurs only with small probability. 

But this then raises another issue on the multiplicity of the spin wave mapping, or in other words, whether one can recover the type from the image of the mapping. This is still non-trivial, and we have to refine the definition of types to this end. Accordingly, when defining our spin wave mapping, we have to modify the classic map by incorporating the differences resulted from different types (see Definitions~\ref{def: type 2 coupling} and~\ref{def: final spin wave global part}), for which the control of the energy costs is also different (see Lemma~\ref{lem: type 2 coupling} and Proposition~\ref{prop: symmetric spin wave}).

\subsection{Good boxes and types}
\label{subsec: Good boxes and types}
In this section, we introduce the definitions that are crucial in the proof of Theorem \ref{thm: thin annulus local time crossing}.  We remind the reader that conditioned on the collection of all loops (without copy labels), each loop is assigned to the $\mathsf z$-copy independently with probability $\frac{\epsilon}{\mathsf N}$.

For notation convenience, we write $\Lambda$ as the annulus $\Lambda_{M,(1+\iota)M}$. We may assume $M_0>0$ is chosen large enough depending on $\iota,\beta_0,\NN$ and consider only $M>M_0$. We fix $\xi_0$ to be the boundary condition on $\partial\Lambda.$ We also assume that both $M$ and $\iota M$ are powers of 2 so that we can partition the annulus $\Lambda$ into dyadic boxes. Fix $\delta, \rho>0$ to be two small constants to be determined later.

\begin{definition}\label{def: def of good box}
    For a dyadic box $\fB$, we write $\fB^{(k)}=(1+k\delta)\fB$ to denote the box concentric with $\fB$ and with side length as $(1+k\delta)$ times that of $\fB$.

    For a (realization of) $\mathsf{z}$-copy local time $L^{\zz}$, we say $\fB$ is a good box if for any boundary condition $\xi $ on $\fB^{(10)}$ the following holds
    \begin{equation}\label{eq: def of good box}
        \mu_{\fB^{(10)},\beta-L^{\zz}}^\xi (\partial\fB\con{\mathsf{xy}}\partial \fB^{(2)})\le \rho.
    \end{equation} Note that in \eqref{eq: def of good box}, the measure $\mu^{\xi}_{\fB^{(10)},\beta-L^{\zz}}$ stands for the XY-measure on $\fB^{(10)}$ with local time conditioning $\beta-L^{\zz}$.
    In addition, whether $ \mathsf{B} $ is a good box depends only on $L^{\zz}$ in $ \mathsf{B}^{(10)} $.
\end{definition}
\begin{lemma}\label{lem: good box probability}
    There exist constants $ C_4>1 $ and $ c_5>0 $ depending on $\epsilon,\beta_0,\delta$ and $\rho$ such that for any $\xx\yy$ boundary condition $\xi$ on $\fB^{(C_4)}$ and any local time $\gamma$ with $\gamma_u\le \beta_0$ for all $ u\in \fB^{(C_4)}$, the following holds
    \begin{equation}\label{eq: good box probability}
        \nu_{\fB^{(C_4)},\gamma}^\xi (\fB\text{ is good})\ge c_5.
    \end{equation}
\end{lemma}
We postpone the proof of Lemma \ref{lem: good box probability} to the end of this section.
For a realization of $\mathsf{z}$-copy local time, we want to construct a collection of good boxes such that it would be difficult for the $\xx,\yy$-copy local time clusters to intersect those boxes. To this end, we define the following adaptive exploration for good boxes.
\begin{definition}\label{def: good box exploration}
    Recall that $\Lambda$ is the annulus $\Lambda_{M,(1+\iota)M}$. Let $ \Lambda^{-1}_*=\Lambda $ and $ \mathfrak{G}^{-1}_*=\varnothing $.
    Let $t\ge -\log_2(40\delta)$ be a constant to be chosen later.
    For each $\ell$, we refer to a box with side length $ \frac{M^{0.1}}{2^{t\ell}} $  as a scale $ \ell $ box. We run the following procedure iteratively for $\ell=\ell_*, \dots, \ell^*$ where $ \ell^*=\min\{ m\in\Z_+: \frac{M^{0.1}}{2^{tm}}\leq10 \} $ {and $ \ell_*=\min\{ m\in\Z_+: \frac{M^{0.01}}{2^{tm}}\leq10 \} $}.
    \begin{enumerate}
        \item Let $\mathfrak B^\ell$ be the collection of all the scale $\ell$ boxes contained in $\Lambda^{\ell-1}_*$ (and we associate these boxes with an arbitrary order), and let $\Lambda^\ell_0 = \Lambda^{\ell-1}_*$ and $\mathfrak{G}^\ell_0 = \mathfrak{G}^{\ell-1}_*$.
        \item Let $\Lambda^{\ell}_0 = \Lambda^{\ell-1}_*$ and $\mathfrak B^\ell_0 = \mathfrak B^\ell$. For $ i\geq 0$, we repeat the following construction until $\mathfrak B^{\ell}_i = \varnothing$. Take the first box $\fB$ in $\mathfrak B^{\ell}_i$ such that $\fB\subset \Lambda^\ell_i$ and let $\mathfrak B^{\ell}_{i+1} = \mathfrak B^\ell_i \setminus \{\fB\}$; if in addition $\fB$ is good, then let $\mathfrak{G}^{\ell}_{i+1} = \mathfrak{G}^\ell_i \cup \{\fB\}$ and $\Lambda^{\ell}_{i+1} = \Lambda^\ell_i \setminus \fB^{(80)}$. We denote by $\Lambda^{\ell}_*$ and $\mathfrak{G}^{\ell}_*$ the values of $\Lambda^{\ell}_i$ and $\mathfrak{G}^{\ell}_i$ upon stopping the construction (and then we move to the iteration for $\ell+1$).
    \end{enumerate}
    Let $\mathfrak{S}^{\ell_{*}}$ denote the collection of good boxes at the end.
\end{definition}
\begin{lemma}\label{lem: covering property for good boxes}
    There exist  constants $C_6,c_7,t_0>0$ where $C_6$ is an absolute constant and $c_7,t_0$ depend on $\delta$ and $\beta_0$ such that  the good boxes in $\mathfrak{S}^{\ell_{*}}$ with parameter $t=t_0$ satisfy the following two properties:
    \begin{enumerate}
        \item For any two good boxes $\fB_1,\fB_2\in \mathfrak{S}^{\ell_{*}}$, we have $\fB_1^{(40)}\cap\fB_2^{(40)}=\emptyset.$\label{item: nonintersecting property for good boxes}
        \item \label{item: covering property for good boxes} For any boundary condition $\xi$ on $\partial\Lambda$, we have (recall that $\iota M$ is the width of $\Lambda$)
        \begin{equation}
            \label{eq: good box covering error}
            \nu_{\Lambda,\beta }^{\xi}\Big( |\Lambda\setminus(\cup_{\fB\in\mathfrak{S}^{\ell_{*}}}\fB)|> C_6\delta|\Lambda|\Big)\le c_7^{-1}\exp(-c_7\iota M).
        \end{equation}
    \end{enumerate}
\end{lemma}
\begin{proof}
Without loss of generality, we assume $\log_2(\delta)\in \mathbb{Z}$ and we choose $t_0> -\log_2(\delta)$. Then, in the construction in Definition \ref{def: good box exploration}, $\Lambda^\ell_*$ can always be divided into scale $\ell+1$ boxes. Let $C_6>0$ be a large constant to be chosen later. Note that \begin{align}
|\Lambda\setminus(\cup_{\fB\in\mathfrak{S}^{\ell_{*}}}\fB)|&\le |\Lambda\setminus(\cup_{\fB\in\mathfrak{S}^{\ell_{*}}}\fB^{(80)})|+\sum_{\fB\in\mathfrak{S}^{\ell_{*}}}|\fB^{(80)}\setminus\fB|\nonumber\\&\le |\Lambda\setminus(\cup_{\fB\in\mathfrak{S}^{\ell_{*}}}\fB^{(80)})|+(80\delta + 4\cdot(80\delta)^2)|\Lambda|. \label{eq: covering estimate 1}
\end{align} Therefore, it suffices to control $|\Lambda\setminus(\cup_{\fB\in\mathfrak{S}^{\ell_{*}}}\fB^{(80)})|$.
    Assertion~\ref{item: nonintersecting property for good boxes} follows directly from the construction in Definition \ref{def: good box exploration} (recall that $\Lambda^{\ell}_{i+1} = \Lambda^\ell_i \setminus \fB^{(80)}$). We next turn to the proof of Assertion~\ref{item: covering property for good boxes}. For notation convenience, we write $M' = \iota M$. Let $\mathfrak{Q}_\ell$ denote the collection of possible realizations for $\Lambda^{\ell}_*$. Note that $\Lambda\setminus\Lambda^\ell_*$ is the union of some scale $\ell$ boxes, and the number of scale $\ell$ boxes is $M'M^{0.8}\cdot 2^{2t\ell}$. Thus, we get that
    \begin{equation}
        \big|\mathfrak{Q}_\ell\big|\le 2^{M'M^{0.8}\cdot 2^{2t\ell}}.\label{eq: enumeration for good box history}
    \end{equation}
    Fix $\kappa\in (0,\frac12)$ as a small constant to be determined later. For any $\ell$ and any realization $\lambda^{\ell}\in \mathfrak{Q}_\ell$, let $\mathcal{E}_{\ell+1}(\lambda^\ell)$ denote the event that the following occurs when we run the $(\ell+1)$-th iteration in Definition \ref{def: good box exploration} with initialization $\Lambda_*^\ell = \lambda^\ell$:
    \begin{equation}
        \big|\Lambda^{\ell+1}_*\big|\le \max\{(1-\kappa)|\lambda^\ell|,\frac{C_6\delta}{2}|\Lambda|\}.\label{eq: large deviation for good boxes event}
    \end{equation}
    Let $\mathcal{E}_{\ell+1}=\cap_{\lambda^{\ell}\in \mathfrak{Q}_\ell}\mathcal{E}_{\ell+1}(\lambda^\ell)$ and let $\mathcal{E}^*=\cap_{\ell}\mathcal{E}_\ell.$ We first prove by contradiction that on the event $\mathcal{E}^*$, we have  $|\Lambda\setminus(\cup_{\fB\in\mathfrak{S}^{\ell_{*}}}\fB^{(80)})|\le C_6\delta|\Lambda|/2$.  Suppose that $|\Lambda\setminus(\cup_{\fB\in\mathfrak{S}^{\ell_{*}}}\fB^{(80)})|> C_6\delta|\Lambda|/2$. Then, by \eqref{eq: large deviation for good boxes event}, we get that 
    \begin{equation}
        \big|\Lambda^{\ell+1}_*\big|\le (1-\kappa)|\Lambda^\ell_*|\le (1-\kappa)^{\ell+1}|\Lambda^0_*|.\label{eq: final region upper bound}
    \end{equation}
    Recalling that $ \ell^*-\ell_*\geq\lceil\frac{\log_2M}{2t}\rceil-20$ and noting that $M'=\iota M>0$ is large enough, we arrive at a contradiction. Therefore, recalling \eqref{eq: covering estimate 1} we obtain that
    \begin{equation}
        \label{eq: good box probability conversion}
        \nu_{\Lambda,\beta}^{\xi}\Big( |\Lambda\setminus(\cup_{\fB\in\mathfrak{S}^{\ell_{*}}}\fB)|> C_6\delta|\Lambda|\Big)\le \nu_{\Lambda,\beta }^{\xi}\Big(\big(\mathcal{E}^*\big)^c\Big).
    \end{equation}
    It remains to control the probability of $\mathcal{E}^*$. For any $\ell$ and $\lambda^\ell\in \mathfrak{Q}_\ell$, we want to upper-bound the probability of $\big(\mathcal{E}_{\ell+1}(\lambda^\ell)\big)^c$. Without loss of generality, we may assume that $|\lambda^\ell|\ge \frac{C_6\delta}{2}|\Lambda|$, since otherwise we get from $\Lambda^{\ell+1}_*\subset\lambda^\ell$ that $\big(\mathcal{E}_{\ell+1}(\lambda^\ell)\big)^c=\varnothing.$ We next divide the boxes in $\mathfrak{B}^{\ell+1}$ (which is now fixed, since it is a function of $\lambda^\ell$) into $m=(C_4)^2$ groups $\mathfrak{B}^{\ell+1,1},\cdots, \mathfrak{B}^{\ell+1,m}$ such that for any two good boxes $\fB_1,\fB_2$ in the same group, we have $\fB_1^{(C_4)}\cap \fB_2^{(C_4)}=\emptyset$.
    Furthermore, in this division we may assume that $\sum_{\fB\in\mathfrak{B}^{\ell+1,1}}|\fB|\ge \frac{|\lambda^\ell|}{m}$ without loss of generality.
Let $\mathfrak{B}^{\ell+1,1}=\{\fB_1,\cdots,\fB_{t}\}$ and denote by $\mathcal{A}_i$ the event that $\fB_i$ is good. Similar to the discussions in \eqref{eq: dmp for three copies} and \eqref{eq: dmp for three copies 2},
 let $\mathcal F_i^{\mathrm{ext}}$ be the $\sigma$-field generated by the $\xx\yy$ boundary condition $\xi_i$ for $\fB_i^{(C_4)}$ and $\beta_i^\zz$ for the $\zz$-copy local times from loops intersecting $(\fB_i^{(C_4)})^c$, where $\xi_i$ encodes  $\xx\yy$ local times on vertices
    outside the interior of $\fB_i^{(C_4)}$  and local time connectivities outside $\fB_i^{(C_4)}$.
    Thus, Lemma \ref{lem: good box probability} gives
    \begin{equation}\label{eq: dmp for three copies 3}
        \nu^{\xi}_{\Lambda,\beta}(\mathcal{A}_i\mid\mathcal F_i^{\mathrm{ext}}) =\nu_{\fB_i^{(C_4)},\beta-\beta^\zz_i}^{\xi_i}(\mathcal{A}_i) \ge c_5.
    \end{equation}
 Note that for $j\neq i$, $\mathcal{A}_j$ is measurable with respect to the edge configuration within $\fB_{j}^{(C_4)}$ and thus $\mathcal{A}_j $ is measurable in $ \mathcal F_i^{\mathrm{ext}}$. Therefore, we get that for any $T: \{0,1\}^{\{1,\cdots,t\}}$ \begin{equation}
        \nu^{\xi}_{\Lambda_{M,2M},\beta}(\mathcal{A}_i\mid \1_{\mathcal{A}_j}=T_i,\forall j\neq i) \ge c_5.\nonumber
    \end{equation}
    Applying the Bernstein inequality, we get that \begin{align}
        \nu_{\Lambda,\beta}^{\xi}\Big(\Big|\big\{\fB\in \mathfrak{B}^{\ell+1,1}: ~\fB \text{ is good}\big\}\Big|\le \frac{c_5}{2}|\mathfrak{B}^{\ell+1,1}|\Big)  &\le \exp(-\bc_1|\mathfrak{B}^{\ell+1,1}|))\nonumber\\&\le \exp(-\bc_2 M'M^{0.8} 2^{2t\ell+2t}).\label{eq: bad event probability in good box}
    \end{align} 
    Next, we will show that on the event $\Big\{\big|\big\{\fB\in \mathfrak{B}^{\ell+1,1}: ~\fB \text{ is good}\big\}\big|> \frac{c_5}{2}|\mathfrak{B}^{\ell+1,1}|\Big\}$, we will get at least $\frac{c_5}{18}|\mathfrak{B}^{\ell+1,1}|$ scale $\ell+1$ good boxes in the $(\ell+1)$-th iteration of Definition~\ref{def: good box exploration}. Let $\mathfrak{C}_{\ell+1}$ denote the collection of scale $\ell$ good boxes in $\mathfrak{S}^{\ell_{*}}$. By Definition \ref{def: good box exploration}, we obtain that for any good box $\fB$ in $\mathfrak{B}^{\ell+1,1}$, it is either contained in $\mathfrak{C}_{\ell+1}$, or there exists a good box $\fC\in \mathfrak{C}_{\ell+1}$ such that $\fB\cap\fC^{(80)}\neq\emptyset$. For any $\fC\in \mathfrak{C}_{\ell+1}$, $\fC^{(80)}$ intersects with at most 9 boxes in $\mathfrak{B}^{\ell+1,1}$.
    Therefore, we obtain that  $|\cup_{\fC\in\mathfrak{C}_{\ell+1}}\fC| \ge \frac{c_5}{18}|\mathfrak{B}^{\ell+1,1}|\cdot\frac{M^{0.2}}{2^{2t(\ell+1)}}\ge \frac{c_5}{18m}|\lambda^\ell|.$ In conclusion, choosing $\kappa=\frac{c_5}{18m}$ we get that
    \begin{equation}
    \Big\{\big|\big\{\fB\in \mathfrak{B}^{\ell+1,1}: ~\fB \text{ is good}\big\}\big|> \frac{c_5}{2}|\mathfrak{B}^{\ell+1,1}|\Big\}\subset\mathcal{E}_{\ell+1}(\lambda^\ell).\nonumber
    \end{equation}
    Combining with \eqref{eq: bad event probability in good box}, we get that
    \begin{equation}
         \nu_{\Lambda,\beta}^{\xi}\Big(\big(\mathcal{E}_{\ell+1}(\lambda^\ell)\big)^c\Big)\le \exp(-\bc_2M'M^{0.8}\cdot 2^{2t\ell+2t}).\label{eq: bad event probability in good box 2}
    \end{equation}
    Summing \eqref{eq: bad event probability in good box 2} over $0\le\ell  \le \lceil\frac{\log_2M}{2t}\rceil-1$ and $\lambda^\ell\in\mathfrak{Q}_\ell$ and combining with \eqref{eq: enumeration for good box history}, we obtain that
    \begin{align}
        \nu_{\Lambda,\beta}^{\xi}(\mathcal{E}^*)
        &\ge 1-\sum_{\ell=0}^{\lceil\frac{\log_2M}{2t}\rceil-1}\sum_{\lambda^\ell\in\mathfrak{Q}_\ell}\nu_{\Lambda,\beta}^{\xi}\Big(\big(\mathcal{E}_{\ell+1}(\lambda^\ell)\big)^c\Big)\nonumber\\
        &\ge 1- \sum_{\ell=0}^{\lceil\frac{\log_2M}{2t}\rceil-1}2^{M'M^{0.8}\cdot 2^{2t\ell}}\exp(-\bc_2 M'M^{0.8} 2^{2t\ell+2t})\nonumber
    \end{align}Choosing $t_0$ large enough such that $\bc_22^{2t}\ge2$, then we that \begin{align}
        \nu_{\Lambda,\beta}^{\xi}(\mathcal{E}^*)\ge 1- \sum_{\ell=0}^{\lceil\frac{\log_2M}{2t}\rceil-1}\exp(-M'M^{0.8}\cdot 2^{2t\ell})\ge1-\exp(-\bc_3 \iota M).\nonumber
    \end{align} Plugging into \eqref{eq: good box probability conversion} finishes the proof of \eqref{eq: good box covering error}.    
\end{proof}
From now on, we will fix a realization of $\mathsf{z}$-copy local time field such that $|\Lambda\setminus(\cup_{\fB\in\mathfrak{S}^{\ell_{*}}}\fB)|\le  c_4\delta|\Lambda|$. Next, we want to determine whether there is an $\mathsf{x}$- or $\yy$-local time cluster crossing through annuli surrounding the good boxes. For those good boxes that do not contain an $\xx$- or $\yy$-local time crossing, we then want to explore all the $\xx$- or $\yy$-loops from the good boxes and consider the measure on the remaining graph.
Note that once loops from a set $A$ are sampled and removed, the remaining loops form (loops of) an XY model conditioned on having local times 0 on $A$ (so equivalently $A$ is removed) and with local times properly adjusted on $A^c$ (i.e., decreased by the amount of the local times contributed by all sampled loops).
By \eqref{eq: def of good box}, we know that the remaining graph has a small volume, and we want to apply a spin wave argument to bound the probability for a local time crossing on the remaining graph. To conduct the intuition, the first difficulty is how to preserve information on the $\xx$- or $\yy$-copy on the good boxes before and after the spin wave. Therefore, we introduce the following definitions.
\begin{definition}\label{def: base spin wave}
    Let $F_0$ denote a base spin wave on $\Lambda=\Lambda_{M, (1+\iota)M}$ such that 
    \begin{equation}\label{eq: def of base spin wave}
        F_0(v)=\frac{\dist_{|\cdot|_{\infty}}(v,\partial\Lambda_{(1+\iota )M})}{\iota M}\cdot \pi, \quad \text{ for } v\in\Lambda.
    \end{equation}
    Then we get that $ F_0(x)=0 $ for $ x\in\partial\Lambda_{(1+\iota)M} $ and $ F_0(x)=\pi $ for $ x\in\partial\Lambda_M $.
\end{definition}

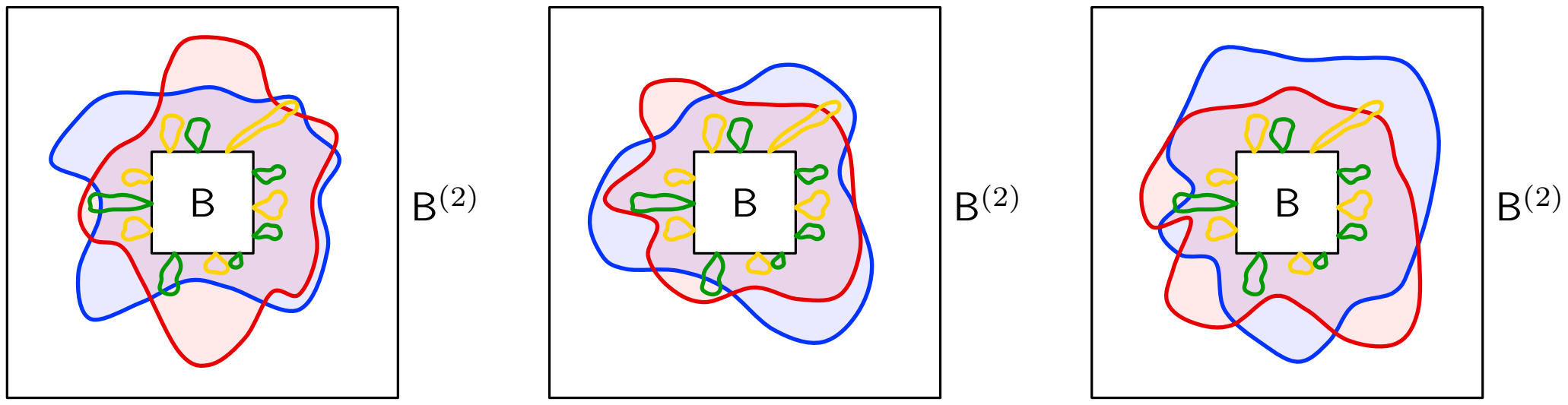
\begin{figure}[htbp]
    \centering
    \begin{tikzpicture}[
  x=1cm, y=1cm,
  scale=.85,
  line cap=round, line join=round,
  every node/.style={font=\normalsize},
  rooted loop/.style={line width=1pt}
]
  \path[use as bounding box] (-2.05,-2.05) rectangle (13.92,2.05);

  \newcommand{\YellowLoops}{%
    \draw[rooted loop,yellow!70!orange]
      (-.34,.52)
      .. controls (-.40,.64) and (-.44,.66) .. (-.41,.74)
      .. controls (-.39,.79) and (-.43,.88) .. (-.33,.89)
      .. controls (-.25,.90) and (-.17,.84) .. (-.22,.77)
      .. controls (-.27,.70) and (-.24,.60) .. (-.34,.52) -- cycle;
    \draw[rooted loop,yellow!70!orange]
      (.25,.52)
      .. controls (.29,.66) and (.47,.76) .. (.60,.84)
      .. controls (.72,.91) and (.78,1.02) .. (.91,1.03)
      .. controls (.96,1.04) and (.99,1.00) .. (.95,.94)
      .. controls (.92,.90) and (.88,.93) .. (.84,.86)
      .. controls (.79,.78) and (.65,.80) .. (.56,.72)
      .. controls (.44,.63) and (.37,.59) .. (.25,.52) -- cycle;
    \draw[rooted loop,yellow!70!orange]
      (-.52,.25)
      .. controls (-.60,.32) and (-.71,.36) .. (-.78,.31)
      .. controls (-.84,.27) and (-.81,.19) .. (-.74,.18)
      .. controls (-.68,.18) and (-.67,.23) .. (-.62,.21)
      .. controls (-.58,.19) and (-.54,.23) .. (-.52,.25) -- cycle;
    \draw[rooted loop,yellow!70!orange]
      (-.52,-.28)
      .. controls (-.61,-.20) and (-.65,-.13) .. (-.75,-.17)
      .. controls (-.82,-.20) and (-.77,-.27) .. (-.81,-.30)
      .. controls (-.84,-.36) and (-.77,-.40) .. (-.71,-.37)
      .. controls (-.64,-.34) and (-.59,-.36) .. (-.52,-.28) -- cycle;
    \draw[rooted loop,yellow!70!orange]
      (.13,-.52)
      .. controls (.09,-.59) and (.02,-.59) .. (.03,-.66)
      .. controls (.04,-.73) and (.11,-.70) .. (.16,-.72)
      .. controls (.24,-.74) and (.29,-.68) .. (.24,-.62)
      .. controls (.20,-.59) and (.17,-.55) .. (.13,-.52) -- cycle;
    \draw[rooted loop,yellow!70!orange]
      (.52,-.05)
      .. controls (.63,.03) and (.67,.00) .. (.70,.07)
      .. controls (.74,.14) and (.84,.12) .. (.85,.04)
      .. controls (.87,-.03) and (.78,-.05) .. (.80,-.11)
      .. controls (.80,-.18) and (.68,-.18) .. (.63,-.12)
      .. controls (.59,-.09) and (.55,-.08) .. (.52,-.05) -- cycle;
  }

  \newcommand{\GreenLoops}{%
    \draw[rooted loop,green!60!black]
      (-.05,.52)
      .. controls (-.09,.61) and (-.19,.65) .. (-.16,.72)
      .. controls (-.12,.78) and (-.15,.84) .. (-.07,.85)
      .. controls (.02,.87) and (.09,.81) .. (.05,.75)
      .. controls (.01,.69) and (.04,.60) .. (-.05,.52) -- cycle;
    \draw[rooted loop,green!60!black]
      (-.52,-.01)
      .. controls (-.66,.04) and (-.76,.10) .. (-.89,.07)
      .. controls (-1.01,.04) and (-1.08,.13) .. (-1.14,.06)
      .. controls (-1.18,.01) and (-1.13,-.02) .. (-1.15,-.06)
      .. controls (-1.16,-.09) and (-1.05,-.10) .. (-.98,-.08)
      .. controls (-.86,-.05) and (-.72,-.11) .. (-.52,-.01) -- cycle;
    \draw[rooted loop,green!60!black]
      (-.29,-.52)
      .. controls (-.33,-.62) and (-.44,-.68) .. (-.42,-.77)
      .. controls (-.41,-.83) and (-.48,-.88) .. (-.40,-.92)
      .. controls (-.33,-.96) and (-.26,-.94) .. (-.25,-.88)
      .. controls (-.24,-.81) and (-.29,-.82) .. (-.25,-.74)
      .. controls (-.20,-.66) and (-.22,-.58) .. (-.29,-.52) -- cycle;
    \draw[rooted loop,green!60!black]
      (.38,-.52)
      .. controls (.34,-.57) and (.26,-.60) .. (.29,-.64)
      .. controls (.32,-.69) and (.40,-.66) .. (.40,-.62)
      .. controls (.39,-.58) and (.40,-.56) .. (.38,-.52) -- cycle;
    \draw[rooted loop,green!60!black]
      (.52,.31)
      .. controls (.59,.37) and (.69,.30) .. (.72,.35)
      .. controls (.76,.40) and (.86,.36) .. (.83,.29)
      .. controls (.80,.23) and (.73,.27) .. (.69,.22)
      .. controls (.63,.16) and (.58,.26) .. (.52,.31) -- cycle;
    \draw[rooted loop,green!60!black]
      (.52,-.34)
      .. controls (.60,-.30) and (.61,-.21) .. (.68,-.23)
      .. controls (.73,-.25) and (.77,-.22) .. (.79,-.28)
      .. controls (.83,-.35) and (.75,-.40) .. (.70,-.36)
      .. controls (.63,-.32) and (.59,-.39) .. (.52,-.34) -- cycle;
  }

  \newcommand{\DrawBoxes}{%
    \draw[black,line width=.6pt] (-2,-2) rectangle (2,2);
    \draw[black,line width=.6pt] (-.52,-.52) rectangle (.52,.52);
    \node at (0,0) {$\mathsf B$};
    \node[anchor=west,inner sep=1.5pt] at (2.06,0) {$\mathsf B^{(2)}$};
  }

  \newcommand{\DrawCombinedPanel}{
    \path[fill=blue,fill opacity=.09,even odd rule]
      \BlueCircuit (-.52,-.52) rectangle (.52,.52);
    \path[fill=red,fill opacity=.09,even odd rule]
      \RedCircuit (-.52,-.52) rectangle (.52,.52);
    \draw[blue!80!cyan,line width=1pt] \BlueCircuit;
    \draw[red!90!black,line width=1pt] \RedCircuit;
    \DrawBoxes
    \YellowLoops
    \GreenLoops
  }

\begin{scope}[shift={(0.00,0)}]
  \def\BlueCircuit{plot[smooth cycle,tension=0.62] coordinates {
    (1.20,0.05) (1.40,0.45) (1.10,0.75) (0.95,1.07)
    (0.52,1.04) (0.10,1.18) (-0.35,1.11) (-0.94,1.11)
    (-1.28,0.86) (-1.55,0.42) (-1.13,0.20) (-1.05,-0.13)
    (-1.27,-0.72) (-1.15,-1.19) (-0.62,-1.02) (-0.15,-0.80)
    (0.27,-0.87) (0.98,-1.10) (1.28,-0.48)
  }}
  \def\RedCircuit{plot[smooth cycle,tension=0.62] coordinates {
    (1.13,0.07) (1.37,0.70) (0.94,0.92) (0.70,1.08)
    (0.50,1.58) (-0.12,1.68) (-0.36,1.37) (-0.45,0.98)
    (-0.62,0.71) (-0.96,0.48) (-1.25,-0.05) (-1.14,-0.37)
    (-0.83,-0.45) (-0.62,-0.80) (-0.24,-1.53) (0.08,-1.66)
    (0.44,-1.38) (0.71,-0.92) (0.99,-0.94) (1.17,-0.50)
  }}
  \DrawCombinedPanel
\end{scope}

\begin{scope}[shift={(5.55,0)}]
  \def\BlueCircuit{plot[smooth cycle,tension=0.62] coordinates {
    (1.12,-0.04) (0.97,0.32) (1.12,0.93) (0.73,1.33)
    (0.33,1.40) (-0.05,1.18) (-0.51,1.03) (-0.79,0.72)
    (-1.28,0.49) (-1.58,-0.11) (-1.47,-0.53) (-1.22,-0.68)
    (-0.68,-0.70) (-0.20,-0.88) (0.29,-1.30) (0.79,-1.43)
    (1.16,-1.12) (1.30,-0.62)
  }}
  \def\RedCircuit{plot[smooth cycle,tension=0.62] coordinates {
    (1.12,0.01) (1.05,0.59) (0.79,0.96) (0.30,1.00)
    (-0.15,1.04) (-0.56,1.24) (-0.98,1.20) (-1.06,0.90)
    (-0.94,0.52) (-1.38,0.14) (-1.37,-0.05) (-0.91,-0.20)
    (-0.99,-0.55) (-0.77,-0.87) (-0.41,-1.02) (0.09,-0.87)
    (0.55,-0.99) (0.98,-0.94) (1.18,-0.48)
  }}

  \DrawCombinedPanel
\end{scope}

\begin{scope}[shift={(11.10,0)}]
  \def\BlueCircuit{plot[smooth cycle,tension=0.62] coordinates {
    (1.47,0.03) (1.54,0.76) (1.25,1.40) (0.68,1.48)
    (0.17,1.46) (-0.29,1.55) (-0.73,1.55) (-1.01,1.04)
    (-1.22,0.48) (-0.96,0.13) (-1.30,-0.24) (-1.10,-0.53)
    (-0.74,-0.60) (-0.64,-1.10) (-0.25,-1.45) (0.12,-1.63)
    (0.43,-1.42) (0.68,-1.05) (1.09,-0.92) (1.26,-0.65)
  }}
  \def\RedCircuit{plot[smooth cycle,tension=0.62] coordinates {
    (1.27,0.06) (1.12,0.39) (0.97,0.85) (0.45,1.00)
    (0.05,1.17) (-0.34,1.01) (-0.98,0.90) (-1.30,0.50)
    (-1.49,0.11) (-1.38,-0.25) (-0.99,-0.17) (-1.16,-0.66)
    (-1.23,-1.01) (-0.99,-1.22) (-0.58,-1.26) (-0.14,-0.96)
    (0.20,-1.10) (0.62,-1.41) (1.10,-1.45) (1.33,-1.14)
    (1.35,-0.56)
  }}
  \DrawCombinedPanel
\end{scope}

\end{tikzpicture}
    \caption{An illustration of a type $ ({\rm I}, -1) $ good box. Green loops: $\mathfrak{L}_\fB^\xx$. Yellow loops: $\mathfrak{L}_\fB^\yy$.   The regions between $\partial\fB$ and the red  circuits are shaded in  red, representing local time clusters from $\fB\cup (\cup \mathfrak{L}_\fB^\xx)$ (i.e., the local time clusters are with respect to spin configurations shown below these pictures as well as the edge marks that stay unchanged; see Definition \ref{def: crossing event decomposition according to loops}). The regions between $\partial\fB$ and the blue  circuits are shaded in blue, representing local time clusters from $\fB\cup (\cup \mathfrak{L}_\fB^\yy)$.  Note that the loops from $\fB$ are all the same in these subfigures.}
    \label{fig: type of good box}
\end{figure}

\begin{definition}\label{def: crossing event decomposition according to loops}
Fix a good box $\fB$.
Let $\mathfrak{L}^\xx_\fB$ denote the collection of $\xx$-loops {intersecting} $\fB$; the same notation applies to $\yy$ and we write $\mathfrak{L}_\fB = \mathfrak{L}_\fB^\xx \cup \mathfrak{L}_\fB^\yy$.
    Let $\Delta_\fB$ denote the collection of loop configurations $\mathfrak{L}_\fB$ such that all loops  in $\mathfrak{L}_\fB$ are in $\fB^{(2)}$.
    Let $\sigma_\fB$ denote the spin configuration (treated as angles) on $\fB^{(2)}\setminus\fB$ under the local time conditioning $\beta-L^{\zz}-\mathfrak{L}_\fB$ ({in what follows} we slightly abuse the notation $\mathfrak{L}_\fB$ to denote its induced local time field)---for convenience, we refer to the corresponding local time clusters as $\sigma_{\mathsf B}$-clusters. In addition, let $\mathfrak{e}$ denote the collection of edge marks that are the random variables determining whether the (metric graph) $\sigma_{\mathsf B}$-clusters percolate through an edge, (see Lemma \ref{lem: FK-cluster-loop}). Furthermore, for an edge $e=\{u,v\}$ where $u\in l\in \mathfrak{L}_\fB^\ii$ and $v\notin \mathfrak{L}_\fB^\ii$, the edge mark also determines whether the $\sigma_{\mathsf B}$-cluster of $v$ intersects with $l$.
    By Lemma \ref{lem: FK-cluster-loop}, we know that the clusters are measurable with respect to $(\sigma_\fB, \mathfrak{e}_\fB)$. Let $\Psi$ denote the law of the edge marks $\mathfrak{e}_\fB$ on $\fB^{(2)}\setminus\fB$.
   
    Let $\mathsf{i}\in\{\xx,\yy\}.$ For any loop configuration $\mathfrak{L}_\fB\in \Delta_\fB$ and any edge mark $\mathfrak{e}_\fB$,
     let $\mathcal{E}^{\mathsf{i}}_*(\mathfrak{L}_\fB,\mathfrak{e}_\fB)$ denote the collection of spin configurations $\sigma_\fB$ such that  $\sigma_{\mathsf B}$-clusters in the $\mathsf i$-copy do not connect $\partial \fB^{(2)}$ and $\fB \cup(\cup  \mathfrak{L}_\fB^\mathsf{i})$ (here we use $\cup \mathfrak{L}_\fB^{\mathsf{i}}$ to denote the union of the traces of loops in $\mathfrak{L}_\fB^{\mathsf{i}}$).
     Moreover, for $k\in \mathbb{Z}$, we define $$\mathcal{E}^{\mathsf{i},k}_*(\mathfrak{L}_\fB,\mathfrak{e}_\fB)=\{\sigma_\fB:\sigma_\fB+k F_0\in \mathcal{E}^{\mathsf{i}}_*(\mathfrak{L}_\fB,\mathfrak{e}_\fB)\}.$$
     For notation convenience, we write $$\mathcal{F}^k_*(\mathfrak{L}_\fB,\mathfrak{e}_\fB)=\cap_{\mathsf{i}\in\{\xx,\yy\},j\in\{0,1,-1\}} \mathcal{E}^{\mathsf{i},j-k}_*(\mathfrak{L}_\fB,\mathfrak{e}_\fB).$$
\end{definition}
\begin{definition}\label{def: type for good boxes}
    Fix a constant $ \kappa\in(0,1) $ to be chosen.
    For each good box $\fB$, we assign an auxiliary variable $\aux_\fB$ where $\{\aux_\fB\}$ is a collection of i.i.d. Bernoulli variables so that they take value 1 with probability $1-\kappa$.
    For a good box $\fB$, we say it is type $({\rm I},k),~k=-1, 0, 1$, if the following properties hold (we allow such types to occur simultaneously for one box):
    \begin{enumerate}
        \item The local time cluster in either $\xx$- or $\yy$-copy in the annulus $\fB^{(2)}\setminus\fB$ does not percolate (we may also interchangeably say the $\mathsf{xy}$-local time clusters do not percolate). We denote this event by $\mathcal{B}=(\mathcal{E}^\xx_0\cup\mathcal{E}^\yy_0)^c$ where $\mathcal{E}^{\mathsf{i}}_0$ denotes the event that the local time cluster in $\mathsf{i}$-copy in the annulus $\fB^{(2)}\setminus\fB$   percolates.
        \item 
        Given $\mathfrak{L}_\fB,\mathfrak{e}_\fB$, the spin configuration $\sigma_\fB$ is in $\mathcal{F}^k_*(\mathfrak{L}_\fB,\mathfrak{e}_\fB)$.
        
        \item The auxiliary variable $\aux_\fB$ takes value $1.$
    \end{enumerate}
     We say it is type $({\rm II},k)$ if it is not type $({\rm I},k)$, for $k=-1, 0, 1$. For notation convenience, we will also say it is type 1 if it is type $({\rm I},0)$ and type 2 if it is type $({\rm II},0).$
     See Figure~\ref{fig: type of good box} for an illustration of a type (I, $k$) good box.
\end{definition}
\begin{lemma}\label{lem: type probability}
Recall $\rho$ from~\eqref{eq: def of good box}.
Fix a good box $\fB$.
For any $ k=-1,0,1 $ and any boundary condition $\xi$ on $\partial\Lambda$, we have 
\begin{equation}\label{eq: type probability without aux}
    \mu_{\Lambda,\beta-L^{\zz}}^{\xi}(\fB \text{ is type } ({\rm II},k)\mid \mathrm{AUX}_\fB=1)\le 6 \rho.
\end{equation}In, addition, we have
\begin{equation}\label{eq: type probability}
    \mu_{\Lambda,\beta-L^{\zz}}^{\xi}(\fB \text{ is type } ({\rm II},k))\le 6 \rho+\kappa.
\end{equation}
\end{lemma}
Before proving Lemma \ref{lem: type probability}, we start with the following decomposition lemma.

\begin{lemma}\label{lem: type 1 loop expansion}
Fix a good box $\fB$ and an expansion $\fB'=\fB^{(C_8)}$.
Recall that $\mathsf{M}(\mathsf{d}\mathfrak{L}_\fB)$ denotes the mass of the  loop configuration $\mathfrak{L}_\fB$ and recall from Definition \ref{def: crossing event decomposition according to loops} that $\Psi(\mathsf{d}\mathfrak{e}_\fB)$ is the law of the edge marks $\mathfrak{e}_\fB$. Recall from \eqref{eq: def--spin O(N) partition function}  that ${Z}^\zeta_{\Omega,\gamma}$ denotes the spin partition function on $\Omega$ with spin boundary condition $\zeta$ and local time condition $\gamma$. Recall $G_{\Omega}(\gamma)$ from Lemma \ref{lem: spin local time factorization}.
Recall the definition of $\mathcal{F}_*^k(\mathfrak{L}_\fB,\mathfrak{e}_\fB)$ from Definition \ref{def: crossing event decomposition according to loops}. 
Let $f$ be a function that coincides with $F_0$ inside $\fB^{(2)}$ and takes constant value on the boundary $\partial\fB'$. Then for any spin boundary condition $\zeta$ on $\partial\fB'$, any nonnegative $\gamma = (\gamma_v)_{v\in \mathsf B'}$ and any $k\in\{-1,0,1\}$, we have
\begin{align}
    \mu_{\fB',\gamma}^{\zeta+kf|_{\partial\fB'}}(\fB \text{ is type } ({\rm I},k))&=(1-\kappa)\cdot q^k_\fB(\zeta),\label{eq: type 1, k loop expansion}
\end{align}
where for $k\in \{-1, 0, 1\}$ 
$$q^k_\fB(\zeta)=\int_{\mathfrak{L}_\fB\in\Delta_\fB}\mathsf{M}(\mathsf{d}\mathfrak{L}_\fB)\cdot  \Psi(\mathsf{d}\mathfrak{e}_\fB)\cdot\frac{G_{\fB'}(\gamma-\mathfrak{L}_\fB)}{G_{\fB'}(\gamma)} \cdot\frac{{Z}^{\zeta+kf|_{\partial{\fB'}}}_{\fB',\gamma-\mathfrak{L}_\fB}\big(\mathcal{F}^k_*(\mathfrak{L}_\fB,\mathfrak{e}_\fB)\big)}{{Z}^{\zeta+kf|_{\partial{\fB'}}}_{\fB',\gamma}}.$$
\end{lemma}
\begin{lemma}
For any boundary partition $\mathcal{P}$, recall the definition of $\mu^{\zeta,\mathcal{P}}_{\fB',\gamma}$ from \eqref{eq: loop relation between spin boundary condition and local time boundary condition}.
    Then for $k=0, -1, 1$, the type $({\rm I},k)$ probability does not depend on the boundary partition $\mathcal{P}$, i.e., \begin{equation}
        \mu_{\fB',\gamma}^{\zeta,\mathcal{P}}(\fB \text{ is type } ({\rm I},k))=\mu_{\fB',\gamma}^{\zeta}(\fB \text{ is type } ({\rm I},k)).
    \end{equation}
\end{lemma}    
\begin{proof}
     Let $\mathcal{A}$ denote the event that $\mathsf{xy}$-local time clusters do not percolate in the annulus $\fB^{(2)}\setminus\fB$. Since  $\mathcal{A}$ is a local time event, the probability of $\mathcal{A}$ does not depend on the partition $\mathcal{P}$ given $\zeta$. Since $
\{\fB \text{ is type } ({\rm I},0)\}\subset\mathcal{A}$, it suffices to consider the conditional law  $\mu^{\zeta,\mathcal{P}}_{\fB',\gamma}(\cdot\mid \mathcal{A}).$ By Lemma \ref{lem: spin boundary condition plus different partitions},  the law of  $\mathfrak{L}_\fB$ under $\mu^{\zeta,\mathcal{P}}_{\fB',\gamma}(\cdot\mid \mathcal{A})$ does not depend on the partition $\mathcal{P}$. In addition, the law of $\mathfrak{e}_\fB$ does not depend on the partition $\mathcal{P}$. Furthermore, given $\mathfrak{L}_\fB,\mathfrak{e}_\fB$, the event $\mathcal{F}^k_*(\mathfrak{L}_\fB,\mathfrak{e}_\fB)$ is a spin configuration event, and thus does not depend on the partition $\mathcal{P}$. In conclusion, the type $({\rm I},0)$ probability does not depend on the boundary partition $\mathcal{P}$, and the same holds for type $({\rm I},k)$ when $k\in\{-1,1\}$.
\end{proof}

In what follows, we decide to drop $\mathcal{P}$ from the superscript when considering the types since it does not change the probability of a good box to be of type $({\rm I},0)$. 

\begin{proof}[Proof of Lemma \ref{lem: type 1 loop expansion}]
{By symmetry, it suffices to prove \eqref{eq: type 1, k loop expansion} for $k=0$.} Recall $\mathcal{Z}^\zeta_{\Omega,\gamma}$ from \eqref{eq: spin boundary local time boundary condition}. By Lemma \ref{lem: spin local time factorization}, it suffices to prove (note that the font of $Z$ below is mathcal instead of plain) \begin{align}
     \mu_{\fB',\gamma}^{\zeta}(\fB \text{ is type } ({\rm I},0))&=(1-\kappa)\cdot \int_{\mathfrak{L}_\fB\in\Delta_\fB,\mathfrak{e}_\fB}\mathsf{M}(\mathsf{d}\mathfrak{L}_\fB)\cdot  \Psi(\mathsf{d}\mathfrak{e}_\fB)\cdot\frac{\mathcal{Z}^{\zeta}_{\fB',\gamma-\mathfrak{L}_\fB}\big(\mathcal{F}^{ 0}_*(\mathfrak{L}_\fB,\mathfrak{e}_\fB)\big)}{\mathcal{Z}^{\zeta}_{\fB',\gamma}}.\label{eq: type 1, 0 loop expansion old}
\end{align}
Fix $\mathfrak{L}_\fB$ and $\mathfrak{e}_\fB$. Let $\mathcal{A}$ be an arbitrary local time event and recall $\mathcal{B}$ from Definition \ref{def: type for good boxes}. By \eqref{eq: def spin b.c.}, in order to handle the boundary condition, we denote by $\mathcal{D}$  the boundary disconnecting event $\{V^+_\xx\ncon[0.5]{\xx}V^-_\xx\}\cap \{V^+_\yy\ncon[0.5]{\yy}V^-_\yy\}$. (Recall the definitions from the paragraph before \eqref{eq: def spin b.c.}.)
Let $\big((\mathfrak L'_{\mathsf B})^\xx , (\mathfrak L'_{\mathsf B})^ \yy\big)$ denotes a loop configuration (including the boundary excursions) with loops (and excursions) not intersecting $\fB$ and let
$$\mathrm{LTC}=\{\mathfrak{L}_\fB+\mathfrak{L}_\fB'=\gamma, ~(\mathfrak{L}_\fB')^\mathsf{i}_u={\gamma_u(\zeta^{\mathsf{i}}_u)^2},\forall u\in\partial\fB', \mathsf i\in \{\xx, \yy\} \}$$
be the local time condition subject to the boundary condition $\zeta$. Let $\Phi^\zeta_{\mathfrak{L}_\fB,\mathfrak{L}_\fB'}(\mathcal{A}\cap \mathcal{B})$ denote the product of the indicator that $\mathcal A \cap \mathcal B$ holds on $(\mathfrak L_{\mathsf B}, \mathfrak L'_{\mathsf B})$ and the probability of sign assignments to the local time clusters in $(\mathfrak{L}_\fB,\mathfrak{L}_\fB')$ are consistent with $\zeta$.
Since on the event $\mathcal{B}$, any loop intersecting $\fB$ cannot intersect $\partial\fB^{(2)}$, we calculate by \eqref{eq: def spin b.c.} and Lemma \ref{lem:cpf-spatial-restriction} that \begin{align}
    \mathcal{Z}^{\zeta}_{\fB',\gamma}\big(\mathcal{A}\cap \mathcal{B}\big)&=\int_{\mathfrak{L}_\fB\in\Delta_\fB,\mathfrak{L}_\fB'}\1_{\mathrm{LTC}}\cdot \mathsf{M}(\mathsf{d}\mathfrak{L}_\fB)\cdot  \mathsf{M}(\mathsf{d}\mathfrak{L}_\fB')\cdot \Phi^\zeta_{\mathfrak{L}_\fB,\mathfrak{L}_\fB'}(\mathcal{A}\cap\mathcal{B})\nonumber\\
    &=\int_{\mathfrak{L}_\fB\in\Delta_\fB} \mathsf{M}(\mathsf{d}\mathfrak{L}_\fB)\cdot \int_{\mathfrak{L}_\fB'}\1_{(\mathfrak{L}_\fB,\mathfrak{L}_\fB')\in\mathrm{LTC}\cap \mathcal{D}\cap\mathcal{A} \cap \mathcal{B}}\cdot \mathsf{M}(\mathsf{d}\mathfrak{L}_\fB')\cdot \Phi^\zeta_{\mathfrak{L}_\fB,\mathfrak{L}_\fB'}(\mathcal{A}\cap\mathcal{B}). \label{eq: local loops expansion with spin bc 1}
\end{align}Note that on the event $\mathcal{B}$, the event $\mathcal{D}$ depends only on $\mathfrak{L}_\fB'$. Denoting by $\Delta(\mathcal D)$ the collection of realizations of $\mathfrak{L}'_\mathsf B$ such that $\mathcal D$ is satisfied if we only consider local time connectivity generated by $\mathfrak{L}'_B$, we get that\begin{align}
&\int_{\mathfrak{L}_\fB'}\1_{(\mathfrak{L}_\fB,\mathfrak{L}_\fB')\in\mathrm{LTC}\cap \mathcal{D}\cap\mathcal{A} \cap \mathcal{B}}\cdot \mathsf{M}(\mathfrak{L}_\fB')\cdot \Phi^\zeta_{\mathfrak{L}_\fB,\mathfrak{L}_\fB'}(\mathcal{A}\cap\mathcal{B})\nonumber\\=~&\int_{\mathfrak{L}_\fB'\in \Delta(\mathcal D)}\1_{(\mathfrak{L}_\fB,\mathfrak{L}_\fB')\in\mathrm{LTC}\cap\mathcal{A} \cap \mathcal{B}}\cdot  \mathsf{M}(\mathfrak{L}_\fB')\cdot \Phi^\zeta_{\mathfrak{L}_\fB,\mathfrak{L}_\fB'}(\mathcal{A}\cap\mathcal{B})\nonumber\\=~&\mathcal{Z}^\zeta_{\fB',\gamma-\mathfrak{L}_\fB}\Big(\mathcal{A}(\mathfrak{L}_\fB) \cap \mathcal{B}(\mathfrak{L}_\fB)\Big)\label{eq: local loops expansion with spin bc 2}
\end{align}where $\mathcal{A}(\mathfrak{L}_\fB)=\{\mathfrak{L}_\fB':(\mathfrak{L}_\fB,\mathfrak{L}_\fB')\in \mathcal{A}\}$ and $\mathcal{B}(\mathfrak{L}_\fB)=\{\mathfrak{L}_\fB':(\mathfrak{L}_\fB,\mathfrak{L}_\fB')\in \mathcal{B}\}$. Combining \eqref{eq: local loops expansion with spin bc 1} with \eqref{eq: local loops expansion with spin bc 2}, we get that \begin{align}
     \mathcal{Z}^{\zeta}_{\fB',\gamma}\big(\mathcal{A}\cap \mathcal{B}\big)&=\int_{{\mathfrak{L}_\fB\in\Delta_\fB}}\mathsf{M}(\mathsf{d}\mathfrak{L}_\fB)\cdot\mathcal{Z}^\zeta_{\fB',\gamma-\mathfrak{L}_\fB}\Big(\mathcal{A}(\mathfrak{L}_\fB) \cap \mathcal{B}(\mathfrak{L}_\fB)\Big).\label{eq: local loops expansion with spin bc 3}
\end{align}
Recall that conditioned on $\mathfrak{L}_\fB$ and $\mathfrak{e}_\fB$ and $\mathrm{Aux}_\fB = 1$, $\{\fB \text{ is type } ({\rm I},0\}$ is equivalent to $\mathcal{F}^{ 0}_*(\mathfrak{L}_\fB, \mathfrak{e}_\fB)$. In addition, since $\{\fB \text{ is type } ({\rm I},0)\}\subset \mathcal{B}$, we get from \eqref{eq: local loops expansion with spin bc 3} (applying to $\mathcal{A}=\{\fB \text{ is type } ({\rm I},0)\}$) \begin{equation}
    \mathcal{Z}^{\zeta}_{\fB',\gamma}\big(\{\fB \text{ is type } ({\rm I},0\}\big)=(1-\kappa)\cdot \int_{\mathfrak{L}_\fB\in\Delta_\fB,\mathfrak{e}_\fB}\mathsf{M}(\mathsf{d}\mathfrak{L}_\fB)\cdot\Psi(\mathsf{d}\mathfrak{e}_\fB)\cdot \mathcal{Z}^\zeta_{\fB',\gamma-\mathfrak{L}_\fB}\Big(\mathcal{F}^{ 0}_*(\mathfrak{L}_\fB,\mathfrak{e}_\fB)\Big).\nonumber
\end{equation}
Thus, we finish the proof of \eqref{eq: type 1, 0 loop expansion old}, thereby establishing \eqref{eq: type 1, k loop expansion}.
\end{proof}
\begin{corollary}\label{cor: loop configuration expansion total measure}
   Recall Definition \ref{def: crossing event decomposition according to loops}.
   Let $\mathrm{NLC}$ denote the event that there exists no loop crossing through the annulus $\fB^{(2)}\setminus\fB$. For any event $\mathcal{A}$ depending only on the (loop and spin) configurations  inside $\fB^{(2)}\setminus\fB$, let $\mathcal{A}_*(\mathfrak{L}_\fB,\mathfrak{e}_{\fB})$ denote the spin configuration $\sigma$ on $\mathsf B'$ (note that what actually matters for $\mathcal A$ is $\sigma|_{\mathsf B^{(2)}}$) such that $\mathcal A$ occurs on $(\sigma, \mathfrak{ L}_\fB, \mathfrak e_{\fB})$. Let $\Sigma$ denote the collection of spin boundary conditions $\zeta$ on $\partial\fB'$.  Then for  boundary condition $\xi$ on $\partial\Lambda$,
   we have
   \begin{align}
     \mu_{\Lambda,\gamma}^{\xi}(\mathcal{A}\cap \mathrm{NLC})&=\int_{\mathfrak{L}_\fB\in\Delta_\fB}\mathsf{M}(\mathsf{d}\mathfrak{L}_\fB)\cdot  \Psi(\mathsf{d}\mathfrak{e}_\fB)\cdot  \frac{\mathcal{Z}_{\Lambda,\gamma-\mathfrak{L}_\fB}^{\xi}\big(\mathcal{A}_*(\mathfrak{L}_\fB,\mathfrak{e}_\fB)\big)}{\mathcal{Z}_{\Lambda,\gamma}^{\xi}}\label{eq: type 1, 0 loop expansion with local time bc}\\
     &= \int_{\zeta\in \Sigma}\mu_{\Lambda,\gamma}^{\xi}(\zeta)\cdot R_{\fB,A}(\zeta),\label{eq: type 1, 0 loop expansion with local time bc 2}
\end{align}
where \begin{equation}
    R_{\fB,\mathcal{A}}(\zeta)=\int_{\mathfrak{L}_\fB\in\Delta_\fB}\mathsf{M}(\mathsf{d}\mathfrak{L}_\fB)\cdot  \Psi(\mathsf{d}\mathfrak{e}_\fB)\cdot\frac{G_{\fB'}(\gamma-\mathfrak{L}_\fB)}{G_{\fB'}(\gamma)} \frac{{Z}_{\fB',\gamma-\mathfrak{L}_\fB}^{\zeta}\big(\mathcal{A}_*(\mathfrak{L}_\fB,\mathfrak{e}_\fB)\big)}{{Z}_{\fB',\gamma}^{\zeta}}.\label{eq: def of RB}
\end{equation}
 In the case that $\mathcal{A}$ is the whole set, we write $R_\fB(\zeta)=R_{\fB,\mathcal{A}}(\zeta)$. Then we have \begin{equation}
    \mu_{\Lambda,\gamma}^{\xi}\big(R_\fB(\zeta)\big)= \mu_{\Lambda,\gamma}^{\xi}(\mathrm{NLC})\le 1.\label{eq: full Rb expectation upper bound}
\end{equation}
\end{corollary}
\begin{proof}
    Note that $\xi$ is a local time boundary condition, and thus \eqref{eq: type 1, 0 loop expansion with local time bc} follows directly from the definition of the partition functions. Let $\mathcal{Z}_{\Lambda\setminus\fB',\gamma-\mathfrak{L}_\fB}^{\xi}\big(\zeta\big)$ denote the partition function on $\Lambda\setminus\fB'$ with local time condition $\gamma-\mathfrak{L}_\fB$, boundary condition $\xi$ on $\partial\Lambda$ and spin boundary condition $\zeta$ on $\partial\fB'.$ To prove \eqref{eq: type 1, 0 loop expansion with local time bc 2}, we calculate that \begin{align}
        \mathcal{Z}_{\Lambda,\gamma-\mathfrak{L}_\fB}^{\xi}\big(\mathcal{A}_*(\mathfrak{L}_\fB,\mathfrak{e}_\fB)\big)&=\int_{\zeta\in\Sigma}\mathcal{Z}_{\fB',\gamma-\mathfrak{L}_\fB}^{\zeta}\big(\mathcal{A}_*(\mathfrak{L}_\fB,\mathfrak{e}_\fB)\big)\cdot \mathcal{Z}_{\Lambda\setminus\fB',\gamma-\mathfrak{L}_\fB}^{\xi}\big(\zeta\big);\nonumber
    \end{align}Plugging into \eqref{eq: type 1, 0 loop expansion with local time bc} and noting that $\mathcal{Z}_{\Lambda\setminus\fB',\gamma-\mathfrak{L}_\fB}^{\xi}\big(\zeta\big)=\mathcal{Z}_{\Lambda\setminus\fB',\gamma}^{\xi}\big(\zeta\big)$ (since $\mathfrak L_B \in \Delta_\mathsf B$ only consists of loops in $\mathsf B^{(2)}$),
 we get that 
    \begin{align}
        \mu_{\Lambda,\gamma}^{\xi}(\mathcal{A}\cap \mathrm{NLC})&=\int_{\mathfrak{L}_\fB\in\Delta_\fB}\mathsf{M}(\mathsf{d}\mathfrak{L}_\fB)\cdot  \Psi(\mathsf{d}\mathfrak{e}_\fB)\cdot\frac{\int_{\zeta\in\Sigma}\mathcal{Z}_{\fB',\gamma-\mathfrak{L}_\fB}^{\zeta}\big(\mathcal{A}_*(\mathfrak{L}_\fB,\mathfrak{e}_\fB)\big)\cdot \mathcal{Z}_{\Lambda\setminus\fB',\gamma}^{\xi}\big(\zeta\big)}{\mathcal{Z}_{\Lambda,\gamma}^{\xi}}\nonumber\\
        &=\int_{\zeta\in\Sigma}\frac{\mathcal{Z}_{\fB',\gamma}^{\zeta}\cdot \mathcal{Z}_{\Lambda\setminus\fB',\gamma}^{\xi}\big(\zeta\big)}{{\mathcal{Z}_{\Lambda,\gamma}^{\xi}}}\int_{\mathfrak{L}_\fB\in\Delta_\fB}\mathsf{M}(\mathsf{d}\mathfrak{L}_\fB)\cdot  \Psi(\mathsf{d}\mathfrak{e}_\fB)\cdot\frac{\mathcal{Z}_{\fB',\gamma-\mathfrak{L}_\fB}^{\zeta}\big(\mathcal{A}_*(\mathfrak{L}_\fB,\mathfrak{e}_\fB)\big)}{\mathcal{Z}_{\fB',\gamma}^{\zeta}}\nonumber\\
        &=\int_{\zeta\in \Sigma}\mu_{\Lambda,\gamma}^{\xi}(\zeta)\cdot R_\fB(\zeta).\nonumber
    \end{align}
    Thus, we finish the proof of \eqref{eq: type 1, 0 loop expansion with local time bc 2} by applying Lemma \ref{lem: spin local time factorization}. Choosing $\mathcal{A}$ to be the whole set in \eqref{eq: type 1, 0 loop expansion with local time bc 2} gives \eqref{eq: full Rb expectation upper bound}.
\end{proof}
\begin{proof}[Proof of Lemma \ref{lem: type probability}]
We focus on the case $k=0$, and the other cases follow similarly. \eqref{eq: type probability} follows directly from \eqref{eq: type probability without aux}. Let $\mathfrak{T}$ denote the event that $\fB \text{ is Type } ({\rm II},0)$ conditioned on $\mathrm{Aux}_\mathsf B = 1$.
Then 
\begin{equation}
        \mu_{\Lambda,\beta-L^{\zz}}^{\xi}\Big(\{\fB \text{ is type } ({\rm II},0)\}\cap \Big(\mathcal{E}^{\mathsf{x}}_0\cup \mathcal{E}^{\mathsf{y}}_0\Big)^c \Big)
        \leq \kappa+(1-\kappa)\cdot  \mu_{\Lambda,\beta-L^{\zz}}^{\xi}\Big(\mathfrak{T}\cap \Big(\mathcal{E}^{\mathsf{x}}_0\cup \mathcal{E}^{\mathsf{y}}_0\Big)^c \Big).\label{eq: type probability expansion wrp to aux}
    \end{equation}
    Let $\bC_1>0$ be a large enough constant to be chosen later and let $\fB'=\fB^{(\bC_1)}$.
    Recall Definition \ref{def: crossing event decomposition according to loops}.
    By \eqref{eq: def of good box},
    we get that \begin{equation}\label{eq: type probability 1}
    \mu_{\Lambda,\beta-L^{\zz}}^{\xi}\Big(\mathcal{E}^{\mathsf{x}}_0\cup \mathcal{E}^{\mathsf{y}}_0\Big)\le \rho.
\end{equation}
Now, on the event $\Big(\mathcal{E}^{\mathsf{x}}_0\cup \mathcal{E}^{\mathsf{y}}_0\Big)^c$, we can explore the loops from $\fB$ and get $\mathfrak{L}_\fB$.
Let $\Sigma$ denote the collection of spin boundary conditions $\zeta$ on $\partial\fB^{(\bC_1)}$. By Lemma \ref{lem: type 1 loop expansion} (with $\fB'=\fB^{(\bC_1)}$), we get that \begin{align}
     &\mu_{\Lambda,\beta-L^{\zz}}^{\xi}\Big(\mathfrak{T}\cap \Big(\mathcal{E}^{\mathsf{x}}_0\cup \mathcal{E}^{\mathsf{y}}_0\Big)^c \Big)=\int_{\zeta\in \Sigma,\mathfrak{L}_\fB,\mathfrak{e}_\fB}\mu_{\Lambda,\beta-L^{\zz}}^{\xi}(\zeta)\cdot \frac{{Z}^{\zeta}_{\fB',\beta-L^{\zz}-\mathfrak{L}_\fB}}{{Z}^{\zeta}_{\fB',\beta-L^{\zz}}}\cdot\mathsf{M}(\mathsf{d}\mathfrak{L}_\fB)\cdot\Psi(\mathsf{d}\mathfrak{e}_\fB)\nonumber\\&\hspace{4cm}\times\frac{G_{\fB'}(\beta-L^{\zz}-\mathfrak{L}_\fB)}{G_{\fB'}(\beta-L^{\zz})} \cdot \mu^{\zeta}_{\fB^{(\bC_1)},\beta-L^{\zz}-\mathfrak{L}_\fB}\Big(\bigcup_{\substack{\mathsf{i}\in\{\xx,\yy\},\\k\in\{-1,1\}}}\big(\mathcal{E}^{\mathsf{i},k}_*( \mathfrak{L}_\fB,\mathfrak{e}_\fB)\big)^c\Big).\label{eq: type probability expansion to edge marks}
    \end{align}

In addition, by \eqref{eq: full Rb expectation upper bound}, we get that \begin{align}
      &\int_{\zeta\in \Sigma,\mathfrak{L}_\fB,\mathfrak{e}_\fB}\mu_{\Lambda,\beta-L^{\zz}}^{\xi}(\zeta)\cdot\frac{G_{\fB'}(\beta-L^{\zz}-\mathfrak{L}_\fB)}{G_{\fB'}(\beta-L^{\zz})}\cdot \frac{{Z}^{\zeta}_{\fB',\beta-L^{\zz}-\mathfrak{L}_\fB}}{{Z}^{\zeta}_{\fB',\beta-L^{\zz}}}\times\mathsf{M}(\mathsf{d}\mathfrak{L}_\fB)\cdot\Psi(\mathsf{d}\mathfrak{e}_\fB)\le 1.\label{eq: loop configuration expansion total measure}
    \end{align}
For a spin wave $f$ and a spin boundary condition $\zeta$ on $\partial\fB^{(\bC_1)}$, let $\mu^{\zeta,f}_{\fB^{(\bC_1)},\beta-L^{\zz}}$ denote the measure of spin configuration on $\fB^{(\bC_1)}$ such that \begin{equation}
    \mu^{\zeta,f}_{\fB^{(\bC_1)},\beta-L^{\zz}}(\sigma)=\frac{1}{{Z}^{\zeta,f}_{\fB^{(\bC_1)},\beta-L^{\zz}}} \exp\Big(H^{\zeta,f}_{\fB^{(\bC_1)},\beta-L^{\zz}-\mathfrak{L}_\fB}(\sigma)\Big)\cdot \mathbf{1}_{\sigma|_{\partial\fB^{(\bC_1)}}=\zeta},\label{eq: tilted spin measure} 
\end{equation} where the Hamiltonian is defined as
\begin{equation}
    H^{\zeta,f}_{\fB^{(\bC_1)},\beta-L^{\zz}}(\sigma)=\sum_{u\sim v}2\sqrt{(\beta_u-L^{\zz}_u)\cdot (\beta_v-L^{\zz}_v) }\cos\Big(\big(\sigma_u+f_u\big)-\big(\sigma_v+f_v\big)\Big)\label{eq: hamiltonian for tilted spin measure}
\end{equation}and the partition function is defined as \begin{equation}
    {Z}^{\zeta,f}_{\fB^{(\bC_1)},\beta-L^{\zz}}=\sum_{\sigma}\exp\Big(H^{\zeta,f}_{\fB^{(\bC_1)},\beta-L^{\zz}-\mathfrak{L}_\fB}(\sigma)\Big)\cdot \mathbf{1}_{\sigma|_{\partial\fB^{(\bC_1)}}=\zeta}.\label{eq: partition function for tilted spin measure}
\end{equation}

Let $f_1$  be the function that is harmonic on $\fB'\setminus\fB^{(2)}$, coincides with $F_0$ inside $\fB^{(2)}$ and vanishes on $\partial\fB^{(\bC_1)}$.  Then for any $\sigma_\fB\in \big(\mathcal{E}^{\mathsf{i},k}_*( \mathfrak{L}_\fB,\mathfrak{e}_\fB)\big)^c$, we have $\sigma_\fB-kf_1\in \big(\mathcal{E}^{\mathsf{i}}_*( \mathfrak{L}_\fB,\mathfrak{e}_\fB)\big)^c$.
Thus, we get that for any $k\in\{1,-1\}$ and $\mathsf{i}\in\{\xx,\yy\}$ 
\begin{equation}
     \mu^{\zeta}_{\fB^{(\bC_1)},\beta-L^{\zz}-\mathfrak{L}_\fB}\Big(\big(\mathcal{E}^{\mathsf{i},k}_*( \mathfrak{L}_\fB,\mathfrak{e}_\fB)\big)^c\Big)=\ \mu^{\zeta,kf_1}_{\fB^{(\bC_1)},\beta-L^{\zz}-\mathfrak{L}_\fB}\Big(\big(\mathcal{E}^{\mathsf{i}}_*( \mathfrak{L}_\fB,\mathfrak{e}_\fB)\big)^c\Big).\label{eq: type 2 + convert to 0}
\end{equation}
By Lemma \ref{lem: push forward probability control} below (with $ \mathsf C = \fB^{(\bC_1)}$ and $f=kf_1$), we get that \begin{align}
     &\mu^{\zeta,kf_1}_{\fB^{(\bC_1)},\beta-L^{\zz}-\mathfrak{L}_\fB}\Big(\big(\mathcal{E}^{\mathsf{i}}_*( \mathfrak{L}_\fB,\mathfrak{e}_\fB)\big)^c\Big)\nonumber\\\le~& \mu^{\zeta}_{\fB^{(\bC_1)},\beta-L^{\zz}-\mathfrak{L}_\fB}\Big(\big(\mathcal{E}^{\mathsf{i}}_*( \mathfrak{L}_\fB,\mathfrak{e}_\fB)\big)^c\Big)+\bC_2\sqrt{\sum_{u\sim v\in \fB^{(\bC_1)}}\big(f_1(u)-f_1(v)\big)^2}.\label{eq: Pinsker application 1}
    \end{align}
By Lemma \ref{lem: harmonic function energy bound for vanishing boundary} below, we get that \begin{equation}
        \sum_{u\sim v\in \fB^{(\bC_1)}} \big(f_1(u)-f_1(v)\big)^2\le \frac{C_{10} }{\iota^2  \log(\bC_1\delta)}+\frac{C_{10}\bC_1^2L^2}{\iota^2\delta^2 M^2},\label{eq: minimal dirichlet energy 1}
    \end{equation}
    where $L$ is the side length of the box $\fB.$ We can choose the constant $\bC_1>0$ to be large enough such that $\frac{C_{10} }{\iota^2  \log(\bC_1\delta)}\le \frac{0.01\rho^2}{\bC_2^2}$. In addition, note that $L\le M^{0.1}$, we can choose $M_0>0$ large enough such that $\frac{C_{10}\bC_1^2L^2}{\iota^2\delta^2 M_0^2}\le \frac{0.01\rho^2}{\bC_2^2}.$ Therefore, we conclude that \begin{equation}
        \sum_{u\sim v\in \fB^{(\bC_1)}} \big(f_1(u)-f_1(v)\big)^2\le \frac{0.02\rho^2}{\bC_2^2}.\label{eq: minimal dirichlet energy}
    \end{equation} Plugging \eqref{eq: minimal dirichlet energy} into \eqref{eq: Pinsker application 1} and recalling $\bC_1=C_7/(\delta\rho)$, we obtain that \begin{equation}\label{eq: no percolation in type 2 tilted 1}
          \mu^{\zeta,kf_1}_{\fB^{(\bC_1)},\beta-L^{\zz}-\mathfrak{L}_\fB}\Big(\big(\mathcal{E}^{\mathsf{i}}_*( \mathfrak{L}_\fB,\mathfrak{e}_\fB)\big)^c\Big)\le \mu^{\zeta}_{\fB^{(\bC_1)},\beta-L^{\zz}-\mathfrak{L}_\fB}\Big(\big(\mathcal{E}^{\mathsf{i}}_*( \mathfrak{L}_\fB,\mathfrak{e}_\fB)\big)^c\Big)+ 0.2\rho.
    \end{equation}
    Summing \eqref{eq: no percolation in type 2 tilted 1} over $k\in\{-1,0,1\}$ and $\mathsf{i}\in\{\xx,\yy\}$ and combining with \eqref{eq: type 2 + convert to 0}, we get that 
    \begin{align}
          &\mu^{\zeta}_{\fB^{(\bC_1)},\beta-L^{\zz}-\mathfrak{L}_\fB}\Big(\bigcup_{\substack{\mathsf{i}\in\{\xx,\yy\},\\k\in\{-1,1\}}}\big(\mathcal{E}^{\mathsf{i},k}_*( \mathfrak{L}_\fB,\mathfrak{e}_\fB)\big)^c\Big)\nonumber\\\le~& 2\mu^{\zeta}_{\fB^{(\bC_1)},\beta-L^{\zz}-\mathfrak{L}_\fB}\Big(\big(\mathcal{E}^{\mathsf{x}}_*( \mathfrak{L}_\fB,\mathfrak{e}_\fB)\big)^c\Big)+2\mu^{\zeta}_{\fB^{(\bC_1)},\beta-L^{\zz}-\mathfrak{L}_\fB}\Big(\big(\mathcal{E}^{\mathsf{y}}_*( \mathfrak{L}_\fB,\mathfrak{e}_\fB)\big)^c\Big)+ 0.8\rho.\label{eq: no percolation in type 2 tilted 2}
    \end{align}
    Plugging \eqref{eq: no percolation in type 2 tilted 2} into \eqref{eq: type probability expansion to edge marks}, we get that
    \begin{align}
     &\mu_{\Lambda,\beta-L^{\zz}}^{\xi}\Big(\mathfrak{T}\cap \Big(\mathcal{E}^{\mathsf{x}}_0\cup \mathcal{E}^{\mathsf{y}}_0\Big)^c \Big)\le\int_{\zeta\in \Sigma,\mathfrak{L}_\fB,\mathfrak{e}_\fB}\mu_{\Lambda,\beta-L^{\zz}}^{\xi}(\zeta)\cdot \frac{{Z}^{\zeta}_{\fB^{(\bC_1)},\beta-L^{\zz}-\mathfrak{L}_\fB}}{{Z}^{\zeta}_{\fB^{(\bC_1)},\beta-L^{\zz}}}\cdot \mathsf{M}(\mathsf{d}\mathfrak{L}_\fB)\cdot\Psi(\mathsf{d}\mathfrak{e}_\fB)\nonumber
     \\&\times \frac{G_{\fB'}(\beta-L^{\zz}-\mathfrak{L}_\fB)}{G_{\fB'}(\beta-L^{\zz})}
     \cdot\Big[2\mu^{\zeta}_{\fB^{(\bC_1)},\beta-L^{\zz}-\mathfrak{L}_\fB}\Big(\big(\mathcal{E}^{\mathsf{x}}_*( \mathfrak{L}_\fB,\mathfrak{e}_\fB)\big)^c\Big)\nonumber\\&\hspace{3.5cm}+2\mu^{\zeta}_{\fB^{(\bC_1)},\beta-L^{\zz}-\mathfrak{L}_\fB}\Big(\big(\mathcal{E}^{\mathsf{y}}_*( \mathfrak{L}_\fB,\mathfrak{e}_\fB)\big)^c\Big)+ 0.8\rho\Big].\label{eq: type probability expansion to edge marks 2}
    \end{align}
    By \eqref{eq: loop configuration expansion total measure}, \eqref{eq: type 1, 0 loop expansion with local time bc 2}, and \eqref{eq: type probability 1}, we get that
    \begin{equation}
    \begin{aligned}
      &\int_{\zeta\in \Sigma,\mathfrak{L}_\fB,\mathfrak{e}_\fB}\mu_{\Lambda,\beta-L^{\zz}}^{\xi}(\zeta)\cdot \frac{{Z}^{\zeta}_{\fB,\beta-L^{\zz}-\mathfrak{L}_\fB}\Big(\big(\mathcal{E}^{\mathsf{x}}_*( \mathfrak{L}_\fB,\mathfrak{e}_\fB)\big)^c\cup \big(\mathcal{E}^{\mathsf{y}}_*( \mathfrak{L}_\fB,\mathfrak{e}_\fB)\big)^c\Big)}{{Z}^{\zeta}_{\fB,\beta-L^{\zz}}}
      \\&\qquad\qquad\qquad
      \times\frac{G_{\fB'}(\beta-L^{\zz}-\mathfrak{L}_\fB)}{G_{\fB'}(\beta-L^{\zz})}\cdot \mathsf{M}(\mathsf{d}\mathfrak{L}_\fB)\cdot\Psi(\mathsf{d}\mathfrak{e}_\fB)  \nonumber
      \\ \stackrel{\eqref{eq: loop configuration expansion total measure} }{\le} & 1- \int_{\zeta\in \Sigma,\mathfrak{L}_\fB,\mathfrak{e}_\fB}\mu_{\Lambda,\beta-L^{\zz}}^{\xi}(\zeta)\cdot \frac{{Z}^{\zeta}_{\fB,\beta-L^{\zz}-\mathfrak{L}_\fB}\Big[\mathcal{E}^{\mathsf{x}}_*( \mathfrak{L}_\fB,\mathfrak{e}_\fB)\cap \mathcal{E}^{\mathsf{y}}_*( \mathfrak{L}_\fB,\mathfrak{e}_\fB)\Big]}{{Z}^{\zeta}_{\fB,\beta-L^{\zz}}}
      \\&\qquad\qquad\qquad
      \times \frac{G_{\fB'}(\beta-L^{\zz}-\mathfrak{L}_\fB)}{G_{\fB'}(\beta-L^{\zz})}\cdot \mathsf{M}(\mathsf{d}\mathfrak{L}_\fB)\cdot\Psi(\mathsf{d}\mathfrak{e}_\fB)\nonumber
      \\\stackrel{\eqref{eq: type 1, 0 loop expansion with local time bc 2}}{=}& 1-\mu_{\Lambda,\beta-L^{\zz}}^{\xi}\Big((\mathcal{E}^\xx_0\cup\mathcal{E}^\yy_0)^c\Big)\stackrel{\eqref{eq: type probability 1}}{\le} \rho.
    \end{aligned}
    \end{equation}
    Plugging into \eqref{eq: type probability expansion to edge marks 2}, we get that 
 \begin{align}
        &\mu_{\Lambda,\beta-L^{\zz}}^{\xi}\Big(\mathfrak{T}\cap \Big(\mathcal{E}^{\mathsf{x}}_0\cup \mathcal{E}^{\mathsf{y}}_0\Big)^c \Big)\le 4\rho+ 0.8\rho.\label{eq: no percolation in type 2 tilted 3}
    \end{align}
    Combining \eqref{eq: no percolation in type 2 tilted 3} with \eqref{eq: type probability 1} and plugging into \eqref{eq: type probability expansion wrp to aux}, we finish the proof of Lemma \ref{lem: type probability}.
\end{proof}
\begin{lemma}\label{lem: harmonic function energy bound for vanishing boundary}
    Let $\fB$ be a box with side length $L$ and let $\fB'=\fB^{(C_9)}$ be an expansion of $\fB$. Recall $F_0$ from Definition \ref{def: base spin wave}. Recall that  $f_1$  be the function that is harmonic on $\fB'\setminus\fB^{(2)}$, coincides with $F_0$ inside $\fB^{(2)}$ and vanishes on $\partial\fB^{(\bC_1)}$.  Then there exists an absolute constant $C_{10}>0$ such that \begin{equation}
        \sum_{u\sim v\in \fB'}\big(f_1(u)-f_1(v)\big)^2\le \frac{C_{10} }{\iota^2  \log(C_9\delta)}+\frac{C_{10}C_9^2L^2}{\iota^2\delta^2 M^2}.\nonumber
    \end{equation}
\end{lemma}
\begin{proof}
For any function $f$, let $\mathcal{E}(f)=\sum_{u\sim v\in \fB'}\big(f(u)-f(v)\big)^2$ be its Dirichlet energy on $\fB'$.  Since $f_1$ is harmonic on $\fB'\setminus\fB^{(2)}$, we get that $f_1$ is the function with minimal Dirichlet energy that coincides with $F_0$ inside $\fB^{(2)}$ and vanishes on $\partial\fB'$.

Let $g_1$ be the harmonic function on $\fB'\setminus\fB^{(2)}$ that takes value 1 on the outer boundary $\partial\fB'$ and takes value 0 inside $\fB^{(2)}$. Then we have (see for example \cite[Equation (9.6)]{Velenik}) \begin{equation}
        \mathcal{E}(g_1)\le \frac{\bC_1}{\log(C_9\delta)}.\label{eq: dirichlet energy for g1}
    \end{equation}    Let $F_\fB=\frac{1}{|\partial\fB^{(2)}|}\sum_{v\in\partial\fB^{(2)}}F_0(v)$ be the average of $f$ on $\partial\fB^{(2)}$. Let $g_2$ be the function that coincides with $F_0$ inside $\fB^{(2)}$,
    takes value $F_\fB$ on $\partial\fB'$ and interpolates from $\partial\fB^{(2)}$ to $\partial\fB^{(\bC_1)}$ linearly in the $ |\cdot|_{\infty} $ distance to $\partial\fB^{(2)}$. Then we have $|g_2(u)-g_2(v)|\le \frac{\bC_2 L}{\iota M}\cdot \frac{1}{C_9\delta L}=\frac{\bC_2 }{C_9\delta\iota M}$ for $u\sim v\in \fB'\setminus\fB^{(2)}$. In addition, for $u\sim v\in\fB^{(2)}$, we have $|g_2(u)-g_2(v)|\le \frac{1}{\iota M}$. Therefore, we conclude  that \begin{equation}
        \mathcal{E}(g_2)\le (C_9\delta L+L)^2\cdot (\frac{\bC_2 }{\delta\iota M})^2\le \frac{C_9^2\bC_2 L^2}{\iota^2\delta^2 M^2}.\label{eq: dirichlet energy for g2}
    \end{equation} Since $g_2-F_\fB\cdot g_1$ coincides with $F_0$ inside $\fB^{(2)}$ and vanishes on $\partial\fB'$, we get from \eqref{eq: dirichlet energy for g1} and \eqref{eq: dirichlet energy for g2} that \begin{equation}
        \mathcal{E}(f_1)\le \mathcal{E}(g_2-F_\fB\cdot g_1)\le 2\mathcal{E}(g_2)+2(F_\fB)^2\cdot \mathcal{E}(g_1)\le \frac{2\bC_{1} }{\iota^2  \log(C_9\delta)}+\frac{2\bC_2 C_9^2L^2}{\iota^2\delta^2 M^2}
    \end{equation}where in the last inequality we used $|F_\fB|\le 1$.
\end{proof}

\begin{lemma}\label{lem: push forward probability control}
Fix a box $\fC$ and fix a spin wave $f$ such that $f$ is a constant on $\partial\fC.$ There exists a constant $C_{11}$ depending only on $\beta_0$ such that for any spin event $\mathcal{E}$ depending only on the spin configuration in $\fC$ and for any spin boundary condition $\zeta$ on $\partial\fC$, we have (recall \eqref{eq: hamiltonian for tilted spin measure})
\begin{equation}\label{eq: probability difference for tilted model}
    |\mu^\zeta_{\fC,\beta}(\mathcal{E})-\mu^{\zeta,f}_{\fC,\beta}(\mathcal{E})|\le C_{11}\sqrt{\sum_{u\sim v}(f_u-f_v)^2}.
\end{equation}
\end{lemma}
\begin{proof}
It suffices to control the total variation distance between the two measures $\mu^\zeta_{\fC,\beta}$ and $\mu^{\zeta,f}_{\fC,\beta}$. To do so, we want to apply Pinsker's inequality, for which we need to first control the relative entropy between the two measures $\mu^\zeta_{\fC,\beta}$ and $\mu^{\zeta,f}_{\fC,\beta}$. Since $f$ is constant on the boundary, we have that 
\begin{equation}
Z^{\zeta,f}_{\fC,\beta}=Z^{\zeta+f}_{\fC,\beta}=Z^{\zeta}_{\fC,\beta}\stackrel{(*)}{=}Z^{\zeta-f}_{\fC,\beta}=Z^{\zeta,-f}_{\fC,\beta}\label{eq: partition function equation for tilted model}
\end{equation}where $(*)$ follows from symmetry.
    Hence, we calculate the relative entropy (recall that $ c^{\beta}_{uv}=2J_{uv}\sqrt{\beta_u\beta_v} $)
    \begin{align}
        &D(\mu^{\zeta}_{\fC,\beta}~||~\mu^{\zeta,f}_{\fC,\beta})= \Big\langle -H^{\zeta,f}_{\fC,\beta}(\sigma)+H^{\zeta}_{\fC,\beta}(\sigma)\Big\rangle^\zeta_{\fC,\beta}\nonumber\\=~&\bigg\langle \sum_{u\sim v}c^{\beta}_{uv}\Big[\cos\Big(\big(\sigma_u-\sigma_v\big)+\big(f_u-f_v\big)\Big)-\cos\Big(\big(\sigma_u-\sigma_v\big)\Big)\Big]\bigg\rangle^\zeta_{\fC,\beta}. \label{eq: relative entropy expansion}
    \end{align}
    Since the relative entropy is non-negative, we get from~\eqref{eq: relative entropy expansion} that
    \begin{align}
        \label{eq: relative entropy symmetrize}
        &D(\mu^{\zeta}_{\fC,\beta}~||~\mu^{\zeta,f}_{\fC,\beta})\le D(\mu^{\zeta}_{\fC,\beta}~||~\mu^{\zeta,f}_{\fC,\beta})+D(\mu^{\zeta}_{\fC,\beta}~||~\mu^{\zeta,-f}_{\fC,\beta})\nonumber\\=~& \Big\langle H^{\zeta,f}_{\fC,\beta}(\sigma)-H^{\zeta}_{\fC,\beta}(\sigma)\Big\rangle^\zeta_{\fC,\beta}+\Big\langle H^{\zeta,-f}_{\fC,\beta}(\sigma)-H^{\zeta}_{\fC,\beta}(\sigma)\Big\rangle^\zeta_{\fC,\beta}\nonumber\\=~&\bigg\langle \sum_{u\sim v}{c^{\beta}_{uv}}\Big[\cos\Big(\big(\sigma_u-\sigma_v\big)+\big(f_u-f_v\big)\Big)+\cos\Big(\big(\sigma_u-\sigma_v\big)-\big(f_u-f_v\big)\Big)\nonumber\\&\hspace{2.5cm}-2\cos\Big(\big(\sigma_u-\sigma_v\big)\Big)\Big]\bigg\rangle^\zeta_{\fC,\beta}\nonumber\\
        =~&\bigg\langle \sum_{u\sim v}{c^{\beta}_{uv}}\Big[2\cos\big(\sigma_u-\sigma_v\big)\cdot\big(\cos(f_u-f_v)-1\big)\Big]\bigg\rangle^\zeta_{\fC,\beta}\nonumber\\\le~&\bigg\langle \sum_{u\sim v}{c^{\beta}_{uv}}\Big[\bC_1(f_u-f_v)^2\Big]\bigg\rangle^\zeta_{\fC,\beta}\le \bC_2\sum_{u\sim v}(f_u-f_v)^2.
    \end{align}
    Combining~\eqref{eq: relative entropy symmetrize} with the Pinsker's inequality, we obtain that
    \begin{equation*}
        ||\mu^{\zeta,f}_{\fC,\beta}-\mu^\zeta_{\fC,\beta}||_{\mathrm{TV}}\le \sqrt{\frac{1}{2}D(\mu^{\zeta}_{\fC,\beta}~||~\mu^{\zeta,f}_{\fC,\beta})}\le \bC_3\sqrt{\sum_{u\sim v}(f_u-f_v)^2}.
    \end{equation*}Thus we finish the proof of \eqref{eq: probability difference for tilted model}.
\end{proof}
\subsection{An empowered spin wave argument}
In this subsection, we implement the spin wave argument. Since type 1 boxes do not contribute to connectivities outside their expansions (by a factor of $1+2\delta$), we can first explore the loops from the type 1 boxes and consider the measure on the remaining graph, motivating the following definition (also recall Definitions \ref{def: crossing event decomposition according to loops} and \ref{def: def of good box}).
\begin{definition}\label{def: final spin wave frozen part}
    Fix a sequence of types $\{T_\fB\}$ for good boxes $\fB$. For a type 1 (good) box $\fB$, we explore the loops $\mathfrak{L}_\fB$ from $\fB$ and we also explore the edge marks $\mathfrak{e}_\fB$ on $\fB^{(2)}$. Let $\Upsilon^0\Big(\{\mathfrak{L}_\fB\},\{T_\fB\}, \{\mathfrak{e}_{\fB}\}\Big)$ be the event that the sequence of types is given by $\{T_\fB\}$, and that we get the loops $\mathfrak{L}_\fB$ and the edge marks $\mathfrak{e}_\fB$ for each type 1 box $\mathsf B$. Let $\mathfrak{B}_{\rm I}$ denote the collection of type 1 boxes and let $\mathfrak{B}_{\rm II}$ denote the collection of type 2 boxes. Let $\Lambda'=\Lambda\setminus(\cup_{\fB\in\mathfrak{B}_{\rm I}}\fB)$ be the remaining region and let $\beta'=\beta-L^{\zz}-\sum_{\fB\in\mathfrak{B}_{\rm I}}\mathfrak{L}_\fB$ be the remaining local time.
\end{definition}
Note that for any $u\in\fB$ and $\fB\in\mathfrak{B}_{\rm I}$, we have $\beta'_u=0$. Then the energy cost when applying the spin wave $F_0$ is restricted to the region $\Lambda'.$
We start with controlling the size of $\Lambda'$.
\begin{lemma}\label{lem: large deviation for type 2 box}
    Recall from Definition \ref{def: good box exploration} that
 $\mathfrak{S}^{\ell_{*}}$ is the collection of good boxes. Recall Definition \ref{def: final spin wave frozen part} from above.
 Then for any boundary condition $\xi$ on $\partial\Lambda$, we have
\begin{equation}\label{eq: large deviation for type 2 box}
        \mu_{\Lambda,\beta-L^{\zz}}^{\xi}\Big(\sum_{\fB\in\mathfrak{B}_{\rm II}}|\fB|\ge 6(\kappa+\rho)^{\frac{2}{3}} |\Lambda|\Big)  \le  (\kappa+\rho)^{\frac{1}{3}}.
    \end{equation}
\end{lemma}
\begin{proof} By Lemma \ref{lem: type probability}, we get that \begin{equation}
        \mu_{\Lambda,\beta-L^{\zz}}^{\xi}\Big(\sum_{\fB\in\mathfrak{B}_{\rm II}}|\fB|\Big)=\sum_{\fB\in \mathfrak{S}^{\ell_{*}}} \mu_{\Lambda,\beta-L^{\zz}}^{\xi}\Big(\fB \text{ is type } ({\rm II},0)\Big)\cdot |\fB|\le 6(\kappa+\rho)\cdot \sum_{\fB\in \mathfrak{S}^{\ell_{*}}}|\fB| \le 6(\kappa+\rho) |\Lambda|.\nonumber
    \end{equation} Combined with Markov's inequality, it yields \eqref{eq: large deviation for type 2 box}. 
\end{proof}
\begin{remark}
    Unlike Lemma \ref{lem: covering property for good boxes}, we are not able to prove a decaying probability with respect to $M$ in \eqref{eq: large deviation for type 2 box}. This is because we are unable to prove Lemma \ref{lem: type probability} (originated from \eqref{eq: def of good box}) under a local spin boundary condition. Due to the lack of independence, we are unable to derive concentration bounds beyond Markov's inequality.
\end{remark}
\begin{definition}\label{def: good sequences and types}
    We say a sequence of good boxes $\{\fB\}$ and a sequence of types $\{T_\fB\}$ are excellent if the following properties  hold:
    \begin{enumerate}
        \item $|\Lambda\setminus(\cup_{\fB\in\mathfrak{S}^{\ell_{*}}}\fB)|\le C_6\delta |\Lambda|.$
        \item $\sum_{\fB\in\mathfrak{B}_{\rm II}}|\fB|\leq 6(\rho+\kappa)^{\frac{2}{3}} |\Lambda|.$
    \end{enumerate}
\end{definition}
Combining Lemmas \ref{lem: covering property for good boxes} and \ref{lem: large deviation for type 2 box}, we see that $\{\fB\}$ and $\{T_\fB\}$ are excellent
with probability close to 1 (by choosing $M>0$ large enough and $\rho,\kappa>0$ small enough) and thus we may assume this. In addition, we may slightly abuse the notation so that boxes (and types) in an excellent sequence $\{\fB\}$ and $\{T_\fB\}$ are still called good (instead of excellent).

We want to apply the spin wave $F_0$ on the remaining graph $\Lambda'= \Lambda\setminus(\cup_{\fB\in\mathfrak{B}_1}\fB) $ and show that there is a positive probability for the two boundaries to have different signs, for which we need to recover the exploration information $\Upsilon^0$ from the image; this is a new challenge that was not present in spin wave arguments in the literature. For a type 1 box $\fB$, by Definition \ref{def: type for good boxes}, we know that it is type $({\rm I},1)$ in the image. However, for a type 2 box $\fB$, we need to make sure it is type $(\rm II,1)$ in the image. Thus, we introduce the following map for type 2 boxes.
\begin{definition}\label{def: type 2 coupling}
    Fix a good box $\fB$, a constant angle $\alpha$ and a spin boundary condition $\zeta$ on $\partial\fB^{(20)}$. 
    {Recall Definition \ref{def: type for good boxes}.
    Let $\Xi(\fB,\zeta)$ denote the collection of spin $\sigma_\fB$, edge mark $\mathfrak{e}_\fB$ and loop $\mathfrak{L}_\fB$ (in the case that $\partial\fB^{(2)}\con{\xx\yy}\partial\fB$ under $(\sigma_\fB,\mathfrak{e}_\fB)$, we use the convention that $\mathfrak{L}_\fB=\emptyset$ since $\partial\fB^{(2)}\con{\xx\yy}\partial\fB$ already implies $\fB$ is  type $({\rm II},0)$)} configurations on $\fB^{(20)}$ with spin boundary condition $\zeta$  such that $\fB$ is type $({\rm II},0)$. In addition, let $\Xi^+(\fB,\zeta, \alpha)$ denote the collection of configurations on $\fB^{(20)}$ with boundary condition $\zeta+\alpha$ such that $\fB$ is type $({\rm II},1)$. Let $\tau^+_{\fB,\zeta,\alpha}$ be a map from $\Xi(\fB,\zeta)$ to $\Xi^+(\fB,\zeta, \alpha)$ such that regardless of the preimage, the image is chosen in $\Xi^+(\fB, \zeta, \alpha)$ uniformly at random. Similarly, we denote by  $\Xi^-(\fB,\zeta, \alpha)$ the collection of configurations on $\fB^{(20)}$ with spin boundary condition $\zeta-\alpha$ such that $\fB$ is type $({\rm II},-1)$. We similarly define $\tau^-_{\fB,\zeta,\alpha}$ to be a random map from $\Xi(\fB,\zeta)$ to $\Xi^-(\fB,\zeta, \alpha)$.
\end{definition}
\begin{remark}
    We remark here that $\tau^+_{\fB,\zeta,\alpha}$ is a fractional mapping (similarly for $ \tau^-_{\mathsf{B}, \zeta, \alpha} $), and the probability density of the image with respect to the pre-image is defined as
    \begin{equation*}
        \frac{(\tau^+_{\fB,\zeta,\alpha})^* \mu^{\zeta+\alpha}_{\fB^{(20)},\beta}\big({\mathrm d}{(\sigma_\fB,\mathfrak{e}_\fB,\mathfrak{L}_\fB)}\big)}{\mu^{\zeta}_{\fB^{(20)},\beta}\big({\mathrm d}{(\sigma_\fB,\mathfrak{e}_\fB,\mathfrak{L}_\fB)}\big)} = \frac{\mu^{\zeta+\alpha}_{\fB^{(20)},\beta}\Big(\Xi(\fB,\zeta+\alpha)\Big)}{\mu^{\zeta}_{\fB^{(20)},\beta}\Big(\Xi(\fB,\zeta)\Big)},
    \end{equation*}
    where $ (\tau^+_{\fB,\zeta,\alpha})^* \mu^{\zeta+\alpha}_{\fB^{(20)},\beta}(\cdot):=\mu^{\zeta+\alpha}_{\fB^{(20)},\beta}(\tau^+_{\fB,\zeta,\alpha}(\cdot)) $ is the induced pull-back measure of the spin configuration $ \sigma $ under the mapping $ \tau^+_{\fB,\zeta,\alpha} $.
\end{remark}
    Recall that $F_0$ is defined in \eqref{eq: def of base spin wave}. To make sure a type $({\rm II},0)$ good box is a type $({\rm II},1)$ good box in the image, we want to apply the map in Definition \ref{def: type 2 coupling}. Thus we need to modify $F_0$ near type $({\rm II},0)$ good boxes $\fB$ to ensure that it is constant on the boundary $\partial\fB^{(20)}$. 
\begin{definition}\label{def: final spin wave global part}
Let $F_2$  be the function that coincides with $F_0$ on $\Lambda'\setminus(\cup_{\fB\in\mathfrak{B}_{\rm II}}\fB^{(40)})$ and takes the constant value $\frac{\sum_{u\in \partial\fB^{(20)}}F_0(u)}{|\partial\fB^{(20)}|}$ (denoted as $F_2(\mathsf B)$) on the boundary $\partial\fB^{(20)}$, and interpolates linearly in the annulus $\fB^{(40)}\setminus\fB^{(20)}$ for each type 2 box $\fB\in\mathfrak{B}_{\rm II}$.
By Lemma \ref{lem: covering property for good boxes}, we know that $\{\fB^{(40)}\}$ is a collection of disjoint expansions of good boxes and thus $F_2$ is well-defined.
Let $F_3^+$ coincide with $F_2$ on ${\Lambda''=}\Lambda'\setminus(\cup_{\fB\in\mathfrak{B}_{\rm II}}\fB^{(20)})$. Furthermore, on each type 2 good box $\fB\in\mathfrak{B}_{\rm II}$, let $F_3^+$ equal the map $\tau^+_{\fB,\xi,F_2(\fB)}$ on $\fB^{(20)}$.
Similarly, we can define $F_3^-$ to be the map that coincides with $-F_2$ on ${\Lambda''}$ and equals the map $\tau^-_{\fB,\xi,-F_2(\fB)}$ on $\fB^{(20)}$.
\end{definition}
By Definition \ref{def: final spin wave global part}, we can recover the type of each good box from the image after applying the spin wave $F_3^+$.
This is because, as we discussed, if in the image we see a type $(\text{I}, 1)$ good box, then the preimage is a type $ (\text{I}, 0) $ good box, and if in the image we see a type $(\text{II}, 1)$ good box, then the preimage is a type $ (\text{II}, 0) $ good box; in addition, these are the only two possibilities. 
Next, we control the energy cost when applying $F_3^+$. The energy cost can be split into two parts: the first part is the cost when applying $F_2$ on ${\Lambda''}$; the second part is the cost when applying the maps $\tau^{\pm}_{\fB,\xi,F_2(\fB)}$ on each type 2 good box $\fB\in\mathfrak{B}_{\rm II}$, which we control below.
\begin{lemma}\label{lem: type 2 coupling}
    Let $\gamma=\beta-L^{\zz}.$
    There exists a constant $C_{12}>0$ such that for any good box $\fB\in\mathfrak{B}_{\rm II}$ such that $L\le M^{0.1}$ where $L$ is the side-length of $\fB$, we have
    \begin{align}
    &\int_{\zeta\in \Sigma}\mu^{\xi}_{\Lambda,\gamma}(\zeta)\cdot\log^-\Big(\frac{ \mu^{\zeta+F_2(\fB)}_{\fB^{(20)},\gamma} \Big(\Xi(\fB,\zeta+F_2(\fB))\Big)\times  \mu^{\zeta-F_2(\fB)}_{\fB^{(20)},\gamma} \Big(\Xi(\fB,\zeta-F_2(\fB))\Big)}{ \Big[\mu^{\zeta}_{\fB^{(20)},\gamma} \big(\Xi(\fB,\zeta)\big)\Big]^2}\Big)\ge -\frac{C_{12}\rho L^2}{\kappa^2\iota^2M^2}\nonumber
\end{align}where $\log^-(x)=\min\{0,\log(x)\}.$
\end{lemma}
The first step in proving Lemma \ref{lem: type 2 coupling} is to apply Lemma \ref{lem: type 1 loop expansion} to expand the Type $\mathrm{II}$ probability in terms of the probability of event $\mathcal{F}^k_*$. As hinted in \eqref{eq: tilted spin measure}, we will compare the measures $\mu^{\zeta,kf(\fB)}_{\fB^{(20)},\gamma}$ for $k=-1,0,1$. Therefore, we start from the following lemma.
\begin{lemma}
    Fix $\mathfrak{L}_\fB$ and $\mathfrak{e}_\fB$. Let $\zeta$ be a spin boundary condition on $\partial\fB'$. Let $f$ be a spin wave such that $f|_{\partial\fB'}$ is constant.
    Let ${c^{(\zz)}_{uv}}=\sqrt{(\beta_u-\mathfrak{L}^\zz_u)\cdot (\beta_v-\mathfrak{L}^\zz_v)}$ and let \begin{equation}
        W(f,\sigma)=\sum_{u\sim v\in \fB'}{c^{(\zz)}_{uv}}\cdot \sin(\sigma_u-\sigma_v)(f_u-f_v).\label{eq: def of W}
    \end{equation} Assume that $|f_u-f_v|\le \frac{C_{13}}{\iota M}$ and the side length of $\mathsf B$ satisfies $L\le M^{0.1}$. Then there exists a constant $C_{14}>0$ depending on $C_{13},\beta_0$ such that \begin{align}
        \langle W(f,\sigma)^2\rangle^\zeta_{\fB',\gamma}\le \frac{C_{14}L^2}{\iota ^2M^2}.\label{eq: second moment for W}
    \end{align}
\end{lemma}
\begin{proof}
    Recall from \eqref{eq: partition function for tilted spin measure} that $Z^\zeta_{\fB',\gamma}$ and $Z^{\zeta,tf}_{\fB',\gamma}$ (for $t\in [-1, 1]$) are the partition functions for the measures $\mu^\zeta_{\fB',\gamma}$ and $\mu^{\zeta,tf}_{\fB',\gamma}$ where $\mu^{\zeta,tf}_{\fB',\gamma}$ is defined as in \eqref{eq: tilted spin measure}.   By \eqref{eq: partition function equation for tilted model}, we get that $Z^{\zeta,tf}_{\fB',\gamma}$ is constant in $t\in \mathbb{R}$.   We compute the second derivative of $Z^{\zeta,tf}_{\fB',\gamma}$ with respect to $t$ and get that \begin{align}
         &\sum_{\sigma}\exp\Big(H^{\zeta,tf}_{\fB',\gamma}(\sigma)\Big)\cdot \Big[\big(\sum_{u\sim v}{c^{(\zz)}_{uv}}(f_u-f_v)\sin(\sigma_u-\sigma_v+tf_u-tf_v)\big)^2\nonumber\\&\hspace{4.5cm}-\sum_{u\sim v}{c^{(\zz)}_{uv}}(f_u-f_v)^2\cos(\sigma_u-\sigma_v+tf_u-tf_v)\Big]=0.\label{eq: second derivative for partition function with respect to t} 
    \end{align}Taking $t=0$ in \eqref{eq: second derivative for partition function with respect to t} and noting that $|f_u-f_v|\le \frac{C_{13}}{\iota M}$, we get that \begin{align}
         \langle W(f,\sigma)^2\rangle^\zeta_{\fB',\gamma}&=\Big\langle\sum_{u\sim v}{c^{(\zz)}_{uv}}(f_u-f_v)^2\cos(\sigma_u-\sigma_v)\Big\rangle^\zeta_{\fB',\gamma} \le \sum_{u\sim v}{c^{(\zz)}_{uv}}(\frac{C_{13}}{\iota M})^2\le \frac{\bC_1L^2}{\iota^2M^2}\label{eq: second moment for first order term} 
         \end{align}where $\bC_{1}>0$ is a constant depending on $C_{13},\beta_0$. Thus, we finish the proof of \eqref{eq: second moment for W}. 
\end{proof}

\begin{proof}[Proof of Lemma~\ref{lem: type 2 coupling}]

Let $f$ be the function with minimal Dirichlet energy that coincides with $F_0$ on $\fB^{(2)}$ and takes value $F_2(\fB)$ on $\partial\fB^{(20)}$. Recall the definition of $\mu^{\zeta,f}_{\fB^{(20)},\beta}$ from \eqref{eq: tilted spin measure}.
By \eqref{eq: type 1, k loop expansion} (with $\fB'=\fB^{(20)}$), we have for $k\in\{-1,0,1\}$,  \begin{equation}
    \mu^{\zeta+kF_2(\fB)}_{\fB^{(20)},\beta} \Big(\Xi(\fB,\zeta+kF_2(\fB))\Big)\leq\kappa+(1-\kappa)\big(1-q^k_\fB(\zeta)\big).
\end{equation}
For notation convenience, we write $(x)^-=\min\{0,x\}.$
Since $\log(\frac{x}{y})\ge \frac{x-y}{y}-\frac{(x-y)^2}{2\kappa^2}$ for any $x\in [\kappa,1]$ and $y\in[\kappa,1]$,  we get that (by applying the preceding inequality twice)
 \begin{align}
    &\log^-\Big(\frac{ \mu^{\zeta+F_2(\fB)}_{\fB^{(20)},\beta} \Big(\Xi(\fB,\zeta+F_2(\fB))\Big)\times  \mu^{\zeta-F_2(\fB)}_{\fB^{(20)},\beta} \Big(\Xi(\fB,\zeta-F_2(\fB))\Big)}{ \Big[\mu^{\zeta}_{\fB^{(20)},\beta} \big(\Xi(\fB,\zeta)\big)\Big]^2}\Big)\nonumber\\\ge ~&\frac{ \Big(\kappa+(1-\kappa)\big(1-q^1_\fB(\zeta)\big)\Big)+  \Big(\kappa+(1-\kappa)\big(1-q^{-1}_\fB(\zeta)\big)\Big)-2\Big(\kappa+(1-\kappa)\big(1-q^0_\fB(\zeta)\big)\Big)^-}{\kappa+(1-\kappa)\big(1-q^0_\fB(\zeta)\big)}\nonumber\\&-\frac{ (1-\kappa)^2\big(q^{-1}_\fB(\zeta)-q^{0}_\fB(\zeta)\big)^2+(1-\kappa)^2\big(q^{1}_\fB(\zeta)-q^{0}_\fB(\zeta)\big)^2}{2\kappa^2}\nonumber\\\ge ~&\frac{ \big(p^1_\fB(\zeta)+p^{-1}_\fB(\zeta)-2p^0_\fB(\zeta)\big)^-}{\kappa+(1-\kappa)\big(1-q^0_\fB(\zeta)\big) }-\frac{\big(p^1_\fB(\zeta)-p^0_\fB(\zeta)\big)^2+\big(p^{-1}_\fB(\zeta)-p^0_\fB(\zeta)\big)^2}{2\kappa^2}\label{eq: type 2 logarithm convert to type 1}
\end{align}
where (recall \eqref{eq: def of RB})
\begin{align}
    p^k_\fB(\zeta)&=R_\fB(\zeta)-q^k_\fB(\zeta)=\int_{\mathfrak{L}_\fB\in\Delta_\fB,\mathfrak{e}_\fB}\mathsf{M}(\mathsf{d}\mathfrak{L}_\fB)\cdot  \Psi(\mathsf{d}\mathfrak{e}_\fB)\cdot\frac{G_{\fB'}(\gamma'-\mathfrak{L}_\fB)}{G_{\fB'}(\gamma')}\nonumber\\&\hspace{5cm}\times\frac{Z^{\zeta+kf|_{\partial\fB^{(20)}}}_{\fB^{(20)},\gamma-\mathfrak{L}_\fB}\Big[\big(\mathcal{F}^k_*(\mathfrak{L}_\fB,\mathfrak{e}_\fB)\big)^c\Big]}{Z^{\zeta}_{\fB^{(20)},\gamma}}.\label{eq: def of pb}
\end{align} We remark that $q^k_\fB(\zeta)\in [0,1]$ but we only have $p_\fB^k(\zeta)\ge 0$ since $R_\fB$ may be larger than 1.

{\bf Bound on the first term in \eqref{eq: type 2 logarithm convert to type 1}.} Recall the definition of $\mathcal{F}^k_*(\mathfrak{L}_\fB,\mathfrak{e}_\fB)$ from Definition \ref{def: crossing event decomposition according to loops}. {For notation convenience, we write $\mathcal{F}_*(\mathfrak{L}_\fB,\mathfrak{e}_\fB)\big)=\mathcal{F}^0_*(\mathfrak{L}_\fB,\mathfrak{e}_\fB)\big)$.} Then we get that (recall that $f$ coincides with $F_0$ on $\mathsf{B}^{(2)}$) $$\mu^{\zeta+kf|_{\partial\fB^{(20)}}}_{\fB^{(20)},\gamma-\mathfrak{L}_\fB}\big(\mathcal{F}^k_*(\mathfrak{L}_\fB,\mathfrak{e}_\fB)\big)=\mu^{\zeta,kf}_{\fB^{(20)},\gamma-\mathfrak{L}_\fB}\big(\mathcal{F}_*(\mathfrak{L}_\fB,\mathfrak{e}_\fB)\big).$$  Thus we compute that (recall \eqref{eq: hamiltonian for tilted spin measure}) \begin{align}
&p^1_\fB(\zeta)+p^{-1}_\fB(\zeta)-2p^0_\fB(\zeta)=\int_{\mathfrak{L}_\fB\in\Delta_\fB,\mathfrak{e}_\fB}\mathsf{M}(\mathsf{d}\mathfrak{L}_\fB)\cdot  \Psi(\mathsf{d}\mathfrak{e}_\fB)\cdot\frac{G_{\fB'}(\gamma'-\mathfrak{L}_\fB)}{G_{\fB'}(\gamma')}\cdot\frac{Z^{\zeta}_{\fB^{(20)},\gamma-\mathfrak{L}_\fB}}{Z^{\zeta}_{\fB^{(20)},\gamma}}\nonumber\\&~\times\sum_{\sigma\notin \mathcal{F}_*(\mathfrak{L}_\fB,\mathfrak{e}_\fB)}\mu^{\zeta}_{\fB^{(20)},\gamma-\mathfrak{L}_\fB}(\sigma)\cdot \Big[\exp\big(H^{\zeta,f}_{\fB^{(20)},\gamma-\mathfrak{L}_\fB}(\sigma)-H^{\zeta}_{\fB^{(20)},\gamma-\mathfrak{L}_\fB}(\sigma)\big)\nonumber\\&\hspace{5cm}+\exp\big(H^{\zeta,-f}_{\fB^{(20)},\gamma-\mathfrak{L}_\fB}(\sigma)-H^{\zeta}_{\fB^{(20)},\gamma-\mathfrak{L}_\fB}(\sigma)\big)-2\Big]. \label{eq: first order term expansion}
\end{align} From now on, we fix $\mathfrak{L}_\fB$ and let $\gamma'=\gamma-\mathfrak{L}_\fB$, $c'_{u,v}=\sqrt{\gamma'_u\gamma_v'}J_{u,v}$.
We calculate that
\begin{align}
    H^{\zeta,kf}_{\fB^{(20)},\gamma'}(\sigma)-H^{\zeta}_{\fB^{(20)},\gamma'}(\sigma)={-}\sum_{u\sim v}c'_{u,v}\cdot 2\sin\big(\sigma_u-\sigma_v+\frac{k}{2}(f_u-f_v)\big)\sin\big(\frac{k}{2}(f_u-f_v)\big).\nonumber
\end{align}
Since $|f_u-f_v|\le \frac{\bC_1}{\iota M}$ and $L\le M^{0.1}$, recalling \eqref{eq: def of W}, we get the following inequalities:\begin{align}
    &\Big|H^{\zeta,kf}_{\fB^{(20)},\gamma'}(\sigma)-H^{\zeta}_{\fB^{(20)},\gamma'}(\sigma)\Big|\le \sum_{u\sim v}2c'_{u,v}\cdot \frac{\bC_1}{\iota M}\le \frac{\bC_2L^2}{M},\label{eq: hamiltonian difference absolute bound}\\
    &\Big|H^{\zeta,kf}_{\fB^{(20)},\gamma'}(\sigma)-H^{\zeta}_{\fB^{(20)},\gamma'}(\sigma){+}kW(f,\sigma)\Big|\nonumber\\
    \le~& \Bigg|\sum_{u\sim v}c'_{u,v}\cdot 2\Big[\sin\big(\sigma_u-\sigma_v+\frac{k}{2}(f_u-f_v)\big)-\sin\big(\sigma_u-\sigma_v\big)\Big]\cdot\sin\big(\frac{k}{2}(f_u-f_v)\big)\Big|\nonumber\\
    &+2\sin\big(\sigma_u-\sigma_v\big)\cdot\Big[\sin\big(\frac{k}{2}(f_u-f_v)\big)-\frac{k}{2}(f_u-f_v)\Big]\Bigg|\nonumber\\
    \le ~&\sum_{u\sim v} 2c'_{u,v}\bC_3(f_u-f_v)^2\le \frac{\bC_4L^2}{\iota^2 M^2}.\label{eq: hamiltonian first order conversion}
\end{align} By \eqref{eq: hamiltonian difference absolute bound} and the fact that $x\le\exp(x)-1$ for $|x|\le 1$, we get that
\begin{align}
    &\exp\big(H^{\zeta,kf}_{\fB^{(20)},\gamma'}(\sigma)-H^{\zeta}_{\fB^{(20)},\gamma'}(\sigma)\big)-1\ge H^{\zeta,kf}_{\fB^{(20)},\gamma'}(\sigma)-H^{\zeta}_{\fB^{(20)},\gamma'}(\sigma).\label{eq: hamiltonian difference absolute bound exponential 2}
\end{align}Combining \eqref{eq: hamiltonian difference absolute bound exponential 2} with \eqref{eq: hamiltonian first order conversion}, we get that \begin{align}
    &\exp\big(H^{\zeta,f}_{\fB^{(20)},\gamma'}(\sigma)-H^{\zeta}_{\fB^{(20)},\gamma'}(\sigma)\big)+\exp\big(H^{\zeta,-f}_{\fB^{(20)},\gamma'}(\sigma)-H^{\zeta}_{\fB^{(20)},\gamma'}(\sigma)\big)-2\nonumber\\&\stackrel{\eqref{eq: hamiltonian difference absolute bound exponential 2}}{\ge} H^{\zeta,f}_{\fB^{(20)},\gamma'}(\sigma)-H^{\zeta}_{\fB^{(20)},\gamma'}(\sigma)+H^{\zeta,-f}_{\fB^{(20)},\gamma'}(\sigma)-H^{\zeta}_{\fB^{(20)},\gamma'}(\sigma)\stackrel{\eqref{eq: hamiltonian first order conversion}}{\ge} -\frac{2\bC_4L^2}{\iota^2 M^2}.\label{eq: hamiltonian first order final}
\end{align}
 Plugging \eqref{eq: hamiltonian first order final} into \eqref{eq: first order term expansion}, we conclude that
\begin{align}
&p^1_\fB(\zeta)+p^{-1}_\fB(\zeta)-2p^0_\fB(\zeta)\nonumber\\\ge~&\int_{\mathfrak{L}_\fB\in\Delta_\fB,\mathfrak{e}_\fB}\mathsf{M}(\mathsf{d}\mathfrak{L}_\fB)\cdot  \Psi(\mathsf{d}\mathfrak{e}_\fB)\cdot\frac{G_{\fB'}(\gamma'-\mathfrak{L}_\fB)}{G_{\fB'}(\gamma')}\cdot\frac{Z^{\zeta}_{\fB^{(20)},\gamma-\mathfrak{L}_\fB}}{Z^{\zeta}_{\fB^{(20)},\gamma}}\times \mu^{\zeta}_{\fB^{(20)},\gamma'}(\mathcal{F}^c_*(\mathfrak{L}_\fB,\mathfrak{e}_\fB))\cdot (-\frac{2\bC_4L^2}{\iota^2 M^2})\nonumber\\\stackrel{\eqref{eq: def of pb}}{=}&-\frac{2\bC_4L^2}{\iota^2 M^2} \big(R_\fB(\zeta)-q^0_\fB(\zeta)\big)=-\frac{2\bC_4L^2}{\iota^2 M^2} p^0_\fB(\zeta). \label{eq: first order term final p}
\end{align} Integrating over $\zeta\in\Sigma$ and combining with Lemma \ref{lem: type probability}, we conclude that 
\begin{align}
   &\int_{\zeta\in \Sigma}\mu^{\xi}_{\Lambda,\gamma}(\zeta)\cdot  \frac{ \big(p^1_\fB(\zeta)+p^{-1}_\fB(\zeta)-2p^0_\fB(\zeta)\big)^-}{\kappa+(1-\kappa)\big(1-q^0_\fB(\zeta)\big) }\ge \int_{\zeta\in \Sigma}\mu^{\xi}_{\Lambda,\gamma}(\zeta)\cdot(-\frac{2\bC_4L^2}{\iota^2 M^2}) \frac{ p^0_\fB(\zeta)}{\kappa+(1-\kappa)\big(1-q^0_\fB(\zeta)\big) }\nonumber\\\ge~&\int_{\zeta\in \Sigma}\mu^{\xi}_{\Lambda,\gamma}(\zeta)\cdot(-\frac{2\bC_4L^2}{\iota^2 M^2}) \frac{ p^0_\fB(\zeta)}{\kappa}\nonumber\\=~&-\frac{2\bC_4L^2}{\kappa\iota^2 M^2}\cdot\mu^{\xi}_{\Lambda,\gamma}(\fB \text{ is type } ({\rm II},k)\mid \mathrm{AUX}=1) \ge -\frac{2\bC_5L^2\rho}{\kappa\iota^2 M^2}.\label{eq: first order term final xi}
\end{align}

{\bf Bound on the second term in \eqref{eq: type 2 logarithm convert to type 1}.}
 Similar to \eqref{eq: first order term expansion}, we calculate that \begin{align}
&p^1_\fB(\zeta)-p^0_\fB(\zeta)=\int_{\mathfrak{L}_\fB\in\Delta_\fB,\mathfrak{e}_\fB}\mathsf{M}(\mathsf{d}\mathfrak{L}_\fB)\cdot  \Psi(\mathsf{d}\mathfrak{e}_\fB)\cdot\frac{G_{\fB'}(\gamma'-\mathfrak{L}_\fB)}{G_{\fB'}(\gamma')}\cdot\frac{Z^{\zeta}_{\fB^{(20)},\gamma'}}{Z^{\zeta}_{\fB^{(20)},\gamma}}\nonumber\\&~\times\sum_{\sigma\notin \mathcal{F}_*(\mathfrak{L}_\fB,\mathfrak{e}_\fB)}\mu^{\zeta}_{\fB^{(20)},\gamma'}(\sigma)\cdot \Big[\exp\big(H^{\zeta,f}_{\fB^{(20)},\gamma'}(\sigma)-H^{\zeta}_{\fB^{(20)},\gamma'}(\sigma)\big)-1\Big]. \label{eq: second order term expansion}
\end{align}
Note that $|e^x-1|\le 2|x|$ for $|x|\le 1$, and thus we get from \eqref{eq: hamiltonian first order conversion} that
\begin{align}
    &\sum_{\sigma\notin \mathcal{F}_*(\mathfrak{L}_\fB,\mathfrak{e}_\fB)}\mu^{\zeta}_{\fB^{(20)},\gamma'}(\sigma)\cdot \Big|\exp\big(H^{\zeta,f}_{\fB^{(20)},\gamma'}(\sigma)-H^{\zeta}_{\fB^{(20)},\gamma'}(\sigma)\big)-1\Big|\nonumber\\&\stackrel{\eqref{eq: hamiltonian first order conversion}}{\le} {2}\sum_{\sigma\notin \mathcal{F}_*(\mathfrak{L}_\fB,\mathfrak{e}_\fB)}\mu^{\zeta}_{\fB^{(20)},\gamma'}(\sigma)\cdot \Big(|W(f,\sigma)|+\frac{\bC_4L^2}{\iota^2M^2}\Big).\label{eq: 5.66}
\end{align}
Plugging \eqref{eq: 5.66} into \eqref{eq: second order term expansion} and combining with the Cauchy-Schwarz inequality, we get that \begin{align}
    &\Big(p^1_\fB(\zeta)-p^0_\fB(\zeta)\Big)^2\le {8}\Bigg[\int_{\mathfrak{L}_\fB\in\Delta_\fB,\mathfrak{e}_\fB}\mathsf{M}(\mathsf{d}\mathfrak{L}_\fB)\cdot  \Psi(\mathsf{d}\mathfrak{e}_\fB)\cdot\frac{G_{\fB'}(\gamma'-\mathfrak{L}_\fB)}{G_{\fB'}(\gamma')}\cdot\frac{Z^{\zeta}_{\fB^{(20)},\gamma'}}{Z^{\zeta}_{\fB^{(20)},\gamma}}\nonumber\\&\hspace{6cm}\times\sum_{\sigma\notin \mathcal{F}_*(\mathfrak{L}_\fB,\mathfrak{e}_\fB)}\mu^{\zeta}_{\fB^{(20)},\gamma'}(\sigma)\cdot|W(f,\sigma)|\Bigg]^2\nonumber\\&\hspace{1cm}+{8}\Bigg[\int_{\mathfrak{L}_\fB\in\Delta_\fB,\mathfrak{e}_\fB}\mathsf{M}(\mathsf{d}\mathfrak{L}_\fB)\cdot  \Psi(\mathsf{d}\mathfrak{e}_\fB)\cdot\frac{G_{\fB'}(\gamma'-\mathfrak{L}_\fB)}{G_{\fB'}(\gamma')}\cdot\frac{Z^{\zeta}_{\fB^{(20)},\gamma'}\big(\mathcal{F}^c_*(\mathfrak{L}_\fB,\mathfrak{e}_\fB)\big)}{Z^{\zeta}_{\fB^{(20)},\gamma}}\cdot \big(\frac{\bC_4L^2}{\iota^2M^2}\big)\Bigg]^2.\label{eq: 5.67}
\end{align}
Applying the Cauchy-Schwarz inequality and combining with \eqref{eq: second moment for W}, we get that the first term on the right hand side in \eqref{eq: 5.67} is upper-bounded by\begin{align}
&{8}\int_{\mathfrak{L}_\fB\in\Delta_\fB,\mathfrak{e}_\fB}\mathsf{M}(\mathsf{d}\mathfrak{L}_\fB)\cdot  \Psi(\mathsf{d}\mathfrak{e}_\fB)\cdot\frac{G_{\fB'}(\gamma'-\mathfrak{L}_\fB)}{G_{\fB'}(\gamma')}\cdot\frac{Z^{\zeta}_{\fB^{(20)},\gamma'}\big(\mathcal{F}^c_*(\mathfrak{L}_\fB,\mathfrak{e}_\fB)\big)}{Z^{\zeta}_{\fB^{(20)},\gamma}}\nonumber\\&\times\int_{\mathfrak{L}_\fB\in\Delta_\fB,\mathfrak{e}_\fB}\mathsf{M}(\mathsf{d}\mathfrak{L}_\fB)\cdot  \Psi(\mathsf{d}\mathfrak{e}_\fB)\cdot\frac{G_{\fB'}(\gamma'-\mathfrak{L}_\fB)}{G_{\fB'}(\gamma')}\cdot\frac{Z^{\zeta}_{\fB^{(20)},\gamma'}}{Z^{\zeta}_{\fB^{(20)},\gamma}}\cdot \langle\big(W(f,\sigma)\big)^2\rangle^{\zeta}_{\fB^{(20)},\gamma'}\nonumber\\&\stackrel{\eqref{eq: second moment for W}}{\le} 2p^0_\fB(\zeta)\cdot R_\fB(\zeta)\cdot \frac{\bC_6L^2}{\iota^2M^2}\label{eq: second order term 1}
\end{align} Noting that $p^0_\fB(\zeta)\le R_\fB(\zeta)$ and plugging \eqref{eq: second order term 1} to \eqref{eq: 5.67}, we get that \begin{align}
    &\Big(p^1_\fB(\zeta)-p^0_\fB(\zeta)\Big)^2\le  p^0_\fB(\zeta)\cdot R_\fB(\zeta)\cdot \frac{\bC_7L^2}{\iota^2M^2}.\label{eq: second order term 000}
\end{align} Substituting the summation $\sum_{\sigma\notin \mathcal{F}_*(\mathfrak{L}_\fB,\mathfrak{e}_\fB)}$ by $\sum_{\sigma\in \mathcal{F}_*(\mathfrak{L}_\fB,\mathfrak{e}_\fB)}$ in \eqref{eq: second order term expansion}, we get the following equality similar to  \eqref{eq: second order term 000}: \begin{align}
    &\Big(q^1_\fB(\zeta)-q^0_\fB(\zeta)\Big)^2\le                      q^0_\fB(\zeta)\cdot R_\fB(\zeta)\cdot \frac{\bC_8L^2}{\iota^2M^2}.\label{eq: second order term 2}
\end{align}Since $q^1_\fB(\zeta)-q^0_\fB(\zeta)=-\big(p^1_\fB(\zeta)-p^0_\fB(\zeta)\big)$, we conclude from \eqref{eq: second order term 000} and \eqref{eq: second order term 2} that \begin{equation}
\Big(p^1_\fB(\zeta)-p^0_\fB(\zeta)\Big)^2\le                      \min\{p^0_\fB(\zeta),q^0_\fB(\zeta)\}\cdot R_\fB(\zeta)\cdot \frac{\bC_8L^2}{\iota^2M^2}.\nonumber
\end{equation}
Let $\Sigma^+=\{\zeta\in\Sigma:p^0_\fB(\zeta)\ge q^0_\fB(\zeta)\}$ and $\Sigma^-=\{\zeta\in\Sigma:p^0_\fB(\zeta)< q^0_\fB(\zeta)\}$.
Integrating over $\zeta$ we get that \begin{align}
&\int_{\zeta\in \Sigma}\mu^{\xi}_{\Lambda,\gamma}(\zeta)\cdot\Big(p^1_\fB(\zeta)-p^0_\fB(\zeta)\Big)^2\nonumber\\ \le~&  \frac{\bC_8L^2}{\iota^2M^2}\times\Big[\int_{\zeta\in \Sigma^+}\mu^{\xi}_{\Lambda,\gamma}(\zeta)q^0_\fB(\zeta)\cdot R_\fB(\zeta) +\int_{\zeta\in \Sigma^-}\mu^{\xi}_{\Lambda,\gamma}(\zeta)p^0_\fB(\zeta)\cdot R_\fB(\zeta)\Big].\label{eq: second order term 3}
\end{align} For the first term in \eqref{eq: second order term 3}, we apply  Lemma \ref{lem: type probability} and the fact that $q^0_\fB(\zeta)\le 1, R_\fB(\zeta)\le 2p^0_\fB(\zeta)$ to get that \begin{align}
    \int_{\zeta\in \Sigma^+}\mu^{\xi}_{\Lambda,\gamma}(\zeta)q^0_\fB(\zeta)\cdot R_\fB(\zeta)\le 2\int_{\zeta\in \Sigma^+}\mu^{\xi}_{\Lambda,\gamma}(\zeta)p^0_\fB(\zeta)\le 12\rho.\label{eq: second order term 4}
\end{align} For the second term in \eqref{eq: second order term 3}, we apply  Lemma \ref{lem: type probability} and the fact that $R_\fB(\zeta)\le 2q^0_\fB(\zeta)\le 2$ to get that \begin{align}
    \int_{\zeta\in \Sigma^-}\mu^{\xi}_{\Lambda,\gamma}(\zeta)p^0_\fB(\zeta)\cdot R_\fB(\zeta)\le 2\int_{\zeta\in \Sigma^-}\mu^{\xi}_{\Lambda,\gamma}(\zeta)p^0_\fB(\zeta)\le 12\rho.\label{eq: second order term 5}
\end{align}Plugging \eqref{eq: second order term 4} and \eqref{eq: second order term 5} into \eqref{eq: second order term 3}, we get that \begin{align}
&\int_{\zeta\in \Sigma}\mu^{\xi}_{\Lambda,\gamma}(\zeta)\cdot\Big(p^1_\fB(\zeta)-p^0_\fB(\zeta)\Big)^2 \le  \frac{\bC_8L^2}{\iota^2M^2}\times 24\rho.\label{eq: second order term final 1}
\end{align}By symmetry, we get the same bound replacing the superscript $1$ by $-1$ in the above. Combining \eqref{eq: second order term final 1} (for both versions) with \eqref{eq: first order term final xi} and plugging into \eqref{eq: type 2 logarithm convert to type 1}, we get that the integration in the lemma statement is at least $-\frac{2\bC_6L^2\rho}{\kappa M^2}-\frac{\bC_9\rho L^2}{\kappa^2\iota^2M^2}$, completing the proof of Lemma \ref{lem: type 2 coupling}.
\end{proof}
Now we are ready to control the energy cost of $F_3^+.$ For notation convenience, we let $\mathfrak{T}$ be a realization of the three sequences $(\{\mathfrak{L}_\fB\},\{T_\fB\}, \{\mathfrak{e}_{\fB}\})$. Recall $\Upsilon^0(\mathfrak{T})$ from Definition \ref{def: final spin wave frozen part}. In addition, for $i\in\{-1,1\}$, let $\Upsilon^{i}(\mathfrak{T})$ denote the realization obtained from $\Upsilon^0$ by replacing $T_\fB$ with $T_\fB^i$, where $T^i_\fB$ is obtained from $T_\fB$ by replacing the second entry with $i$.

Recall from Definition \ref{def: best x bc for crossing} that $ \xi^* $ denotes the all-$\mathsf{x}$-copy on $\partial\Lambda$ and recall Remark~\ref{rmk: xi^* largest boundary condition}.
Thus, it suffices to consider the measure $\mu^{\xi^*}_{\Lambda,\beta-L^{\zz}}.$
Let $ \mathsf{SAME} $ denote the event that the outer boundary takes the same sign as the inner boundary in the $\mathsf{x}$-copy. Similarly, let $\mathsf{OPP}$  denote the event that the outer boundary takes the opposite sign from the inner boundary in the $\mathsf{x}$-copy. Finally, we remark that, while each $\Upsilon^i$ is a realization of explored objects, each of them corresponds to a collection of configurations in the full space.
\begin{proposition}
\label{prop: symmetric spin wave}
Recall $\mathfrak{T}=\Big(\{\mathfrak{L}_\fB\},\{T_\fB\}, \{\mathfrak{e}_{\fB}\}\Big)$.
Suppose that $\{\fB\}$ and $\{T_\fB\}$ are excellent.
Let $\Gamma$ denote the collection of spin configurations $\sigma$ such that \begin{equation}
    \sum_{\fB\in\mathfrak{B}_{\rm II}}\log\Big(\mu^{\sigma|_{ \fB^{( 20)}}}_{\fB^{(20)},\beta'}(F_3^{+}(\sigma))\Big)+\log\Big(\mu^{\sigma|_{ \fB^{( 20)}}}_{\fB^{(20)},\beta'}(F_3^{-}(\sigma))\Big)-2\log\Big(\mu^{\sigma|_{ \fB^{( 20)}}}_{\fB^{(20)},\beta'}(\sigma)\Big)\ge -1.\label{eq: def of good config for type 2 map}
\end{equation}
    Then, letting $\xi^*$ be the all-$\mathsf x$-copy boundary condition on $\partial \Lambda$ (as in Definition \ref{def: best x bc for crossing}), we have \begin{equation}
        \mu^{\xi^*}_{\Lambda,\beta-L^{\zz}}(\Upsilon^0(\mathfrak{T})\cap\mathsf{SAME}\cap\Gamma)
        \le  10\mu^{\xi^*}_{\Lambda,\beta-L^{\zz}}(\Upsilon^{1}(\mathfrak{T})\cap\mathsf{OPP})+10  \mu^{\xi^*}_{\Lambda,\beta-L^{\zz}}(\Upsilon^{-1}(\mathfrak{T})\cap\mathsf{OPP}).\label{eq: spin wave comparison}
    \end{equation}
\end{proposition}
\begin{proof}
    For any configuration $\theta \in\Upsilon^0\cap\mathsf{SAME}$, it can be decomposed into the explored configuration $(\mathfrak{B}_{\rm I}, \{\mathfrak{L}_\mathsf B, \mathfrak{e}_\mathsf B\})_{B\in \mathfrak{B}_{\rm I}}$, {$(\mathfrak{B}_{\rm II}, \{\mathfrak{L}_\mathsf B, \mathfrak{e}_\mathsf B,\sigma_\fB\})_{B\in \mathfrak{B}_{\rm II}}$ (recall the convention in Definition \ref{def: type 2 coupling})} before applying the spin wave and the spin configuration $\sigma$ on the remaining graph ${\Lambda''}$. 
Then the probability for this configuration $\theta$ can be written as (note that below we treat $\mathfrak{L}_\mathsf B, \mathfrak{e}_\mathsf B, \mathfrak{B}_{\rm I}, \sigma$ as an realization)
    \[{
    \mu^{\xi^*}_{\Lambda,\beta-L^{\zz}}(\{\mathfrak{L}_\fB,\mathfrak{e}_\fB\}_{\fB\in\mathfrak{B}_{\rm I}})\cdot \mu^{\xi^*}_{\Lambda'',\beta'}(\sigma)\cdot \prod_{\fB\in\mathfrak{B}_{\rm II}}\mu^{\sigma|_{\partial\fB^{(20)}}}_{\fB^{(20)},\beta-L^\zz}(\mathfrak{L}_\mathsf B, \mathfrak{e}_\mathsf B,\sigma_\fB),}
    \]
    where we denote by $\beta' = \beta -  L^\mathsf z -  \sum_{\mathsf B}\mathfrak{ L}_{\mathsf B}$ the remaining local time condition after the exploration and we use the convention that $\mu^{\sigma|_{\partial \mathsf B^{(20)}}}_{\fB^{(20)},\beta-L^\zz}(\mathfrak{L}_\mathsf B, \mathfrak{e}_\mathsf B,\sigma_\fB) = \mu^{\sigma|_{\partial \mathsf B^{(20)}}}_{\fB^{(20)},\beta-L^\zz}( \mathfrak{e}_\mathsf B,\sigma_\fB)$ if $\mathcal L_\mathsf B = \emptyset$.

    Note that $\{T_\mathsf B\}$ is determined by $(\mathfrak{B}_{\rm I}, \{\mathfrak L_\mathsf B, \mathfrak e_\mathsf B\})_{\fB\in \mathfrak{B}_{\rm I}}$, {$(\mathfrak{B}_{\rm II}, \{\mathfrak{L}_\mathsf B, \mathfrak{e}_\mathsf B,\sigma_\fB\})_{B\in \mathfrak{B}_{\rm II}}$} and $\sigma$; this is why $\{T_\mathsf B\}$ did not appear explicitly above except through $\mathfrak{B}_{\rm I}$ {and $\mathfrak{B}_{\rm II}$}. In addition, we define $\Theta_{\{\mathfrak{L}_\fB,\mathfrak{e}_\fB\}}(\{T_\fB\})$ to be the collection of spin configurations $\sigma$ {and $(\mathfrak{B}_{\rm II}, \{\mathfrak{L}_\mathsf B, \mathfrak{e}_\mathsf B,\sigma_\fB\})_{B\in \mathfrak{B}_{\rm II}}$ which we referred to as $\theta_{\rm II}$}  
    such that $(\{\mathfrak{L}_\fB,\mathfrak{e}_\fB\}, \sigma)$ satisfies the type condition $\{T_\fB\}$ {(note that $\sigma$ only influences the boxes in $\mathfrak{B}_{\rm I}$)}. Then we calculate that \begin{align}
        \mu^{\xi^*}_{\Lambda,\beta-L^{\zz}}(\Upsilon^0(\mathfrak{T})\cap\text{SAME}\cap \Gamma)&=\mu^{\xi^*}_{\Lambda,\beta-L^{\zz}}(\{\mathfrak{L}_\fB,\mathfrak{e}_\fB\}_{\fB\in\mathfrak{B}_{\rm I}})\cdot\sum_{(\sigma{,\theta_{\rm II}})\in \Theta_{\{\mathfrak{L}_\fB,\mathfrak{e}_\fB\}}(\{T_\fB\})\cap \text{SAME}\cap \Gamma} \mu^{\xi^*}_{\Lambda',\beta'}(\sigma{,\theta_{\rm II}})\nonumber
    \end{align}
    where $$
    \mu^{\xi^*}_{\Lambda',\beta'}(\sigma{,\theta_{\rm II}})=\mu^{\xi^*}_{\Lambda'',\beta'}(\sigma)\times \prod_{\fB\in\mathfrak{B}_{\rm II}}\mu^{\sigma|_{\partial\fB^{(20)}}}_{\fB^{(20)},\beta-L^\zz}(\mathfrak{L}_\mathsf B, \mathfrak{e}_\mathsf B,\sigma_\fB).
    $$
    Similarly, we also have for any $i\in \{1,-1\}$ \begin{align}
        \mu^{\xi^*}_{\Lambda,\beta-L^{\zz}}(\Upsilon^i(\mathfrak{T})\cap\mathsf{OPP})&=\mu^{\xi^*}_{\Lambda,\beta-L^{\zz}}(\{\mathfrak{L}_\fB,\mathfrak{e}_\fB\}_{\fB\in\mathfrak{B}_{\rm I}})\cdot\sum_{(\sigma{,\theta_{\rm II}})\in \Theta_{\{\mathfrak{L}_\fB,\mathfrak{e}_\fB\}}(\{T_\fB^i\})\cap \mathsf{OPP}} \mu^{\xi^*}_{{\Lambda''},\beta'}\mu^{\xi^*}_{\Lambda',\beta'}(\sigma{,\theta_{\rm II}}).\nonumber
    \end{align}
    It suffices to consider the probability of $\Theta_{\{\mathfrak{L}_\fB,\mathfrak{e}_\fB\}}(\{T_\fB\})\cap \text{SAME}{\cap\Gamma}$, $\Theta_{\{\mathfrak{L}_\fB,\mathfrak{e}_\fB\}}(\{T_\fB^1\})\cap \mathsf{OPP}$ and $\Theta_{\{\mathfrak{L}_\fB,\mathfrak{e}_\fB\}}(\{T_\fB^{-1}\})\cap \mathsf{OPP}$ under the measure $\mu^{\xi^*}_{{\Lambda''},\beta'}$ {and the corresponding Type $\rm II$ probability}.
    For any configuration $(\sigma{,\theta_{\rm II}})\in \Theta_{\{\mathfrak{L}_\fB,\mathfrak{e}_\fB\}}(\{T_\fB\})\cap \text{SAME}$, we consider the two configurations $F_3^+(\sigma{,\theta_{\rm II}})$ and $F_3^-(\sigma{,\theta_{\rm II}})$. By Definition \ref{def: final spin wave global part}, we get that $F_3^+(\sigma{,\theta_{\rm II}})\in \Theta_{\{\mathfrak{L}_\fB,\mathfrak{e}_\fB\}}(\{T_\fB^1\})\cap \mathsf{OPP}$ and $F_3^-(\sigma{,\theta_{\rm II}})\in \Theta_{\{\mathfrak{L}_\fB,\mathfrak{e}_\fB\}}(\{T_\fB^{-1}\})\cap \mathsf{OPP}$. Therefore, it suffices to compare the probability $\mu^{\xi^*}_{{\Lambda'},\beta'}(\sigma{,\theta_{\rm II}})$ with $\mu^{\xi^*}_{{\Lambda'},\beta'}\big(F_3^+(\sigma{,\theta_{\rm II}})\big)+ \mu^{\xi^*}_{{\Lambda'},\beta'}\big(F_3^-(\sigma{,\theta_{\rm II}})\big).$

For any $k$, write $\Lambda^{(\setminus k)} = \Lambda'\setminus(\cup_{\fB\in\mathfrak{B}_{\rm II}}\fB^{(k)})$ {and note that ${\Lambda''}=\Lambda^{(\setminus 20)}$}.
We calculate that \begin{equation}\label{eq: expansion for pre-image}
    \mu^{\xi^*}_{\Lambda',\beta'}(\sigma{,\theta_{\rm II}})=\frac{1}{\mathcal{Z}^{\xi^*}_{\Lambda',\beta'}}\cdot \exp\big(H^{\xi^*}_{\Lambda^{(\setminus 20)},\beta'}(\sigma)\big)\cdot\prod_{\fB\in\mathfrak{B}_{\rm II}}\Big(\mathcal{Z}^{\sigma|_{ \partial\fB^{( 20)}}}_{\fB^{(20)},\beta'}\cdot \mu^{\sigma|_{\partial \fB^{( 20)}}}_{\fB^{(20)},\beta'}(\sigma|_{ \fB^{( 20)}})\Big).
\end{equation}
Similarly, we also have \begin{align}
    \mu^{\xi^*}_{\Lambda',\beta'}(F_3^{\pm}(\sigma{,\theta_{\rm II}}))&=\frac{1}{\mathcal{Z}^{\xi^*}_{\Lambda',\beta'}}\cdot \exp\big(H^{\xi^*}_{\Lambda^{(\setminus 20)},\beta'}(F_3^\pm{}(\sigma))\big)\nonumber\\&~\cdot\prod_{\fB\in\mathfrak{B}_{\rm II}}\Big(\mathcal{Z}^{\sigma|_{ \partial\fB^{( 20)}}\pm F_2(\fB)}_{\fB^{(20)},\beta'}\cdot \mu^{\sigma|_{\partial \fB^{( 20)}}}_{\fB^{(20)},\beta'}(F_3^{\pm}(\sigma)|_{ \fB^{( 20)}})\Big).\label{eq: expansion for image}
\end{align}

Recalling  that $F_2(B)$ is the constant value $F_2$ takes on $\partial\fB^{(20)}$, we get that $\mathcal{Z}^{\sigma|_{ \fB^{( 20)}}\pm F_2(\fB)}_{\fB^{(20)},\beta'}=\mathcal{Z}^{\sigma|_{ \fB^{( 20)}}}_{\fB^{(20)},\beta'}$.
    Thus, combining \eqref{eq: expansion for pre-image} and \eqref{eq: expansion for image} with \eqref{eq: def of good config for type 2 map}, we get that for $\sigma\in \Gamma$,
    \begin{align}
        &\frac{\mu^{\xi^*}_{\Lambda',\beta'}\big(F_3^{+}(\sigma{,\theta_{\rm II}})\big)\cdot\mu^{\xi^*}_{\Lambda',\beta'}\big(F_3^{-}(\sigma{,\theta_{\rm II}})\big)}{\Big(\mu^{\xi^*}_{\Lambda',\beta'}(\sigma{,\theta_{\rm II}})\Big)^2}\nonumber\\\ge~& e^{-1}\cdot \exp\Big(H^{\xi^*}_{\Lambda^{(\setminus 20)},\beta'}\big(F_3^+(\sigma)\big)+ H^{\xi^*}_{\Lambda^{(\setminus 20)},\beta'}\big(F_3^-(\sigma)\big)-2H^{\xi^*}_{\Lambda^{(\setminus 20)},\beta'}(\sigma)\Big). \label{eq: energy cost for pm spin wave}
    \end{align}
    In the above and thereon, we slightly abuse the notation $F_3^\pm(\sigma)$ to denote the spin configuration on $\Lambda^{(\setminus 20)}$ in $F_3^\pm(\sigma{,\theta_{\rm II}})$ since the spin configuration in $ F^{\pm}_3(\sigma, \theta_{\rm II}) $ does not depend on $ \theta_{\rm II} $.
    We now bound \eqref{eq: energy cost for pm spin wave}.
    Similarly to \eqref{eq: relative entropy symmetrize}, we will consider the contributions of $F_3^+(\sigma)$ and $F_3^-(\sigma)$ together. Thus we calculate that \begin{align}
        &H^{\xi^*}_{\Lambda^{(\setminus 20)},\beta'}\big(F_3^+(\sigma)\big)+ H^{\xi^*}_{\Lambda^{(\setminus 20)},\beta'}\big(F_3^-(\sigma)\big)-2H^{\xi^*}_{\Lambda^{(\setminus 20)},\beta'}(\sigma)\nonumber\\=~&\sum_{u\sim v\in \Lambda^{(\setminus 20)}}\sqrt{\beta'_u\beta'_v}\Big[\cos\Big((\sigma_u-\sigma_v)+\big(F_2(u)-F_2(v)\big)\Big)+\cos\Big((\sigma_u-\sigma_v)-\big(F_2(u)-F_2(v)\big)\Big)\nonumber\\
        &\hspace{2.8cm}-2\cos(\sigma_u-\sigma_v)\Big]\nonumber\\
        =~&\sum_{u\sim v\in \Lambda^{(\setminus 20)}}2\sqrt{\beta'_u\beta'_v}\cos(\sigma_u-\sigma_v)\cdot \Big(\cos\big(F_2(u)-F_2(v)\big)-1\Big)\nonumber\\
        \ge ~&-\sum_{u\sim v\in \Lambda^{(\setminus 20)}}2\sqrt{\beta'_u\beta'_v}\cdot \bC_1\big(F_2(u)-F_2(v)\big)^2\ge -\bC_2\sum_{u\sim v\in \Lambda^{(\setminus 20)}}\big(F_2(u)-F_2(v)\big)^2.\label{eq: cosine approximation + symmetrize}
        \end{align}
    Recalling Definition \ref{def: final spin wave global part}, $F_2$ coincides with $F_0$ in the region $\Lambda^{(\setminus 40)}$. Thus by \eqref{eq: def of base spin wave}, we can control the Dirichlet energy on $\Lambda^{(\setminus 40)}$ as follows: \begin{align}
       \sum_{u\sim v\in \Lambda^{(\setminus 40)}}\big(F_2(u)-F_2(v)\big)^2\le \sum_{u\sim v\in \Lambda^{(\setminus 40)}}(\frac{\pi}{\iota M})^2\le \frac{\bC_3|\Lambda^{(\setminus 40)}|}{\iota^2 M^2}.\label{eq:  hamiltonian in spin wave region 1}
    \end{align}For each $\mathsf B \in \mathfrak B_{\rm II}$, note that $F_2$  interpolates linearly in the annulus $\fB^{(40)}\setminus\fB^{(20)}$.
    Then we calculate that \begin{equation}
        \sum_{u\sim v\in \fB^{(40)}\setminus\fB^{(20)}}\big(F_2(u)-F_2(v)\big)^2\le \sum_{u\sim v\in \fB^{(40)}\setminus\fB^{(20)}}(\frac{\bC_4}{\delta \iota M})^2\le \bC_5\frac{|\fB|}{\delta  \iota^2 M^2}.\label{eq:  hamiltonian in spin wave region 2}
    \end{equation}Plugging \eqref{eq:  hamiltonian in spin wave region 1} and \eqref{eq:  hamiltonian in spin wave region 2} into \eqref{eq: cosine approximation + symmetrize}, we get that  for $\sigma\in \Gamma$, \begin{equation}
        H^{\xi^*}_{\Lambda^{(\setminus 20)},\beta'}\big(F_3^+(\sigma)\big)+ H^{\xi^*}_{\Lambda^{(\setminus 20)},\beta'}\big(F_3^-(\sigma)\big)-2H^{\xi^*}_{\Lambda^{(\setminus 20)},\beta'}(\sigma)\ge -\frac{\bC_2\bC_3|\Lambda^{(\setminus 40)}|}{\iota^2 M^2}-\sum_{\fB\in\mathfrak{B}_{\rm II}}\frac{\bC_2\bC_5|\fB|}{\delta \iota^2 M^2}.\label{eq: hamiltonian difference 1}
    \end{equation} Since $\Lambda^{(\setminus 40)}\subset\Lambda\setminus(\cup_{\fB\in \mathfrak{S}^{\ell_{*}}}\fB)$, we get from Definition \ref{def: good sequences and types} that $|\Lambda^{(\setminus 40)}|\le C_6\delta M^2.$ In addition, by Definition \ref{def: good sequences and types}, we get that $\sum_{\fB\in\mathfrak{B}_{\rm II}}\frac{C_2C_5|\fB|}{\delta \iota^2 M^2}\le \frac{\bC_6(\kappa+\rho)^{\frac23}}{\iota^2\delta}.$ Plugging this into \eqref{eq: hamiltonian difference 1} and combining with \eqref{eq: energy cost for pm spin wave}, we get that \begin{equation}
        \frac{\mu^{\xi^*}_{\Lambda',\beta'}\big(F_3^+(\sigma{,\theta_{\rm II}})\big)}{\mu^{\xi^*}_{\Lambda',\beta'}(\sigma{,\theta_{\rm II}})}\cdot \frac{\mu^{\xi^*}_{\Lambda',\beta'}\big(F_3^-(\sigma{,\theta_{\rm II}})\big)}{\mu^{\xi^*}_{\Lambda',\beta'}(\sigma{,\theta_{\rm II}})}
        \ge \exp(-1- \frac{\bC_6(\kappa+\rho)^{\frac23}}{\iota^2\delta}-\frac{\bC_7\delta}{\iota^2}).\label{eq: spin wave comparison frozen 1}
    \end{equation} Choosing $\delta=(\kappa+\rho)^{\frac{1}{3}}$ and choosing $\rho,\kappa$ small enough in \eqref{eq: spin wave comparison frozen 1} such that $\frac{\bC_6(\kappa+\rho)^{\frac23}}{\iota^2\delta}+\frac{\bC_7\delta}{\iota^2}\le 2$, we get from Lemma \ref{lem: push forward probability control} that  for $\sigma\in \Gamma$, \begin{align}
        \mu^{\xi^*}_{\Lambda',\beta'}\big(F_3^+(\sigma{,\theta_{\rm II}})\big)+\mu^{\xi^*}_{\Lambda',\beta'}\big(F_3^-(\sigma{,\theta_{\rm II}})\big)&\ge 2\sqrt{\mu^{\xi^*}_{\Lambda',\beta'}\big(F_3^+(\sigma{,\theta_{\rm II}})\big)\cdot \mu^{\xi^*}_{\Lambda',\beta'}\big(F_3^-(\sigma{,\theta_{\rm II}})\big)}\nonumber\\&\ge 2e^{-3}\mu^{\xi^*}_{\Lambda',\beta'}(\sigma{,\theta_{\rm II}}).\label{eq: spin wave comparison frozen 2}
    \end{align}The desired result follows from integrating \eqref{eq: spin wave comparison frozen 2} over $(\sigma{,\theta_{\rm II}})\in \Theta_{\{\mathfrak{L}_\fB,\mathfrak{e}_\fB\}}(\{T_\fB\})\cap \text{SAME}\cap \Gamma$ and observing $2e^{-3}\ge 0.1.$
\end{proof}

Now we are ready to prove Theorem~\ref{thm: thin annulus local time crossing}.
\begin{proof}[Proof of Theorem~\ref{thm: thin annulus local time crossing}]
We first establish a uniform bound for the crossing probability in the $\mathsf{x}$-copy.

Fix $\iota>0$ and write $\Lambda=\Lambda_{M,(1+\iota)M}$.
Let $\mathcal C^{\xx}$ denote the event that an $\xx$-cluster crosses $\Lambda$.
Throughout this proof, crossing paths use cable edges of the annulus.
Let $\lambda_{\Lambda,\beta}^{\xi}$ be the distribution of $ L^{\mathsf z}$ under $\nu_{\Lambda,\beta}^{\xi}$.

For a realization $\mathfrak{z}$ of $ L^\mathsf z$, let $\mu_{\mathfrak{z}}^{\xi^*}$ be the XY measure with local time condition $\beta-\mathfrak{z}$, and with all-$\mathsf x$-copy boundary condition on $\partial \Lambda$ as in Definition \ref{def: best x bc for crossing}.
Write $\mathsf{SAME}$ and $\mathsf{OPP}$ for the events that the $\mathsf{x}$-signs on the outer and inner boundaries agree and differ, respectively.

By Corollary~\ref{cor:fkg-boundary-domination}, applied after conditioning on $ L^{\mathsf z}=\mathfrak{z}$,
\begin{equation}\label{eq: z copy condition boundary comparison}
    \nu_{\Lambda,\beta}^{\xi}\bigl((\mathcal C^{\xx})^c
\mid  L^{\mathsf z}=\mathfrak{z}\bigr)
\geq
\mu_{\mathfrak{z}}^{\xi^*}\bigl((\mathcal C^{\xx})^c\bigr).
\end{equation}
Under $\mu_{\mathfrak{z}}^{\xi^*}$, on $(\mathcal C^{\xx})^c$, the two boundary clusters have independent uniform signs.
Consequently,
\begin{equation}\label{eq: spin wave boundary symmetry 1}
    \mu_{\mathfrak{z}}^{\xi^*}\bigl((\mathcal C^{\xx})^c
\cap\mathsf{SAME}\bigr)
=
\mu_{\mathfrak{z}}^{\xi^*}(\mathsf{OPP}).
\end{equation}
Consider the auxiliary joint measure
\[
\mathbb Q_{\Lambda,\beta}^{\xi}(d\mathfrak{z},d\omega)=\lambda_{\Lambda,\beta}^{\xi}(d\mathfrak{z})\mu_{\mathfrak{z}}^{\xi^*}(d\omega),
\]
where $\omega$ includes the auxiliary variables used to define the box types. Integrating \eqref{eq: z copy condition boundary comparison} and \eqref{eq: spin wave boundary symmetry 1} over $\mathfrak{z}$ gives \begin{equation}\label{eq: spin wave boundary symmetry 2}
    \nu_{\Lambda,\beta}^{\xi}\bigl((\mathcal C^{\xx})^c
\bigr)
\geq
\mathbb Q_{\Lambda,\beta}^{\xi}\bigl((\mathcal C^{\xx})^c\bigr)\ge\mathbb Q_{\Lambda,\beta}^{\xi}\bigl((\mathcal C^{\xx})^c
\cap\mathsf{SAME}\bigr)
=
\mathbb Q_{\Lambda,\beta}^{\xi}(\mathsf{OPP}).
\end{equation}
Recalling Definition~\ref{def: good sequences and types}, let $\mathcal E$ be the event that the sequence of good boxes and their types are excellent. Recall the definition of $\Gamma$ from \eqref{eq: def of good config for type 2 map}. By Lemma~\ref{lem: type 2 coupling}, we get that \begin{align}
    &\sum_{\sigma}\mu^{\xi^*}_{\mathfrak{z}}(\sigma)\times\sum_{\fB\in\mathfrak{B}_{\rm II}}\Big[\log\Big(\mu^{\sigma|_{ \fB^{( 20)}}}_{\fB^{(20)},\beta'}(F_3^{+}(\sigma))\Big)+\log\Big(\mu^{\sigma|_{ \fB^{( 20)}}}_{\fB^{(20)},\beta'}(F_3^{-}(\sigma))\Big)-2\log\big(\mu^{\sigma|_{ \fB^{( 20)}}}_{\fB^{(20)},\beta'}(\sigma)\big)\Big]\nonumber\\\ge ~&\sum_{\fB\in\mathfrak{B}}-\frac{C_{12}\rho L^2}{\kappa^2\iota^2M^2}\ge -\frac{\bar{C}_1\rho }{\kappa^2\iota}.\label{eq: total energy cost for type 2 map}
\end{align}Applying the Markov inequality to \eqref{eq: def of good config for type 2 map} and combining with \eqref{eq: total energy cost for type 2 map}, we get from Lemma \ref{lem: type 2 coupling} that \begin{equation}
    \int_{\mathfrak{z}}\mu^{\xi^*}_{\mathfrak{z}}(\Gamma^c)\lambda_{\Lambda,\beta}^{\xi}(\mathrm{d}\mathfrak{z})\le \frac{\bar{C}_1\rho }{\kappa^2\iota}.\nonumber
\end{equation}
Combined with Lemmas~\ref{lem: covering property for good boxes} and~\ref{lem: large deviation for type 2 box}, it yields that
\begin{align*}
    &\mathbb Q_{\Lambda,\beta}^{\xi}(\mathcal E^c\cup \Gamma^c)\nonumber\\ \leq~& \lambda_{\Lambda,\beta}^{\xi}\Big( |\Lambda\setminus(\cup_{\fB\in\mathfrak{S}^{\ell_{*}}}\fB)|> c_4\delta|\Lambda|\Big)+ \int_{\mathfrak{z}}\mu^{\xi^*}_{\mathfrak{z}}\Big(\sum_{\fB\in\mathfrak{B}_{\rm II}}|\fB|\ge 6(\kappa+\rho)^{\frac{2}{3}} M^2\Big)\lambda_{\Lambda,\beta}^{\xi}(\mathrm{d}\mathfrak{z})+\frac{\bar{C}_1\rho }{\kappa^2\iota}  \nonumber\\ \leq~& c_7^{-1}\exp(-c_7\iota M)+(\kappa+\rho)^{\frac{1}{3}}+\frac{\bar{C}_1\rho }{\kappa^2\iota}.
\end{align*}

Applying Proposition~\ref{prop: symmetric spin wave} at  $\mathfrak{z}$ and the excellent exploration records $\mathfrak{T}=\Big(\{\mathfrak{L}_\fB\},\{T_\fB\}, \{\mathfrak{e}_{\fB}\}\Big)$, we get that
\begin{equation}\label{eq: spin wave comparison with exploration}
    \mu_{\mathfrak{z}}^{\xi^*}(\Upsilon^0(\mathfrak{T})\cap\mathsf{SAME}\cap \Gamma)
\leq10
\sum_{s\in\{-1,1\}}
\mu_{\mathfrak{z}}^{\xi^*}(\Upsilon^s(\mathfrak{T})\cap\mathsf{OPP}).
\end{equation}
Since $\Upsilon^s(\mathfrak{T})$ are disjoint for different $\mathfrak{T}$, integrating \eqref{eq: spin wave comparison with exploration} over $\mathfrak{T}\in \mathcal{E}$  gives
\[
\mathbb Q_{\Lambda,\beta}^{\xi}(\mathcal E\cap\mathsf{SAME}\cap \Gamma)
\leq 20\mathbb Q_{\Lambda,\beta}^{\xi}(\mathsf{OPP}).
\]
It follows that
\begin{equation}
    1
=
\mathbb Q_{\Lambda,\beta}^{\xi}(\mathsf{SAME})+\mathbb Q_{\Lambda,\beta}^{\xi}(\mathsf{OPP})
\leq
21\mathbb Q_{\Lambda,\beta}^{\xi}(\mathsf{OPP})+c_7^{-1}\exp(-c_7\iota M)+(\kappa+\rho)^{\frac{1}{3}}+\frac{\bar{C}_1\rho }{\kappa^2\iota}.\label{eq: 5.100}
\end{equation}
Take $M$ sufficiently large, $\kappa=\rho^{0.1}>0$ and $\rho>0$ sufficiently small that $c_7^{-1}\exp(-c_7\iota M)+(\kappa+\rho)^{\frac{1}{3}}+\frac{\bar{C}_1\rho }{\kappa^2\iota}\leq1/2$. Then, by \eqref{eq: 5.100} and \eqref{eq: spin wave boundary symmetry 2}, we obtain
\begin{equation}\label{eq: no crossing probability lower bound for one copy}
    \nu_{\Lambda,\beta}^{\xi}\bigl((\mathcal C^{\xx})^c\bigr)
\geq \mathbb Q_{\Lambda,\beta}^{\xi}(\mathsf{OPP})\geq \frac{1}{42}.
\end{equation}
By symmetry, the same bound holds for the $\yy$-copy.

Next, we wish to show that no crossing for the $\xx$ and $\yy$-copies simultaneously with positive probability. To this end, we get around considering the correlation between these two non-crossing events, and instead, we consider two disjoint subannuli (where each annulus blocks an $\xx$- or $\yy$-crossing). To be precise, set $r=(1+\iota)^{1/3}$ and consider
\[
A_1=\Lambda_{M,rM},
\qquad
A_2=\Lambda_{r^2M,r^3M}.
\]
Let $\mathcal{D}_1$ be the event that $A_1$ has no $\xx$-crossing, and let $\mathcal{D}_2$ be the event that $A_2$ has no $\yy$-crossing.
Increasing $M$ if necessary, such that the preceding one-copy estimates apply to both subannuli.

By \eqref{eq: no crossing probability lower bound for one copy}, we see that $\mathcal D_2$ has probability at least $\frac1{42}$. In addition, $\mathcal D_2$ is measurable with respect to configurations outside $A_2$, and thus conditioned on $\mathcal D_2$ one can apply DMP and \eqref{eq: no crossing probability lower bound for one copy} to deduce that $\mathbb P(\mathcal D_1 \mid \mathcal D_2) \geq 1/42$ (note that \eqref{eq: no crossing probability lower bound for one copy} holds for all possible boundary conditions by Remark \ref{rmk: xi^* largest boundary condition}), leading to $\mathbb P(\mathcal D_1 \cap \mathcal D_2) \geq 1/42^2$.
Since $\xx$-crossing of $\Lambda$ contains an $\xx$-crossing of $A_1$ and every $\yy$-crossing of $\Lambda$ contains a $\yy$-crossing of $A_2$,
\[
\nu_{\Lambda_{M,(1+\iota)M},\beta}^{\xi}
\bigl(
\partial\Lambda_M
\con{\mathsf{xy}}
\partial\Lambda_{(1+\iota)M}
\bigr)
\leq 1-\frac{1}{42^2},
\]
which proves the theorem.
\end{proof}

\subsection{Proof of Lemma~\ref{lem: good box probability}}

In this subsection, we prove Lemma~\ref{lem: good box probability}.
Our starting point is the following: for any $\kappa>0$, by Theorem~\ref{thm: apriori--Gaussian tail for number of crossing loops}, there exists some $k_{*}$ such that: (i) the probability for the total number of crossing loops in all copies to be exactly $ k_* $ is at least $ \delta $ and (ii) the probability for the total number of crossing loops being strictly larger than $k_*$ is upper-bounded by $\kappa \delta/k_*$. By Lemma~\ref{lem:yyy},
this indicates that with probability at least $ (\frac{\epsilon}{2+\epsilon})^{k_*}\delta $, there exist $k_*$ crossing loops in the $\zz$-copy, and conditioned on this, the probability for the existence of a crossing loop in the $\xx$- or $\yy$- copy is at most $O( \epsilon)$. Note that in the latter part of the assertion, there is a subtlety in the enumeration of choosing which loop to be in the $\xx$- or $\yy$- copy when we have $k_*+1$ loops as opposed to only $k_*$ loops in the $\zz$-copy (see Lemma~\ref{lem: super exponential decay analysis} for an elaborative version of this computation). This subtle point is why the factor of $1/k_*$ disappears in the conditional probability, and thus is exactly why we need the decay rate in Theorem~\ref{thm: apriori--Gaussian tail for number of crossing loops} to be substantially better than $ e^{-k\log k} $ (a simple argument can yield an upper bound of $ e^{-ck\log k} $ where $c$ is small for large $\beta$).

To summarize, roughly speaking, we may derive the version for the crossing loop of Lemma~\ref{lem: good box probability} relatively easily from Corollary~\ref{cor: Gaussian tail for number of macroscopic loops}.
Then the main difficulty is to derive a small loop cluster crossing probability, taking the small loop crossing probability as input, as incorporated in the next proposition.
Recall from Remark~\ref{rmk: RSW, CBC, FKG} the terminologies of $ \mathsf{RSW} $, $ \mathsf{CBC} $ and $ \mathsf{FKG} $.

\begin{figure}[htb]
    \centering
    \begin{tikzpicture}[
  scale=.85,
  line cap=round,
  line join=round,
  box line/.style={black,line width=.9pt},
  guide/.style={black,line width=.7pt,
    dash pattern=on 6pt off 6pt},
  blue region/.style={fill=regionBlue,draw=none},
  purple loop/.style={draw=loopPurple,line width=.85pt},
  contact/.style={draw=red!85!black,line width=2.2pt},
  yellow connection/.style={draw=yellow!85!orange,line width=1.7pt},
  leader/.style={draw=black,line width=.45pt},
  contact label/.style={font=\small,inner sep=2pt}
]

\definecolor{regionBlue}{RGB}{70,125,235}
\definecolor{loopPurple}{RGB}{145,55,165}

\newcommand{\purpleLoop}[5]{%
  \draw[purple loop,rotate around={#5:(#1,#2)}]
    (#1,#2) ellipse[x radius=#3,y radius=#4];
}


\def\outerCollarInnerEdge{%
  (.30,.34)
    .. controls (1.10,.24) and (1.75,.39) .. (2.55,.30)
    .. controls (3.30,.22) and (3.88,.40) .. (4.65,.31)
    .. controls (5.40,.23) and (6.12,.39) .. (6.85,.30)
    .. controls (7.70,.22) and (8.55,.37) .. (9.66,.30)
    .. controls (9.78,1.20) and (9.60,1.90) .. (9.70,2.70)
    .. controls (9.78,3.50) and (9.60,4.20) .. (9.69,5.05)
    .. controls (9.77,5.90) and (9.62,6.60) .. (9.70,7.40)
    .. controls (9.78,8.30) and (9.60,9.20) .. (9.68,10.68)
    .. controls (8.70,10.77) and (8.10,10.60) .. (7.25,10.69)
    .. controls (6.40,10.78) and (5.72,10.61) .. (4.90,10.70)
    .. controls (4.10,10.78) and (3.38,10.61) .. (2.55,10.69)
    .. controls (1.70,10.77) and (.90,10.62) .. (.31,10.68)
    .. controls (.23,9.60) and (.40,8.80) .. (.30,7.95)
    .. controls (.22,7.10) and (.39,6.40) .. (.31,5.55)
    .. controls (.23,4.70) and (.40,3.90) .. (.30,3.05)
    .. controls (.22,2.20) and (.39,1.20) .. (.30,.34)
  -- cycle
}

\def\innerCollarOuterEdge{%
  (3.91,4.94)
    .. controls (4.35,4.84) and (4.64,4.99) .. (5.02,4.90)
    .. controls (5.44,4.83) and (5.88,4.99) .. (6.32,4.93)
    .. controls (6.42,5.40) and (6.28,5.72) .. (6.36,6.13)
    .. controls (6.43,6.55) and (6.28,6.92) .. (6.32,7.37)
    .. controls (5.86,7.46) and (5.52,7.29) .. (5.12,7.39)
    .. controls (4.70,7.47) and (4.30,7.29) .. (3.92,7.38)
    .. controls (3.83,6.94) and (3.98,6.56) .. (3.90,6.15)
    .. controls (3.82,5.73) and (3.98,5.34) .. (3.91,4.94)
  -- cycle
}

\path[fill=gray!22,even odd rule]
  (0,0) rectangle (10,11)
  (4.20,5.20) rectangle (6.05,7.10);

\path[fill=regionBlue,even odd rule]
  (0,0) rectangle (10,11)
  \outerCollarInnerEdge;

\path[fill=regionBlue,even odd rule]
  \innerCollarOuterEdge
  (4.20,5.20) rectangle (6.05,7.10);

\draw[guide] (1,1.15) rectangle (8.85,10.15);
\draw[guide] (3.05,4.10) rectangle (7.00,8.00);


\path[blue region]
  (4.12,11)
    .. controls (4.00,10.78) and (4.29,10.70) .. (4.18,10.48)
    .. controls (4.08,10.29) and (4.25,10.28) .. (4.22,10.15)
  -- (4.84,10.15)
    .. controls (4.71,10.33) and (4.94,10.45) .. (4.78,10.63)
    .. controls (4.67,10.79) and (4.87,10.90) .. (4.80,11)
  -- cycle;

\path[blue region]
  (7.65,11)
  -- (8.48,11)
    .. controls (8.37,10.80) and (8.72,10.76) .. (8.62,10.55)
    .. controls (8.54,10.38) and (8.77,10.30) .. (8.66,10.15)
  -- (7.63,10.15)
    .. controls (7.75,10.35) and (7.52,10.50) .. (7.64,10.69)
    .. controls (7.73,10.83) and (7.56,10.92) .. (7.65,11)
  -- cycle;

\path[blue region]
  (0,7.18)
    .. controls (.26,7.35) and (.44,7.25) .. (.60,7.48)
    .. controls (.73,7.66) and (.87,7.55) .. (1,7.66)
  -- (1,7.05)
    .. controls (.73,7.12) and (.60,7.05) .. (.44,6.99)
    .. controls (.27,6.94) and (.17,6.90) .. (0,6.94)
  -- cycle;

\path[blue region]
  (0,5.47)
    .. controls (.25,5.51) and (.32,5.28) .. (.55,5.26)
    .. controls (.72,5.24) and (.82,5.01) .. (1,5.01)
  -- (1,4.66)
    .. controls (.75,4.64) and (.61,4.90) .. (.39,4.92)
    .. controls (.22,4.94) and (.12,5.11) .. (0,5.08)
  -- cycle;

\path[blue region]
  (8.85,7.72)
    .. controls (9.12,7.81) and (9.24,7.60) .. (9.46,7.70)
    .. controls (9.67,7.79) and (9.83,7.65) .. (10,7.70)
  -- (10,6.64)
    .. controls (9.78,6.75) and (9.63,6.55) .. (9.43,6.68)
    .. controls (9.19,6.82) and (9.04,6.61) .. (8.85,6.75)
  -- cycle;

\path[blue region]
  (1.60,1.15)
  -- (2.55,1.15)
    .. controls (2.43,.94) and (2.69,.82) .. (2.57,.60)
    .. controls (2.47,.42) and (2.62,.22) .. (2.58,0)
  -- (1.66,0)
    .. controls (1.71,.25) and (1.50,.42) .. (1.64,.63)
    .. controls (1.78,.84) and (1.50,1.02) .. (1.60,1.15)
  -- cycle;

\path[blue region]
  (8.48,1.15)
    .. controls (8.75,1.25) and (8.76,.94) .. (8.99,1.04)
    .. controls (9.17,1.14) and (9.39,1.22) .. (9.50,1.01)
    .. controls (9.60,.81) and (9.43,.64) .. (9.61,.46)
    .. controls (9.76,.31) and (9.70,.11) .. (9.82,0)
  -- (8.50,0)
    .. controls (8.36,.20) and (8.51,.37) .. (8.40,.57)
    .. controls (8.31,.77) and (8.48,.96) .. (8.48,1.15)
  -- cycle;


\path[blue region]
  (5.32,7.10)
    .. controls (5.29,7.35) and (5.47,7.45) .. (5.70,7.48)
    .. controls (6.01,7.50) and (6.15,7.88) .. (6.42,7.81)
    .. controls (6.62,7.75) and (6.77,8.05) .. (7.00,7.96)
  -- (7.00,7.18)
    .. controls (6.99,6.93) and (6.75,6.78) .. (6.50,6.77)
  -- (6.05,6.77)
  -- (6.05,7.10)
  -- cycle;

\path[blue region]
  (4.20,6.56)
    .. controls (3.98,6.61) and (4.04,6.33) .. (4.01,6.15)
    .. controls (3.99,5.95) and (4.03,5.72) .. (3.85,5.65)
    .. controls (3.63,5.55) and (3.47,5.74) .. (3.30,5.58)
    .. controls (3.18,5.47) and (3.13,5.52) .. (3.05,5.48)
  -- (3.05,5.10)
    .. controls (3.27,5.26) and (3.44,5.19) .. (3.64,5.25)
    .. controls (3.87,5.34) and (4.01,5.17) .. (4.20,5.20)
  -- cycle;

\path[blue region]
  (4.70,5.20)
  -- (5.16,5.20)
    .. controls (5.15,4.99) and (5.25,4.91) .. (5.42,4.82)
    .. controls (5.60,4.73) and (5.32,4.58) .. (5.25,4.45)
    .. controls (5.20,4.31) and (5.25,4.21) .. (5.16,4.10)
  -- (4.42,4.10)
    .. controls (4.52,4.29) and (4.82,4.30) .. (4.72,4.48)
    .. controls (4.61,4.65) and (4.58,4.76) .. (4.70,4.92)
    .. controls (4.79,5.05) and (4.65,5.12) .. (4.70,5.20)
  -- cycle;

\path[blue region]
  (6.05,5.32)
  -- (6.05,6.18)
    .. controls (6.22,6.25) and (6.23,6.58) .. (6.43,6.43)
    .. controls (6.61,6.30) and (6.54,6.11) .. (6.76,6.03)
    .. controls (6.94,5.96) and (7.00,5.92) .. (7.00,5.85)
  -- (7.00,5.48)
    .. controls (6.79,5.57) and (6.64,5.75) .. (6.48,5.83)
    .. controls (6.38,5.64) and (6.33,5.34) .. (6.05,5.32)
  -- cycle;

\path[fill=gray!22]
  (7.98,10.15)
  arc[start angle=180,end angle=0,x radius=.17,y radius=.24]
  -- cycle;
\begin{scope}
  \clip (8.15,10.15) ellipse[x radius=.17,y radius=.24];
  \draw[guide] (1,1.15) rectangle (8.85,10.15);
\end{scope}


\path[fill=gray!22,draw=none]
  (10,8.86)--(9.43,8.95)--(9.43,9.35)--cycle
  (10,3.92)--(9.43,3.55)--(9.43,3.86)--cycle
  (5.91,0)--(5.98,.57)--(6.40,.57)--cycle
  (0,2.56)--(.57,2.62)--(.57,3.00)--cycle
  (0,6.16)--(.57,5.80)--(.57,6.12)--cycle
  (2.62,11)--(2.67,10.43)--(3.05,10.43)--cycle
  (6.22,11)--(5.78,10.43)--(6.12,10.43)--cycle;

\path[fill=gray!22,draw=none]
  (6.05,6.58)--(6.65,6.59)--(6.65,6.75)--cycle
  (5.81,5.20)--(5.47,4.65)--(5.79,4.65)--cycle
  (4.40,5.20)--(4.22,4.70)--(4.48,4.70)--cycle
  (4.66,7.10)--(4.72,7.65)--(5.08,7.65)--cycle;


\draw[box line] (0,0) rectangle (10,11);

\draw[yellow connection]
  (2.32,1.15)
    .. controls (2.45,1.56) and (3.05,1.62) .. (3.30,1.10)
    .. controls (3.53,.61) and (4.21,.61) .. (4.48,1.15)
    .. controls (4.85,1.89) and (5.78,2.45) .. (5.65,3.25)
    .. controls (5.55,3.86) and (5.08,3.65) .. (5.08,4.10);


\draw[contact] (4.22,10.15)--(4.84,10.15);
\draw[contact] (7.63,10.15)--(7.98,10.15);
\draw[contact] (8.32,10.15)--(8.66,10.15);

\draw[contact] (1,7.05)--(1,7.66);
\draw[contact] (1,4.66)--(1,5.01);
\draw[contact] (8.85,6.75)--(8.85,7.72);

\draw[contact] (1.60,1.15)--(2.55,1.15);
\draw[contact] (8.48,1.15)--(8.66883501,1.15);

\draw[contact] (3.05,5.10)--(3.05,5.48);
\draw[contact] (4.42,4.10)--(5.16,4.10);
\draw[contact] (7.00,7.18)--(7.00,7.96);
\draw[contact] (7.00,5.48)--(7.00,5.85);

\node[contact label] (Aone) at (8.15,9.05) {$A_1$};
\draw[leader] (7.805,10.15)--(Aone.north west);
\draw[leader] (8.49,10.15)--(Aone.north east);

\node[contact label] (Atwo) at (8.00,6.60) {$A_2$};
\draw[leader] (8.85,7.235)--(Atwo.east);

\node[contact label] (Athree) at (8.40,1.85) {$A_3$};
\draw[leader] (8.574417505,1.15)--(Athree.south);

\node[contact label] (Afour) at (2.075,1.85) {$A_4$};
\draw[leader] (2.075,1.15)--(Afour.south);

\node[contact label] (Afive) at (2.00,4.60) {$A_5$};
\draw[leader] (1,4.835)--(Afive.west);

\node[contact label] (Asix) at (1.70,7.355) {$A_6$};
\draw[leader] (1,7.355)--(Asix.west);

\node[contact label] (Aseven) at (4.53,9.45) {$A_7$};
\draw[leader] (4.53,10.15)--(Aseven.north);

\node[contact label] (Bone) at (7.80,8.12) {$B_1$};
\draw[leader] (7.00,7.57)--(Bone.west);

\node[contact label] (Btwo) at (7.95,5.40) {$B_2$};
\draw[leader] (7.00,5.665)--(Btwo.west);

\node[contact label] (Bthree) at (4.79,3.22) {$B_3$};
\draw[leader] (4.79,4.10)--(Bthree.north);

\node[contact label] (Bfour) at (2.20,5.80) {$B_4$};
\draw[leader] (3.05,5.29)--(Bfour.east);

\draw[box line] (4.20,5.20) rectangle (6.05,7.10);
\node[font=\large] at (5.125,6.15) {$\fB$};
\node[anchor=south west,inner sep=3pt,font=\small]
  at (10,11) {$\fB^{(2)}$};


\purpleLoop{4.18}{11.13}{.25}{.55}{-35}
\purpleLoop{4.83}{11.10}{.47}{.29}{30}
\purpleLoop{4.02}{10.13}{.38}{.24}{40}
\purpleLoop{4.84}{10.28}{.43}{.21}{8}

\purpleLoop{7.96}{11.15}{.53}{.25}{40}
\purpleLoop{7.57}{10.38}{.49}{.27}{25}
\purpleLoop{7.89}{10.10}{.22}{.48}{10}
\purpleLoop{8.39}{10.24}{.25}{.47}{18}
\purpleLoop{8.70}{10.66}{.45}{.23}{15}

\purpleLoop{.04}{7.58}{.35}{.52}{32}
\purpleLoop{.39}{7.27}{.36}{.22}{20}

\purpleLoop{.03}{5.57}{.42}{.31}{-15}
\purpleLoop{.44}{5.31}{.24}{.43}{-38}
\purpleLoop{.03}{4.83}{.40}{.26}{20}
\purpleLoop{1.14}{4.86}{.44}{.32}{18}

\purpleLoop{8.99}{7.73}{.25}{.42}{18}
\purpleLoop{9.57}{7.71}{.26}{.43}{-28}
\purpleLoop{8.63}{6.98}{.49}{.23}{12}
\purpleLoop{9.22}{6.70}{.23}{.48}{10}
\purpleLoop{9.74}{6.73}{.26}{.46}{25}

\purpleLoop{1.69}{1.17}{.29}{.47}{24}
\purpleLoop{2.54}{1.10}{.35}{.49}{-12}
\purpleLoop{1.53}{.29}{.61}{.24}{15}
\purpleLoop{2.28}{-.13}{.34}{.63}{-6}

\purpleLoop{8.36}{1.03}{.34}{.55}{14}
\purpleLoop{9.10}{1.27}{.34}{.58}{-35}
\purpleLoop{9.87}{.76}{.67}{.28}{32}
\purpleLoop{8.30}{.08}{.62}{.24}{12}
\purpleLoop{9.05}{-.10}{.32}{.52}{12}
\purpleLoop{9.73}{-.08}{.29}{.60}{24}


\purpleLoop{5.48}{7.81}{.27}{.56}{15}
\purpleLoop{6.15}{7.94}{.29}{.58}{-20}
\purpleLoop{6.89}{7.77}{.57}{.19}{2}
\purpleLoop{7.08}{6.97}{.50}{.24}{-18}
\purpleLoop{6.26}{6.94}{.19}{.15}{18}

\purpleLoop{3.97}{6.50}{.23}{.37}{30}
\purpleLoop{3.61}{5.66}{.24}{.42}{-15}
\purpleLoop{3.63}{5.24}{.25}{.43}{12}
\purpleLoop{2.89}{5.18}{.57}{.24}{30}

\purpleLoop{4.86}{4.95}{.20}{.25}{-12}
\purpleLoop{4.44}{4.63}{.52}{.24}{28}
\purpleLoop{5.48}{4.86}{.49}{.20}{12}
\purpleLoop{5.50}{4.46}{.46}{.19}{-8}
\purpleLoop{4.70}{3.98}{.27}{.50}{10}

\purpleLoop{6.66}{6.34}{.44}{.22}{28}
\purpleLoop{6.49}{5.91}{.22}{.47}{-28}
\purpleLoop{7.22}{5.96}{.52}{.28}{-8}
\purpleLoop{6.38}{5.35}{.22}{.46}{8}

\purpleLoop{.95}{10.82}{.48}{.22}{12}
\purpleLoop{3.43}{10.86}{.43}{.25}{15}
\purpleLoop{5.58}{10.84}{.44}{.24}{-12}
\purpleLoop{6.90}{10.83}{.43}{.23}{18}
\purpleLoop{9.45}{10.84}{.40}{.22}{-15}

\purpleLoop{.65}{.16}{.39}{.21}{20}
\purpleLoop{3.40}{.16}{.46}{.23}{-12}
\purpleLoop{4.70}{.15}{.50}{.22}{12}
\purpleLoop{6.80}{.14}{.43}{.23}{-14}
\purpleLoop{7.60}{.18}{.35}{.20}{18}

\purpleLoop{4.28}{7.30}{.29}{.18}{12}
\purpleLoop{5.06}{7.25}{.28}{.18}{-12}

\end{tikzpicture}
    \caption{ Illustration for the exploration $\mathcal{F}$. Blue region: the local time clusters. Red parts:  $A_1,\cdots,A_t,B_1,\cdots,B_s$. Purple loops: $\mathfrak{L}^\yy_1$ and $\mathfrak{L}^\yy_2$. Yellow line: local time connection between $A_4$ and $B_3$ in $\Omega.$}
    \label{fig: local time crossing exploration}
\end{figure}
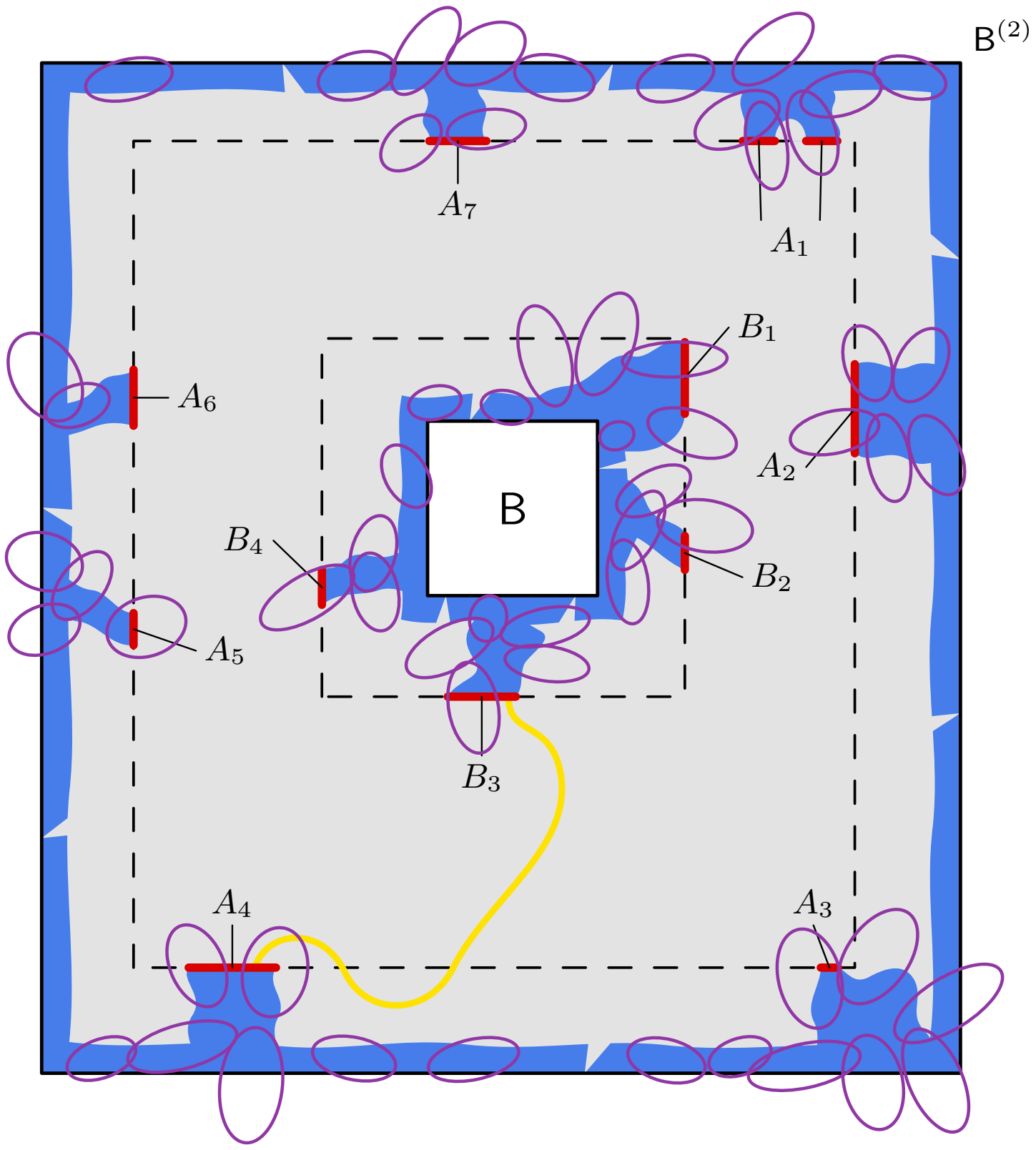

\begin{proposition}\label{prop: small loop crossing implies small cluster crossing}
    Let $\fB$ be a dyadic box and let $\fB'$ be a box concentric with $\fB$ {such that $\fB^{(2)}\subset \fB'$}. Let $ \gamma=(\gamma_u) $ be a local time condition with $0\leq \gamma_u\leq\beta_0$ for all vertex $u$.
    Assume that for some boundary condition $\eta$ on $\partial\fB'$, we have \begin{equation}
        \mu^\eta_{\fB', \gamma}\Big(\text{there exists an $\mathsf{x}$-loop intersecting $\fB^{(2)}\setminus\fB$ with diameter at least ${{\frac{\delta}{2}}\cdot \diam(\fB)}$} \Big)\le \kappa.\label{eq: small loop crossing probability for good box}
    \end{equation}Then there exists a constant $C_{15}>0$ depending only on $\beta_0$ such that
    \begin{equation}
        \mu^\eta_{\fB', \gamma}\Big(\partial\fB\con{\mathsf{x}}\partial \fB^{(2)}\Big)\le C_{15}\kappa^{0.1}.\label{eq: small cluster crossing probability for good box}
    \end{equation}
\end{proposition}
\begin{proof}
We first sample the $\xx$-local time clusters from $\partial\fB^{(2)}$ restricted to the complement of $ \mathsf{B}^{(1.5)} $.
Let $\mathfrak{C}_1$ denote the collection of $\xx$-copy clusters and let $ \mathcal{C}_{1,1}, \cdots, \mathcal{C}_{1,t} \in \mathfrak{C}_1 $ denote the $\mathsf{x}$-copy clusters intersecting  $\partial\fB^{(1.5)}$. In addition, let $A_i=\mathcal{C}_{1,i}\cap \partial\fB^{(1.5)}$.
See Figure \ref{fig: local time crossing exploration} for an illustration. We further sample the $\yy$-loops $\mathfrak{L}_1^\yy$ intersecting $\cup \mathfrak{C}_1$ where $\cup \mathfrak{C}_1$ denotes the union of the traces of the clusters in $ \mathfrak{C}_1$. Similarly, we sample the $\xx$-local time clusters from $\partial\fB$ restricted to $ \mathsf{B}^{(0.5)}\setminus\mathsf{B} $. Let $\mathfrak{C}_2$ denote the collection of these $\xx$-copy clusters and let $\mathcal{C}_{2,1},\cdots,\mathcal{C}_{2,s}\in\mathfrak{C}_2$ denote those intersecting  $\partial\fB^{(0.5)}$. In addition, let $B_j=\mathcal{C}_{2,j}\cap \partial\fB^{(0.5)}$. We further sample the $\yy$-loops $\mathfrak{L}_2^\yy$ intersecting $\cup \mathfrak{C}_2$.

Let $\Xi$ be the collection of possible realizations $\omega$ of $(\mathfrak{C}_1,\mathfrak{C}_2,\mathfrak{L}_1^\yy,\mathfrak{L}_2^\yy)$. We write $\Omega=\Omega(\omega)=\fB^{(2)}\setminus\big(\fB\cup(\cup \mathfrak{C}_1)\cup(\cup \mathfrak{C}_2)\big)$ and write $ \gamma^{\prime}=\gamma^{\prime}(\omega)=(\gamma^{\prime}_u) $ such that
\begin{align*}
    &\gamma'_u=\gamma_u-(\mathfrak{L}_1^\yy\cup\mathfrak{L}_2^\yy)_u,\quad \text{for } u\in \Omega\setminus\partial\Omega
\end{align*}
where $(\mathfrak{L}_i^\yy)_u$ denotes the local time accumulated at $u$ from loops in $\mathfrak{L}_i^\yy$.
Let $\xi=\xi(\omega)$ denote the boundary condition such that
\begin{equation}
    \xi^\mathsf{x}_u=
        \begin{cases}
            \gamma_u-(\mathfrak{L}_1^\yy \cup \mathfrak{L}_2^\yy)_u, & u\in \cup_{i=1}^tA_i\cup \cup_{j=1}^sB_j,\\
            0, & u\in \partial\Omega\setminus\big(\cup_{i=1}^tA_i\cup \cup_{j=1}^sB_j\big),
        \end{cases}
    ~~{  \rm and }~\xi^\mathsf{y}_u=0,\quad\forall u\in\partial\Omega.
\end{equation}
and that each $A_i$ for $1\leq i\leq t$ and each $B_j$ for $1\leq j\leq s$ are wired into single points correspondingly.
In conclusion, we get that $ \mu^\eta_{\fB', \gamma}(\cdot|_{\Omega}\mid \omega)=\mu^\xi_{\Omega,\gamma'}(\cdot)$. Let $\mathcal{A}$ denote the event that there exists an $\mathsf{x}$-loop intersecting $\Omega$ with diameter at least ${\delta\cdot \diam(\fB)}$.
By Lemma \ref{lem:cpf-spatial-restriction} (see also Remark \ref{rmk: spatial restriction for excursion}) and Corollary \ref{cor: DMP for loop}, we get that
\begin{align}
    \mu^\eta_{\fB', \gamma}\Big(\partial\fB\con{\mathsf{x}}\partial \fB^{(2)}\Big)&=\int_{\omega\in \Xi}\mu^\eta_{\fB', \gamma}(\omega)\cdot\mu^\xi_{\Omega,\gamma^{\prime}} \Big(\partial\fB\con{\mathsf{x}}\partial \fB^{(2)}\Big),\label{eq: annulus config exploration in loop+local time relation}\\
     \mu^\eta_{\fB', \gamma}\Big(\mathcal{A}\Big)&=\int_{\omega\in \Xi}\mu^\eta_{\fB', \gamma}(\omega)\cdot\mu^\xi_{\Omega,\gamma^{\prime}} \Big(\mathcal{A}\Big).\label{eq: annulus config exploration in loop+local time relation 2}
\end{align}
Fix $\omega\in\Xi$.
On the event $\{\partial\fB\con{\mathsf{x}}\partial \fB^{(2)}\}$, there exist $1\le i\le t$ and $1\le j\le s$ such that $A_i\con{\mathsf{x}}B_j$.
Then by Lemma~\ref{lem: cluster-loop comparison, all-x local time b.c.} (see also Remark \ref{rmk: boundary two point function}), we have
    \begin{align}
            \mu^{\xi}_{\Omega,\gamma'}(\partial\fB\con{\mathsf{x}}\partial \fB^{(2)})
            & \leq \sum_{1\leq i\leq t,1\leq j\leq s}{\mu}^{\xi}_{\Omega, \gamma'}(A_i\stackrel{\xx}{\longleftrightarrow}B_j)\nonumber \\
            & \leq 2^{3/4}\sum_{1\leq i\leq t,1\leq j\leq s}{\mu}^{\xi}_{\Omega, \gamma'}(\mathcal{A})^{1/4} \quad \text{(by Lemma \ref{lem: cluster-loop comparison, all-x local time b.c.})}\nonumber \\
            & \leq 2st\cdot {\mu}^{\xi}_{\Omega, \gamma'}(\mathcal{A})^{1/4}.\label{eq: 5.107} 
    \end{align}
    Applying Jensen's inequality and \eqref{eq: annulus config exploration in loop+local time relation 2} in order, we get that 
    \begin{align}
        &\int_{\omega\in \Xi}\mu^\eta_{\fB', \gamma}(\omega)\cdot  \big[\mu^\xi_{\Omega,\gamma'}(\mathcal{A})\big]^{\frac14}\le \Big[\int_{\omega\in \Xi}\mu^\eta_{\fB', \gamma}(\omega)\cdot  \mu^{\xi}_{\Omega',\gamma'}(\mathcal{A})\Big]^{\frac14}= \Big[\mu^\eta_{\fB',\gamma}(\mathcal{A})\Big]^{\frac14}.\label{eq: connecting prob after exploration split 2}
    \end{align}
    Now we turn to control $s,t$. We claim that
    \begin{equation}
        \mu^\eta_{\fB', \gamma}\Big(\big\{\omega\in \Xi: s\ge k\big\}\Big)\le \exp(-\bc_1k).\label{eq: disjoint cluster crossing exp decay}
    \end{equation}
    We can explore the local time clusters crossing from $\partial\fB^{(2)}$ through $\fB^{(2)}\setminus\fB^{(1.5)}$ one by one. Since $\mathsf{RSW}$ holds uniformly for any region, we get that conditioned on the exploration history, the probability of getting a new local time cluster crossing through $\fB^{(2)}\setminus\fB^{(1.5)}$ is at most $\bc_2<1$. Multiplying the probability over the $k$ explorations, we arrive at \eqref{eq: disjoint cluster crossing exp decay} (note that  with some more work following~\cite{AizenmanPercolationSuperExpDecay}, one can prove a Gaussian decay here). By symmetry, the same bound holds replacing $s$ in \eqref{eq: disjoint cluster crossing exp decay} by $t$.
    Let $\mathcal{B}(k)\subset\Xi$ denote the collection of configurations $\omega\in \Xi$ such that $s\le k $ and $t\le k$. From (both versions) of \eqref{eq: disjoint cluster crossing exp decay}, we get that
    \begin{equation}
        \mu^\eta_{\fB', \gamma}\big(\mathcal{B}(k)\big)\ge1- 2\exp(-\bc_1k).\label{eq: disjoint cluster crossing exp decay final}
    \end{equation}
    
    Plugging \eqref{eq: small loop crossing probability for good box} into \eqref{eq: connecting prob after exploration split 2} and combining with \eqref{eq: 5.107}, we obtain that \begin{align}
        \int_{\omega\in \mathcal{B}(k)}\mu^\eta_{\fB', \gamma}(\omega)\cdot\mu^\xi_{\Omega,\gamma'} \Big(\partial\fB\con{\mathsf{x}}\partial \fB^{(2)} \Big)&\le k^2\int_{\omega\in \Xi}\mu^\eta_{\fB', \gamma}(\omega)\cdot  \big[\mu^\xi_{\Omega,\gamma'}(\mathcal{A})\big]^{\frac14}\nonumber\\&\le k^2\Big[\mu^\eta_{\fB',\gamma}(\mathcal{A})\Big]^{\frac14}\le k^2 \kappa^{\frac{1}{4}}.\label{eq: small cluster crossing probability for good box step 1}
    \end{align}Combining \eqref{eq: small cluster crossing probability for good box step 1} and \eqref{eq: disjoint cluster crossing exp decay final} and plugging into \eqref{eq: annulus config exploration in loop+local time relation}, we obtain that \begin{align}
        \mu^\eta_{\fB',\gamma} \Big(\partial\fB\con{\mathsf{x}}\partial \fB^{(2)}\Big)&\le \int_{\omega\in \mathcal{B}(k)}\mu^\eta_{\fB', \gamma}(\omega)\cdot\mu^\xi_{\Omega,\gamma'} \Big(\partial\fB\con{\mathsf{x}}\partial \fB^{(2)}\Big)+\mu^\xi_{\fB',\gamma}\Big(\mathcal{B}(k)\big)^c\Big)\nonumber\\&\le k^2 \kappa^{\frac{1}{4}}+ 2\exp(-\bc_1k).\nonumber
    \end{align}Choosing $k=\kappa^{-0.01}$ and taking $C_{15}>0$ large enough finishes the proof of \eqref{eq: small cluster crossing probability for good box}.
\end{proof}

Next, we prove that the bound in~\eqref{eq: small cluster crossing probability for good box} is not sensitive to the boundary condition.

\begin{lemma}\label{lem: weak strong cluster crossing prob}
    Let $\fB$ be a dyadic box and let $\fB'$ be a box concentric with $\fB$ such that $\fB^{(10)}\subset \fB'$. Let $\gamma=(\gamma_u)$ be a local time condition with $0\leq \gamma_u\le \beta_0 $ for all vertex $ u $. Assume that for some boundary condition $\eta$ on $\partial\fB'$, we have \begin{equation}
        \mu^\eta_{\fB', \gamma}\Big(\partial\fB\con{\mathsf{x}}\partial \fB^{(2)}\Big)\le \kappa.\label{eq: small cluster crossing probability for good box weak}
    \end{equation}Then there exists a constant $C_{16}>0$ depending only on $\beta_0$ and $\delta$ such that for any boundary condition $\xi$ on $\partial\fB'$ we have \begin{equation}
        \mu^\xi_{\fB', \gamma}\Big(\partial\fB\con{\mathsf{x}}\partial \fB^{(2)}\Big)\le C_{16}\kappa.\label{eq: small cluster crossing probability for good box strong}
    \end{equation}
\end{lemma}
\begin{proof}
    For notation convenience, we write $\mathcal{E}_{\mathsf{x}}:=\{\partial\fB\con{\mathsf{x}}\partial \fB^{(2)}\}$ and write $\mathcal{E}_{\mathsf{y}}:=\{\partial\fB\con{\mathsf{y}}\partial \fB^{(2)}\}$. By $ \mathsf{CBC} $, it suffices to consider the case that $\eta$ is the worst boundary condition for $\mathcal{E}_{\mathsf{x}}$ and $\xi$ is the best boundary condition for $\mathcal{E}_{\mathsf{x}}$ and show that \eqref{eq: small cluster crossing probability for good box weak} implies \eqref{eq: small cluster crossing probability for good box strong}. By symmetry and \eqref{eq: small cluster crossing probability for good box weak}, we get that
    \begin{equation}
        \mu^\xi_{\fB', \gamma}\Big(\mathcal{E}_{\mathsf{y}}\Big)\le \kappa.\label{eq: y cannot cluster cross}
    \end{equation}
To control the boundary influence to $\fB^{(2)}$, define an event $\mathcal{A}:=\{\partial\fB'\ncon[0.5]{\mathsf{x}}\partial \fB^{(2)}\}$. 
Let $\mathfrak{C}$ denote the collection of realizations for boundary $\xx$-clusters $\mathcal{C}$ from $\partial\fB'$ that satisfies $\mathcal{A}$. Note that $\mathcal{C}$ stands not only for the vertex set $V_{\mathcal{C}}$ but also for the local time $\gamma^{\mathcal{C}}$ on those vertices. Let $\mu^{\xi,\mathcal{C}}_{\fB',\gamma}=\mu^{\xi}_{\fB',\gamma}(\cdot\mid \mathcal{C})$ denote the XY-measure on $\fB'$ with local time condition $\gamma$ on $\fB'\setminus V_{\mathcal{C}}$ and boundary condition
\begin{equation}
L^{\xx}_v=0,\quad L^{\yy}_v=\gamma_v- \gamma^{\mathcal{C}}_v   ,\quad\forall v\in  V_{\mathcal{C}}.
\end{equation}
By $ \mathsf{CBC} $, we get that \begin{equation}
\mu^{\xi,\mathcal{C}}_{\fB',\gamma}\Big(\mathcal{E}_{\mathsf{x}}\Big)\le \mu^{\xi,\mathcal{C}}_{\fB',\gamma}\Big(\mathcal{E}_{\mathsf{y}}\Big).\nonumber\label{eq: no crossing is stronger than boundary condition step 1}
\end{equation}Summing over $\mathcal{C}\in \mathfrak{C}$, we obtain that \begin{align}
\mu^\xi_{\fB', \gamma}\Big(\mathcal{A}\cap \mathcal{E}_{\mathsf{x}}\Big)=\sum_{\mathcal{C}\in \mathfrak{C}} \mu^\xi_{\fB', \gamma}(\mathcal{C})\cdot    \mu^{\xi,\mathcal{C}}_{\fB',\gamma}\Big(\mathcal{E}_{\mathsf{x}}\Big)\le \sum_{\mathcal{C}\in \mathfrak{C}} \mu^\xi_{\fB', \gamma}(\mathcal{C})\cdot    \mu^{\xi,\mathcal{C}}_{\fB',\gamma}\Big(\mathcal{E}_{\mathsf{y}}\Big)=\mu^{\xi}_{\fB',\gamma}\Big(\mathcal{A}\cap \mathcal{E}_{\mathsf{y}}\Big)\le \mu^{\xi}_{\fB',\gamma}\Big(\mathcal{E}_{\mathsf{y}}\Big).\nonumber\label{eq: no crossing is stronger than boundary condition step 2}
\end{align}
Recalling \eqref{eq: y cannot cluster cross}, to prove \eqref{eq: small cluster crossing probability for good box weak} it remains to show that
\begin{equation}
\mu^{\xi}_{\fB',\gamma}\Big(\mathcal{E}_{\mathsf{x}}\Big)\le \bC_1 \mu^{\xi}_{\fB',\gamma}\Big(\mathcal{A}\cap \mathcal{E}_{\mathsf{x}}\Big)\label{eq: up to constant reverse FKG}
\end{equation}
for some constant $\bC_1>0$ depending only on $\beta_0$. This is equivalent to lower-bounding the conditional probability of $\mathcal A$ given $\mathcal E_\mathsf x$. Note that the event $\mathcal A$ is measurable with respect to the local time field restricted in $\mathsf B^{(2)} \setminus \mathsf B$. Conditioned on local times thereon, by $\mathsf{CBC}$ the $\mathsf x$-copy local time field on $\mathsf B'\setminus \mathsf B^{(2)}$ is dominated by $\mu^{\xi^*}_{\mathsf B'\setminus \mathsf B^{(2)}, \gamma}$ where $\xi^*$ is the all-$\mathsf x$-copy local time boundary condition on the annulus $\mathsf B'\setminus \mathsf B^{(2)}$ (recall Definition \ref{def: best x bc for crossing}). Combined with $\mathsf{RSW}$, we get that \begin{equation}
    \mu^{\xi}_{\fB',\gamma}\Big(\mathcal{A}\mid \mathcal{E}_{\xx}\Big)\ge  \mu^{\xi^*}_{\mathsf B'\setminus \mathsf B^{(2)}, \gamma}\Big(\mathcal{A}\Big)\ge \bc_2.
\end{equation}
Hence, we finish the proof of \eqref{eq: up to constant reverse FKG}, thereby completing the proof of the lemma.
\end{proof}
Combining Proposition \ref{prop: small loop crossing implies small cluster crossing} with Lemma \ref{lem: weak strong cluster crossing prob}, we arrive at the following corollary.
\begin{corollary}\label{cor: small loop crossing implies strong cluster crossing}
     Let $\fB$ be a dyadic box and let $\fB'$ be a box concentric with $\fB$ such that $\fB^{(10)}\subset \fB'$. Let $ \gamma=(\gamma_u) $ be a local time condition with $0\leq \gamma_u\leq\beta_0$ for all vertex $u$. Assume that for \textbf{some} boundary condition $\eta$ on $\partial\fB'$, we have \begin{equation}
        \mu^\eta_{\fB', \gamma}\Big(\text{there exists an $\mathsf{x}$-loop intersecting $\fB^{(2)}\setminus\fB$ with diameter at least ${\delta\cdot \diam(\fB)}$} \Big)\le \kappa.\label{eq: small loop crossing probability for good box 2}
    \end{equation}Then there exists a constant $C_{17}>0$ depending only on $\beta_0$ such that for \textbf{all} boundary condition $\xi$ on $\partial\fB'$ we have \begin{equation}
        \mu^\xi_{\fB', \gamma}\Big(\partial\fB\con{\mathsf{x}}\partial \fB^{(2)}\Big)\le C_{17}\kappa^{0.1}.\label{eq: small cluster crossing probability for good box strong 2}
    \end{equation}
\end{corollary}

We also need the following lemma of calculus.
\begin{lemma}\label{lem: super exponential decay analysis}
For any $\rho,c_{18}>0$, there exist $c_{19},\, K_0>0$ depending only on $ (c_{18}, \rho, \epsilon)$ such that the following holds. 
    For any sequence of  non-negative numbers $ (q_k)_{k\ge 0}$ such that $\sum_{k=0}^{\infty}q_k=1$, $q_k\le c_{18}^{-1}\exp(-c_{18}k^2)$ for $k\le K_0$ and $\sum_{k=K_0+1}^{\infty}q_{k}\le c_{18}^{-1}\exp(-c_{18}(K_0+1)^2)$, there exists $k_*>0$ with $k_*\le \frac{1}{c_{19}}$,  $q_{k_*}\ge c_{19}$ and
       \begin{equation}
    \label{eq: proof of lemma 4.5 - 2}
    \sum_{m\geq k_*+1} q_m\binom{m}{k_*}\big(\frac{2}{2+\epsilon}\big)^{m-k_*}\leq \rho q_{k_*}.
\end{equation}
\end{lemma}
\begin{proof}
{Let $K_0>0$ be a large enough constant to be determined later. Let $\bar{q}_k=q_k\1_{k\le K_0}$, then we have $\bar{q}_k\le c_{18}^{-1}\exp(-c_{18}k^2)$ for any $q\ge0$.}
Let
$$s=\frac{2+\epsilon}{2}\log(1+\rho/2)\qquad{\rm and}
\qquad
b_k=\frac{k!\,\bar{q}_k}{s^k}\quad{\rm for}
\quad k\geq0.$$
For $k\geq1$, the upper bound of $ \bar{q}_k $ and the inequality $k!\leq k^k$ give
$$
0\leq b_k
\leq c_{18}^{-1}\exp\!\left\{-c_{18} k^{2}+k\log k-k\log s\right\}
\longrightarrow0\quad{\rm as}\quad k\to\infty.
$$
Thus, the sequence $(b_k)_{k\geq0}$ attains a positive maximum at some $k_*\geq0$.
In particular, $\bar{q}_{k_*}>0$.
By maximality, for every $j\geq1$,
\[
\frac{(k_*+j)!\,\bar{q}_{k_*+j}}{s^{k_*+j}}
\leq \frac{k_*!\,\bar{q}_{k_*}}{s^{k_*}},
\]
and hence
\[
\bar{q}_{k_*+j}\binom{k_*+j}{k_*}
\leq q_{k_*}\frac{s^j}{j!}.
\]
Consequently,
\begin{equation}\label{eq: eq: 5.116}
    \sum_{m\geq k_*+1}\bar{q}_m\binom{m}{k_*}
\big(\frac{2}{2+\epsilon}\big)^{m-k_*}\leq \bar{q}_{k_*}\sum_{j\geq1}\frac{(\frac{2s}{2+\epsilon})^j}{j!}=\bar{q}_{k_*}\big(e^{\frac{2s}{2+\epsilon}}-1\big)=\frac\rho2 \bar{q}_{k_*}.
\end{equation}
Next, we turn to lower-bound $q_{k_*}$, and the upper bound of $k_*$ will follow directly from the lower bound of $q_{k_*}$ and the assumption $q_k\le c_{18}^{-1}\exp(-c_{18}k^2).$  Let $k_0=\lceil\frac{1}{c_{18}}\rceil+10$. Then we calculate \begin{align}
    \max_{1\le i\le k_0}\bar{q}_i&\ge\frac{1}{k_0}\sum_{k=0}^{k_0-1}\bar{q}_k=\frac{1}{k_0}(1-\sum_{k=k_0}^{\infty}q_k)\ge \frac{1}{k_0}\big(1-\sum_{k=k_0}^{K_0+1}c_{18}^{-1}\exp(-c_{18}k^2)\big)\nonumber\\& \ge \frac{1}{k_0}\big(1-\sum_{k=k_0}^{\infty}c_{18}^{-1}\exp(-k)\big)
    \ge \frac{1}{k_0}\cdot\big(1-\frac{c_{18}^{-1}\exp(-k_0)}{1-e^{-1}}\big)\ge \frac{1}{2k_0}.\label{eq: lower bouond for 10}
\end{align}Since $b_k\rightarrow 0$ as $k\rightarrow\infty$, there exists an integer $k_1$ such that for any $k\ge k_1$, we have \begin{equation}
    b_{k}\le \frac{1}{2k_0}\le \max_{1\le i\le k_0}q_i\le  \max_{1\le i\le k_0}b_i.\label{eq: b0 lower bound}
\end{equation} Therefore, we get that $k_*\le k_1$. By maximality of $b_{k_*}$, we conclude that \begin{equation}\label{eq: qk lower bound}
    q_{k_*}=\frac{b_{k_*}s^{k_*}}{(k_*)!}\ge \frac{\max_{1\le i\le k_0}\{b_i\}\times s^{k_*}}{(k_*)!}\stackrel{\eqref{eq: b0 lower bound}}{\ge}\frac{b_0\min\{1,s\}^{k_1}}{(k_1)!}\ge \frac{\min\{1,s\}^{k_1}}{2k_0(k_1)!}.
\end{equation}
Noting that both $k_0$ and $k_1$ depend only on $c,\epsilon$ and $\rho$, {we get the desired bound for $k_*$ and $q_{k_*}$. In addition, we can choose $K_0$ large enough such that for any $m\ge K_0$ \begin{equation}
    \binom{m}{k_*}(\frac{2}{2+\epsilon})^{m-k_*}\le \frac{\rho}{2}\frac{\min\{1,s\}^{k_1}}{2k_0(k_1)!}.\label{eq: choose of K_0}
\end{equation}Plugging \eqref{eq: choose of K_0} into \eqref{eq: eq: 5.116}, we get that \begin{align}\label{eq: cutoff add back}
    \sum_{m\geq k_*+1}\bar{q}_m\binom{m}{k_*}
\big(\frac{2}{2+\epsilon}\big)^{m-k_*}&\le \sum_{m\geq k_*+1}\bar{q}_m\binom{m}{k_*}
\big(\frac{2}{2+\epsilon}\big)^{m-k_*}+\sum_{m\geq K_0+1}\binom{m}{k_*}
\big(\frac{2}{2+\epsilon}\big)^{m-k_*}\nonumber\\&\le \frac\rho2 \bar{q}_{k_*}+\frac{\rho}{2}\frac{\min\{1,s\}^{k_1}}{2k_0(k_1)!}.
\end{align}\eqref{eq: proof of lemma 4.5 - 2} then follows from combining \eqref{eq: cutoff add back} and \eqref{eq: qk lower bound}.
}
\end{proof}
We are now ready to prove Lemma~\ref{lem: good box probability}.
\begin{proof}[Proof of Lemma~\ref{lem: good box probability}]
Let $\mathcal{A}_k$ denote the event that there are exactly $k$ loops with diameters at least $\delta\cdot\diam(\fB)$ intersecting $\fB^{(2)}\setminus\fB$. In addition, let $\mathcal{A}_k^\zz$ denote the event that there are exactly $k$ $\zz$-copy loops with diameters at least $\delta\cdot\diam(\fB)$ intersecting $\fB^{(2)}\setminus\fB$. Let $\bC_1>0$ be a large enough constant to be chosen later and write $\fB'=\fB^{(\bC_1)}$. Let $\mathcal{B}$ denote the event $\partial\fB'\ncon[1]{\xx\yy\zz\text{-loop}}\partial\fB^{(2)}$ (here we slightly abuse the terminology of loop by letting it denote loop and excursion trajectories) and let $\Gamma$ denote the collection of realizations for loops   from $ \partial \mathsf{B}^{\prime} $ such that the event $\mathcal{B}$ occurs. For $ \gamma\in\Gamma $, let $ L_{\gamma} $ denote the copywise local time field generated by $\gamma$.
Then by Lemma \ref{lem:cpf-spatial-restriction} {(see also Remark \ref{rmk: spatial restriction for excursion})}, we get that \begin{equation}
    \nu_{\fB',\beta}^\xi \big(\mathcal{A}_k\cap \mathcal{B}\big)=\int_{\gamma\in \Gamma}\nu_{\fB',\beta}^\xi(\gamma)\cdot \nu^{\mathsf{f}}_{\fB',\beta-L_{\gamma}}(\mathcal{A}_k),\quad \nu_{\fB',\beta}^\xi \big(\mathcal{A}_k^\zz\cap \mathcal{B}\big)=\int_{\gamma\in \Gamma}\nu_{\fB',\beta}^\xi(\gamma)\cdot \nu^{\mathsf{f}}_{\fB',\beta-L_{\gamma}}(\mathcal{A}_k^\zz).\nonumber
\end{equation}
Applying the preceding two identities in reverse order and using Lemma \ref{lem:yyy}, we get that \begin{align}
    \nu_{\fB',\beta}^\xi \big(\mathcal{A}_k^\zz\cap \mathcal{B}\big)&= \int_{\gamma\in \Gamma}\nu_{\fB',\beta}^\xi(\gamma)\cdot\sum_{m\ge k} \nu^{\mathsf{f}}_{\fB',\beta-L_{\gamma}}(\mathcal{A}_m)\cdot \binom{m}{k} (\frac{\epsilon}{2+\epsilon})^{k}(\frac{2}{2+\epsilon})^{m-k}\nonumber\\&=\sum_{m\ge k} \nu_{\fB',\beta}^\xi \big(\mathcal{A}_m\cap \mathcal{B}\big)\cdot \binom{m}{k} (\frac{\epsilon}{2+\epsilon})^{k}(\frac{2}{2+\epsilon})^{m-k}.\label{eq: boundary exploration for k loops}
\end{align}
Let $\rho'\in(0,\rho)$ be determined later.
By Corollary~\ref{cor: DMP for loop} and $ \mathsf{RSW} $, we get that for any local time $\beta'$ such that $\beta_u'\le \beta_0$,
\begin{equation}
   \nu_{\fB',\beta'}^\xi \big( \mathcal{B}^c\big) \le (\bC_1\delta)^{-\bc_2}.\label{eq: boundary crossing prob in good box existence}
\end{equation} We slightly abuse the notation that \eqref{eq: boundary crossing prob in good box existence} holds for both the XY-measure and the XYZ-measure (the result for the XYZ measure can be derived by conditioning on the $ \zz $-copy first).
We choose $\bC_1>0$ large enough (depending on $\beta_0$) such that $(\bC_1\delta)^{-\bc_2}\le \frac{\rho'}{2}$. Let $q_k=\nu_{\fB',\beta}^\xi \big(\mathcal{A}_k\mid\mathcal{B}\big)$. By Corollary~\ref{cor: Gaussian tail for number of macroscopic loops} (applied to the $ \xx\yy $-copies after conditioning on the $ \zz $-copy and by the fact that the number of loops in the $\zz$-copy is dominated by that in the $\xx$-copy), we get that $ \nu^{\xi}_{\mathsf{B}^{\prime}, \beta}(\mathcal{A}_k)\leq 2c^{-1}\exp(-ck^2) $ {for $k\le \diam(\fB)$. Recall Definition \ref{def: good box exploration}, we get that $\diam(\fB)\ge M^{0.01}$ and we can without loss of generality assume $M$ is large enough such that $M^{0.01}\ge K_0$ in Lemma \ref{lem: super exponential decay analysis}.} Combined with~\eqref{eq: boundary crossing prob in good box existence}, it yields that
$q_k\le 2\bc_3^{-1}\exp(-\bc_3k^2)$ and thus $(q_k)_{k\geq0}$ satisfies the condition in Lemma \ref{lem: super exponential decay analysis}. Then we get $k_*$ such that $q_{k_*}\ge c_{19}$, $k_*\le \frac{1}{c_{19}}$ and satisfies \eqref{eq: proof of lemma 4.5 - 2} with $ \rho $ replaced by $ \rho^{\prime} $.
We say a loop is macroscopic, if it intersects $ \mathsf{B}^{(2)}\setminus \mathsf{B} $ and has diameter at least $ \delta\cdot\mathrm{diam}(\mathsf{B}) $, and let $ \mathcal{C}^{\xx} $ denote the event that there exists a macroscopic $\mathsf{x}$-loop.
We claim that \begin{equation}
    \nu^\xi_{\fB', \beta}(\mathcal{C}^{\mathsf{x}}\mid \mathcal{A}_{k_*}^\zz\cap \mathcal{B})\le \rho'.\label{eq: small loop crossing prob conditioned on z}
\end{equation}
To prove the claim, we calculate that
    \begin{align}
        &\quad\,\nu^\xi_{\fB', \beta}(\mathcal{C}^{\xx}\mid \mathcal{A}_{k_*}^\zz\cap \mathcal{B})\nonumber\\
        &\le \frac{\nu^\xi_{\fB', \beta}\bigg(
        \begin{gathered}
            \exists\,\text{at least ${k_*}+1$  macroscopic loops, and}\nonumber\\
            \exists\,\text{exactly ${k_*}$ them are in $\zz$-copy}
        \end{gathered}
        \bigg)}{\nu^\xi_{\fB', \beta}(\mathcal{A}_{k_*}^\zz\cap \mathcal{B})}\nonumber\\
        &= \frac{\sum_{m\geq {k_*}+1} q_m{m \choose {k_*}}(\frac{\epsilon}{2+\epsilon})^{k_*}(\frac{2}{2+\epsilon})^{m-{k_*}} }{\sum_{\ell\geq{k_*}}q_{\ell}{\ell \choose {k_*}}(\frac{\epsilon}{2+\epsilon})^{k_*}(\frac{2}{2+\epsilon})^{\ell-{k_*}}}\leq \frac{\sum_{m\geq {k_*}+1} q_m{m \choose {k_*}}(\frac{\epsilon}{2+\epsilon})^{k_*}(\frac{2}{2+\epsilon})^{m-{k_*}}}{q_{k_*}\cdot(\frac{\epsilon}{2+\epsilon})^{k_*}}\nonumber\\
        &= \frac{\sum_{m\geq {k_*}+1} q_m{m \choose {k_*}}(\frac{2}{2+\epsilon})^{m-{k_*}v}}{q_{k_*}}\stackrel{\eqref{eq: proof of lemma 4.5 - 2}}{\le}\rho'.\label{eq: proof of lemma 4.5 - 1}
    \end{align}
This completes the proof of \eqref{eq: small loop crossing prob conditioned on z}.
Define $ \mathcal{E} $ to be the collection of realizations for the $\zz$-copy local time field $ L^{\zz} $ on $\fB'$ such that for $ L^{\zz}\in\mathcal{E} $,
\begin{equation}
    \mu^\xi_{\fB', \beta-L^{\zz}}(\mathcal{C}^{\xx}\mid \mathcal{B})\le 2\rho'.\label{eq: def of good z local time}
\end{equation} For any $\zz$-copy local time field $ L^{\zz} $, we write $\mu^\f_{\fB', L^{\zz} }$ to be the loop measure on $\zz$-copy with local time condition $ L^{\zz} $ and free boundary condition (this echoes the assumption that $\xi$ is an $\mathsf xy$-boundary condition).
Then we have
\begin{align}
    \nu^\xi_{\fB', \beta}(\mathcal{C}^{\xx}\cap\mathcal{A}_{k_*}^\zz\cap \mathcal{B})&=\sum_{L^{\zz}}\nu^\xi_{\fB', \beta}(L^{\zz})\cdot\mu^{\f}_{\fB', L^{\zz} }(\mathcal{A}_{k_*}^\zz)\cdot  \mu^\xi_{\fB', \beta-L^{\zz}}(\mathcal{C}^{\xx}\cap \mathcal{B}).\label{eq: expansion for good z box 2}
\end{align} Plugging \eqref{eq: expansion for good z box 2} into \eqref{eq: proof of lemma 4.5 - 1} and recalling the definition of $\mathcal{E}$ from \eqref{eq: def of good z local time}, we get that \begin{equation}
    \nu^\xi_{\fB', \beta}(\mathcal{C}^{\xx}\cap\mathcal{A}_{k_*}^\zz\cap \mathcal{B})\ge\sum_{L^{\zz}\in\mathcal{E}^c}\nu^\xi_{\fB', \beta}(L^{\zz})\cdot\mu^{\f}_{\fB', L^{\zz} }(\mathcal{A}_{k_*}^\zz)\cdot  \mu^\xi_{\fB', \beta-L^{\zz}}(\mathcal{B})\cdot 2\rho'.\label{eq: expansion for good z box 3}
\end{equation}
Combining with \eqref{eq: proof of lemma 4.5 - 1}, we get that
\begin{align}
   \nu^\xi_{\fB', \beta}(\mathcal{A}_{k_*}^\zz\cap \mathcal{B})&\stackrel{\eqref{eq: proof of lemma 4.5 - 1}}{\ge}\frac{1}{\rho'} \nu^\xi_{\fB', \beta} (\mathcal{C}^{\xx}\cap\mathcal{A}_{k_*}^\zz\cap \mathcal{B})\nonumber\\ &\stackrel{\eqref{eq: expansion for good z box 3}}{\ge}2 \sum_{L^\zz\in\mathcal{E}^c}\nu^\xi_{\fB', \beta}(L^{\zz})\cdot  \mu^{\f}_{\fB', L^{\zz} }(\mathcal{A}_{k_*}^\zz) \mu^\xi_{\fB', \beta-L^{\zz}}(\mathcal{B})\nonumber\\ &~=2\nu^\xi_{\fB', \beta}(\mathcal{E}^c\cap \mathcal{A}_{k_*}^\zz\cap \mathcal{B}).\nonumber
\end{align} Therefore, we conclude that
\begin{align}
    \nu^\xi_{\fB', \beta}(\mathcal{E})&~\ge \nu^\xi_{\fB', \beta}(\mathcal{E}\cap \mathcal{A}_{k_*}^\zz\cap \mathcal{B})\ge \frac{1}{2}\nu^\xi_{\fB', \beta}(\mathcal{A}_{k_*}^\zz\cap \mathcal{B}).\nonumber
\end{align}
Combined with \eqref{eq: boundary exploration for k loops}, it yields that  \begin{equation}
     \nu^\xi_{\fB', \beta}(\mathcal{E})\ge\frac{1}{2} q_{k_*}(\frac{\epsilon}{2+\epsilon})^{k_*}\times\nu^\xi_{\fB', \beta}(\mathcal{B})\ge\frac{c_{19}}{4}\cdot (\frac{\epsilon}{2+\epsilon})^{\frac{1}{c_{19}}}\times\nu^\xi_{\fB', \beta}(\mathcal{B}).\nonumber
\end{equation}
Combining with \eqref{eq: boundary crossing prob in good box existence}, we obtain that 
\begin{align}
     \nu^\xi_{\fB', \beta}(\mathcal{E})\ge\frac{c_{19}}{8}\cdot (\frac{\epsilon}{2+\epsilon})^{\frac{1}{c_{19}}}.\label{eq: weak z good box prob lower bound}
\end{align}
Combining \eqref{eq: boundary crossing prob in good box existence} with \eqref{eq: def of good z local time}, we get for any $ L^{\zz}\in\mathcal{E} $,
\begin{equation}
    \mu^\xi_{\fB', \beta-L^{\zz}}(\mathcal{C}^{\xx})\le\mu^\xi_{\fB', \beta-L^{\zz}}(\mathcal{C}^{\xx}\mid \mathcal{B})+\mu^\xi_{\fB', \beta-L^{\zz}}( \mathcal{B}^c)\le 3\rho'.\label{eq: small loop crossing condition on z}
\end{equation} By Corollary \ref{cor: small loop crossing implies strong cluster crossing}, we get that for any $ L^\zz\in\mathcal{E} $ any boundary condition $\eta$ on {$\partial \mathsf{B}^{(10)} $}
\begin{equation}
     \mu^\eta_{\fB', \beta-L^{\zz}}(\partial\fB^{(2)}\con{\xx}\partial\fB)\le \bC_{4} (\rho')^{0.1}.\nonumber
\end{equation}
Choosing $\rho'>0$ small such that $\bC_{4} (\rho')^{0.1}\le \rho/2$ and combining with \eqref{eq: weak z good box prob lower bound} finishes the proof of Lemma~\ref{lem: good box probability} with $ C_4=\bC_1 $ (since the bound on the $ \yy $-copy can be derived in the same way).
\end{proof}

\bibliographystyle{alpha}
\bibliography{ref}
\end{document}